\newif \ifSp
\Spfalse     

\ifSp
  \documentclass[graybox,envcountchap,envcountsame,smartqed]{svmono}
\else
  \documentclass[smartqed]{book}
  \usepackage{geometry}
  \usepackage{fancyhdr}
\fi

\usepackage{tikz}

\usepackage{type1cm}

\usepackage{makeidx}         
\usepackage{graphicx}        
\usepackage{multicol}        
\usepackage[bottom]{footmisc}

\usepackage{amsfonts,amssymb,amsmath}

\usepackage{newtxtext}       %
\usepackage{newtxmath}       

\usepackage{ifthen}
\usepackage{ifmtarg}
\usepackage{algpseudocode}
\usepackage{algorithm}
\usepackage{algorithmicx}
\usepackage{bm} 
\usepackage{bbm}
\usepackage{psfrag}
\usepackage{epsfig}
\usepackage{xspace}
\usepackage[caption=false]{subfig}

\usepackage{hho_package}

\usepackage[frozencache]{minted}

\ifSp

\spdefaulttheorem{Exp}{Example}{\bfseries}{\rmfamily}
\spdefaulttheorem{Exo}{Exercise}{\bfseries}{\rmfamily}
\spdefaulttheorem{hypothesis}{Hypothesis}{\bfseries}{\rmfamily}

\newcommand*{\bTheo}[1][]{
\ifthenelse{\equal{#1}{}}{\begin{theorem} \ \par}{\begin{theorem}[#1] } \ \par}
\def\eTheo{\end{theorem}}
\newcommand*{\bLem}[1][]{
\ifthenelse{\equal{#1}{}}{\begin{lemma} \ \par}{\begin{lemma}[#1]} \ \par}
\def\eLem{\end{lemma}}
\newcommand*{\bProp}[1][]{
\ifthenelse{\equal{#1}{}}{\begin{proposition} \ \par }{\begin{proposition}[#1] } \ \par}
\def\eProp{\end{proposition}}

\else
\usepackage{amsthm}

\makeatletter
\def\th@plain{%
  \thm@notefont{}
  \itshape 
}
\def\th@definition{%
  \thm@notefont{}
  \normalfont 
}
\makeatother

\newtheorem{theorem}{Theorem}[chapter]{\bfseries}{\it}
\newtheorem{proposition}[theorem]{Proposition}{\bfseries}{\it}
\newtheorem{lemma}[theorem]{Lemma}{\bfseries}{\it}
\newtheorem{corollary}[theorem]{Corollary}{\bfseries}{\it}
\theoremstyle{definition}
\newtheorem{definition}[theorem]{Definition}{\bfseries}{\rmfamily}
\newtheorem{remark}[theorem]{Remark}{\bfseries}{\rmfamily}
\newtheorem{Exp}[theorem]{Example}{\bfseries}{\rmfamily}
\newtheorem{hypothesis}[theorem]{Hypothesis}{\bfseries}{\rmfamily}

\def\bTheo{\begin{theorem}}
\def\eTheo{\end{theorem}}

\def\bLem{\begin{lemma}}
\def\eLem{\end{lemma}}

\def\bProp{\begin{proposition}}
\def\eProp{\end{proposition}}

\def\bCor{\begin{corollary}}
\def\eCor{\end{corollary}}

\fi

\def\bProof{\begin{proof}}
\def\eProof{\end{proof}}

\def\bDef{\begin{definition}}
\def\eDef{\hfill$\square$\end{definition}}

\def\bRem{\begin{remark}}
\def\eRem{\hfill$\square$\end{remark}}

\def\bExp{\begin{Exp}}
\def\eExp{\hfill$\square$\end{Exp}}

\def\bHypo{\begin{hypothesis}}
\def\eHypo{\end{hypothesis}}

\def\bproof{\begin{proof}}
\def\eproof{\end{proof}}

\newcommand{\eq}[1]{\begin{equation} #1 \end{equation}}

\makeatletter

\@ifundefined{bigtimes}{\def\bigtimes{\bigcross}}{}

\makeatother

\AtBeginEnvironment{minted}{%
	}
\expandafter\def\csname PY@tok@err\endcsname{}

\makeindex             

\begin{document}

\ifSp
\author{Matteo Cicuttin, Alexandre Ern, Nicolas Pignet}
\title{Hybrid high-order methods. \\ A primer with applications to solid mechanics}
\subtitle{-- Monograph --} 
\else
\renewcommand{\thefootnote}{\fnsymbol{footnote}} 
\author{Matteo Cicuttin\footnotemark[2], Alexandre Ern\footnotemark[3], Nicolas Pignet\footnotemark[4]}
\title{Hybrid high-order methods. \\ A primer with applications to solid mechanics\footnotemark[1]}
\fi

\maketitle

\ifSp \else
\renewcommand{\thefootnote}{\fnsymbol{footnote}} 
\footnotetext[1]{This is a preprint of the following work: M. Cicuttin, A. Ern, N. Pignet, \textit{Hybrid high-order methods. A primer with applications to solid mechanics}, Springer, (in press) 2021, reproduced with the permission of the publisher.}
\footnotetext[2]{University of Li\`ege (Montefiore Institute), All\'ee de la d\'ecouverte 10, B-4000 Li\`ege, Belgium}
\footnotetext[3]{CERMICS, \'Ecole des Ponts, 6 \& 8 avenue Blaise Pascal, F-77455 Marne-la-vall\'ee cedex 2, France and INRIA Paris, F-75589 France}
\footnotetext[4]{EDF R\&D, 7 Boulevard Gaspard Monge, F-91120 Palaiseau, France}
\renewcommand{\thefootnote}{\arabic{footnote}}
\thispagestyle{empty}
\phantom{a}
\pagebreak
\fi

\frontmatter

\ifSp \else
\fancyhead[LE]{}
\fancyhead[RE]{Preface}
\fancyhead[LO]{}
\fancyhead[RO]{}
\fi
%
%

\ifSp
  \preface
\else
  \chapter*{Preface}
  \addcontentsline{toc}{chapter}{Preface}
\fi

Hybrid high-order (HHO) methods attach discrete unknowns to the cells and to the faces of the mesh. 
At the heart of
their devising lie two intuitive ideas: \textup{(i)} a local operator reconstructing
in every mesh cell a gradient (and possibly a potential for the gradient) from
the local cell and face unknowns and \textup{(ii)} a local stabilization operator
weakly enforcing in every mesh cell 
the matching of the trace of the cell unknowns with the face unknowns. These two
local operators are then combined into a local discrete bilinear form, and the
global problem is assembled cellwise as in standard finite element methods. HHO methods
offer many attractive features: support of polyhedral meshes, 
optimal convergence rates,
local conservation principles, a dimension-independent formulation,
and robustness in various regimes (e.g.,
no volume-locking in linear elasticity). Moreover, their
computational efficiency hinges on the possibility of locally eliminating
the cell unknowns by static condensation, leading to a global transmission
problem coupling only the face unknowns.

HHO methods were introduced in \cite{DiPEL:14,DiPEr:15}
for linear diffusion and quasi-incom\-pressible linear elasticity. 
A high-order method in mixed form sharing the same devising principles was introduced
in \cite{DiPEr:17}, and shown in \cite{AgBDP:15} to lead after hybridization to a 
HHO method with a slightly different, yet equivalent, writing of
the stabilization.
The realm of applications of HHO methods has been substantially expanded over the last few years.
Developments in solid mechanics include 
nonlinear elasticity \cite{BoDPS:17}, hyperelasticity \cite{AbErPi:18}, plasticity \cite{AbErPi:19,AbErPi:19a}, poroelasticity \cite{BoBDP:16,BoDPS:20}, Kirchhoff--Love plates \cite{BDPGK:18}, 
the Signorini \cite{CaChE:20}, obstacle \cite{CiErG:20} and
two-membrane contact \cite{DabDe:20} problems, 
Tresca friction \cite{ChErPi:20}, and acoustic and elastic wave propagation \cite{BDE:21,BDES:21}. Those related to fluid mechanics include 
convection-diffusion in various regimes \cite{DPDEr:15}, 
Stokes \cite{AgBDP:15,DPELS:16}, 
Navier--Stokes \cite{DiPKr:18,BoDiPD:19,CaDiP:20}, Bingham \cite{CaBCE:18}, 
creeping non-Newtonian \cite{BoCDH:21}, 
and Brinkman \cite{BDiPD:18} flows,
flows in fractured porous media \cite{ChDPF:18,HePiE:21}, 
single-phase miscible flows \cite{AndDr:18}, and
elliptic \cite{BurEr:18} and Stokes \cite{BuDeE:21} interface problems.
Other interesting applications include 
the Cahn--Hilliard problem \cite{ChDMP:16}, Leray--Lions equations \cite{DiPDr:17}, 
elliptic multiscale problems \cite{CiErL:19}, $H^{-1}$ loads \cite{ErnZa:20}, spectral problems \cite{CaCDE:19,CaErP:21}, domains with curved boundary \cite{BoDiP:18,BurEr:18,BurEr:19},
and magnetostatics \cite{ChDPL:20}.

Bridges and unifying viewpoints emerged progressively between HHO methods and several other
discretization methods which also attach unknowns to the mesh cells and faces. Already
in the seminal work \cite{DiPEL:14}, a connection was established between the lowest-order
HHO method and the hybrid finite volume method from \cite{Eymard2010} (and, thus, to the
broader setting of hybrid mimetic mixed methods in \cite{DrEGH:10}). Perhaps the most 
salient connection was made in \cite{CoDPE:16} where HHO methods were embedded into the broad
setting of hybridizable discontinuous Galerkin (HDG) methods \cite{CoGoL:09}. One originality of equal-order HHO methods is the use of the (potential) reconstruction operator in the stabilization. 
Moreover, the analyses of HHO and HDG methods follow somewhat different paths, since the former relies on orthogonal projections, whereas the latter often invokes a more specific approximation operator \cite{CoGoS:10}. We believe that the links between HHO and HDG methods are mutually beneficial, as, for instance, recent HHO developments can be transposed to the HDG setting.
Weak Galerkin (WG) methods \cite{WangYe:13,WangYe:14}, which were embedded into the 
HDG setting in \cite[Sect.~6.6]{Cockburn:16}, are, thus, also closely related to HHO.
WG and HHO were developed independently and share a common devising 
viewpoint combining reconstruction (called weak gradient in WG) 
and stabilization. Yet, the WG stabilization often relies on plain least-squares penalties, whereas the more sophisticated HHO stabilization is key to a higher-order consistency property. Furthermore, the work
\cite{CoDPE:16} also bridged HHO methods to the nonconforming virtual element method \cite{LipMa:14,AyLiM:16}. 
Finally, the connection to the multiscale hybrid mixed method from \cite{HarParVal13} was uncovered in \cite{ChELV:21}.

A detailed monograph on HHO methods appeared this year \cite{DiPDr:20}. The present text is shorter and does not cover as many aspects of the analysis and applications of HHO methods. Its originality lies in targetting the material to computational mechanics without sacrificing mathematical rigor, while including on the one hand some mathematical results with their own specific twist and on the other hand numerical illustrations drawn from industrial examples. Moreover, several topics not covered in \cite{DiPDr:20} are treated here: domains with curved boundary, hyperelasticity, plasticity, contact, friction, and wave propagation. The present material is organized into eight chapters: the first three gently introduce the basic principles of HHO methods on a linear diffusion problem, the following four present various challenging applications to solid mechanics, and the last one reviews implementation aspects.

This book is primarily intended for graduate students, researchers (in applied mathematics, numerical analysis, and computational mechanics), and engineers working in related fields of application. Basic knowledge of the devising and analysis of finite element methods is assumed. Special effort was made to streamline the presentation so as to pinpoint the essential ideas, address key mathematical aspects, present examples, and provide bibliographic pointers. This book can also be used as a support for lectures. As a matter of fact, its idea originated from a series of lectures given by one of the authors during the Workshop on Computational Modeling and Numerical Analysis (Petr\'opolis, Brasil, 2019).

We are thankful to many colleagues for stimulating discussions at various occasions. Special thanks go to G. Delay (Sorbonne University) and S. Lemaire (INRIA) for their careful reading of parts of this manuscript. 

\vspace{\baselineskip}
\begin{flushright}\noindent
Namur and Paris, December 2020\hfill \phantom{a}\\
\hfill {\it Matteo Cicuttin, Alexandre Ern and Nicolas Pignet}
\end{flushright}

\ifSp \else
\fancyhead[LE]{}
\fancyhead[RE]{Table of contents}
\fancyhead[LO]{}
\fancyhead[RO]{}
\fi

\setcounter{tocdepth}{1}
\tableofcontents

\mainmatter

\ifSp \else
\renewcommand{\chaptermark}[1]%
{\markboth{\chaptername\ \thechapter. #1}{}}
\renewcommand{\sectionmark}[1]%
{\markright{\thesection\ #1}}
\fancyhead[LE]{\leftmark}
\fancyhead[RE]{}
\fancyhead[RO]{\rightmark}
\fancyhead[LO]{}
\fi

\chapter{Getting started: Linear diffusion}
\label{chap:diffusion}

The objective of this chapter is to gently introduce the hybrid high-order (HHO)
method on one of the simplest model problems: the Poisson problem
with homogeneous Dirichlet boundary conditions.
Our goal is to present the key ideas underlying the devising of the method
and state its main properties (most of them without proof). The keywords of this
chapter are cell and face unknowns,
local reconstruction and stabilization operators,
elementwise assembly, static condensation, energy minimization, and
equilibrated fluxes.

\section{Model problem}
\label{sec:model_diff}

Let $\Dom$ be an open, bounded, connected, Lipschitz
subset of $\Real^d$ in space dimension $d\ge2$. The one-dimensional case $d=1$ can also be covered, and we refer the reader to Sect.~\ref{sec:1D} for an outline of HHO methods in this setting. Vectors in $\Real^d$ and
vector-valued functions are denoted in bold font,
$\ba\SCAL\bb$ denotes the Euclidean inner product
between two vectors $\ba,\bb\in \Real^d$ and $\|\SCAL\|_{\ell^2}$ 
the Euclidean norm in $\Real^d$.
Moreover, $\#S$ denotes the cardinality of a finite set $S$.

We use standard notation for the Lebesgue and Sobolev spaces; see, \eg \cite[Chap.~4 \& 8]{Brezis:11}, \cite[Chap.~1-4]{ErnGu:21a}, and \cite{AdaFo:03,Evans:98}. In particular, $\Ldeux$ is the Lebesgue space composed of square-integrable functions over $\Dom$, and $\Hun$ is the Sobolev space composed of those functions in $\Ldeux$ whose (weak) partial derivatives are square-integrable functions over $\Dom$. Moreover, $\Hunz$ is the subspace of $\Hun$ composed of functions with zero trace on the boundary $\front$. Inner products and norms in these spaces are denoted by $(\cdot,\cdot)_{\Ldeux}$, $\|\SCAL\|_{\Ldeux}$, $(\cdot,\cdot)_{\Hun}$, and $\|\SCAL\|_{\Hun}$. Recall that for a real-valued function $v$:
\begin{equation}
\|v\|_{\Ldeux}^2\eqq \int_\Dom v^2\dif \bx, \qquad
\|v\|_{\Hun}^2\eqq \|v\|_{\Ldeux}^2+\ell_\Dom^2\|\GRAD v\|_{\Ldeuxd}^2,
\end{equation}
where the length scale $\ell_\Dom\eqq\diam(\Dom)$ (the diameter of $\Omega$) is introduced to be
dimensionally consistent. Owing to the Poincar\'e--Steklov inequality (a.k.a.~Poincar\'e inequality; see \cite[Rem.~3.32]{ErnGu:21a} for a discussion on the terminology), there is
$C_{\textsc{ps}}>0$ such that $C_{\textsc{ps}}\|v\|_{L^2(\Omega)} \le \ell_\Dom \|\GRAD v\|_{\Ldeuxd}$
for all $v\in \Hunz$. 

The model problem we want to approximate in this chapter
is the Poisson problem with source
term $f\in\Ldeux$ and homogeneous
Dirichlet boundary conditions, \ie $-\Delta u =f$ in $\Omega$ and $u=0$ on $\front$.
The weak formulation of this problem reads as follows:
Seek $u \in V\eqq H^1_0(\Dom)$ such that
\begin{equation}\label{eq:weak_diff}
a(u,w)=\ell(w), \quad \forall w\in V,
\end{equation}
with the following bounded bilinear and linear forms:
\begin{equation}
a(v,w)\eqq \psv[\Dom]{\GRAD v}{\GRAD w}, \qquad
\ell(w) \eqq \pss[\Dom]{f}{w},
\end{equation}
for all $v,w\in V$. Since we have $a(v,v)=\|\GRAD v\|_{\Ldeuxd}^2$, the 
Poincar\'e--Steklov inequality implies that the bilinear form
$a$ is coercive on $V$. Hence, the model problem~\eqref{eq:weak_diff} is
well-posed owing to the Lax--Milgram lemma.

\section{Discrete setting}
\label{sec:discrete_diff}

In this section, we present the setting to formulate the HHO discretization of the model problem~\eqref{eq:weak_diff}.

\subsection{The mesh}\label{sec:mesh}

For simplicity, we assume in what follows that the domain $\Dom$ is a polyhedron in $\Real^d$, so that its boundary is composed of a finite union of portions of affine hyperplanes with mutually disjoint interiors. The case of domains with a curved boundary is discussed in Sect.~\ref{sec:curved}.

Since $\Omega$ is a polyhedron, it can be covered exactly by a mesh $\calT$ composed of a finite collection of (open) polyhedral mesh cells $T$, all mutually disjoint, \ie we have
$\overline\Dom=\bigcup_{T\in\calT}\overline T$.
Notice that by definition of a polyhedron,
the mesh cells have straight edges if $d=2$ and planar faces if $d=3$. For a generic mesh cell $T\in\calT$, its boundary is denoted by $\dT$, its unit outward normal by $\nT$, and its diameter by $h_T$. The mesh size is defined as the largest cell diameter in the mesh and is denoted by $h_\calT$, and more simply by $h$ when there is no ambiguity.
When establishing error estimates, one is interested in the process $h\to0$ corresponding to a sequence of successively refined meshes. In this case, one needs to introduce a notion of shape-regularity for the mesh sequence. This notion is detailed in Sect.~\ref{sec:basic_tools}.

\begin{figure}
    \centering
       \includegraphics[scale=0.5]{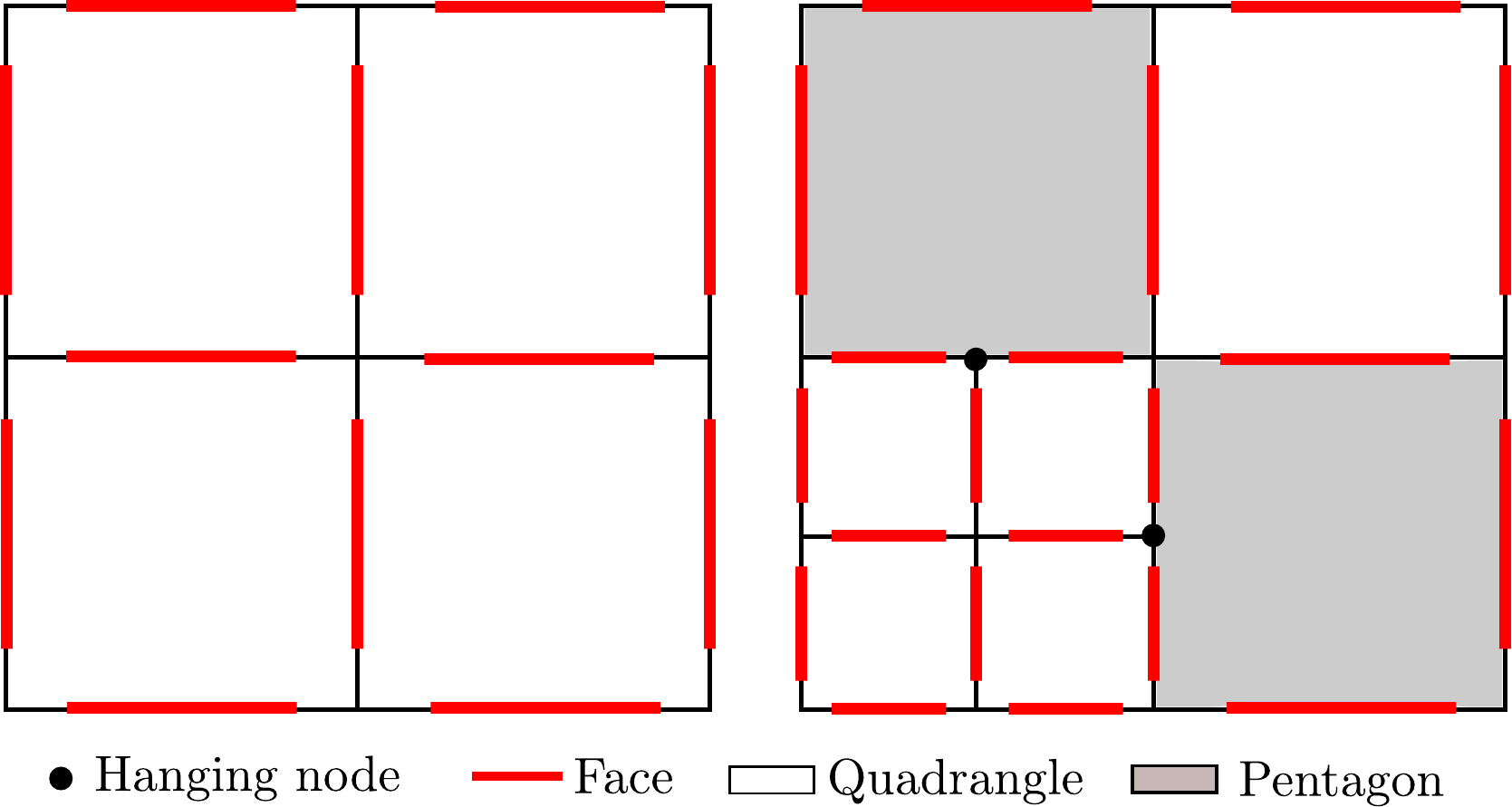}
    \caption{Local refinement of a quadrilateral mesh; the mesh cells containing hanging nodes are treated as polygons (here, pentagons).} \label{fig::hanging_node}
\end{figure}

The possibility of handling meshes composed of polyhedral mesh cells is an attractive feature of HHO methods. For instance, it allows one to treat quite naturally the presence of hanging nodes arising from local mesh refinement; see Figure~\ref{fig::hanging_node} 
for an illustration. However, the reader can assume for simplicity that the mesh is composed of cells with a single shape, such as simplices (triangles in 2D, tetrahedra in 3D) or (rectangular) cuboids, without loosing anything essential in the understanding of the devising and analysis of HHO methods.

Besides the mesh cells, the mesh faces also play an important role in HHO methods. We say that the $(d-1)$-dimensional subset $F\subset \overline\Omega$ is a mesh face if $F$ is a subset of an affine hyperplane, say $H_F$, such that the following holds: \textup{(i)} either there are two distinct mesh cells $T_-,T_+\in\calT$ such that
\begin{equation} \label{eq:def_interface}
F=\dT_-\cap \dT_+ \cap H_F,
\end{equation}
and $F$ is called a (mesh) interface; \textup{(ii)} or there is one mesh cell $T_-\in\calT$ such that
\begin{equation} \label{eq:def_bndy_face}
F=\dT_-\cap \front \cap H_F,
\end{equation}
and $F$ is called a (mesh) boundary face. The interfaces are collected in the set $\calFi$, the boundary faces in the set $\calFb$, so that the set
\begin{equation}
\calF\eqq \calFi \cup \calFb
\end{equation}
collects all the mesh faces. For a mesh cell $T\in\calT$, $\FT$ denotes the collection of the mesh faces composing its boundary $\dT$. Notice that the above definition of the mesh faces implies that each mesh face is straight in 2D and planar in 3D. Hence, for every mesh cell $T\in\calT$, $\nTF$ is a constant vector on every face $F\in\FT$. Notice also that the definitions \eqref{eq:def_interface}
and \eqref{eq:def_bndy_face} do not allow for the case of several coplanar faces that could be shared by two cells or a cell and the boundary, respectively; this choice is only made for simplicity.

\subsection{Discrete unknowns}

The discrete unknowns in HHO methods are
polynomials attached to the mesh cells and to the mesh faces. The idea is that the cell
polynomials approximate the exact solution in the mesh cells, and that the face polynomials
approximate the trace of the exact solution on the mesh faces (although they are not the trace of the cell polynomials). To ease the exposition,
we consider here the equal-order HHO method where the cell and face polynomials have the same
degree. Variants are considered in Sect.~\ref{sec:variants_cell}.

Let $k\ge0$ be the polynomial degree.
Let $\Pkd$ be the space composed of $d$-variate (real-valued) polynomials
of total degree at most $k$. For every mesh cell $T\in\calT$, $\Pkd(T)$ denotes
the space composed of the restriction to $T$ of the polynomials in $\Pkd$.
To define the $(d-1)$-variate polynomial space attached to a mesh face $F\in\calF$
(which is a subset of $\Real^d$),
we consider an affine geometric mapping $\bT_F:\Real^{d-1}\to H_F$ (recall that $H_F$
is the affine hyperplane in $\Real^d$ supporting $F$). Then we set
\begin{equation} \label{eqn:face_affine_mapping}
\PkF(F)\eqq \PkF\circ (\bT_F^{-1})_{|F}.
\end{equation}
It is easy to see that the definition of
$\PkF(F)$ is independent of the choice of the affine geometric mapping $\bT_F$.
(Notice that defining polynomials on the mesh faces is meaningful since we are assuming $d\ge2$.)

\begin{figure}[htb]
\centering
\includegraphics[scale=0.65]{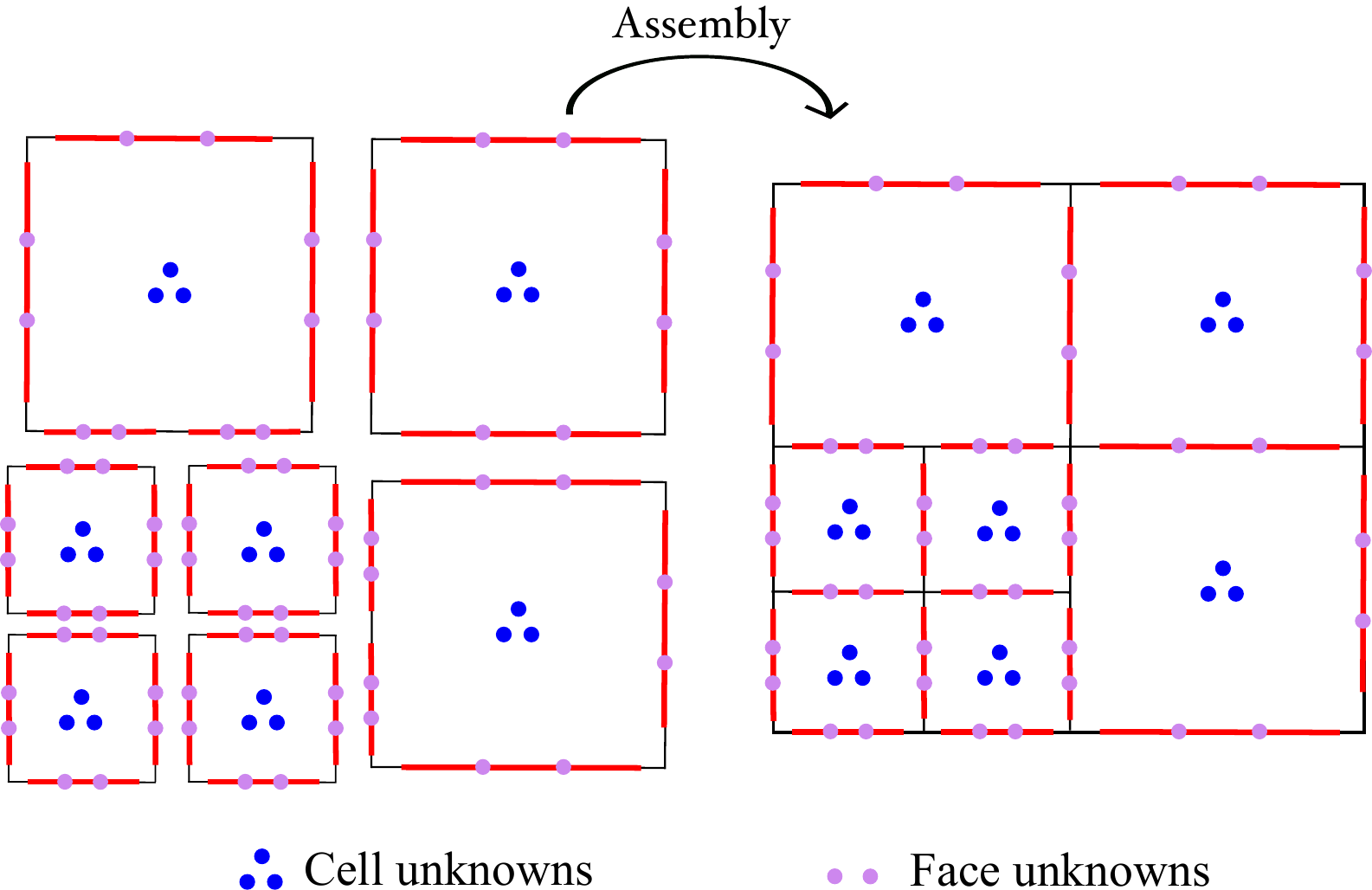}
\caption{Local (left) and global (right) unknowns for the HHO method ($d=2$, $k=1$). Each bullet on the faces and in the cells conventionally represents one basis
function.}
\label{fig:dofs_hho}
\end{figure}

Let us first consider a local viewpoint. For every mesh cell $T\in\calT$, we set
\begin{equation}
\shVkT \eqq \Pkd(T) \times \sVkdT, \qquad
\sVkdT \eqq \bigtimes_{F\in\FT} \PkF(F).
\label{def_of_VkK_VpartialKk}
\end{equation}
A generic element in $\shVkT$ is denoted by $\shvT\eqq (\svT,\svdT)$. We shall
systematically employ the hat notation to indicate a pair of (piecewise) functions, one attached
to the mesh cell(s) and one to the mesh face(s). Notice that the trace of $v_T$ on $\dT$ differs from $v_{\dT}$; in particular, the former is a smooth function over $\dT$, whereas the latter generally exhibits jumps from one face in $\FT$ to an adjacent one.
To define the global discrete HHO unknowns, we follow a similar paradigm; see
Figure~\ref{fig:dofs_hho}.

\bDef[HHO space]
The equal-order HHO space is defined as follows:
\begin{equation}
\shVh \eqq V_{\calT}^k \times V_{\calF}^k,
\qquad
V_{\calT}^k \eqq \bigtimes_{T\in\calT} \Pkd(T),
\qquad
V_{\calF}^k \eqq \bigtimes_{F\in\calF} \PkF(F).
\end{equation}
We have $\dim(\shVh) = {k+d\choose d}\#\calT + {k+d-1\choose d-1}\#\calF$.
\eDef

A generic element in $\shVh$ is denoted by $\shvh\eqq (\vcT,\vcF)$ with
$\vcT\eqq (v_T)_{T\in\calT}$ and $\vcF\eqq (v_F)_{F\in\calF}$. 
Notice that in general $\vcT$ is only piecewise smooth, \ie it can jump 
across the mesh interfaces, and similarly $\vcF$ can jump from one mesh face 
to an adjacent one. Moreover, for all
$\shvh \in \shVh$ and all $T\in\calT$, it is convenient to localize the components
of $\shvh$ associated with $T$ and its faces by using the notation
\begin{equation}
\shvT\eqq \big(v_T,v_{\dT}\eqq (v_F)_{F\in\FT}\big)\in \shVkT.
\end{equation}

At this stage, a natural question that arises is how to reduce a generic function
$v\in\Hun$ (think of the weak solution to~\eqref{eq:weak_diff}) to some member
of the discrete space $\shVh$. In the context of finite elements, this task is usually realized by means of the interpolation operator associated with 
the finite element. In the context of HHO methods, this task is realized in a simple way by considering $L^2$-orthogonal projections. Let $T\in\calT$. Let $\PikTs:L^2(T)\to \Pkd(T)$ and $\PikdTs:L^2(\dT)\to \sVkdT$ be the $L^2$-orthogonal projections defined such that for all $v\in L^2(T)$ and all $w\in L^2(\dT)$,
\begin{alignat}{2}
&(\PikTs(v)-v,p)_{L^2(T)} = 0&\quad&\forall p\in \Pkd(T),\label{eq:def_loc_proj_a}\\
&(\PikdTs(w)-w,q)_{L^2(\dT)} = 0&\quad&\forall q\in \sVkdT.\label{eq:def_loc_proj_b}
\end{alignat}
Notice that for all $F\in\FT$, $(\PikdTs(w))_{|F}=\PikFs(w_{|F})$, with the $L^2$-orthogonal projection $\PikFs:L^2(F)\to \PkF(F)$. The global $L^2$-orthogonal projections $\PikcTs:\Ldeux\to V_{\calT}^k$ and $\PikcFs:L^2(\bigcup_{F\in\calF}F)\to V_{\calF}^k$ are defined similarly to~\eqref{eq:def_loc_proj_a}-\eqref{eq:def_loc_proj_b}.

\bDef[HHO reduction operator] 
For all $T\in\calT$, the local HHO reduction operator $\shIkT:H^1(T)\to \shVkT$ is defined such that for all $v\in H^1(T)$, 
\begin{equation} \label{eq:def_IkT}
\shIkT(v)\eqq (\PikTs(v),\PikdTs(v_{|\dT}))\in\shVkT.
\end{equation}
Similarly, the global HHO reduction operator $\shIkh:\Hun\to \shVh$ is defined such that for all $v\in \Hun$,
\begin{equation} \label{eq:def_Ikh}
\shIkh(v)\eqq (\PikcTs(v),\PikcFs(v_{|\calF}))\in\shVh.
\end{equation}
Since $v\in\Hun$, $v_{|\calF}$ is well-defined on all the mesh faces composing $\calF$.
\eDef

\section{Local reconstruction and stabilization}
\label{sec:heart}

Local reconstruction and stabilization operators associated with
each mesh cell lie at the heart of HHO methods. The goal of this section
is to present these two operators and their main properties. In the whole section,
$T\in\Th$ denotes a generic mesh cell.

\subsection{Local reconstruction}
\label{sec:local_rec}

The main purpose of the reconstruction operator is to compute a gradient
in the mesh cell $T\in\Th$ given a pair of discrete unknowns
$\shvT:=(\svT,\svdT)\in\shVkT$. Obviously, a simple possibility is
to take the gradient of the cell unknown. However, as we shall now see,
taking also into account the face unknowns leads to a reconstruction operator
with better approximation properties.

To stay simple, we consider for the time being a local reconstruction
operator $\opRec : \shVkT \to \Pkpd(T)$, so that the gradient is reconstructed
locally as $\GRAD\opRec(\shvT) \in \GRAD \Pkpd(T) \subset \Pkd(T;\Real^d)$
(see Sect.~\ref{sec:variants_rec} for some variants).

\bDef[Reconstruction]
The local reconstruction operator
$\opRec : \shVkT \to \Pkpd(T)$ is such that for all
$\shvT:=(\svT,\svdT)\in\shVkT$, the function
$\opRec(\shvT) \in \Pkpd(T)$ is uniquely defined by the following equations:
\begin{align}
&\psv[T]{\GRAD \opRec(\shvT)}{\GRAD q} =
- \pss[T]{\svT}{\Delta q} + \pss[\dT]{\svdT}{\nT\SCAL\GRAD q},
\label{eq:def_Rec_HHO} \\
&\pss[T]{\opRec(\shvT)}{1} = \pss[T]{\svT}{1},\label{eq:def_Rec_HHO_mean}
\end{align}
where \eqref{eq:def_Rec_HHO} holds for all 
$q\in \Pkpd(T)^\perp\eqq \{q\in \Pkpd(T)\tq (q,1)_{L^2(T)}=0\}$.
\eDef

Integrating by parts
in~\eqref{eq:def_Rec_HHO} readily yields for all $q\in \Pkpd(T)^\perp$,
\begin{equation} \label{eq:def_Rec_HHO_bis}
\psv[T]{\GRAD \opRec(\shvT)}{\GRAD q} = \psv[T]{\GRAD \svT}{\GRAD q} - \pss[\dT]{\svT-\svdT}{\nT\SCAL\GRAD q}.
\end{equation}
Moreover, we notice that $\opRec(\shvT)=\svT$ if
$\svdT=v_{T|\dT}$, \ie $\opRec(\shvT)$ is in $\Pkpd(T)$ and not just in $\Pkd(T)$ only
if $\svdT\ne v_{T|\dT}$.
In practice, computing $\opRec(\shvT)$
requires choosing a basis of $\Pkpd(T)^\perp$,
inverting the corresponding local stiffness matrix of size
${k+d+1\choose d}-1$, and adjusting the mean-value of $\opRec(\shvT)$
in $T$ using~\eqref{eq:def_Rec_HHO_mean}.

To motivate the above definition of $\opRec$, we show that the composed operator
$\opRec\circ\shIkT$ enjoys a higher-order approximation property.

\bLem[Elliptic projection] \label{lem:ell_proj_HHO}
We have $\calE_T^{k+1} = \opRec\circ\shIkT$ where $\calE_T^{k+1}  : H^1(T) \to \Pkpd(T)$
is the elliptic projection uniquely defined such that for all $v\in H^1(T)$,
\begin{alignat}{2}
&\psv[T]{\GRAD\calE_T^{k+1}(v)}{\GRAD q}=\psv[T]{\GRAD v}{\GRAD q},
&\quad&\forall q\in \Pkpd(T)^\perp,  \label{eq:def_ellp_a}\\
&\pss[T]{\calE_T^{k+1}(v)}{1}=\pss[T]{v}{1}.\label{eq:def_ellp_b}
\end{alignat}
\eLem

\bProof
Consider an arbitrary function $v\in H^1(T)$ and to alleviate the notation, let us
set $\phi\eqq \opRec(\shIkT(v))=\opRec(\PikTs(v),
\PikdTs(v_{|\dT}))$.
Using the definition~\eqref{eq:def_Rec_HHO} of $\opRec$, we infer that
\begin{align*}
\psv[T]{\GRAD\phi}{\GRAD q}
&=- \pss[T]{\PikTs(v)}{\Delta q}
+ \pss[\dT]{\PikdTs(v_{|\dT})}{\nT\SCAL\GRAD q} \\
&=
- \pss[T]{v}{\Delta q} + \pss[\dT]{v}{\nT\SCAL\GRAD q}
= \psv[T]{\GRAD v}{\GRAD q},
\end{align*}
for all $q\in \Pkpd(T)^\perp$,
since $\Delta q\in \Pkmd(T)\subset \Pkd(T)$ and
$\nT\SCAL\GRAD q \in \sVkdT$ (here, we use that all the faces are planar
so that $\nT$ is piecewise constant; projectors were removed owing to \eqref{eq:def_loc_proj_a} and \eqref{eq:def_loc_proj_b}).
Moreover, we have
\[
\pss[T]{\phi}{1}=\pss[T]{\opRec(\shIkT(v))}{1} = \pss[T]{\PikTs(v)}{1} =
\pss[T]{v}{1},
\]
owing to the definition of $\opRec$ and $\shIkT$.
These two identities prove that $\phi$ satisfies
\eqref{eq:def_ellp_a}-\eqref{eq:def_ellp_b}, so that
$\phi=\calE_T^{k+1}(v)$ for all $v\in H^1(T)$. Hence, $\calE_T^{k+1} = \opRec\circ\shIkT$.
\eProof

\subsection{Local stabilization}
\label{sec:stab_diff}

The main issue with the reconstruction operator is that
$\GRAD \opRec(\shvT)=\vzero$ does not imply that $\svT$ and
$\svdT$ are constant functions taking the same value.
Indeed, since $\#\FT\ge d+1$, 
we have $\dim(\shVkT)=\dim(\Pkd)+\dim(\PkF)\#\FT
\ge \dim(\Pkd)+\dim(\PkF)(d+1)$,
and the rank theorem together with Lemma~\ref{lem:ell_proj_HHO}
give $\dim(\KER(\opRec)) = \dim(\shVkT)-\dim(\IM(\opRec)) = \dim(
\shVkT)-\dim(\Pkpd)$. Combining these two inequalities and since
$\dim(\polP_{d'}^l)={l+d'\choose d'}$, this shows that
$\dim(\KER(\opRec))\ge1$.

To fix this
issue, a local stabilization operator is introduced. Among various
possibilities, we focus on an operator that maps $\shVkT$ to
face-based functions $\opSt : \shVkT \to \sVkdT$ such that
for all $\shvT \in \shVkT$,
\begin{equation} \label{eq:def_SKk_HHO}
\opSt(\shvT) \eqq \PikdTs\Big(v_{T|\dT}
-\svdT+\big((I-\PikTs)\opRec(\shvT)\big)_{|\dT}\Big),
\end{equation}
where $I$ is the identity operator.
Letting $\delta_{\dT}\eqq v_{T|\dT}
- \svdT$ be the difference between the trace of the cell
component and the face component on $\dT$, we observe that
\begin{equation}
\opRec(\shvT)=\opRec(\svT,v_{T|\dT})-\opRec(0,\delta_{\dT})=v_T-\opRec(0,\delta_{\dT}).
\end{equation}
Since $v_T\in \Pkd(T)$, the operator
$\opSt$ in~\eqref{eq:def_SKk_HHO} can be rewritten as follows:
\begin{equation}\label{eq:def_SKk_HHO_bis}
\opSt(\shvT) = \PikdTs\Big(\delta_{\dT}-
\big((I-\PikTs)\opRec(0,\delta_{\dT})\big)_{|\dT}\Big).
\end{equation}
This shows that $\opSt(\shvT)$ only
depends (linearly) on the difference
$(v_{T|\dT}-\svdT)$.  The role of $\opSt$ is to help
enforce the matching between the trace of the cell component and the
face component. In the discrete problem, this matching is enforced in a
least-squares manner (see~Sect.\ref{sec:assembling_HHO}).  In practice,
computing $\opSt(\shvT)$ requires to evaluate $L^2$-orthogonal
projections in the cell and on its faces, which entails inverting the
mass matrix in $\T$, which is of size ${k+d\choose d}$, and inverting the mass
matrix in each face $F\in\FT$, which is of size
${k+d-1\choose d-1}$.

Let us finally state an important stability result motivating
the introduction of the operator $\opSt$. To this purpose,
we equip the space $\shVkT$
with the following $H^1$-like seminorm: For all $\shvT \in \shVkT$,
\begin{equation} \label{eq:def_norme_locale_HHO}
\snorme[\shVkT]{\shvT}^2 \eqq \normev[T]{\GRAD \svT}^2
+ h_T^{-1}\normes[\dT]{\svT-\svdT}^2.
\end{equation}
Notice that $\snorme[\shVkT]{\shvT}=0$ implies that $\svT$ and
$\svdT$ are constant functions taking the same value. Then, as shown in
Sect.~\ref{sec:stability_proof},
there are $0<\alpha\le \omega < + \infty$, independent of the mesh size $h$, such that
for all $T\in\calT$ and all $\shvT\in \shVkT$,
\begin{align}
\alpha \snorme[\shVkT]{\shvT}^2 &\le \normev[T]{\GRAD\opRec(\shvT)}^2 +
h_T^{-1}\normes[\dT]{\opSt(\shvT)}^2
\le \omega \snorme[\shVkT]{\shvT}^2.\label{eq:stab_HHO}
\end{align}

\subsection{Example: lowest-order case}

Let us briefly illustrate the above reconstruction and stabilization operators
in the lowest-order case where $k=0$.
Then, for all $\shvT\eqq(\svT,\svdT)\in \hat V_T^0$, $\svT$ is constant on $\T$
and $\svdT$ is piecewise constant on $\dT$. Moreover, $\GRAD \opRec(\shvT)$ is a constant vector in $T$, $\opRec(\shvT)\in \polP_d^1(T)$, and $\opSt(\shvT)$
is piecewise constant on $\dT$.

\bProp[Lowest-order realization]
Assume $k=0$. Let $T\in\calT$. Let $\bx_T$ be the barycenter of $T$
and $\bx_F$ that of the face $F\in\FT$.
For all $\shvT\eqq(\svT,\svdT)\in \hat V_T^0$, setting
$v_F\eqq v_{\dT|F}$ for all $F\in\FT$, we have
\begin{alignat}{2}
&\GRAD \opRec(\shvT) = \sum_{F\in\FT} \frac{\mes{F}}{\mes{T}} (v_F-\svT) \bn_{T|F},&&\label{eq:HFV_a}\\
&\opRec(\shvT)(\bx) = \svT + \GRAD \opRec(\shvT)\SCAL(\bx-\bx_T),&\quad&\forall \bx\in T,\label{eq:HFV_b}\\
&\opSt(\shvT)_{|F} = \svT-v_F-\GRAD\opRec(\shvT)\SCAL(\bx_T-\bx_F),&\quad&\forall F\in\FT.\label{eq:HFV_c}
\end{alignat}
\eProp

\begin{proof}
The proof revolves around the fact that
any polynomial $q\in \polP_d^1(T)$ is such that $q(\bx)=\bar{q}_T
+\bG_q\SCAL(\bx-\bx_T)$ for all $\bx\in T$, where $\bar{q}_T$ is the mean-value
of $q$ in $T$ and $\bG_q\eqq \GRAD q$ is a constant vector in $T$.
Using~\eqref{eq:def_Rec_HHO_bis} and $\GRAD\svT=\vzero$ gives for all $q\in\polP_d^1(T)$,
\begin{align*}
\mes{T} \GRAD \opRec(\shvT)\SCAL \bG_q = (\GRAD \opRec(\shvT),\GRAD q)_{\bL^2(T)} &= (\GRAD \svT,\GRAD q)_{\bL^2(T)} \!- (\svT-\svdT,\nT\SCAL\GRAD q)_{L^2(\dT)} \\
&= \sum_{F\in\FT} \mes{F} (v_F-\svT) \bn_{T|F}\SCAL\bG_q.
\end{align*}
Since $\bG_q$ can be chosen arbitrarily in $\Real^d$, this
proves~\eqref{eq:HFV_a}. The expression~\eqref{eq:HFV_b} then follows from
the above characterization of polynomials in $\polP_d^1(T)$ and~\eqref{eq:def_Rec_HHO_mean}.
Finally, since $\Pi_T^0(\opRec(\shvT)) = \svT$ and $\Pi_{F}^0(\opRec(\shvT))=\svT+\GRAD\opRec(\shvT)\SCAL(\bx_F-\bx_T)$ for all $F\in\FT$, inserting these expressions into~\eqref{eq:def_SKk_HHO} yields~\eqref{eq:HFV_c}.
\eProof

The idea of reconstructing a gradient in each mesh cell by means of~\eqref{eq:HFV_a} and adding a stabilization proportional to~\eqref{eq:HFV_c} has been considered
in the hybrid finite volume (HFV) method from \cite{Eymard2010}.

\section{Assembly and static condensation}
\label{sec:disc_pb_diff}

In this section, we present the discrete problem resulting from the HHO
approximation of the weak problem~\eqref{eq:weak_diff}. We then
highlight the algebraic realization of the discrete problem and show that
the cell unknowns can be eliminated locally by a Schur complement technique often
called static condensation. 

\subsection{The discrete problem}
\label{sec:assembling_HHO}

The discrete problem is formulated by means of a discrete bilinear form
$a_h:\shVh\times\shVh\to\Real$ which is assembled cellwise in the same spirit
as in the finite element method. Thus, for all $\shvh,\shwh\in \shVh$, we set
\begin{equation}
a_h(\shvh,\shwh) \eqq
\sum_{T\in\Th} a_T(\shvT,\shwT),
\end{equation}
where we recall that $\shvT\in \shVkT$ (resp., $\shwT\in \shVkT$) collects the components
of $\shvh$ (resp., $\shwh$) associated with the cell $T\in\Th$ and the faces $F\in\FT$
composing its boundary.
The local bilinear form $a_T:\shVkT\times \shVkT\to\Real$ is devised by
using the local reconstruction and stabilization operators introduced in the
previous section by setting
\begin{align} \label{eq:local_aT}
a_T(\shvT,\shwT) &\eqq
\psv[T]{\GRAD\opRec(\shvT)}{\GRAD\opRec(\shwT)}
+ h_T^{-1}\pss[\dT]{ \opSt(\shvT)}{\opSt(\shwT)}.
\end{align}
The first term on the right-hand side is the counterpart of the local term
$(\GRAD v,\GRAD w)_{\bL^2(T)}$ in the exact bilinear form $a$, whereas the second
term acts as a stabilization that weakly enforces the matching between the trace of
the cell unknowns and the face unknowns. Notice that the scaling by $h_T^{-1}$ makes both terms in~\eqref{eq:local_aT} dimensionally consistent and, at the same time, ensures optimally-decaying error estimates (see Chapter~\ref{chap:math}).
Defining the piecewise polynomial space $\Pkpd(\calT)\eqq \{v\in\Ldeux\tq v_{|T}\in\Pkpd(T),\ \forall T\in\calT\}$,
the global reconstruction operator $\opRecT:\shVh\to \Pkpd(\calT)$ is such that
\begin{equation} \label{eq:def_opRecT}
\opRecT(\shvh)_{|T}\eqq \opRec(\shvT), \quad \forall \shvh\in \shVh,\quad \forall T\in\calT.
\end{equation}
We also define the global stabilization bilinear form
$s_h:\shVh\times\shVh\to\Real$ such that
\begin{equation}
s_h(\shvh,\shwh) \eqq \sum_{T\in\Th} h_T^{-1}\pss[\dT]{ \opSt(\shvT)}{\opSt(\shwT)}.
\end{equation}
The discrete bilinear form $a_h$ can then be rewritten as follows:
\begin{equation} \label{eq:def_ah_globale}
a_h(\shvh,\shwh) = (\GRAD_\calT \opRecT(\shvh),\GRAD_\calT
\opRecT(\shwh))_{\Ldeuxd} + s_h(\shvh,\shwh),
\end{equation}
with the broken gradient operator $\GRAD_\calT$ acting locally in every mesh cell.

We enforce strongly the homogeneous Dirichlet boundary condition by zeroing out the
discrete unknowns associated with the boundary faces, \ie we consider the subspace
\begin{equation} \label{eq:def_HHO_Vh0k}
\shVkhz\eqq V_{\calT}^k \times V_{\calF,0}^k,
\quad
V_{\calF,0}^k \eqq \{\svFh\in V_{\calF}^k \st v_{F}=0,\ \forall F\in\Fhb\}.
\end{equation}
The discrete problem is as follows:
\begin{equation}
\left\{
\begin{array}{l}
\text{Find $\shuh\in \shVkhz$ such that}\\[2pt]
a_h(\shuh,\shwh) = \ell(\wcT),
\quad \forall \shwh\eqq(\wcT,\wcF) \in \shVkhz.
\end{array}
\right.
\label{eq:weak_HHO}
\end{equation}
Notice that only the cell component of the test function $\shwh$
is used to evaluate the load term since $\ell(\wcT)\eqq (f,\wcT)_{L^2(\Omega)}
= \sum_{T\in\Th}(f,w_T)_{L^2(T)}$ (we keep the same symbol $\ell$ for simplicity).
A more subtle treatment of the load term is needed if one works with
loads in the dual Sobolev space $H^{-1}(\Omega)$ (see
\cite{ErnZa:20} for further insight).

To establish the well-posedness of~\eqref{eq:weak_HHO}, we prove that the bilinear form
$a_h$ is coercive on $\shVkhz$. To this purpose, we equip this space with a
suitable norm. Recall the $H^1$-like seminorm $|\SCAL|_{\shVkT}$ defined in~\eqref{eq:def_norme_locale_HHO}.

\bLem[Norm]
The following map defines a norm on $\shVkhz$:
\begin{equation} \label{eq:def_norm_Vh0}
\shVkhz\ni \shvh \longmapsto
\norme[\shVkhz]{\shvh}\eqq \bigg(\sum_{T\in\Th}
\snorme[\shVkT]{\shvT}^2\bigg)^{\frac12} \in [0,+\infty).
\end{equation}
\eLem

\begin{proof}
The only nontrivial property to verify is the definiteness of the map.
Let $\shvh\in \shVkhz$ be such that $\norme[\shVkhz]{\shvh}=0$,
\ie $\snorme[\shVkT]{\shvT}=0$ for all
$T\in\Th$. Owing to~\eqref{eq:def_norme_locale_HHO}, we
infer that $v_{T}$ and $\svdT$ are constant functions
taking the same value in each mesh cell. On cells having a boundary
face, this value must be zero since $\svFh$ vanishes on the
boundary faces. We can repeat the argument for the cells sharing an
interface with those cells, and we can move inward and reach all the
cells in $\Th$ by repeating this process a finite number of times.
Thus, $\shvh=(0,0)\in \shVkhz$.
\eProof

\bLem[Coercivity and well-posedness]  \label{lem:coer_glob_HHO}
The bilinear form $a_h$ is coercive on $\shVkhz$, and the
discrete problem~\eqref{eq:weak_HHO} is well-posed.
\eLem

\bProof
The coercivity of $a_h$ follows by summing the lower bound in
\eqref{eq:stab_HHO} over the mesh cells, which yields
\begin{equation} \label{eq:stab_HHO_glob}
a_h(\shvh,\shvh) \ge \alpha \, \norme[\shVkhz]{\shvh}^2,
\quad \forall \shvh\in \shVkhz.
\end{equation}
Well-posedness is
a consequence of the Lax--Milgram lemma.
\eProof

Standard convexity arguments show that 
the weak solution $u\in V= \Hunz$ to \eqref{eq:weak_diff} is the
unique minimizer in $V$ of the energy functional
\begin{equation}
\engy : V \ni v \longmapsto \frac12 \|\GRAD v\|_{\Ldeuxd}^2 - \ell(v) \in \Real.
\end{equation}
The HHO solution $\shuh\in \shVkhz$ of~\eqref{eq:weak_HHO} can also be characterized
as the unique minimizer in $\shVkhz$ of a suitable energy functional, namely
\begin{equation} \label{eq:def_engy_diff}
\engy_h : \shVkhz \ni \shvh \longmapsto \frac12 \|\GRAD_\calT \opRecT(\shvh)\|_{\bL^2(\Dom)}^2
+ \frac12 s_h(\shvh,\shvh)- \ell(\vcT) \in \Real.
\end{equation}

\bProp[HHO energy minimization] \label{prop:HHO_min_diff}
Let $\engy_h:\shVkhz\to \Real$ be defined in~\eqref{eq:def_engy_diff}. Then
$\shuh\in \shVkhz$ solves~\eqref{eq:weak_HHO} if and only if $\shuh$ minimizes $\engy_h$ in
$\shVkhz$.
\eProp

\begin{proof}
Owing to the coercivity of the discrete bilinear form $a_h$ established in
Lemma~\ref{lem:coer_glob_HHO}, the discrete energy functional $\engy_h$ is
strongly convex in $\shVkhz$. Moreover, this functional is
Fr\'echet-differentiable at any $\shvh\in \shVkhz$,
and a straightforward calculation shows that for all $\shwh\in \shVkhz$,
$D\engy_h(\shvh)[\shwh] = a_h(\shvh,\shwh)-\ell(\wcT).$
This proves the claimed equivalence.
\eProof

To streamline the presentation, we postpone the statement and proof of
the main error estimates regarding the HHO method to the next chapter.
At this stage, we merely announce that, under reasonable assumptions,
the (broken) $H^1$-seminorm of the error decays as $\calO(h^{k+1})$
and the $L^2$-norm of the error decays as $\calO(h^{k+2})$ where $h$ denotes the mesh size.
More precise statements can be found in Sect.~\ref{sec:error_anal_HHO}-\ref{sec:error_anal_HHO_L2}. A residual-based a posteriori error analysis can be found in \cite{DiPSp:16}.

\bRem[Face unknowns] 
Consider the energy functional $\engy_{\calF,0}(\vcT,\SCAL):
V_{\calF,0}^k\to \Real$ such that
$\engy_{\calF,0}(\vcT,\SCAL) = \frac12 \|\GRAD \opRecT(\vcT,\SCAL)\|_{\Ldeuxd}^2
+ \frac12 s_h((\vcT,\SCAL),(\vcT,\SCAL))$ for all $\vcT\in V_\calT^k$.
Elementary arguments show that $\engy_{\calF,0}(\vcT,\SCAL)$ admits a unique minimizer in
$V_{\calF,0}^k$ which we denote $\vcF^*(\vcT)\in V_{\calF,0}^k$ for all $\vcT\in V_\calT^k$.
Let $\engy_\calT:V_\calT^k\to \Real$  be the energy functional such that
$\engy_\calT(\vcT)\eqq \engy_h(\vcT,\vcF^*(\vcT))$. Then,
$\shuh=(\ucT,\ucF)\in \shVkhz$ solves~\eqref{eq:weak_HHO} if and only if
$\ucF=\vcF^*(\ucT)$ and $\ucT$ is the unique minimizer of $\engy_\calT$ in $V_\calT^k$.\label{rem:elim_face_unkn}%
\eRem

\subsection{Algebraic realization}
\label{sec:static_HHO}

Let $N_\calT^{k}\eqq \dim(V_\calT^{k})={k+d\choose d}\#\calT$ and let 
$N_{\calF,0}^{k}\eqq \dim(V_{\calF,0}^{k})={k+d-1\choose d-1}\#\calFi$.
Let $(\sfU_\calT,\sfU_\calF)\in \Real^{N_\calT^{k}\times N_{\calF,0}^k}$
be the component vectors of the discrete solution
$\shuh\eqq (\ucT,\ucF)\in \shVkhz$ once bases
$\{\varphi_i\}_{1\le i\le N_\calT^{k}}$
and $\{\psi_j\}_{1\le j\le N_{\calF,0}^k}$ for $V_\calT^{k}$ and $V_{\calF,0}^k$, respectively,
have been chosen. (Notice that the components of $\sfU_\calF$ are attached only to
the mesh interfaces.)
Let $\sfF_\calT\in \Real^{N_\calT^{k}}$ have components
given by $\sfF_i\eqq (f,\varphi_i)_{L^2(\Omega)}$ for all
$1\le i\le N_\calT^{k}$.
The algebraic realization of~\eqref{eq:weak_HHO} is
\begin{equation}
\begin{bmatrix}
\sfA_{\mathcal{T}\mathcal{T}} & \sfA_{\mathcal{T}\mathcal{F}}\\
\sfA_{\mathcal{F}\mathcal{T}} & \sfA_{\mathcal{F}\mathcal{F}}
\end{bmatrix} \begin{bmatrix}
\sfU_\calT\\
\sfU_\calF
\end{bmatrix} = \begin{bmatrix}
\sfF_\calT\\
0
\end{bmatrix},
\label{eq:HHO_alg}
\end{equation}
where the symmetric positive-definite stiffness matrix $\sfA$ is of size $N_\calT^k+N_{\calF,0}^k$ and is composed of the blocks
$\sfA_{\mathcal{T}\mathcal{T}}$, $\sfA_{\mathcal{T}\mathcal{F}}$,
$\sfA_{\mathcal{F}\mathcal{T}}$, $\sfA_{\mathcal{F}\mathcal{F}}$
associated with the bilinear form $a_h$ and the cell and face
basis functions. Assume that the basis functions associated with a given cell
or face are ordered consecutively. Then
the submatrix
$\sfA_{\calT\calT}$ is block-diagonal, whereas this is not the case for the submatrix
$\sfA_{\calF\calF}$ since the entries attached to faces belonging to the same cell
are coupled together. A computationally effective way to solve the linear system~\eqref{eq:HHO_alg} is to eliminate locally the cell unknowns and solve first for the face unknowns. Defining the Schur complement matrix
\begin{equation}
\sfA_{\mathcal{F}\mathcal{F}}^{\textsc{s}}\eqq \sfA_{\mathcal{F}\mathcal{F}}-
\sfA_{\mathcal{F}\mathcal{T}}\sfA_{\mathcal{T}\mathcal{T}}^{-1}\sfA_{\mathcal{T}\mathcal{F}},
\end{equation}
the global transmission problem coupling all the face unknowns is
\begin{equation}
\sfA_{\mathcal{F}\mathcal{F}}^{\textsc{s}}\sfU_{\calF} = - \sfA_{\mathcal{F}\mathcal{T}}\sfA_{\mathcal{T}\mathcal{T}}^{-1}\sfF_\calT.
\end{equation}
This linear system is only of size $N_{\calF,0}^k$.
Once it is solved, one recovers locally the cell unknowns by using that
$\sfU_\calT = \sfA_{\mathcal{T}\mathcal{T}}^{-1}(\sfF_\calT-\sfA_{\mathcal{T}\mathcal{F}}\sfU_\calF).$
This procedure is called static condensation.

It can be instructive to reformulate the above manipulations by working
directly on the discrete bilinear forms $a_T$ and the discrete HHO unknowns.
To this purpose, for every mesh cell $T\in\calT$, we define
$U_\mu \in \Pkd(T)$ for all $\mu \in \sVkdT$, and we define $U_{r}\in \Pkd(T)$
for all $r\in L^2(T)$ as follows:
\begin{alignat}{2}\label{eq:def_Tsk_HHO_1}
&a_T((U_\mu,0),(q,0)) \eqq -a_T((0,\mu),(q,0)), &\quad &\forall q\in \Pkd(T), \\
&a_T((U_{r},0),(q,0)) \eqq \pss[T]{r}{q}, &\quad &\forall q\in \Pkd(T).
\label{eq:def_Tsk_HHO_2}
\end{alignat}
These problems are well-posed since
$a_T$ is coercive on $\Pkd(T)\CROSS\{0\}$
owing to~\eqref{eq:stab_HHO}.

\bProp[Transmission problem] \label{prop:trans_HHO}
The pair $\shuh:=(\ucT,\ucF)\in \shVkhz$ solves the HHO problem  \eqref{eq:weak_HHO} if
and only if the cell component satisfies $u_{T}=U_{\sudT} + U_{f_{|T}}$ for all
$T\in\calT$, and the face component $\ucF\in V_{\calF,0}^k$ solves the following
global transmission problem:
\begin{equation} \label{eq:glob_trans_pb_HHO}
\sum_{T\in\Th} a_T((U_{\sudT},\sudT),(U_{\swdT},\swdT))
= \sum_{T\in\Th} (f,U_{\swdT})_{L^2(T)}, \quad \forall \wcF \in V_{\calF,0}^k.
\end{equation}
\eProp

\bProof
\textup{(i)} Assume that $\shuh$ solves~\eqref{eq:weak_HHO}.
Let $T\in\Th$ and $\swT\in \Pkd(T)$. Since
$a_T((\suT,\sudT),(\swT,0)) = (f,\swT)_{L^2(T)} = a_T((U_{f_{|T}},0),(\swT,0))$,
we infer that
\begin{align*}
a_T((\suT-U_{f_{|T}},\sudT),(\swT,0))
&=0 = a_T((U_{\sudT},\sudT),(\swT,0)),
\end{align*}
showing that $\suT-U_{f_{|T}}=U_{\sudT}$. This implies that
for all $\swdT\in \sVkdT$,
\begin{align*}
\ifSp &a_T((U_{\sudT},\sudT),(U_{\swdT},\swdT)) \\
\else a_T((U_{\sudT},\sudT),(U_{\swdT},\swdT)) \fi
&= a_T((\suT,\sudT),(U_{\swdT},\swdT))
- a_T((U_{f_{|T}},0),(U_{\swdT},\swdT)) \\
&= a_T((\suT,\sudT),(U_{\swdT},\swdT)),
\end{align*}
where we used the symmetry of $a_T$ and $a_T(U_{\swdT},\swdT),(q,0))=0$ for all $q\in
\Pkd(T)$. Summing over $T\in\Th$ and using~\eqref{eq:weak_HHO}
shows that $\ucF$ solves~\eqref{eq:glob_trans_pb_HHO}.
\\
\textup{(ii)} Assume that $\ucF$ solves~\eqref{eq:glob_trans_pb_HHO}. 
Let $\shwh\in \shVkhz$.
Setting $\shuT\eqq (\suT,\sudT)\eqq (U_{f_{|T}},0)+(U_{\sudT},\sudT)$
for all $T\in\calT$, we infer that
\begin{align*}
a_T(\shuT,\shwT) ={}&
a_T((U_{f_{|T}},0)+(U_{\sudT},\sudT),(\swT-U_{\swdT},0)) 
\\& + a_T((U_{f_{|T}},0)+(U_{\sudT},\sudT),(U_{\swdT},\swdT)) \\
={}& a_T((U_{f_{|T}},0),(\swT-U_{\swdT},0)) + (f,U_{\swdT})_{L^2(T)} \\
&+ a_T((U_{\sudT},\sudT),(U_{\swdT},\swdT)) - (f,U_{\swdT})_{L^2(T)} \\
={}& (f,\swT)_{L^2(T)} + a_T((U_{\sudT},\sudT),(U_{\swdT},\swdT)) - (f,U_{\swdT})_{L^2(T)},
\end{align*}
using that $a_T((U_{\sudT},\sudT),(y_T,0))=0$ for all
$y_T\in \Pkd(T)$, a similar argument for $(U_{\swdT},\swdT)$
together with the symmetry of $a_T$, and the definition of $U_{f_{|T}}$.
Summing over $T\in\calT$ and using \eqref{eq:glob_trans_pb_HHO}
shows that $\shuh$ solves~\eqref{eq:weak_HHO}.
\eProof


\section{Flux recovery and embedding into HDG methods}
\label{sec:flux_recovery_HDG}

In this section, following \cite{CoDPE:16}, we uncover equilibrated fluxes in 
the HHO method. These fluxes, which are associated with
all the faces of every mesh cell, are
in equilibrium at every mesh interface and are balanced in every mesh cell with
the source term. With these fluxes in hand, we embed HHO methods 
into the broad class of hybridizable discontinuous Galerkin (HDG) methods.
\enlargethispage{\baselineskip}
\subsection{Flux recovery}
\label{sec:flux_recovery}

Let $\optSt : \sVkdT\to \sVkdT$ for all $T\in\Th$ be defined such that
\begin{equation}
\optSt(\mu) \eqq \PikdTs\Big(\mu-\big((I-\PikTs)\opRec(0,\mu)\big)_{|\dT}\Big),
\end{equation}
so that the stabilization operator satisfies
$\opSt(\shvT) = \optSt(v_{T|\dT}-\svdT)$ (see~\eqref{eq:def_SKk_HHO_bis}).
By definition, the adjoint of $\optSt$, say
$\optSt^{*}: \sVkdT\to \sVkdT$, is such that
$\pss[\dT]{\optSt^{*}(\lambda)}{\mu} = \pss[\dT]{\lambda}{\optSt(\mu)}$
for all $\lambda,\mu\in \sVkdT$.
The numerical fluxes of a pair $\shvh\in \shVh$ at the boundary of every mesh cell $T\in\calT$ are defined as 
\begin{equation} \label{eq:def_num_flux}
\phi_{\dT}(\shvT) \eqq -\nT\SCAL \GRAD \opRec(\shvT)_{|\dT}
+ h_T^{-1}(\optSt^*\circ \optSt) (v_{T|\dT}-\svdT) \in \sVkdT.
\end{equation}

\bProp[HHO rewriting with fluxes]  \label{prop:flux_recovery}
Let $\shuh\in \shVkhz$ solve~\eqref{eq:weak_HHO} and let the numerical fluxes
$\phi_{\dT}(\shuT)\in \sVkdT$
be defined as in~\eqref{eq:def_num_flux} for all $T\in\calT$. The following holds:\\
\textup{(i)} Equilibrium at every mesh interface $F=\partial T_-\cap \partial T_+ \cap H_F\in \calFi$:
\begin{equation} \label{eq:equil_flux}
\phi_{\dT_-}(\shuTm)_{|_F} + \phi_{\dT_+}(\shuTp)_{|_F} = 0.
\end{equation}
\textup{(ii)} Balance with the source term in every mesh cell $T\in\calT$:
\begin{equation} \label{eq:balance_flux}
\psv[T]{\GRAD \opRec(\shuT)}{\GRAD q} + \pss[\dT]{\phi_{\dT}(\shuT)}{q} = \pss[T]{f}{q}, \quad \forall q\in \Pkd(T).
\end{equation}
\textup{(iii)}
\eqref{eq:equil_flux}-\eqref{eq:balance_flux} are an equivalent rewriting of~\eqref{eq:weak_HHO} that fully characterizes the HHO solution $\shuh\in \shVkhz$.
\eProp

\begin{proof}
\textup{(i)} Let $F\in\calFi$. The identity
\eqref{eq:equil_flux} is proved by taking a test function $\shwh$ in~\eqref{eq:weak_HHO} whose only nonzero component is attached to the interface $F$. Let $w_F\in \PkF(F)$ and take $\shwh\eqq (0,w_\calF)$ with $w_\calF\eqq (\delta_{F,F'}w_F)_{F'\in\calF}$, where $\delta_{F,F'}$ is the Kronecker delta. This is a legitimate test function, \ie $\shwh\in \shVkhz$. Letting $\calT_F\eqq\{T_-,T_+\}$, and using the definitions of $a_h$ and $a_T$, we infer that
\begin{align*}
0 &= \!\sum_{T\in\calT_F}a_T(\shuT,(0,w_{\dT})) \\ &=
\!\sum_{T\in\calT_F} \!\psv[T]{\GRAD\opRec(\shuT)}{\GRAD\opRec(0,w_{\dT})} - h_{T}^{-1}\pss[\dT]{\optSt(u_{T|\dT}-\sudT)}{\optSt(w_{\dT})}.
\end{align*}
Using that $\psv[T]{\GRAD\opRec(\shuT)}{\GRAD\opRec(0,w_{\dT})}
= \pss[F]{\nT\SCAL\GRAD\opRec(\shuT)}{w_F}$ and the definition of the adjoint operator
$\optSt^*$ then gives
\begin{align*}
0&= \sum_{T\in\calT_F} \pss[F]{\nT\SCAL\GRAD\opRec(\shuT)}{w_F} -
h_{T}^{-1}\pss[F]{(\optSt^*\circ \optSt)(u_{T|\dT}-\sudT)}{w_F} \\
&= \sum_{T\in\calT_F} \pss[F]{\phi_{\dT}(\shuT)}{w_F}.
\end{align*}
Since $\phi_{\dT}(\shuT)_{|F}\in \PkF(F)$
for all $T\in\calT_F$ and $w_F$ is arbitrary in $\PkF(F)$, we conclude that \eqref{eq:equil_flux} holds true.
\\
\textup{(ii)} Let $T\in\calT$. The identity
\eqref{eq:balance_flux} is proved by taking a test function $\shwh$ in~\eqref{eq:weak_HHO} whose only nonzero component is attached to the mesh cell $T$. Let $q\in \Pkd(T)$
and take $\shwh\eqq (w_\calT,0)$ with $w_\calT\eqq (\delta_{T,T'}q)_{T'\in\calT}$, so that $w_T=q$. This is a legitimate test function, \ie $\shwh\in \shVkhz$.
Since \eqref{eq:def_Rec_HHO} implies that
$\psv[T]{\GRAD\opRec(\shuT)}{\GRAD \opRec(\shwT)} = \psv[T]{\GRAD\opRec(\shuT)}{\GRAD
q} - \pss[\dT]{\nT\SCAL\GRAD\opRec(\shuT)}{q}$, we have
\begin{align*}
\ifSp &\pss[T]{f}{q} = a_T(\shuT,\shwT) \\ \else
\pss[T]{f}{q} &= a_T(\shuT,\shwT) \\ \fi
&= \psv[T]{\GRAD\opRec(\shuT)}{\GRAD \opRec(\shwT)} + h_{T}^{-1}\pss[\dT]{\optSt(u_{T|\dT}-\sudT)}{\optSt(q)} \\
&=\psv[T]{\GRAD\opRec(\shuT)}{\GRAD q} + \pss[\dT]{-\nT\SCAL\GRAD\opRec(\shuT)+h_{T}^{-1}(\optSt^*\circ \optSt)(u_{T|\dT}-\sudT)}{q} \\
&=\psv[T]{\GRAD\opRec(\shuT)}{\GRAD q}+\pss[\dT]{\phi_{\dT}(\shuT)}{q}.
\end{align*}
\textup{(iii)} The last assertion is a direct consequence of the above two proofs since the considered test functions span $\shVkhz$.
\eProof

\subsection{Embedding into HDG methods}
\label{sec:HHO_HDG}

HDG methods were introduced in~\cite{CoGoL:09} (see also \cite{Cockburn:16} for an overview). In such methods, one
approximates a triple, whereas one approximates a
pair in HHO methods.
Let us consider the dual variable $\vecteur{\sigma}\eqq -\GRAD u$
(sometimes called flux), the primal variable $u$, and its trace
$\lambda\eqq u_{|\calF}$ on the mesh faces.
HDG methods approximate the triple
$(\vecteur{\sigma},u,\lambda)$ by introducing some local
spaces $\vecteur{S}_T$, $V_T$, and $V_F$ for all $T\in\calT$ and all
$F\in\calF$, and by defining a numerical flux trace that
includes a stabilization operator.
Defining the global spaces
\begin{align}
\vecteur{S}_{\calT} &\eqq \{\vecteur{\tau}_{\calT}\eqq(\vecteur{\tau}_T)_{T\in\calT}\in\Ldeuxd\tq \vecteur{\tau}_T\in\vecteur{S}_T,\,\forall
T\in\calT\}, \\
V_{\calT} &\eqq \{\svTh\eqq (v_T)_{T\in\calT}\in \Ldeux\tq v_{T}\in V_T,\,\forall
T\in\calT\}, \label{HDG_V_Th}\\
V_{\calF} &\eqq \{\mu_{\calF}\eqq (\mu_F)_{F\in\calF}\in L^2(\calF)\tq \mu_{F}\in V_F,\,\forall
F\in\calF\}, \label{HDG_V_Fh}
\end{align}
as well as $V_{\calF,0}\eqq \{\mu_{\calF}\in V_{\calF}\tq \mu_{F}=0,\, \forall F\in\Fhb\}$,
the HDG method
consists in seeking the triple
$(\vecteur{\sigma}_{\calT},u_{\calT},\lambda_{\calF})\in \vecteur{S}_{\calT}{\times} V_{\calT}{\times}
V_{\calF,0}$ such that the following holds true:
\begin{align}
&\psv[T]{\vecteur{\sigma}_T}{\vecteur{\tau}_T} - \pss[T]{u_T}{\DIV\vecteur{\tau}_T} + \pss[\dT]{\lambda_{\dT}}{\vecteur{\tau}_T\SCAL\nT} = 0, \label{eq:HDG1} \\
&-\psv[T]{\vecteur{\sigma}_T}{\GRAD \swT} + \pss[\dT]{\vecteur{\phi}_{\dT}\SCAL\nT}{\swT} = \pss[T]{f}{\swT}, \label{eq:HDG2} \\
&\pss[F]{\jump{\vecteur{\phi}_{\partial\calT}}\SCAL\no_F}{\mu_F} =0, \label{eq:HDG3}
\end{align}
for all
$(\vecteur{\tau}_T,\swT,\mu_F) \in \vecteur{S}_T{\times} V_T{\times} V_F$, all $T\in\Th$, and
all $F\in\calFi$,
with $\lambda_{\dT}\eqq (\lambda_F)_{F\in\FT}$, the HDG
numerical flux trace $\vecteur{\phi}_{\partial\calT}\eqq(\vecteur{\phi}_{\dT})_{T\in\calT}$
such that
\begin{equation} \label{eq:bphi_HDG}
\vecteur{\phi}_{\dT} \eqq \vecteur{\sigma}_{T|\dT}
+ s_\dT\upHDG(u_{T|\dT}-\lambda_{\dT}) \nT, \quad \forall T\in\calT,
\end{equation}
the normal jump across the interface $F=\partial T_-\cap \partial T_+ \cap H_F\in\calFi$ defined by
\begin{equation}
\jump{\vecteur{\phi}_{\partial\calT}}\SCAL\no_F \eqq \big(\vecteur{\phi}_{\dT_-|F}-\vecteur{\phi}_{\dT_+|F}\big) \SCAL\no_F = \big(\vecteur{\phi}_{\dT_-}\SCAL\no_{T_-}\big)_{|F}+\big(\vecteur{\phi}_{\dT_+}\SCAL\no_{T_+}\big)_{|F},
\end{equation}
(i.e., $\no_F\eqq \no_{T_-|F}=-\no_{T_+|F}$),
and finally, $s_\dT\upHDG$ is a linear
stabilization operator (to be specified).
The equation~\eqref{eq:HDG1} is the discrete
counterpart of $\vecteur{\sigma}=-\GRAD u$,
the equation \eqref{eq:HDG2} that of $\DIV\vecteur{\sigma}=f$, and the
equation~\eqref{eq:HDG3} weakly
enforces the continuity of the normal component of the numerical flux
trace across the mesh interfaces.

Within the above setting, HDG methods are realized by choosing the local spaces
$\vecteur{S}_T$, $V_T$, $V_F$, and the HDG stabilization operator
$s_{\dT}$.
Following \cite{CoDPE:16}, let us apply this paradigm to the HHO method.

\bProp[HHO as HDG method] \label{prop:HHO_HDG}
The HHO method studied above is rewritten as an HDG method by taking
\begin{equation} \label{eq:HDG_local_spaces}
\vecteur{S}_T \eqq \GRAD \Pkpd(T), \qquad
V_T \eqq \Pkd(T), \qquad
V_F \eqq \PkF(F),
\end{equation}
and the HDG stabilization operator
\begin{equation} \label{eq:tau_HHO}
s_{\dT}\upHDG(v) \eqq h_T^{-1} (\optSt^{*} \circ \optSt)(v),
\quad \forall v\in V_{\dT}\eqq \bigtimes_{F\in\FT} V_F.
\end{equation}
The HDG dual variable is then
$\vecteur{\sigma}_T= -\GRAD\opRec(\shuT)$, the HDG trace variable
is $\lambda_{\dT}=u_{\dT}$,
and the HDG numerical flux trace satisfies
$\vecteur{\phi}_{\dT}\SCAL\nT=\phi_{\dT}(\shuT)$ for all $T\in\calT$.
\eProp

\bProof
Owing to the choice \eqref{eq:HDG_local_spaces} for the local spaces and the
definition of the reconstruction operator, \eqref{eq:HDG1} can be rewritten
$\psv[T]{\vecteur{\sigma}_T+\GRAD\opRec(u_T,\lambda_{\dT})}{\vecteur{\tau}_T}=0$
for all $\vecteur{\tau}_T\in \vecteur{S}_T$. Since
$\vecteur{\sigma}_T+\GRAD\opRec(u_T,\lambda_{\dT})\in \vecteur{S}_T$,
this implies that $\vecteur{\sigma}_T=-\GRAD\opRec(u_T,\lambda_{\dT})$.
The rest of the proof is a direct consequence of the identities derived in
Proposition~\ref{prop:flux_recovery}.
\eProof

HHO methods were devised 
independently of HDG methods by adopting the primal
viewpoint outlined in Sect.~\ref{sec:heart}-\ref{sec:disc_pb_diff}, \ie without
introducing a dual variable explicitly. 
The analysis of HHO methods (see the next chapter)
relies on the approximation properties of $L^2$-orthogonal and
elliptic projections, whereas the analysis of HDG methods generally invokes
a specific projection operator using
Raviart--Thomas finite elements \cite{CoGoS:10}
(see also \cite{DuSay:19}).
Furthermore, the HHO stabilization operator from
\cite{DiPEr:15, DiPEL:14} did not have, at the time of its introduction,
a counterpart in the setting of HDG methods. Indeed, this operator uses
the reconstruction operator, so that at any point $\bx\in \dT$,
$s_{\dT}(v)(\bx)$ depends on the values taken by $v$ over the whole boundary $\dT$.
Instead, the HDG stabilization
operator often acts pointwise, that is, 
$s_{\dT}(v)(\bx)$ only depends on the value taken by $v$ at $\bx$. 
The HHO stabilization operator delivers optimal error estimates
for all $k\ge0$ even on polyhedral meshes.
Achieving this result for HDG methods
with a stabilization operator acting pointwise
requires a subtle design of the
local spaces, as explored for instance in~\cite{CocQS:12}.
The Lehrenfeld--Sch\"oberl stabilization~\cite{Lehrenfeld:10,LehSc:16}
for HDG+ methods (where the cell unknowns are one degree higher than the face
unknowns) is of different nature since $s_{\dT}(v)(\bx)$ depends on
the values taken by $v$ on the face containing $\bx$. This operator
is considered in the context of HHO methods in Sect.~\ref{sec:variants_cell}.

\bRem[Weak Galerkin]
The weak Galerkin (WG) method introduced in \cite{WangYe:13,WangYe:14} can also be embedded
into the setting of HDG methods, as shown in \cite[Sect.~6.6]{Cockburn:16}. The gradient of
the HHO reconstruction operator is called weak gradient in the WG method (not to be
confused with the weak gradient in functional analysis). 
HHO and WG methods were developed independently. In WG methods,
the stabilization operator is often based on plain least-squares penalties. A WG method with
Lehrenfeld--Sch\"oberl stabilization was considered in \cite{MuWaY:15}.
\eRem

\section{One-dimensional setting}
\label{sec:1D}

This section briefly outlines the HHO method in 1D. The model problem is then $-u''=f$ in $\Omega\eqq (a,b)$ with the boundary conditions $u(a)=u(b)=0$. We enumerate the mesh vertices as $(x_i)_{0\le i\le N+1}$ with $x_0\eqq a$, $x_{N+1}\eqq b$. Let $T_i \eqq (x_i,x_{i+1})$ be a generic mesh cell of size $h_i$ for all $0\le i\le N$. 
In 1D, the HHO method simplifies since the face unknowns reduce to one real number attached to every mesh vertex. Thus, the choice of the polynomial degree is only relevant to the cell unknowns which are denoted by
$u_\calT\eqq (u_i\eqq u_{T_i})_{0\le i\le N}$ with $u_i\in \polP^k_1(T_i)$ for all $0\le i\le N$.
The face unknowns are denoted by $u_\calF\eqq \lambda\eqq (\lambda_i)_{0\le i\le N+1}$ with $\lambda_i\in\Real$ for all $0\le i\le N+1$, and $\lambda_0=\lambda_{N+1}=0$ owing to the homogeneous Dirichlet boundary condition. We use the (obvious) notation $\lambda\in \Real^{N+2}_{0,0}$ for the face unknowns.
It is convenient to define the piecewise affine polynomial 
$\pi^1_\lambda:\Dom\to\Real$ such that $\pi^1_\lambda(x_i)=\lambda_i$ for all $0\le i\le N+1$.

Let us first consider the case $k=0$. Then, on the cell $T_i$, the discrete unknowns are the real number $u_i$ attached to the cell and the two real numbers $(\lambda_i,\lambda_{i+1})$ attached to the two endpoints of the cell. A direct computation shows that $R'_{T_i}(u_i,(\lambda_i,\lambda_{i+1})) = h_i^{-1}(\lambda_{i+1}-\lambda_i)$ and $S_{\dT_i}(u_i,(\lambda_i,\lambda_{i+1})) = u_i-\frac12(\lambda_i+\lambda_{i+1})$ at both endpoints of $T_i$. The local discrete equations are for all $(v_i)_{0\le i\le N}$ and all $(\mu_i)_{0\le i\le N+1}$, 
\begin{multline} \label{eq:HHO_1d_full}
\sum_{0\le i\le N} \bigg( h_i^{-1}(\lambda_{i+1}-\lambda_i)(\mu_{i+1}-\mu_i) 
\ifSp \\ \fi + 2h_i^{-1}
\big( u_i-\tfrac12(\lambda_i+\lambda_{i+1})\big)
\big( v_i-\tfrac12(\mu_i+\mu_{i+1})\big)\bigg) = \sum_{0\le i\le N} h_i\bar f_i v_i,
\end{multline}
where $\bar f_i$ denotes the mean-value of $f$ over $T_i$. 
Taking first $v_i=\tfrac12(\mu_i+\mu_{i+1})$ for all $0\le i\le N$ leads to
\begin{equation} \label{eq:1d_k=0}
\sum_{0\le i\le N} 
h_i^{-1}(\lambda_{i+1}-\lambda_i)(\mu_{i+1}-\mu_i) = \sum_{0\le i\le N} 
h_i\bar f_i \frac12(\mu_i+\mu_{i+1}),
\end{equation}
which is nothing but the transmission problem identified in Proposition~\ref{prop:trans_HHO}.
Using the piecewise affine polynomials $\pi_\lambda^1$ and $\pi_\mu^1$, \eqref{eq:1d_k=0} can be rewritten as
\begin{equation} \label{eq:1d_k=0_pi}
\int_\Dom (\pi_\lambda^1)'(\pi_\mu^1)'\dxs = \int_\Dom \Pi_\calT^0(f)\pi_\mu^1\dxs,
\quad \forall \mu\in\Real^{N+2}_{0,0},
\end{equation}
where we recall that 
$\Pi_\calT^0$ is the $L^2$-orthogonal projection onto piecewise constant functions. 
We recognize in \eqref{eq:1d_k=0_pi} the usual finite element discretization of the 1D model problem, up to the projection of the source term.
The algebraic realization of~\eqref{eq:1d_k=0_pi} is $\sfA \Lambda = \sfF$, where $\sfA$ is the tridiagonal matrix of size $N$ with entries $(-h_{i-1}^{-1},h_{i-1}^{-1}+h_i^{-1},-h_i^{-1})$, $\Lambda\in \Real^{N}$ is the vector formed by the $\lambda_i$'s at the interior vertices, and $\sfF\in\Real^{N}$ has components given by $\sfF_i\eqq\frac12(h_{i-1}\bar f_{i-1}+h_{i}\bar f_{i})$ for all $1\le i\le N$. Once the $\lambda_i$'s have been computed, the cell unknowns are recovered from~\eqref{eq:HHO_1d_full} by taking arbitrary cell test functions and zero face test functions. This gives
$u_i = \frac12 h_i^2 \bar f_i + \frac12(\lambda_i+\lambda_{i+1})$ for all $0\le i\le N$.
 
A remarkable fact for the HHO method in 1D is that the global transmission problem is the same for all $k\ge1$. Thus, only the way to post-process locally the face unknowns in order to compute the cell unknowns changes if one modifies the polynomial degree.

\bProp[Transmission problem in 1D, $k\ge1$]
For all $k\ge1$, the global transmission problem is: Find 
$\lambda \in\Real^{N+2}_{0,0}$ such that
\begin{equation} \label{eq:1d_k_ge_1}
\int_\Dom (\pi_\lambda^1)'(\pi_\mu^1)'\dxs = \int_\Dom f\pi_\mu^1\dxs,
\quad \forall \mu\in\Real^{N+2}_{0,0}.
\end{equation}
\eProp

\bproof
Since $k\ge1$, we can consider the pair $(\pi^1_\mu,\mu)$
as a test function in the HHO method for all $\mu\in\Real^{N+2}_{0,0}$. 
Since the trace of the cell component equals the face component at every mesh vertex, we infer that $R_{T_i}(\pi^1_{\mu|T_i},(\mu_i,\mu_{i+1}))=\pi^1_{\mu|T_i} \in \polP_1^1(T_i)$, and hence (recall that $k\ge1$), $S_{\dT_i}(\pi^1_{\mu|T_i},(\mu_i,\mu_{i+1}))=0$ for all $0\le i\le N$. Using these identities in the HHO method and letting $r_i\eqq R_{T_i}(u_i,(\lambda_i,\lambda_{i+1}))\in \polP_1^{k+1}(T_i)$ gives
\[
\sum_{0\le i\le N} (r_i',(\pi_{\mu|T_i}^1)')_{L^2(T_i)} = \int_\Dom f\pi_\mu^1\dxs.
\]
Since $(\pi_{\mu|T_i}^1)'$ is constant, it only remains to show that 
$r_i(x_i)=\lambda_i$ and $r_i(x_{i+1})=\lambda_{i+1}$. This is a remarkable property of the reconstruction in 1D for $k\ge1$. To prove this fact, we observe that the definition of the
reconstruction implies that for all $q\in \polP_1^{k+1}(T_i)^\perp$,
\[
\big( r_i(x_{i+1})-\lambda_{i+1}\big) q'(x_{i+1}) -
\big( r_i(x_{i})-\lambda_{i}\big) q'(x_{i}) = \int_{T_i} (r_i-u_i)q''\dxs.
\]
Since $k\ge1$, we can take any polynomial $q\in \polP_1^{2}(T_i)^\perp$. Recalling that $\int_{T_i} (r_i-u_i)\dxs=0$ by definition, the claim follows by taking $q\in \polP_1^{2}(T_i)^\perp$ such that $q'(x)=h_i^{-1}(x-x_i)$ and then such that $q'(x)=h_i^{-1}(x_{i+1}-x)$ for all $x\in T_i$. 
\eproof

\bRem[Comparison with FEM]
The HHO method with cell unknowns of degree at most $k$
has as many discrete unknowns as the finite element method based
on continuous, piecewise polynomials of degree at most $(k+2)$. This latter method
is more efficient to use since it delivers error estimates with one-order higher convergence 
rate while it is also amenable to static condensation. 
\eRem


\chapter{Mathematical aspects}
\label{chap:math}

The objective of this chapter is to put the HHO method presented in the previous
chapter on a firm mathematical ground. In particular, we prove the key stability and
convergence results announced in the previous chapter.

\section{Mesh regularity and basic analysis tools}
\label{sec:basic_tools}

In this section, we give an overview of the basic mathematical notions
underlying the analysis of HHO methods: mesh regularity, functional
and discrete inverse inequalities, and polynomial approximation
properties in Sobolev spaces.

\subsection{Mesh regularity}
\label{sec:mesh_reg}

Recall that a mesh $\Th$ is composed of polyhedral mesh cells, and $h_\calT$
denotes the mesh size, \ie the largest diameter of the cells in $\Th$.
For simplicity, we assume that $\Dom$ is a polyhedron in $\Real^d$, $d\ge2$, so that any
mesh covers $\Dom$ exactly, \ie there is no error in the geometric representation of the computational domain. We address the case of a domain with a curved boundary in Sect.~\ref{sec:curved}. We recall that, by assumption, the mesh faces are planar.
This property is used to assert that the normal derivative of a $d$-variate polynomial at a cell boundary is a piecewise $(d-1)$-variate polynomial.

Since we are interested in a convergence process where the meshes are successively refined,
we consider a mesh sequence $\meshfam$, that is, a countable family of meshes
such that $0$ is the unique accumulation point of $\{h_\calT\}_{\TinT}$.
The notion of shape-regularity of a mesh sequence is crucial when performing the
convergence analysis of any discretization method, since it is instrumental to
derive fundamental results on polynomial approximation in the mesh cells, as well
as various discrete inverse and functional inequalities. In the simple case where
every mesh $\Th\in\meshfam$ is composed of simplices (without hanging nodes), 
the notion of regularity goes back to
Ciarlet \cite{Ciarlet1978}: the mesh sequence is said to be shape-regular if
there exists a shape-regularity parameter
$\rho > 0$ such that for all $\calT\in\meshfam$ and all $T \in \calT$ with 
diameter $h_T$,
$\rho h_T \leq r_T$, where $r_T$ denotes the inradius of the simplex $T$.
In the more general case of meshes composed of polyhedral cells,
the mesh sequence is said to be shape-regular if
\textup{(i)} any mesh $\calT\in\meshfam$ admits a matching simplicial submesh $\calS_\calT$
such that any cell (or face) of $\calS_\calT$ 
is a subset of a cell (or at most one face) of
$\calT$ and \textup{(ii)} there exists a shape-regularity parameter
$\rho > 0$ such that for all $\calT\in\meshfam$ and all $T \in \calT$
and all $S \in \calS_\calT$ such that $S \subset T$,
we have $\rho h_S \leq r_S$ and $\rho h_T \leq h_S$.
The idea of considering a simplicial submesh to define the regularity of a polyhedral
mesh sequence is rather natural. It was considered, \eg in
\cite{Brenn:03} and in \cite{DiPEr:12} in the context of discontinuous Galerkin
methods. This is also the approach followed in the seminal works on
HHO methods \cite{DiPEL:14,DiPEr:15}. We notice that more general approaches are available,
for instance to handle meshes with cells having some very small
faces \cite{CaDGH:17,CaDoG:21}.

In what follows, it is implicitly understood that any mesh belongs to a
shape-regular mesh sequence, and we do not mention explicitly the mesh sequence.
For simplicity, the mesh size is then denoted by $h$.
Moreover, we use the symbol $C$
to denote a generic constant whose value can change at each occurrence as long as it
is uniform in the mesh sequence, so that it is, in particular,
independent of the mesh size $h$. The value of $C$ can depend on
the domain $\Omega$ and the regularity assumptions on the exact solution, the
underlying polynomial
degree (\eg the one used in the HHO method), and the shape-regularity parameter
$\rho$ of the mesh sequence.

\subsection{Functional and discrete inverse inequalities}

Let $S$ be a subset of $\Omega$ 
(typically, $S$ is composed of a collection of mesh cells).
For any locally integrable function $v:S\to\Real$,
$\partial^\alpha v$ denotes the weak partial derivative of $v$
with multi-index $\alpha\eqq(\alpha_1,\ldots,\alpha_d)\in\polN^d$
of length $|\alpha|\eqq \alpha_1+\ldots+\alpha_d$.
Let $m\in\polN$ and recall the Sobolev space
$H^m(S)\eqq\{v\in L^2(S)\tq \partial^\alpha v\in L^2(S),\, \forall \alpha\in\polN^d,\,
|\alpha|\le m\}$ equipped with the following norm and seminorm:
\begin{equation}
\norme[H^{m}(S)]{v}\eqq \Big(\sum_{|\alpha|\le m}
\ell_S^{2|\alpha|}\|\partial^\alpha v\|_{L^2(S)}^2\Big)^{\frac12},
\quad
\snorme[H^{m}(S)]{v}\eqq \Big(\sum_{|\alpha|=m} \|\partial^\alpha v\|_{L^2(S)}^2\Big)^{\frac12},
\end{equation}
where the length scale $\ell_S\eqq\diam(S)$ is introduced to be
dimensionally consistent (notice that the norm and seminorm have different scalings). 
In some cases, we shall also consider Sobolev spaces
of fractional order.
Let $s=m+\sigma\in (0,+\infty)\setminus\polN$
with $m\eqq\floor{s}\in \polN$ and $\sigma\eqq s-m\in(0,1).$
We define $H^{s}(S)\eqq\{ v\in H^{m}(S)\tq
|\partial^\alpha v|_{H^{\sigma}(S)}<+\infty, \, \forall
\alpha\in\polN^d, \, |\alpha|=m\}$ with the Sobolev--Slobodeckij seminorm
\begin{equation}
|w|_{H^{\sigma}(S)}\eqq \bigg(\int_S \int_S
\frac{|w(\bx)-w(\by)|^2}{\|\bx-\by\|_{\ell^2}^{2\sigma+d}}\dif \bx
\dif \by\bigg)^{\frac12}.
\end{equation}
We equip $H^s(S)$ with the seminorm $|v|_{H^{s}(S)}\eqq
\big(\sum_{|\alpha|=m}|\partial^\alpha v|_{H^{\sigma}(S)}^2\big)^{\frac12}$
and the norm $\|v\|_{H^{s}(\Dom)} \eqq \big(\|v\|_{H^{m}(S)}^2 +
\ell_S^{2s} |v|_{H^{s}(S)}^2\big)^{\frac12}$. 

Let us now state two important functional inequalities valid on every mesh cell $T\in\Th$:
the Poincar\'e--Steklov inequality (a.k.a.~Poincar\'e inequality; see \cite[Rem.~3.32]{ErnGu:21a} for a discussion on the terminology) and the multiplicative trace inequality.

\bLem[Poincar\'e--Steklov inequality] \label{lem:PS}
There is $C_{\textsc{ps}}$ such that for all $T\in\Th$ and all $v\in H^1(T)$,
\begin{equation}\label{eq:Poincare}
\|v-\Pi_T^0(v)\|_{L^2(T)} \le C_{\textsc{ps}} h_T\|\GRAD v\|_{\bL^2(T)},
\end{equation}
where $\Pi_T^0(v)$ is the mean-value of $v$ over $T$, \ie the $L^2$-orthogonal projection
of $v$ onto $\polP_d^0(T)$. Moreover, there is $C$ such that for all $s\in (0,1)$, 
all $T\in\Th$ and all $v\in H^s(T)$,
\begin{equation}\label{eq:Poincare_frac}
\|v-\Pi_T^0(v)\|_{L^2(T)} \le C h_T^s|v|_{H^s(T)}.
\end{equation}
\eLem

\bLem[Multiplicative trace inequality]  \label{lem:mtr}
There is $C$ such that for all $T\in\Th$ and all $v\in H^1(T)$,
\begin{equation}\label{eq:mtr}
\|v\|_{L^2(\dT)} \le C\big( h_T^{-\frac12}\|v\|_{L^2(T)} + \|v\|_{L^2(T)}^{\frac12}
\|\GRAD v\|_{\bL^2(T)}^{\frac12}\big).
\end{equation}
Moreover, for all $s\in(\frac12,1)$, there is $C$ such that for all $T\in\Th$ and all $v\in H^s(T)$,
\begin{equation}\label{eq:mtr_frac}
\|v\|_{L^2(\dT)} \le C\big( h_T^{-\frac12}\|v\|_{L^2(T)} + h_T^{s-\frac12}
|v|_{H^s(T)}\big).
\end{equation}
The constant $C$ is uniform with respect to $s$ as long as $s$ is bounded away from $\frac12$.
\eLem

\bRem[Literature]
If the mesh cell $T$ is a convex set, the Poincar\'e--Steklov 
inequality~\eqref{eq:Poincare} holds true
with constant $C_{\textsc{ps}}=\frac{1}{\pi}$ \cite{PayWe:60,Beben:03}. 
In the general case, one
decomposes $T$ into the subsimplices resulting from the shape-regularity 
assumption on the mesh. We refer the reader to \cite[Sect.~2.3]{VeeVe:12}
and \cite[Lem.~5.7]{ErnGu:17} for proofs of Poincar\'e--Steklov inequalities on
composite elements and to \cite[Lem.~7.1]{ErnGu:17} for the fractional Poincar\'e--Steklov inequality~\eqref{eq:Poincare_frac}.
The idea behind the proof of the multiplicative trace inequality~\eqref{eq:mtr} 
in a simplex is
to lift the trace using the lowest-order Raviart--Thomas polynomial associated with
the face in question (see \cite[App.~B]{MonSu:99} and \cite[Thm.~4.1]{CarFu:00}).
In a polyhedral cell, for each subface composing
$\partial T$, one carves a subsimplex inside $T$ having equivalent height
(see, \eg \cite[Lem.~1.49]{DiPEr:12}). For the fractional multiplicative trace inequality,  one considers a pullback to the reference simplex if $T$ is a simplex
\cite[Lem.~7.2]{ErnGu:17}, and one considers a subsimplex as above if $T$ is polyhedral.
\eRem

In contrast to functional inequalities, discrete inverse inequalities are only
valid in polynomial spaces, and their proof hinges on norm equivalence in
a finite-dimensional space. For this reason, discrete inverse inequalities are proven
first on a reference simplex and then a geometric mapping is invoked to pass to
a generic mesh simplex. In the case of a polyhedral mesh cell, one exploits its
decomposition into a finite number of subsimplices; we refer the reader, \eg
to \cite[Lem.~1.44 \& 1.46]{DiPEr:12} for more details.

\bLem[Discrete inverse inequalities] \label{lem:disc_inv}
Let $l\in\polN$ be the polynomial degree.
There is $C$ such that for all $T\in\Th$ and all $q\in \polP_d^l(T)$,
\begin{align}
&\|\GRAD q\|_{\bL^2(T)} \le C h_T^{-1}\|q\|_{L^2(T)}, \label{eq:disc_inv}\\
&\|q\|_{L^2(\dT)} \le Ch_T^{-\frac12}\|q\|_{L^2(T)}.\label{eq:disc_trace}
\end{align}
\eLem

\subsection{Polynomial approximation}

The last question we need to address is how well it is possible to approximate a
given function in some Sobolev space by a polynomial. In the context of HHO methods,
it is sufficient to consider the approximation by the $L^2$-orthogonal projection.

\bLem[Approximation by $L^2$-projection]  \label{lem:poly_approx}
Let $l\in \polN$ be the polynomial degree. Let $\Pi_T^l$ be the
$L^2$-orthogonal projection onto $\polP_d^l(T)$.
There is $C$ such that for all $r\in [0,l+1]$, all $m\in \{0,\ldots,\lfloor r\rfloor\}$,
all $T\in\Th$, and all $v\in H^r(T)$,
\begin{equation} \label{eq:pol_app_L2}
|v-\Pi_T^l(v)|_{H^m(T)} \le Ch_T^{r-m}|v|_{H^r(T)}.
\end{equation}
Moreover, if $r\in(\frac12,l+1]$ and $r\in(\frac32,l+1]$, $l\ge1$, respectively, we have
\begin{equation} \label{eq:pol_app_face}
\|v-\Pi_T^l(v)\|_{L^2(\dT)} \le Ch_T^{r-\frac12}|v|_{H^r(T)},
\quad
\|\GRAD(v-\Pi_T^l(v))\|_{\bL^2(\dT)} \le Ch_T^{r-\frac32}|v|_{H^r(T)}.
\end{equation}
\eLem

\bproof
The estimate \eqref{eq:pol_app_L2} can be proved by standard arguments in a simplicial cell and by using the arguments from the proof of \cite[Lem.~5.6]{ErnGu:17} in a polyhedral cell (the proof combines the Poincar\'e--Steklov inequalities from Lemma~\ref{lem:PS}
with a polynomial built using mean-values of the derivatives of $v$ in $T$). Let us prove the first bound in~\eqref{eq:pol_app_face}. If $r\in[1,l+1]$, we invoke
the multiplicative trace inequality~\eqref{eq:mtr} which yields
\[
\normes[\dT]{v-\Pi_T^l(v)} \le  C\big(h_T^{-\frac12}\normes[T]{v-\Pi_T^l(v)}+h_T^{\frac12}|v-\Pi_T^l(v)|_{H^1(T)}\big).
\]
The first term is bounded using~\eqref{eq:pol_app_L2} with $m\eqq0$ and the second
term using~\eqref{eq:pol_app_L2} with $m\eqq1$ (this is possible since 
$\lfloor r\rfloor\ge1$). If instead $r\in (\frac12,1)$,
the triangle inequality, the discrete trace 
inequality~\eqref{eq:disc_trace}, and the bound $\normes[T]{\Pi_T^l(v)-\Pi_T^0(v)}
\le 2\normes[T]{v-\Pi_T^0(v)}$ imply that
\begin{align*}
\normes[\dT]{v-\Pi_T^l(v)} &\le \normes[\dT]{v-\Pi_T^0(v)} + \normes[\dT]{\Pi_T^l(v)-\Pi_T^0(v)} \\
&\le \normes[\dT]{v-\Pi_T^0(v)} + Ch_T^{-\frac12}\normes[T]{v-\Pi_T^0(v)}.
\end{align*}
Invoking the fractional multiplicative trace inequality~\eqref{eq:mtr_frac}
to bound the first term on the right-hand side, the bound~\eqref{eq:pol_app_face}
follows from~\eqref{eq:pol_app_L2} with $l=m\eqq0$
and $|v-\Pi_T^0(v)|_{H^r(T)}=|v|_{H^r(T)}$.
Finally, the proof of the second bound in~\eqref{eq:pol_app_face} is similar, up to the use of the discrete inverse inequality~\eqref{eq:disc_inv} together with~\eqref{eq:disc_trace}. 
\eproof


\section{Stability}
\label{sec:stability_proof}

Recall that for all $T\in\Th$, $\shVkT$ is equipped with the $H^1$-like seminorm
$\snorme[\shVkT]{\shvT}^2 \eqq \normev[T]{\GRAD \svT}^2
+ h_T^{-1}\normes[\dT]{\svT-\svdT}^2$ for all $\shvT\eqq(\svT,\svdT) \in \shVkT$
(see~\eqref{eq:def_norme_locale_HHO}).

\bLem[Stability]  \label{lem:stab_HHO}
There are $0<\alpha\le \omega<+\infty$ such that
for all $T\in\Th$ and all $\shvT\in \shVkT$,
\begin{align}
\alpha \snorme[\shVkT]{\shvT}^2 &\le a_T(\shvT,\shvT)
\le \omega \snorme[\shVkT]{\shvT}^2,\label{eq:stab_HHO_bis}
\end{align}
recalling that $a_T(\shvT,\shvT)
= \normev[T]{\GRAD\opRec(\shvT)}^2 +
h_T^{-1}\normes[\dT]{\opSt(\shvT)}^2$.
\eLem

\bProof
Let $\shvT\in \shVkT$ and set $r_T\eqq \opRec(\shvT)$. \\
\textup{(i)} Lower bound. Let us first bound $\normev[T]{\GRAD \svT}$.
Taking $q\eqq \svT-\Pi_T^0(\svT)$ in the definition~\eqref{eq:def_Rec_HHO_bis} of $r_T$
and using the Cauchy--Schwarz inequality leads to
\begin{align*}
\ifSp \alhere \fi \normev[T]{\GRAD \svT}^2 
\ifSp \else \alhere \fi = \psv[T]{\GRAD
r_T}{\GRAD \svT} +
\pss[\dT]{\svT-\svdT}{\nT\SCAL\GRAD \svT}\\
&\le\normev[T]{\GRAD r_T} \normev[T]{\GRAD \svT}+ h_{T}^{-\frac12}
\normes[\dT]{\svT-\svdT} h_{T}^{\frac12}\normes[\dT]{\nT\SCAL\GRAD \svT}.
\end{align*}
Invoking the discrete trace inequality \eqref{eq:disc_trace} to bound $\normes[\dT]{\nT\SCAL\GRAD \svT}$ gives
\begin{equation} \label{eq:estim_stab_HHO1}
\normev[T]{\GRAD \svT} \le C \big(\normev[T]{\GRAD r_T} + h_T^{-\frac12}
\normes[\dT]{\svT-\svdT}\big).
\end{equation}
Let us now bound $h_T^{-1}\normes[\dT]{\svT-\svdT}$.
We have
\begin{align}
&\normes[\dT]{\PikdTs(((I-\PikTs)r_T)_{|\dT})}
\le \normes[\dT]{(I-\PikTs)r_T} \nonumber \\
&\le C h_T^{-\frac12} \normes[T]{(I-\PikTs)r_T}
\le C h_T^{-\frac12} \normes[T]{(I-\pizT)r_T}
\le C h_T^{\frac12} \normev[T]{\GRAD r_T}, \label{eq:bnd_on_trace_rT}
\end{align}
owing to the $L^2$-stability of $\PikdTs$,
the discrete trace inequality \eqref{eq:disc_trace}, and the
Poincar\'e--Steklov inequality~\eqref{eq:Poincare} (recall that the value of $C$
can change at each occurrence).
Using the definition \eqref{eq:def_SKk_HHO} of $\opSt$
and the fact that
$v_{T|\dT}-\svdT$ is in $\sVkdT$, we infer that
$v_{T|\dT}-\svdT = \opSt(\shvT) -
\PikdTs(((I-\PikTs)r_T)_{|\dT})$.
The triangle inequality and \eqref{eq:bnd_on_trace_rT} imply that
\[
h_{T}^{-\frac12}\normes[\dT]{\svT-\svdT} \le
h_{T}^{-\frac12}\normes[\dT]{\opSt(\shvT)} +
C\normev[T]{\GRAD r_T}.\]
Combining this estimate with~\eqref{eq:estim_stab_HHO1} proves the lower bound
in~\eqref{eq:stab_HHO_bis}.
\\
\textup{(ii)} Upper bound. Using the definition~\eqref{eq:def_Rec_HHO_bis} of $r_T$ with $q\eqq r_T-\Pi_T^0(r_T)$ leads to
$\normev[T]{\GRAD r_T}^2 = \psv[T]{\GRAD \svT}{\GRAD r_T}
-\pss[\dT]{\svT-\svdT}{\nT\SCAL\GRAD r_T}$.
Invoking the Cauchy--Schwarz inequality and
the discrete trace inequality \eqref{eq:disc_trace} gives
\[
\normev[T]{\GRAD r_T} \le \normev[T]{\GRAD \svT}
+ C h_{T}^{-\frac12}\normes[\dT]{\svT-\svdT}.
\]
Moreover, the triangle inequality and the bound~\eqref{eq:bnd_on_trace_rT} imply that
\begin{align*}
h_{T}^{-\frac12}\normes[\dT]{\opSt(\shvT)} &\le h_{T}^{-\frac12}\normes[\dT]{\svT-\svdT}
+ h_{T}^{-\frac12}\normes[\dT]{\PikdTs(((I-\PikTs)r_T)_{|\dT})}\\
&\le h_{T}^{-\frac12}\normes[\dT]{\svT-\svdT} + C \normev[T]{\GRAD r_T}.
\end{align*}
Combining the above bounds proves the upper bound in~\eqref{eq:stab_HHO_bis}.
\eProof


\section{Consistency}

Let us first prove that the stabilization operator leads to
optimal approximation properties when combined with the reduction
operator $\shIkT$ defined in~\eqref{eq:def_IkT}. 
Recall that $\calE_T^{k+1}$ denotes the elliptic projection
operator onto $\Pkpd(T)$ (see Lemma~\ref{lem:ell_proj_HHO}).

\bLem[Approximation property of $\opSt\circ\shIkT$] 
\label{lem:approx_HHO_St}
There is $C$ such that for all $T\in\Th$ and
all $v\in H^1(T)$,
\begin{equation} \label{eq:approx_HHO_St}
h_{T}^{-\frac12}\normes[\dT]{\opSt(\shIkT(v))}
\le C \normev[T]{\GRAD(v-\calE_T^{k+1}(v))},
\end{equation}
\ie we have $h_{T}^{-\frac12}\normes[\dT]{\opSt(\shIkT(v))}
\le C\min_{q\in \Pkpd(T)}\normev[T]{\GRAD(v-q)}$.
\eLem

\bProof
Let $v\in H^1(T)$ and set $\eta\eqq v-\calE_T^{k+1}(v)$.
Owing to the definition~\eqref{eq:def_SKk_HHO} of $\opSt$,
the definition~\eqref{eq:def_IkT} of $\shIkT$,
and since $\opRec\circ \shIkT = \calE_T^{k+1}$
(see Lemma~\ref{lem:ell_proj_HHO}), we have
\begin{align*}
\opSt(\shIkT(v)) &= \PikdTs\left(\PikTs(v)_{|\dT}
-\PikdTs(v_{|\dT})+((I-\PikTs)\calE_T^{k+1}(v))_{|\dT}\right) \\
&= \PikTs(\eta)_{|\dT} - \PikdTs(\eta_{|\dT}),
\end{align*}
since $\PikdTs(\PikTs(\eta)_{|\dT})=
\PikTs(\eta)_{|\dT}$ and $\PikdTs\circ \PikdTs=\PikdTs$.
Invoking the triangle inequality, the $L^2$-stability of $\PikdTs$,
the discrete trace inequality~\eqref{eq:disc_trace}, and the $L^2$-stability of $\PikTs$ 
leads to (recall that the value of $C$ can change at each occurrence)
\begin{align*}
\ifSp &\normes[\dT]{\opSt(\shIkT(v))} \\ \else
\normes[\dT]{\opSt(\shIkT(v))} \fi &\le 
\normes[\dT]{\PikTs(\eta)} + \normes[\dT]{\PikdTs(\eta_{|\dT})}
\le \normes[\dT]{\PikTs(\eta)} + \normes[\dT]{\eta}\\
&\le Ch_T^{-\frac12}\normes[T]{\PikTs(\eta)}
+ \normes[\dT]{\eta} \le Ch_T^{-\frac12}\normes[T]{\eta}
+ \normes[\dT]{\eta} \le Ch_T^{\frac12}\normev[T]{\GRAD\eta},
\end{align*}
where the last bound follows from the multiplicative
trace inequality~\eqref{eq:mtr}
and the Poincar\'e--Steklov inequality~\eqref{eq:Poincare}
(since $\pss[T]{\eta}{1}=0$). This proves the 
bound~\eqref{eq:approx_HHO_St}, and the bound using the
minimum over $q\in \Pkpd(T)$ readily follows from
the definition of the elliptic projection.
\eProof

Loosely speaking, the consistency error is measured by inserting the exact
solution into the discrete equations and bounding the resulting truncation error.
To realize this operation within HHO methods, the idea is to
insert $\shIkh(u)$ into the discrete equations, where $\shIkh$ is the global reduction
operator defined in~\eqref{eq:def_Ikh}. Notice that $\shIkh(u)\in \shVkhz$
since $u\in\Hunz$. With this tool in hand,
we define the consistency error $\delta_{h} \in (\shVkhz)'$ as the linear form
such that for all $\shwh\in\shVkhz$,
\begin{equation} \label{eq:delta_calI_HHO}
\langle \delta_{h},\shwh\rangle
\eqq \ell(\wcT)- a_h(\shIkh(u),\shwh).
\end{equation}
Bounding the consistency error then amounts to bounding the dual norm
\begin{equation}
\norme[*]{\delta_{h}} \eqq \sup_{\shwh\in\shVkhz}
\frac{|\langle \delta_{h},\shwh\rangle|}{\norme[\shVkhz]{\shwh}},
\end{equation} 
where the
norm $\norme[\shVkhz]{\SCAL}$ is defined in~\eqref{eq:def_norm_Vh0}.
It is implicitly understood here and in what follows that the 
argument is nonzero when evaluating
the dual norm by means of the supremum.
To avoid distracting technicalities, we henceforth assume
that the exact solution satisfies
$u\in H^{1+r}(\Dom)$, $r>\frac12$. This assumption actually follows from
elliptic regularity theory (see, \eg \cite[p.~158]{Dauge_1988}). It
implies that $\GRAD u$ can be localized as a single-valued function at every mesh face 
(see \cite[Rmk.~18.4]{ErnGu:21a} and also \cite[Sect.~41.5]{ErnGu:21b} on how to go beyond this assumption for heterogeneous diffusion problems).
For all $T\in\Th$ and all $v\in H^{1+r}(T)$, $r>\frac12$, we define the local seminorm
\begin{equation} \label{eq:def_norme_sharp}
\snorme[\sharp,T]{v} \eqq\normev[T]{\GRAD v}
+ h_T^{\frac12}\normev[\dT]{\GRAD v},
\end{equation}
as well as the global counterpart 
\begin{equation} \label{eq:def_norme_sharp_global}
\snorme[\sharp,\calT]{v}\eqq \Big(\sum_{T\in\Th}
\snorme[\sharp,T]{v}^2\Big)^{\frac12},
\end{equation} 
for all $v\in H^{1+r}(\calT)\eqq
\{v\in \Ldeux\tq v_{|T}\in H^{1+r}(T),\ \forall T\in\Th\}$.
Let $\calE_\calT^{k+1}:\Hun\to\Pkpd(\Th)$ be the global elliptic projection operator
such that for all $v\in\Hun$,
\begin{equation}
\calE_\calT^{k+1}(v)_{|T}\eqq \calE_T^{k+1}(v_{|T}), \quad \forall T\in\Th.
\end{equation}

\bLem[Bound on consistency error]  \label{lem:bnd_HHO}
Assume that the exact solution satisfies $u\in H^{1+r}(\Dom)$, $r>\frac12$.
There is $C$ such that
\begin{equation} \label{eq:bnd_HHO}
\norme[*]{\delta_{h}}
\le C |u-\calE_\calT^{k+1}(u)|_{\sharp,\calT}.
\end{equation}
\eLem

\bProof
Let $\shwh\in\shVkhz$.
Integrating by parts in every mesh cell $T\in\Th$, recalling that $f=-\Delta u$, and since $\nT\SCAL\GRAD u$ is meaningful on every face $F\in\FT$, we obtain
\begin{align*}
\ell(\wcT) &=\sum_{T\in\Th} \pss[T]{f}{\swT} = \sum_{T\in\Th} -\pss[T]{\Delta u}{\swT} \\ &= \sum_{T\in\Th} \big(\psv[T]{\GRAD u}{\GRAD \swT} - \pss[\dT]{\nT\SCAL\GRAD u}{\swT}\big) \\
&= \sum_{T\in\Th} \big(\psv[T]{\GRAD u}{\GRAD \swT} - \pss[\dT]{\nT\SCAL\GRAD u}{\swT-\swdT}\big),
\end{align*}
where we used that
$\sum_{T\in\Th} \pss[\dT]{\nT\SCAL\GRAD u}{\swdT} = 0$
since $\GRAD u$ and $\wcF$ are single-valued on the mesh
interfaces and $\wcF$ vanishes on the boundary faces. Moreover,
since $\calE_T^{k+1}=\opRec\circ\shIkT$, using the
definition of
$\opRec(\shwT)$ (with $q \eqq \calE_T^{k+1}(u)$) leads to
\begin{align*}
\ifSp \alhere \fi 
\psv[T]{\GRAD \opRec (\shIkT(u))}{\GRAD\opRec(\shwT)} 
\ifSp \else \alhere \fi
= \psv[T]{\GRAD \calE_T^{k+1}(u)}{\GRAD\opRec(\shwT)} \\
&= \psv[T]{\GRAD \calE_T^{k+1}(u)}{\GRAD \swT}
- \pss[\dT]{\nT\SCAL\GRAD \calE_T^{k+1}(u)}{\swT-\swdT}.
\end{align*}
Let us set $\eta_T \eqq u_{|T}-\calE_{T}^{k+1}(u)$.
Using the definition of $a_T$ and since
$\psv[T]{\GRAD \eta_T}{\GRAD \swT}=0$ owing to~\eqref{eq:def_ellp_a}, we have
$\langle \delta_{h},\shwh\rangle
= -\sum_{T\in\Th}(\term_{1,T}+\term_{2,T})$ with
\begin{align*}
\term_{1,T}&\eqq \pss[\dT]{\nT\SCAL\GRAD\eta_T}{\swT-\swdT}, \\
\term_{2,T}&\eqq h_{T}^{-1} \pss[\dT]{\opSt(\shIkT(u))}{\opSt(\shwT)}.
\end{align*}
The Cauchy--Schwarz
inequality and the definition of $\snorme[\shVkT]{\shwT}$ imply that
\begin{align*}
|\term_{1,T}| &\le \normev[\dT]{\GRAD\eta_T} \norme[L^2(\dT)]{\swT-\swdT}
\le h_T^{\frac12}\normev[\dT]{\GRAD \eta_T} \snorme[\shVkT]{\shwT}.
\end{align*}
Moreover, we have
\[
|\term_{2,T}| \le h_{T}^{-\frac12}\normes[\dT]{\opSt(\shIkT(u))}
h_{T}^{-\frac12}\normes[\dT]{\opSt(\shwT)}.
\]
The first factor is bounded in Lemma~\ref{lem:approx_HHO_St},
and the second one in Lemma~\ref{lem:stab_HHO}.
\ifSp Hence, \else This implies that \fi $|\term_{2,T}|\le
C\normev[T]{\GRAD\eta_T} \snorme[\shVkT]{\shwT}$.
Collecting these bounds and summing over the mesh cells proves~\eqref{eq:bnd_HHO}.
\eProof

\section{$H^1$-error estimate}
\label{sec:error_anal_HHO}

To allow for a more compact notation, we
consider the broken gradient operator $\GRAD_\calT$
and the global reconstruction operator $\opRecT:\shVkhz\to\Pkpd(\Th)$
(see~\eqref{eq:def_opRecT}). 
For nonnegative real numbers $\theta,\beta$
and a function $\phi\in H^\beta(\calT)\eqq \{\phi\in \Ldeux\tq \phi_{|T}\in H^\beta(T),
\ \forall T\in\Th\}$,
we use the shorthand notation 
\begin{equation}
|h^\theta\phi|_{H^\beta(\calT)}\eqq \Big( \sum_{T\in\Th}
h_T^{2\theta}|\phi|_{H^\beta(T)}^2\Big)^{\frac12}.
\end{equation} 
Let us introduce the discrete error
\begin{equation}
\hat e_h \eqq (e_\calT,e_\calF) \eqq \shuh-\shIkh(u)\in \shVkhz,
\end{equation}
so that $e_\calT=u_\calT-\Pi_\calT^k(u)$ and $e_\calF=u_\calF-\Pi_\calF^k(u_{|\calF})$.

\bLem[Discrete $H^1$-error estimate]  \label{lem:esterr_HHO}
Let $u\in\Hunz$ be the exact solution and let $\shuh\in \shVkhz$
be the HHO solution solving~\eqref{eq:weak_HHO}. 
Assume that $u\in H^{1+r}(\Dom)$, $r>\frac12$.
There is $C$ such that
\begin{equation}
\label{eq:bnd_heh}
\norme[\shVkhz]{\hat e_h} + \|\GRAD_\calT\opRecT(\hat e_h)\|_{\Ldeuxd} +
s_h(\hat e_h,\hat e_h)^{\frac12} \le C\snorme[\sharp,\calT]{u-\calE_\calT^{k+1}(u)}.
\end{equation}
\eLem

\bProof
We have $a_h(\hat e_h,\hat e_h)=\langle \delta_h,\hat e_h
\rangle$, where $\delta_h$ is the consistency error defined in~\eqref{eq:delta_calI_HHO}.
Summing the lower bound in Lemma~\ref{lem:stab_HHO} over all the mesh cells
yields
\[
\alpha \norme[\shVkhz]{\hat e_h}^2\le a_h(\hat e_h,\hat e_h)=\langle \delta_h,\hat e_h
\rangle \le \|\delta_h\|_*\norme[\shVkhz]{\hat e_h}.
\]
Hence, $\norme[\shVkhz]{\hat e_h} \le \frac{1}{\alpha}\|\delta_h\|_*$, and
Lemma~\ref{lem:bnd_HHO} yields $\norme[\shVkhz]{\hat e_h}
\le C\snorme[\sharp,\calT]{u-\calE_\calT^{k+1}(u)}$. Finally,
\eqref{eq:bnd_heh} follows from
$a_h(\hat e_h,\hat e_h) = \|\GRAD_\calT\opRecT(\hat e_h)\|_{\Ldeuxd}^2 +
s_h(\hat e_h,\hat e_h)$
and the above bounds on $a_h(\hat e_h,\hat e_h)$, $\|\delta_h\|_*$, and
$\norme[\shVkhz]{\hat e_h}$.
\eProof

\bTheo[$H^1$-error estimate]  \label{th:esterr_HHO}
Under the assumptions of Lemma~\ref{lem:esterr_HHO},
there is $C$ such that
\begin{equation}
\label{eq:est_err_HHO1}
\normev[\Dom]{\GRAD_\calT (u-\opRecT(\shuh))} + s_h(\shuh,\shuh)^{\frac12}
\le C\snorme[\sharp,\calT]{u-\calE_\calT^{k+1}(u)}.
\end{equation}
Moreover, if $u\in H^{t+1}(\calT)$ for some $t \in (\frac12,k+1]$, we have
\begin{equation} \label{eq:est_err_HHO2}
\normev[\Dom]{\GRAD_\calT (u-\opRecT(\shuh))} + s_h(\shuh,\shuh)^{\frac12}
\le C\snorme[H^{t+1}(\calT)]{h^tu}.
\end{equation}
This estimate is optimal when $t=k+1$ and converges at rate $\calO(h^{k+1})$.
\eTheo

\bProof
\textup{(i)} The estimate~\eqref{eq:est_err_HHO1}
follows from~\eqref{eq:bnd_heh} and the triangle inequality.
Indeed, since $\opRecT\circ \shIkh=\calE_\calT^{k+1}$, we have
\begin{align*}
&\|\GRAD_\calT(u-\opRecT(\shuh))\|_{\Ldeuxd} \le \|\GRAD_\calT(u-\calE_\calT^{k+1}(u))\|_{\Ldeuxd} +
\|\GRAD_\calT \opRecT(\hat e_h)\|_{\Ldeuxd}, \\
&s_h(\shuh,\shuh)^{\frac12}\le s_h(\shIkh(u),\shIkh(u))^{\frac12} +
s_h(\hat e_h,\hat e_h)^{\frac12},
\end{align*}
$\|\GRAD_\calT \opRecT(\hat e_h)\|_{\Ldeuxd}$ and $s_h(\hat e_h,\hat e_h)^{\frac12}$
are bounded in Lemma~\ref{lem:esterr_HHO},
and $s_h(\shIkh(u),\shIkh(u))^{\frac12}$ is bounded in Lemma~\ref{lem:approx_HHO_St}.
\\
\textup{(ii)} The estimate~\eqref{eq:est_err_HHO2}
results from~\eqref{eq:est_err_HHO1} and
the approximation properties of the local elliptic projection.
Indeed, let us set $\eta \eqq u-\calE_\calT^{k+1}(u)$.
Owing to the optimality property of the local elliptic projection in the
$H^1$-seminorm
and to the approximation property~\eqref{eq:pol_app_L2}
of $\Pi_T^{k+1}$ (with $l\eqq k+1$, $r\eqq 1+t$, $m\eqq1$, so that $r\le l+1$ since
$t\le k+1$), we have for all $T\in\Th$,
\[
\|\GRAD\eta\|_{\bL^2(T)} \le \|\GRAD(u-\Pi_T^{k+1}(u))\|_{\bL^2(T)} \le C h_T^{t}|u|_{H^{t+1}(T)}.
\]
Using the same arguments together with the triangle inequality, the approximation
property~\eqref{eq:pol_app_face}, and the
discrete trace inequality~\eqref{eq:disc_trace}, we infer that
\begin{align*}
h_T^{\frac12} \|\GRAD\eta\|_{\bL^2(\dT)} \le{}&
h_T^{\frac12} \|\GRAD(u-\Pi_T^{k+1}(u))\|_{\bL^2(\dT)}
+ h_T^{\frac12} \|\GRAD(\calE_T^{k+1}(u)-\Pi_T^{k+1}(u))\|_{\bL^2(\dT)} \\
\le{}& C\big(h_T^{t}|u|_{H^{t+1}(T)}
+ \|\GRAD(\calE_T(u)-\Pi_T^{k+1}(u))\|_{\bL^2(T)}\big) \\
\le{}& C\big(h_T^{t}|u|_{H^{t+1}(T)}
+ 2\|\GRAD(u-\Pi_T^{k+1}(u))\|_{\bL^2(T)} \big) \le C h_T^{t}|u|_{H^{t+1}(T)}.
\end{align*}
We conclude by squaring and summing over the mesh cells.
\eProof

\section{Improved $L^2$-error estimate}
\label{sec:error_anal_HHO_L2}

As is classical with elliptic problems, an error estimate with a higher-order
convergence rate can be established on the $L^2$-norm of the error. To this purpose,
one uses that there are a constant $C_{\mathrm{ell}}$
and a regularity pickup index $s\in(\frac12,1]$ such that for all $g\in\Ldeux$,
the unique function $\zeta_g\in \Hunz$ such that $a(v,\zeta_g)=(v,g)_{\Ldeux}$ for all $v\in \Hunz$ satisfies the regularity estimate 
\begin{equation} \label{eq:reg_dual}
\|\zeta_g\|_{H^{1+s}(\Dom)}\le C_{\mathrm{ell}} \ell_\Dom^{2}
\|g\|_{\Ldeux},
\end{equation} 
where the scaling factor $\ell_\Dom\eqq \diam(\Dom)$ is introduced to make the constant $C_{\mathrm{ell}}$ dimensionless (recall that $\|\SCAL\|_{H^{1+s}(\Dom)}$ and $\|\SCAL\|_{\Ldeux}$ have the same scaling and that $-\Delta \zeta_g=g$). The elliptic regularity property 
\eqref{eq:reg_dual} holds true for the
Poisson model problem posed in a polyhedron (see \cite[Chap.~4]{Gr85},
\cite[p.~158]{Dauge_1988}).

\bLem[Discrete $L^2$-error estimate]  \label{lem:L2_est_HHO}
Let $u\in\Hunz$ be the exact solution and let $\shuh\in \shVkhz$
be the HHO solution. 
Assume that $u\in H^{1+r}(\Dom)$, $r>\frac12$.
Let $s\in (\frac12,1]$ be the pickup index in
the elliptic regularity property.
Let $\delta\eqq s$ if $k=0$ and $\delta\eqq0$ if $k\ge1$.
There is $C$ such that
\begin{equation} \label{eq:est_L2_disc}
\|e_\calT\|_{\Ldeux} \le C\ell_\Dom^{1-s}h^s\Big(
|u-\calE_\calT^{k+1}(u)|_{\sharp,\calT}
+ \ell_\Dom^\delta \|h^{1-\delta}(f-\Pi_\calT^k(f))\|_{\Ldeux} \Big).
\end{equation}
\eLem

\bProof
Let $\zeta_e\in\Hunz$
be such that $a(v,\zeta_e)=(v,e_\calT)_{\Ldeux}$ for all $v\in \Hunz$. 
Since $-\Delta \zeta_e=e_\calT$, we have
\[
\|e_\calT\|_{\Ldeux}^2 = -(e_\calT,\Delta\zeta_e)_{\Ldeux}
= \sum_{T\in\Th} \big((\GRAD e_T,\GRAD \zeta_e)_{\bL^2(T)} +
(e_{\dT}-e_T,\nT\SCAL\GRAD\zeta_e)_{L^2(\dT)} \big),
\]
where we used that $\zeta_e\in H^{1+s}(\Omega)$, $s>\frac12$, and $e_F=0$
for all $F\in\calFb$
to infer that $\sum_{T\in\Th}(e_{\dT},\nT\SCAL\GRAD\zeta_e)_{L^2(\dT)}=0$.
Let us set $\xi\eqq \zeta_e-\calE_\calT^{k+1}(\zeta_e)$.
Adding and subtracting $\opRec(\shIkT(\zeta_e))=\calE_T^{k+1}(\zeta_e)$ for all
$T\in\Th$ in the above expression, using the
definition of $\opRec(\hat e_T)$ and since $(\GRAD e_T,\GRAD \xi)_{\bL^2(T)}=0$, 
we infer that
\[
\|e_\calT\|_{\Ldeux}^2 =
\sum_{T\in\Th} (e_{\dT}-e_T,\nT\SCAL\GRAD\xi)_{L^2(\dT)}+\term_1,
\]
with
\begin{align*}
\term_1&\eqq \sum_{T\in\Th} (\GRAD\opRec(\hat e_T),\GRAD\opRec(\shIkT(\zeta_e)))_{\bL^2(T)}
=-s_h(\hat e_h,\shIkh(\zeta_e))+a_h(\hat e_h,\shIkh(\zeta_e))\\
&=-s_h(\hat e_h,\shIkh(\zeta_e))+(f,\Pi_\calT^k(\zeta_e))_{\Ldeux}
-a_h(\shIkh(u),\shIkh(\zeta_e)),
\end{align*}
where we used the definition of $\hat e_h$ and the fact that $\shuh$ solves the
HHO problem. Since
$(f,\zeta_e)_{\Ldeux}=(\GRAD u,\GRAD\zeta_e)_{\Ldeuxd}$, re-arranging the terms leads
to $\|e_\calT\|_{\Ldeux}^2=\term_2+\term_3-\term_4$ with
\begin{align*}
\term_2 &\eqq \sum_{T\in\Th} 
(e_{\dT}-e_T,\nT\SCAL\GRAD\xi)_{L^2(\dT)} -s_h(\hat e_h,\shIkh(\zeta_e)),\\
\term_3 &\eqq (\GRAD u,\GRAD\zeta_e)_{\Ldeuxd}-a_h(\shIkh(u),\shIkh(\zeta_e)),\\
\term_4 &\eqq (f,\zeta_e-\Pi_\calT^k(\zeta_e))_{\Ldeux} =
(f-\Pi_\calT^k(f),\zeta_e-\Pi_\calT^k(\zeta_e))_{\Ldeux}.
\end{align*}
It remains to bound these three terms. The Cauchy--Schwarz inequality,
Lemma~\ref{lem:stab_HHO} (to bound $s_h(\hat e_h,\hat e_h)^{\frac12}$),
and Lemma~\ref{lem:approx_HHO_St} (to bound $s_h(\shIkh(\zeta_e),\shIkh(\zeta_e))^{\frac12}$) give
\[
|\term_2| \le C \norme[\shVkhz]{\hat e_h}
|\zeta_e-\calE_\calT^{k+1}(\zeta_e)|_{\sharp,\calT}.
\]
The approximation property of the elliptic projection gives
$|\zeta_e-\calE_\calT^{k+1}(\zeta_e)|_{\sharp,\calT}\le Ch^s|\zeta_e|_{H^{1+s}(\Dom)}$,
and $\ell_\Dom^{1+s}|\zeta_e|_{H^{1+s}(\Dom)}\le
\|\zeta_e\|_{H^{1+s}(\Dom)}\le C_{\mathrm{ell}} \ell_\Dom^{2}
\|e_\calT\|_{\Ldeux}$ by the elliptic regularity property. 
Using~\eqref{eq:bnd_heh} to bound $\norme[\shVkhz]{\hat e_h}$, we infer that
\[
|\term_2| \le C |u-\calE_\calT^{k+1}(u)|_{\sharp,\calT} \ell_\Dom^{1-s}h^s\|e_\calT\|_{\Ldeux}.
\]
Furthermore, using the definition of $a_h$, the identity $\opRecT\circ\shIkh=\calE_\calT^{k+1}$, and the orthogonality property of the elliptic projection yields
\begin{align*}
\term_3&=(\GRAD u,\GRAD\zeta_e)_{\Ldeuxd}-(\GRAD_\calT \calE_\calT^{k+1}(u),\GRAD_\calT\calE_\calT^{k+1}(\zeta_e))_{\Ldeuxd}- s_h(\shIkh(u),\shIkh(\zeta_e)) \\
&=(\GRAD_\calT (u-\calE_\calT^{k+1}(u)),\GRAD_\calT(\zeta_e-\calE_\calT^{k+1}(\zeta_e)))_{\Ldeuxd}
- s_h(\shIkh(u),\shIkh(\zeta_e)).
\end{align*}
Hence, $|\term_3| \le C\|\GRAD_\calT (u-\calE_\calT^{k+1}(u))\|_{\Ldeuxd}\|\GRAD_\calT (\zeta_e-\calE_\calT^{k+1}(\zeta_e))\|_{\Ldeuxd}$ by
the Cauchy--Schwarz inequality and Lemma~\ref{lem:approx_HHO_St}.
Invoking again the approximation property of the elliptic projection and
the elliptic regularity property yields
\[
|\term_3| \le C\|\GRAD_\calT (u-\calE_\calT^{k+1}(u))\|_{\Ldeuxd}
\ell_\Dom^{1-s}h^s\|e_\calT\|_{\Ldeux}.
\]
Finally, we have $\|\zeta_e-\Pi_T^k(\zeta_e)\|_{L^2(T)}\le 
Ch_T^{1+\gamma}|\zeta_e|_{H^{1+\gamma}(T)}$ with $\gamma\eqq s-\delta$ (\ie $\gamma=0$ if
$k=0$ and $\gamma=s$ if $k\ge1$). The Cauchy--Schwarz inequality implies that
$|\term_4|\le C\|h^{1+\gamma}(f-\Pi_\calT^k(f))\|_{\Ldeux}|\zeta_e|_{H^{1+\gamma}(\Omega)}$.
Invoking the elliptic regularity property yields
$|\term_4|\le C\ell_\Dom^{1-\gamma}\|h^{1+\gamma}(f-\Pi_\calT^k(f))\|_{\Ldeux}\|e_\calT\|_{\Ldeux}$, so that
\[
|\term_4| \le C\ell_\Dom^{\delta} 
\|h^{1-\delta}(f-\Pi_\calT^k(f))\|_{\Ldeux}\ell_\Dom^{1-s} h^{s}\|e_\calT\|_{\Ldeux}.
\]
Putting together the bounds on $\term_2$, $\term_3$, and $\term_4$ completes the proof.
\end{proof}

\bTheo[$L^2$-error estimate]  \label{th:L2_est_HHO}
Under the assumptions of Lemma~\ref{lem:L2_est_HHO}, there is $C$ such that
\begin{equation}
\label{eq:est_L2_cont}
\|u-\opRecT(\shuh)\|_{\Ldeux} \le
Ch|u-\calE_\calT^{k+1}(u)|_{\sharp,\calT}
+\|e_\calT\|_{\Ldeux}.
\end{equation}
Moreover, if $u\in H^{t+1}(\calT)$
and $f\in H^{\tau}(\calT)$ for some $t\in (\frac12,k+1]$, we have
\begin{equation} \label{eq:conv_L2_rate}
\|e_\calT\|_{\Ldeux} + \|u-\opRecT(\shuh)\|_{\Ldeux}
\le C\ell_\Dom^{1-s}h^s\big( \snorme[H^{t+1}(\calT)]{h^tu}
+ \ell_\Dom^\delta |h^{1-\delta+\tau}f|_{H^{\tau}(\calT)}\big),
\end{equation}
with $\tau\eqq \max(t-1+\delta,0)$
(recall that $\delta\eqq s$ if $k=0$ and $\delta\eqq0$ if $k\ge1$).
This estimate is optimal when $t=k+1$ and $s=1$ and converges at rate $\calO(h^{k+2})$.
\eTheo

\bProof
\textup{(i)} The triangle inequality and the Poincar\'e--Steklov inequality~\eqref{eq:Poincare} give
\begin{align*}
\|u-\opRecT(\shuh)\|_{\Ldeux}
&\le \|u-\calE_\calT^{k+1}(u)\|_{\Ldeux}
+ \|\opRecT(\shuh)-\calE_\calT^{k+1}(u)\|_{\Ldeux}\\
&\le Ch\|\GRAD_\calT(u-\calE_\calT^{k+1}(u))\|_{\Ldeuxd}
+ \|\opRecT(\shuh)-\calE_\calT^{k+1}(u)\|_{\Ldeux}.
\end{align*}
Moreover, we have $\|v\|_{\Ldeux}\le Ch\|\GRAD_\calT v\|_{\Ldeuxd}+\|\Pi_\calT^0(v)\|_{\Ldeux}$
for all $v\in H^1(\calT)$, owing to the triangle inequality and the Poincar\'e--Steklov inequality~\eqref{eq:Poincare}.
Applying this bound to $v\eqq \opRecT(\shuh)-\calE_\calT^{k+1}(u)=\opRecT(\hat e_h)$ and using that $\Pi_\calT^0(v)=\Pi_\calT^0(u_{\calT}-u)$, we infer that
\[
\|\opRecT(\shuh)-\calE_\calT^{k+1}(u)\|_{\Ldeux} \le
Ch\|\GRAD_\calT \opRecT(\hat e_h)\|_{\Ldeuxd}
+ \|e_\calT\|_{\Ldeux},
\]
since $\|\Pi_\calT^0(u_\calT-u)\|_{\Ldeux}\le \|\Pi_\calT^k(u_\calT-u)\|_{\Ldeux}$
and $\Pi_\calT^k(u_{\calT}-u)=e_\calT$. Finally, the estimate~\eqref{eq:est_L2_cont}
follows by combining the above inequalities, using~\eqref{eq:bnd_heh}
to bound $\|\GRAD_\calT \opRecT(\hat e_h)\|_{\Ldeuxd}$, and since
$\|\GRAD_\calT(u-\calE_\calT^{k+1}(u))\|_{\Ldeuxd}\le |u-\calE_\calT^{k+1}(u)|_{\sharp,\calT}$.
\\
\textup{(ii)} \eqref{eq:conv_L2_rate} follows from~\eqref{eq:est_L2_disc}-\eqref{eq:est_L2_cont} and the approximation properties of the elliptic projection and the $L^2$-orthogonal projection (notice that in all cases, $\tau\in [0,k+1]$).
\eProof

\bRem[Regularity assumption]
If $k\ge1$, the regularity assumption on $f$ in Theorem~\ref{th:L2_est_HHO}
is $f\in H^{t-1}(\calT)$ which is consistent with the assumption $u\in H^{t+1}(\calT)$
and the fact that $-\Delta u=f$. If $k=0$ and, say, $t=1$, the assumptions become $u\in
H^2(\calT)$ and $f\in H^s(\calT)$, so that some extra regularity on $f$ is required.
\eRem


\chapter{Some variants}
\label{chap:variants}

The goal of this chapter is to explore some variants of the HHO method
devised in Chapter~\ref{chap:diffusion} and analyzed in Chapter~\ref{chap:math}.
We first study two variants of the gradient reconstruction operator that will turn useful,
for instance, when dealing with nonlinear problems in Chapters~\ref{Chap:Elas}
and \ref{chap::plashpp}.
Then, we explore a mixed-order variant of the HHO method that is useful, for instance,
to treat domains with a curved boundary. Finally, we
bridge the HHO method to the finite element and virtual element viewpoints.

\section{Variants on gradient reconstruction}
\label{sec:variants_rec}

In this section, we discuss two variants of the gradient
reconstruction operator defined in Sect.~\ref{sec:local_rec}.
Let $k\ge0$ be the polynomial degree. 
Recall that for every mesh cell $T\in\Th$, letting
$\shVkT \eqq \Pkd(T) \times \sVkdT$, the local reconstruction
operator $\opRec : \shVkT \to \Pkpd(T)$ is defined such that
for all $\shvT\in\shVkT$,
\begin{align}
&\psv[T]{\GRAD \opRec(\shvT)}{\GRAD q} =
- \pss[T]{\svT}{\Delta q} + \pss[\dT]{\svdT}{\nT\SCAL\GRAD q},
\label{eq:def_Rec_HHO_var} \\
&\pss[T]{\opRec(\shvT)}{1} = \pss[T]{\svT}{1},\label{eq:def_Rec_HHO_mean_var}
\end{align}
where \eqref{eq:def_Rec_HHO_var} holds for all 
$q\in \Pkpd(T)^\perp\eqq \{q\in \Pkpd(T)\tq (q,1)_{L^2(T)}=0\}$.
The gradient is then reconstructed
locally as $\GRAD\opRec(\shvT) \in \GRAD \Pkpd(T)$. 

A first variant is to reconstruct the gradient in the larger space
$\vPkd(T)\eqq \Pkd(T;\Real^d)$. Notice that $\GRAD\Pkpd(T)
\subsetneq \vPkd(T)$ for all $k\ge1$, whereas $\GRAD\polP_d^1(T)
=\bpolP_d^0(T)$. Although it may be surprising at first sight
to reconstruct a gradient in a space that is not composed of curl-free fields, this choice 
is relevant in the context of nonlinear problems, as highlighted in
\cite{DiPDr:17} for Leray--Lions problems and in \cite{BoDPS:17,AbErPi:18}
for nonlinear elasticity. Indeed, looking at the consistency
proof in Lemma~\ref{lem:bnd_HHO}, one sees that one exploits locally the definition of 
the reconstructed gradient of the test function, $\GRAD\opRec(\shwT)$, 
acting against the
reconstructed gradient of some interpolate of the exact solution,
$\GRAD\opRec(\shIkT(u))$. However, in the nonlinear case, $\GRAD\opRec(\shwT)$
acts against some nonlinear transformation of $\GRAD\opRec(\shIkT(u))$,
and there is no reason that this transformation preserves curl-free fields. 
For further mathematical insight using the notion of limit-conformity, 
we refer the reader to \cite[Sect.~4.1]{DiPDr:17}. 

The devising of the gradient reconstruction operator $\opGRec:\shVkT \to \vPkd(T)$
follows the same principle as the one for $\opRec$: it is based on integration by parts.
Here, $\opGRec(\shvT)\in \vPkd(T)$ is defined such that for all $\shvT\in\shVkT$,
\begin{equation} \label{eq:def_opGRec}
(\opGRec(\shvT),\bq)_{\bL^2(T)}=-(\svT,\DIV\bq)_{L^2(T)}+(\svdT,
\nT\SCAL\bq)_{L^2(\dT)}, \quad \forall \bq\in \vPkd(T).
\end{equation}
To compute $\opGRec(\shvT)$, it suffices to invert the mass matrix associated with the scalar-valued polynomial space $\Pkd(T)$ since only the right-hand side changes when computing each Cartesian component of $\opGRec(\shvT)$.

\bLem[Gradient reconstruction] 
\textup{(i)} $\vecteur{\Pi}_{\GRAD \Pkpd}(\opGRec(\shvT)) =\GRAD\opRec(\shvT)$ for all
$\shvT\in \shVkT$, where $\vecteur{\Pi}_{\GRAD \Pkpd}$ is the
$\bL^2$-orthogonal projection onto $\GRAD \Pkpd(T)$. 
\textup{(ii)} $\opGRec(\shIkT(v))=\PikTv(\GRAD v)$ for all $v\in H^1(T)$,
where $\PikTv$ is the $\bL^2$-orthogonal projection onto $\vPkd(T)$.
\label{lem:grad_rec}
\eLem

\bproof
\textup{(i)} Let $\shvT\in \shVkT$. 
For all $q\in \Pkpd(T)^\perp$, since $\GRAD q\in \vPkd(T)$, \eqref{eq:def_opGRec} yields
\begin{align*}
(\opGRec(\shvT),\GRAD q)_{\bL^2(T)}&=-(\svT,\Delta q)_{L^2(T)}+(\svdT,\nT\SCAL\GRAD q)_{L^2(\dT)} = (\GRAD\opRec(\shvT),\GRAD q)_{\bL^2(T)},
\end{align*}
where the second equality follows from \eqref{eq:def_Rec_HHO_var}. Since $\GRAD q$
is arbitrary in $\GRAD\Pkpd(T)$ and $\GRAD\opRec(\shvT)\in\GRAD\Pkpd(T)$, this proves
that $\vecteur{\Pi}_{\GRAD \Pkpd}(\opGRec(\shvT)) =\GRAD\opRec(\shvT)$.\\
\textup{(ii)} Let $v\in H^1(T)$. Since $\shIkT(v)= (\PikTs(v),\PikdTs(v_{|\dT}))$, \eqref{eq:def_opGRec} yields for all $\bq\in \vPkd(T)$,
\begin{align*}
(\opGRec(\shIkT(v)),\bq)_{\bL^2(T)}&=-(\PikTs(v),\DIV\bq)_{L^2(T)}+(\PikdTs(v_{|\dT}),
\nT\SCAL\bq)_{L^2(\dT)} \\
&=-(v,\DIV\bq)_{L^2(T)}+(v,\nT\SCAL\bq)_{L^2(\dT)} = (\GRAD v,\bq)_{\bL^2(T)},
\end{align*}
where we used $\DIV\bq\in\Pkmd(T)\subset\Pkd(T)$, $\nT\SCAL\bq_{|\dT}\in \sVkdT$, 
and integration by parts. Since $\bq$ is arbitrary in $\vPkd(T)$, this proves 
that $\opGRec(\shIkT(v))=\PikTv(\GRAD v)$.
\eproof

The property \textup{(ii)} from Lemma~\ref{lem:grad_rec} is the counterpart of the
identity $\opRec\circ\shIkT=\calE_T^{k+1}$ (see Lemma~\ref{lem:ell_proj_HHO}).
By inspecting the proofs of Lemma~\ref{lem:bnd_HHO} and Theorem~\ref{th:esterr_HHO},
one readily sees that devising the HHO method with the local bilinear form 
\begin{equation}
a_T(\shvT,\shwT) \eqq
\psv[T]{\opGRec(\shvT)}{\opGRec(\shwT)}
+ h_T^{-1}\pss[\dT]{ \opSt(\shvT)}{\opSt(\shwT)},
\end{equation}
again leads to optimal $H^1$- and $L^2$-error estimates.

Another interesting variant on simplicial meshes
is to reconstruct the gradient in the
even larger Raviart--Thomas 
space $\RT^k_d(T)\eqq \vPkd(T)\oplus\bx\tilde{\polP}_d^k(T)$, where 
$\tilde{\polP}_d^k(T)$ is composed of the restriction to $T$ of
the homogeneous $d$-variate polynomials of degree $k$. Notice that
$\vPkd(T)\subsetneq \RT^k_d(T) \subsetneq \bpolP_d^{k+1}(T)$.
Similarly to~\eqref{eq:def_opGRec}, 
$\opGRecRT:\shVkT \to \RT_d^k(T)$ is defined such that for all $\shvT\in\shVkT$,
$\opGRecRT(\shvT)\in \RT_d^k(T)$ satisfies
\begin{equation} \label{eq:def_opGRecRT}
(\opGRecRT(\shvT),\bq)_{\bL^2(T)}=-(\svT,\DIV\bq)_{L^2(T)}+(\svdT,
\nT\SCAL\bq)_{L^2(\dT)}, \quad \forall \bq\in \RT_d^k(T).
\end{equation}
In practice, $\opGRecRT(\shvT)$ is computed by inverting the mass matrix
associated with the space $\RT_d^k(T)$ (it is not possible here to
compute the Cartesian components of $\opGRecRT(\shvT)$ separately).
Following the seminal idea from \cite{JoNeS:16} in the context of penalty-free
discontinuous Galerkin methods, the motivation for reconstructing a gradient
using Raviart--Thomas polynomials 
is that it allows one to discard the stabilization operator
in the HHO method on simplicial meshes \cite{DiPDM:18,AbErPi:18}.
Recall the $H^1$-like seminorm such that
$\snorme[\shVkT]{\shvT}^2 \eqq \normev[T]{\GRAD \svT}^2
+ h_T^{-1}\normes[\dT]{\svT-\svdT}^2$ for all $\shvT\in \shVkT$.

\bLem[Raviart--Thomas gradient reconstruction] 
\textup{(i)} $\PikTv(\opGRecRT(\shvT)) =\opGRec(\shvT)$ for all
$\shvT\in \shVkT$. 
\textup{(ii)} $\opGRecRT(\shIkT(v))=\PikTv(\GRAD v)$ for all $v\in H^1(T)$.
\textup{(iii)} Assuming that the mesh belongs to a shape-regular sequence of simplicial meshes, there is $C>0$ such that 
$\|\opGRecRT(\shvT)\|_{\bL^2(T)} \ge C\snorme[\shVkT]{\shvT}$
for all $T\in\Th$ and all $\shvT\in \shVkT$. 
\label{lem:grad_rec_RT}
\eLem

\bproof
\textup{(i)} follows from $\vPkd(T)\subset \RT^k_d(T)$, and 
\textup{(ii)} is proved by proceeding as in the proof of Lemma~\ref{lem:grad_rec}
and observing that $\DIV\bq\in\Pkd(T)$ and $\nT\SCAL\bq_{|\dT}\in \sVkdT$
for all $\bq\in \RT_d^k(T)$ (even if $T$ is not a simplex). 
Finally, on a simplex, using classical properties of
Raviart--Thomas polynomials (see, \eg \cite{BoBrF:13,ErnGu:21a}), one can
show that for all $\shvT\in \shVkT$, there is $\bq_v\in\RT_d^k(T)$ such that
\[
\PikdTs(\nT\SCAL\bq_{v|\dT})=h_T^{-1}(\svdT-v_{T|\dT}),\quad
\vecteur{\Pi}_T^{k-1}(\bq_v)=\GRAD v_T,\quad
\|\bq_v\|_{\bL^2(T)} \le C\snorme[\shVkT]{\shvT}.
\]
Using the test function $\bq_v$ in \eqref{eq:def_opGRecRT} and integrating by parts
gives
\begin{align*}
(\opGRecRT(\shvT),\bq_v)_{\bL^2(T)} &= -(\svT,\DIV\bq_v)_{L^2(T)}+(\svdT,\nT\SCAL\bq_v)_{L^2(\dT)} \\
&= (\GRAD \svT,\bq_v)_{\bL^2(T)}-(\svT-\svdT,\nT\SCAL\bq_v)_{L^2(\dT)} \\
&= \normev[T]{\GRAD \svT}^2 + h_T^{-1}\normes[\dT]{\svT-\svdT}^2 = \snorme[\shVkT]{\shvT}^2,
\end{align*}
since $\GRAD\svT \in \bpolP_d^{k-1}(T)$ and $v_{T|\dT}-\svdT\in \sVkdT$.
The Cauchy--Schwarz inequality and the above bound on $\|\bq_v\|_{\bL^2(T)}$
finally imply that
\begin{align*}
\snorme[\shVkT]{\shvT}^2= (\opGRecRT(\shvT),\bq_v)_{\bL^2(T)}
&\le \|\opGRecRT(\shvT)\|_{\bL^2(T)} \|\bq_v\|_{\bL^2(T)} 
\le C\|\opGRecRT(\shvT)\|_{\bL^2(T)}\snorme[\shVkT]{\shvT},
\end{align*}
which proves the assertion \textup{(iii)}.
\eproof

The property \textup{(iii)} from Lemma~\ref{lem:grad_rec_RT} is the cornerstone ensuring the stability of the HHO method on simplicial meshes using the unstabilized bilinear form
\begin{equation}
a_T(\shvT,\shwT) \eqq
\psv[T]{\opGRecRT(\shvT)}{\opGRecRT(\shwT)},
\end{equation}
and the property \textup{(ii)} is key to deliver optimal error estimates. Notice that
the property~\textup{(ii)} fails if the gradient is reconstructed locally in the even larger space $\bpolP_d^{k+1}(T)$ since the normal component on $\dT$ of polynomials in this space does not necessarily belong to $\sVkdT$. Notice also that the property~(iii) can be achieved on polyhedral meshes by considering Raviart--Thomas polynomials on the simplicial submesh of each mesh cell (see \cite{DiPDM:18}). Another possibility pursued in the context
of weak Galerkin methods is to reconstruct the gradient in $\bpolP_d^{k+n-1}(T)$ where
$n$ is the number of faces of $T$ \cite{YeZha:20}; 
however, the energy-error estimate only decays as $O(h^k)$.

\section{Mixed-order variant and application to curved boundaries}
\label{sec:variants}

In this section, we briefly discuss the possibility of considering cell and face
unknowns that are polynomials of different degrees. As an example of application,
we show how a mixed-order variant of the HHO method 
lends itself to the approximation of problems
posed on a domain with a curved boundary.

\subsection{Mixed-order variant with higher cell degree}
\label{sec:variants_cell}

Let $k\ge 0$ be the polynomial degree for the face unknowns. The degree of the cell
unknowns is now set to $k'\eqq k+1$, leading to a mixed-order HHO method
(a mixed-order variant with lower cell degree is briefly addressed below).
The mixed-order HHO space is then defined as follows:
\begin{equation}
\shVkph \eqq V_{\calT}^{k'} \times V_{\calF}^k,
\qquad
V_{\calT}^{k'} \eqq \bigtimes_{T\in\calT} \polP^{k'}_d(T),
\qquad
V_{\calF}^k \eqq \bigtimes_{F\in\calF} \PkF(F),
\end{equation}
and the local components of a generic member $\shvh\in \shVkph$ associated with a mesh cell
$T\in\Th$ and its faces are denoted by $\shvT\eqq(v_T,v_{\dT})\in \shVkpT\eqq
\polP^{k'}_d(T) \times \PkF(\FT)$ with $\PkF(\FT)\eqq\bigtimes_{F\in\FT}\PkF(F)$. 
The HHO reduction operators
$\shIkpT:H^1(T)\to \shVkpT$ and $\shIkph:\Hun\to \shVkph$ are defined such that
\begin{equation} \label{eqn:reduction_mixed}
\shIkpT(v)\eqq (\Pi_T^{k'}(v),\PikdTs(v_{|\dT})), \qquad
\shIkph(v)\eqq (\Pi^{k'}_\calT(v),\PikcFs(v_{|\calF})).
\end{equation}
The local reconstruction operator $\opRec^+ : \shVkpT \to \Pkpd(T)$ is defined exactly as
in~\eqref{eq:def_Rec_HHO}-\eqref{eq:def_Rec_HHO_mean}, and one readily verifies that
the identity from Lemma~\ref{lem:ell_proj_HHO} can be extended to the mixed-order case, \ie we
have $\calE_T^{k+1} = \opRec^+\circ\shIkpT$ on $H^1(T)$.

The main difference between the equal-order and mixed-order versions of the HHO method lies in the stabilization
operator. Indeed, its expression is simpler in the mixed-order case and reads for all $T\in\Th$
(compare with~\eqref{eq:def_SKk_HHO}),
\begin{equation} \label{eq:stab_HDG}
\opSt^+(\shvT) \eqq
\PikdTs(v_{T|\dT}-\svdT)=\PikdTs(v_{T|\dT})-\svdT.
\end{equation}
The local bilinear form $a_T^+:\shVkpT\times\shVkpT\to\Real$ is defined as
\begin{equation}
a_T^+(\shvT,\shwT)\eqq (\GRAD\opRec^+(\shvT),\GRAD\opRec^+(\shwT))_{\bL^2(T)}
+ h_T^{-1}(\opSt^+(\shvT),\opSt^+(\shwT))_{L^2(\dT)},
\end{equation}
and the global bilinear form $a_h^+:\shVkph\times\shVkph\to\Real$ is still assembled
by summing the local contributions cellwise.
The discrete problem takes a similar form to~\eqref{eq:weak_HHO}:
\begin{equation}
\left\{
\begin{array}{l}
\text{Find $\shuh\in \shVkphz\eqq V_{\calT}^{k'} \times V_{\calF,0}^k$ such that}\\[2pt]
a_h^+(\shuh,\shwh) = \ell(\wcT),
\quad \forall \shwh \in \shVkphz.
\end{array}
\right.
\end{equation}
The cell unknowns can be eliminated locally by static condensation
(see Sect.~\ref{sec:static_HHO}), and by proceeding as in Sect.~\ref{sec:flux_recovery},
one can recover equilibrated fluxes. Recalling
Sect.~\ref{sec:HHO_HDG}, we observe that the HDG rewriting of the above mixed-order
HHO method has been considered by Lehrenfeld and Sch\"oberl~\cite{Lehrenfeld:10,LehSc:16}
(see also \cite{Oikawa:15}) and is often called HDG+ method.

The analysis of the mixed-order HHO method is quite similar to that of the equal-order version, and we only outline the few changes in the analysis.

\bLem[Stability]
Let $\snorme[\shVkpT]{\SCAL}$ denote the extension to $\shVkpT$ of
the $H^1$-like seminorm defined in \eqref{eq:def_norme_locale_HHO}.
There are $0<\alpha\le \omega<+\infty$ such that
for all $T\in\Th$ and all $\shvT\in \shVkpT$,
\begin{align}
\alpha \snorme[\shVkpT]{\shvT}^2 &\le a_T^+(\shvT,\shvT)
\le \omega \snorme[\shVkpT]{\shvT}^2.\label{eq:stab_HHO_bis_p}
\end{align}
 \label{lem:stability_p}
\eLem

\begin{proof}
Only a few adaptations are needed from the proof of Lemma~\ref{lem:stab_HHO}.
For the lower bound, setting $r_T\eqq\opRec^+(\shvT)$, 
a slightly sharper version of \eqref{eq:estim_stab_HHO1} is
\[
\normev[T]{\GRAD \svT} \le C \big(\normev[T]{\GRAD r_T} + h_T^{-\frac12}
\normes[\dT]{\Pi^k_{\dT}(v_{T|\dT})-\svdT}\big)
\le C a_T^+(\shvT,\shvT)^{\frac12}.
\]
(Recall that the value of $C$ can change at each occurrence.)
Moreover, the triangle inequality, the $L^2$-optimality of $\Pi^k_{\dT}$, the discrete trace inequality~\eqref{eq:disc_trace}, and the Poincar\'e--Steklov inequality~\eqref{eq:Poincare} imply that
\begin{align*}
h_T^{-\frac12}\|v_T-v_{\dT}\|_{L^2(\dT)} &\le h_T^{-\frac12}
\|v_T-\PikdTs(v_{T|\dT})\|_{L^2(\dT)} + h_T^{-\frac12}\|\opSt^+(\shvT)\|_{L^2(\dT)} \\
&\le h_T^{-\frac12}
\|v_T'-\PikdTs(v'_{T|\dT})\|_{L^2(\dT)} + h_T^{-\frac12}\|\opSt^+(\shvT)\|_{L^2(\dT)}\\
&\le h_T^{-\frac12}
\|v_T'\|_{L^2(\dT)} + h_T^{-\frac12}\|\opSt^+(\shvT)\|_{L^2(\dT)} \\
&\le
C\|\GRAD v_T\|_{\bL^2(T)} + h_T^{-\frac12}\|\opSt^+(\shvT)\|_{L^2(\dT)},
\end{align*}
where $v'_T\eqq v_T-\Pi^0_T(v_T)$. Combining these estimates proves the lower bound in~\eqref{eq:stab_HHO_bis_p}.
The proof of the upper bound is similar to the one of Step (ii) in 
Lemma~\ref{lem:stab_HHO}, \ie it combines the bounds
$\normev[T]{\GRAD r_T} \le \normev[T]{\GRAD \svT}
+ C h_{T}^{-\frac12}\normes[\dT]{\svT-\svdT}$ and 
$h_T^{-\frac12}\|\opSt^+(\shvT)\|_{L^2(\dT)}\le C\|\GRAD v_T\|_{\bL^2(T)}+
h_T^{-\frac12}\|v_T-v_{\dT}\|_{L^2(\dT)}$.
\eProof

\bLem[Consistency] 
There is $C$ such that for all $T\in\Th$ and all $v\in H^1(T)$,
\begin{equation} \label{eq:approx_S_p}
h_{T}^{-\frac12}\normes[\dT]{\opSt^+(\shIkpT(v))}
\le C \normev[T]{\GRAD(v-\Pi_T^{k+1}(v))}.
\end{equation}
Moreover, defining the consistency error as 
$\langle \delta_{h}^+,\shwh\rangle
\eqq \ell(\wcT)- a_h(\shIkph(u),\shwh)$ for all $\shwh \in \shVkphz$,
with the dual norm $\|\delta_h^+\|_{*}\eqq \sup_{\shwh\in \shVkphz} \frac{|\langle \delta_{h}^+,\shwh\rangle|}{\|\shwh\|_{\shVkphz}}$ and the norm $\|\shwh\|_{\shVkphz}\eqq\big(\sum_{T\in\Th}
|\shwT|_{\shVkpT}^2\big)^{\frac12}$, letting the seminorm $|\SCAL|_{\sharp,\calT}$ be defined as in~\eqref{eq:def_norme_sharp}-\eqref{eq:def_norme_sharp_global}, and assuming
that the exact solution satisfies $u\in H^{1+r}(\Dom)$, $r>\frac12$, we have
\begin{equation} \label{eq:bnd_HHO_p}
\norme[*]{\delta_{h}^+}
\le C \big( |u-\calE_\calT^{k+1}(u)|_{\sharp,\calT} + \normev[\Dom]{\GRAD_\calT(u-\Pi_\calT^{k+1}(u))}\big).
\end{equation}
\label{lem:consistency_p}
\eLem

\begin{proof}
By definition, we have $\opSt^+(\shIkpT(v)) = \PikdTs(\Pi^{k+1}_T(v)_{|\dT})-\PikdTs(v_{|\dT})
= \PikdTs((\Pi^{k+1}_T(v)-v)_{|\dT})$. Using the $L^2$-stability of $\PikdTs$, the multiplicative trace
inequality~\eqref{eq:mtr}, and the Poincar\'e--Steklov inequality~\eqref{eq:Poincare}, we infer that
\begin{align*}
h_{T}^{-\frac12}\normes[\dT]{\opSt^+(\shIkpT(v))} &=
h_{T}^{-\frac12}\normes[\dT]{\PikdTs((\Pi^{k+1}_T(v)-v)_{|\dT})} \\
&\le h_{T}^{-\frac12}\normes[\dT]{\Pi^{k+1}_T(v)-v} \le C \normev[T]{\GRAD(v-\Pi_T^{k+1}(v))}.
\end{align*}
This proves~\eqref{eq:approx_S_p}. Finally, the proof of~\eqref{eq:bnd_HHO_p} is identical to that of
Lemma~\ref{lem:bnd_HHO} except that we now invoke~\eqref{eq:approx_S_p} instead of Lemma~\ref{lem:approx_HHO_St}.
\eProof

Using the above stability and consistency results and reasoning as in the proofs of
Lemmas~\ref{lem:esterr_HHO} and \ref{lem:L2_est_HHO} and of
Theorems~\ref{th:esterr_HHO} and \ref{th:L2_est_HHO} leads to optimally converging
$H^1$- and $L^2$-error estimates. The statements and proofs are omitted for brevity
(the $L^2$-error estimate does not require a further regularity assumption on $f$
when $k=0$ if $k'=k+1$).

\bRem[Mixed-order variant with lower cell degree] 
As observed in~\cite{CoDPE:16}, if the face polynomial degree satisfies $k\ge1$, the
cell polynomial degree can also be set to $k'\eqq k-1$.
One advantage is that there are less cell unknowns to eliminate locally by
static condensation. However, the stabilization operator 
must include a correction depending on the local
reconstruction operator as in~\eqref{eq:def_SKk_HHO} (as in the equal-order case). The stability, consistency, and convergence analysis
presented in Chapter~\ref{chap:math} 
can be adapted to the mixed-order HHO method with $k'= k-1$ as well. The only salient difference
is that the improved $L^2$-norm error estimate requires $k\ge2$ (this fact was not stated in \cite{CoDPE:16}).
Interestingly, as shown in~\cite{CoDPE:16},
the mixed-order HHO method with $k'= k-1$ can be bridged to the nonconforming virtual element method
introduced in~\cite{LipMa:14} and analyzed in~\cite{AyLiM:16}.
\eRem 

\subsection{Domains with a curved boundary}
\label{sec:curved}

The main idea to treat domains with a curved boundary is to consider the mixed-order HHO
method with a higher cell degree and to avoid placing unknowns on the 
boundary faces. Instead, all the terms involving the boundary are evaluated locally 
by means of the trace of the corresponding cell unknown. Moreover, the 
boundary condition is enforced weakly by means of a consistent penalty technique
inspired by the seminal work of Nitsche \cite{nit71}. One novelty is that here
the consistency term is directly incorporated into the reconstruction operator, thereby avoiding the need for a penalty parameter that has to be large enough. 
The main ideas behind the HHO method presented in this section were introduced
in \cite{BurEr:18,BurEr:19} with a different reconstruction operator 
and later simplified in \cite{BCDE:21}, where the 
presentation dealt with the more general case of an elliptic interface problem.

One way to mesh a domain $\Dom$ with a curved boundary consists in embedding
it into a larger polyhedral domain $\Dom'$ and considering a shape-regular sequence
of meshes of $\Dom'$. Notice that these meshes are built without bothering about the location
of $\front$ inside $\Dom'$. Then, from every mesh $\Th'$ of $\Dom'$, 
one generates a mesh $\Th$ of
$\Dom$ by dropping the cells in $\Th'$ outside $\Dom$, keeping those inside $\Dom$
(called the interior cells), and keeping only the part inside $\Dom$
of those cells that are cut by the boundary $\front$ (producing the so-called boundary
cells). With this process, the cells composing $\Th$ cover $\Dom$ exactly,
the interior cells have planar faces, whereas the 
boundary cells have one curved face lying on $\front$ and planar faces lying inside $\Dom$;
see Figure~\ref{fig:curved_mesh}.

\begin{figure}
\begin{center}
   \includegraphics[scale=0.5]{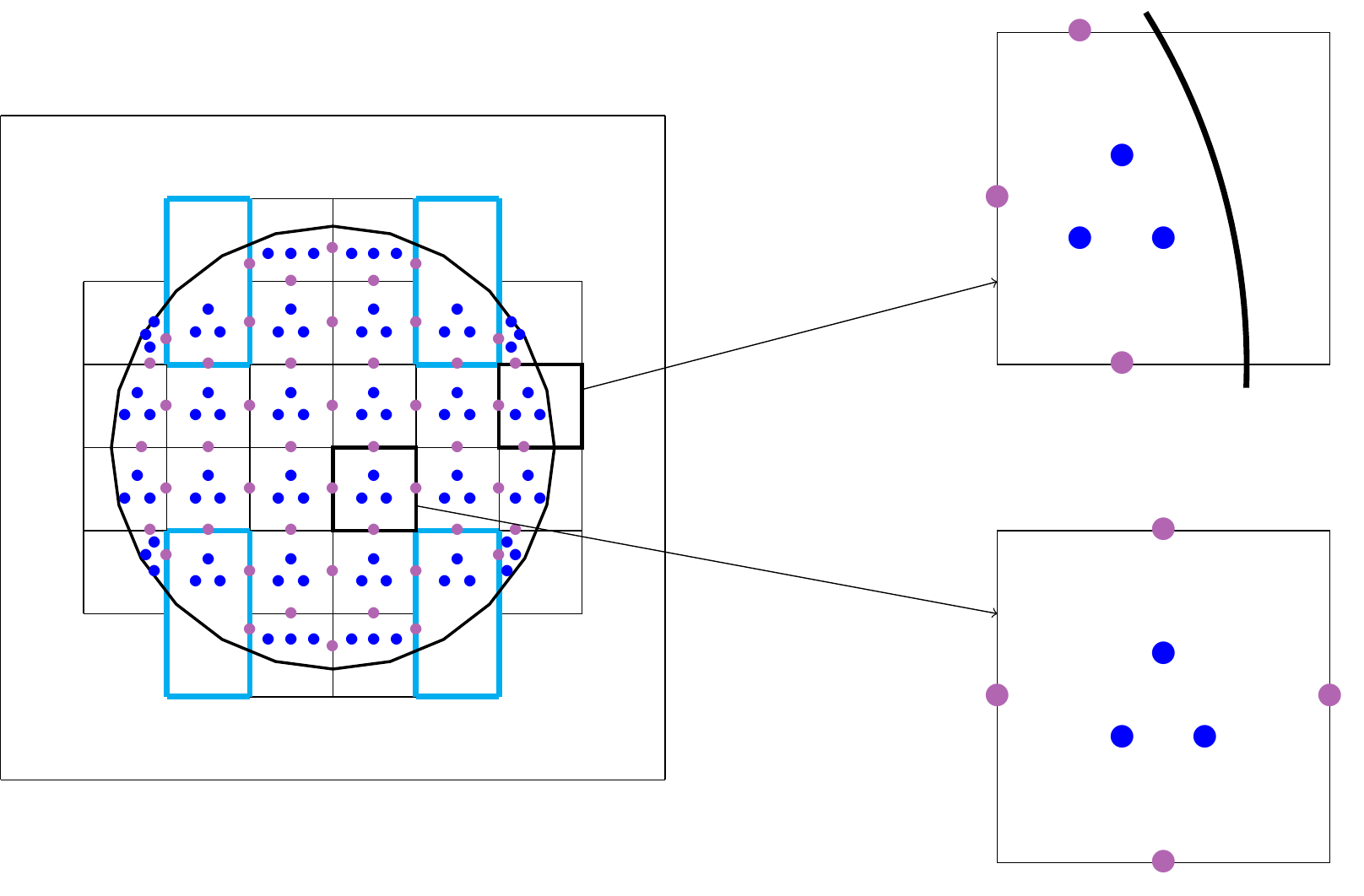}
\end{center}
\caption{Circular domain embedded into a square domain that is meshed by a quadrangular mesh; HHO unknowns associated with the interior and the boundary cells are shown for $k=0$; four agglomerated cells have been created} to avoid bad cuts.
\label{fig:curved_mesh}
\end{figure}

Some adjustments are still necessary to ensure that the basic analysis tools 
outlined in Sect.~\ref{sec:basic_tools} are available on the mesh $\Th$ that has been
constructed this way. The difficulty lies in the fact that the original mesh $\Th'$ 
was deployed without taking into account the position of the boundary $\front$ which can
therefore cut the cells in an arbitrary way. In particular, some boundary cells of $\calT$
can be very small, very flat or have an irregular shape. One possible remedy inspired from
the work \cite{JohLa:13} on discontinuous Galerkin methods
is to use a local cell-agglomeration procedure for 
badly cut cells. This procedure essentially ensures that each mesh cell, possibly after 
local agglomeration, contains a ball with diameter equivalent to its own diameter.
It is shown in \cite[Lem.~6.4]{BurEr:18}, \cite[Sect.~4.3]{BCDE:21} that this is
indeed possible if the mesh is fine enough, and \cite[Lem.~3.4]{BurEr:18},
\cite[Lem.~3.4]{BCDE:21} establish that the discrete inverse 
inequalities~\eqref{eq:disc_inv} and~\eqref{eq:disc_trace} then hold true,
together with a Poincar\'e--Steklov inequality on discrete functions. 
Moreover, it is shown in \cite[Lem.~6.1]{BurEr:18} that if the mesh size is small enough
with respect to the curvature of the boundary, every boundary cell $T\in\Th$ can be embedded
into a ball $T^\dagger$ with equivalent diameter 
such that the following multiplicative trace inequality holds
true: There is $C$ such that for all $T\in\Th$ and all $v\in H^1(T^\dagger)$ (setting $T^\dagger\eqq T$
for interior cells),
\begin{equation}\label{eq:mtr_curved}
\|v\|_{L^2(\dT)} \le C\big( h_T^{-\frac12}\|v\|_{L^2(T^\dagger)} + \|v\|_{L^2(T^\dagger)}^{\frac12}
\|\GRAD v\|_{\bL^2(T^\dagger)}^{\frac12}\big).
\end{equation}
Notice that part of the ball $T^\dagger$ may lie outside $\Omega$. 
Assuming that the exact solution is in $H^{t+1}(\Dom)$ with $t\in[1,k+1]$,
polynomial approximation
is realized by considering the $L^2$-orthogonal projection on $T^\dagger$ composed 
with the stable Calder\'on--Stein extension operator $E:H^{t+1}(\Dom)\to H^{t+1}(\Real^d)$. 
Thus we set
\begin{equation} \label{eq:def_JkpT}
J_T^{k+1}(v) \eqq \big( \Pi_{T^\dagger}^{k+1}(E(v)_{|T^\dagger}) \big)_{|T} \in \Pkpd(T),
\end{equation}
for all $T\in\Th$ and all $v\in H^{t+1}(\Dom)$. 
Notice that $J_T^{k+1}(v)=\Pi_T^{k+1}(v)$ if $T$ is an interior cell.
Reasoning as in
\cite[Lem.~5.6]{BurEr:18} (\ie using the approximation properties 
of $\Pi_{T^\dagger}^{k+1}$ and the multiplicative trace inequality~\eqref{eq:mtr_curved}), 
one can show that there is $C$ such that for all 
$T\in\Th$ and all $v\in H^{t+1}(\Dom)$, setting $\eta\eqq v-J_T^{k+1}(v)$,
\begin{equation}
\label{eq:pol_app_curved}
\|\eta\|_{L^2(T)}+h_T^{\frac12}\|\eta\|_{L^2(\dT)}+h_T\|\GRAD\eta\|_{\bL^2(T)}
+h_T^{\frac32}\|\GRAD\eta\|_{\bL^2(\dT)} \le Ch_T^{t+1}|E(v)|_{H^{t+1}(T^\dagger)}.
\end{equation}

With these tools in hand, we can devise the mixed-order HHO method for domains
with curved boundary.
For all $T\in\Th$, we consider the partition $\dT=\dTi\cup\dTb$ with
$\dTi\eqq\dT\cap\Dom$, $\dTb\eqq\dT\cap\front$, and the partition
$\FT=\FTi\cup\FTb$ with $\FTi\eqq\FT\cap\calFi$, $\FTb\eqq\FT\cap\calFb$ (the sets $\calFi$, $\calFb$ refer to the faces of $\Th$ and recall that $T\subset\Dom$).
Referring to Figure~\ref{fig:curved_mesh},
the mixed-order HHO space is redefined as follows (we keep the same notation for simplicity):
\begin{equation}
\shVkph \eqq V_{\calT}^{k'} \times V_{\calFi}^k,
\qquad
V_{\calT}^{k'} \eqq \bigtimes_{T\in\calT} \polP^{k'}_d(T),
\qquad
V_{\calF}^k \eqq \bigtimes_{F\in\calFi} \PkF(F),
\end{equation}
\ie $\calFi$ is now used in place of $\calF$ (no unknowns are attached to the mesh boundary faces in $\calFb$). 
The local components of $\shvh\in \shVkph$ associated with a mesh cell
$T\in\Th$ and its interior faces are denoted by $\shvT\eqq(v_T,v_{\dT})\in \shVkpT\eqq
\polP^{k'}_d(T) \times \PkF(\FTi)$ with $\PkF(\FTi)\eqq\bigtimes_{F\in\FTi}\PkF(F)$. 
Adapting~\eqref{eq:def_Rec_HHO}, the reconstruction operator 
$\opRec^+:\shVkpT\to \Pkpd(T)$ is such that for all $\shvT\in \shVkpT$,
\begin{align}
&\psv[T]{\GRAD \opRec^+(\shvT)}{\GRAD q} =
- \pss[T]{\svT}{\Delta q} + \pss[\dTi]{v_{\dT}}{\nT\SCAL\GRAD q},
\label{eq:def_Rec_HHO_curved} \\
&\pss[T]{\opRec^+(\shvT)}{1} = \pss[T]{\svT}{1},\label{eq:def_Rec_HHO_mean_curved}
\end{align}
where \eqref{eq:def_Rec_HHO_curved} holds for all $q\in \Pkpd(T)^\perp$.
Notice that~\eqref{eq:def_Rec_HHO_curved} is equivalent to
\begin{align} 
\psv[T]{\GRAD \opRec^+(\shvT)}{\GRAD q} = {}&(\GRAD v_T,\GRAD q)_{\bL^2(T)} 
- \pss[\dTi]{\svT-v_{\dT}}{\nT\SCAL\GRAD q} \nonumber \\
& - \pss[\dTb]{\svT}{\nT\SCAL\GRAD q}.\label{eq:def_Rec_HHO_curved_ipp}
\end{align}
Furthermore, the stabilization operator $\opSti^+:\shVkpT\to \PkF(\FTi)$ is such that
\begin{equation} \label{eq:stab_HDG_curved}
\opSti^+(\shvT) \eqq
\PikdTis(v_{T|\dTi}-\svdT)=\PikdTis(v_{T|\dTi})-\svdT.
\end{equation}
The local bilinear form $a_T^+:\shVkpT\times\shVkpT\to\Real$ is defined as
\begin{align}
a_T^+(\shvT,\shwT)\eqq {}& (\GRAD\opRec^+(\shvT),\GRAD\opRec^+(\shwT))_{\bL^2(T)}
+ h_T^{-1}(\opSti^+(\shvT),\opSti^+(\shwT))_{L^2(\dTi)} \nonumber \\
&+ h_T^{-1}(v_T,w_T)_{L^2(\dTb)},
\end{align}
where the last term results from the Nitsche's boundary penalty technique.
The global bilinear form $a_h^+:\shVkph\times\shVkph\to\Real$ is still assembled
by summing the local contributions cellwise,
and the discrete problem seeks $\shuh\in \shVkph$ such that
\begin{equation}\label{eq:disc_pb_curved}
a_h^+(\shuh,\shwh) = \ell(\wcT),
\quad \forall \shwh \in \shVkph.
\end{equation}
The cell unknowns can still be eliminated locally by means of static condensation,
and one can again recover equilibrated numerical fluxes. If the Dirichlet condition
is non-homogeneous, the right-hand side~\eqref{eq:disc_pb_curved} has to be modified
to preserve consistency (see \cite{BurEr:19,BCDE:21}).

The error analysis hinges, as usual, on stability and consistency properties.
Adapting~\eqref{eq:def_norme_locale_HHO}, let us equip $\shVkpT$ with
the $H^1$-like seminorm 
\begin{equation}
\snorme[\shVkpT]{\shvT}^2 \eqq \normev[T]{\GRAD \svT}^2
+ h_T^{-1}\normes[\dTi]{\svT-\svdT}^2 + h_T^{-1}
\normes[\dTb]{\svT}^2, 
\end{equation}
for all $\shvT\in \shVkpT$.
Proceeding as in the proof of Lemma~\ref{lem:stability_p}, one readily establishes the
following stability result.

\bLem[Stability]
There are $0<\alpha\le \omega<+\infty$ such that
$\alpha \snorme[\shVkpT]{\shvT}^2 \le a_T^+(\shvT,\shvT)
\le \omega \snorme[\shVkpT]{\shvT}^2$ for all $T\in\Th$ and all $\shvT\in \shVkpT$.
\eLem

To measure the consistency error, we need to adapt the HHO reduction operator. 
We now define $\shIkph(u)\in \shVkph$ such that its local components associated
with a mesh cell $T\in\Th$ are
\begin{equation} \label{eq:def_shIkpT}
\shIkpT(u)\eqq (J_T^{k+1}(u),\PikdTis(u_{|\dTi})) \in \shVkpT.
\end{equation} 
Notice that the cell component of $\shIkpT(u)$ is now defined using the operator $J_T^{k+1}$
(see~\eqref{eq:def_JkpT}), whereas the face component (which is now restricted to
interior faces) is still defined using an $L^2$-orthogonal projection on each face.
The use of $J_T^{k+1}$ is motivated by the approximation property~\eqref{eq:pol_app_curved},
whereas the use of $\PikdTis$ is instrumental in the following proofs. 
Before bounding the consistency error, we need to study the approximation properties
of the interpolation operator $\calJ_T^{k+1} \eqq \opRec^+\circ \shIkpT$.
This operator no longer coincides with the elliptic projection on $T$ 
because the boundary term in~\eqref{eq:def_Rec_HHO_curved} is integrated only over
$\dTi$ and because the operator $J_T^{k+1}$ differs from $\Pi_T^{k+1}$. Nonetheless,
the operator $\calJ_T^{k+1}$ still enjoys an optimal approximation property.

\bLem[Approximation property for $\calJ_T^{k+1} \eqq\opRec^+\circ \shIkpT$] 
Assume $u\in H^{t+1}(\Omega)$ with $t\in[1,k+1]$.
There is $C$ such that for all $T\in\Th$,
\begin{equation}
\snorme[\sharp,T]{u-\calJ_T^{k+1}(u)} \le Ch_T^t|E(u)|_{H^{t+1}(T^\dagger)},
\end{equation}
recalling that $\snorme[\sharp,T]{v} \eqq\normev[T]{\GRAD v}
+ h_T^{\frac12}\normev[\dT]{\GRAD v}$.
\label{lem:approx_curved}
\eLem

\bproof
Owing to the triangle inequality and the approximation property~\eqref{eq:pol_app_curved},
it suffices to bound $\snorme[\sharp,T]{q}$ with $q\eqq J_T^{k+1}(u)
-\opRec^+(\shIkpT(u))\in \Pkpd(T)$. Owing to~\eqref{eq:def_Rec_HHO_curved_ipp} and
the definition~\eqref{eq:def_shIkpT} of $\shIkpT$, we infer that
\begin{align*}
\|\GRAD q\|_{\bL^2(T)}^2 &= (\GRAD(J_T^{k+1}(u)
-\opRec^+(\shIkpT(u))),\GRAD q)_{\bL^2(T)} \\
&= -(\PikdTis(u)-J_T^{k+1}(u),\nT\SCAL\GRAD q)_{L^2(\dTi)} 
+ (J_T^{k+1}(u),\nT\SCAL\GRAD q)_{L^2(\dTb)} \\
&= -(u-J_T^{k+1}(u),\nT\SCAL\GRAD q)_{L^2(\dTi)} 
- (u-J_T^{k+1}(u),\nT\SCAL\GRAD q)_{L^2(\dTb)},
\end{align*}
where we used the $L^2$-orthogonality property of $\PikdTis$ and the fact that
the exact solution vanishes on $\front$. Invoking the Cauchy--Schwarz inequality,
the discrete trace inequality~\eqref{eq:disc_trace}, and the 
approximation property~\eqref{eq:pol_app_curved} shows that $\|\GRAD q\|_{\bL^2(T)}
\le Ch_T^t|E(u)|_{H^{t+1}(T^\dagger)}$. Finally, we bound $h_T^{\frac12}\|\GRAD q\|_{\bL^2(\dT)}$
by using the discrete trace inequality~\eqref{eq:disc_trace}.
\eproof 

As above, we define the consistency error $\delta_{h}^+ \in (\shVkph)'$ as the linear form
such that 
$\langle \delta_{h}^+,\shwh\rangle
\eqq \ell(\wcT)- a_h^+(\shIkph(u),\shwh)$ for all $\shwh\in\shVkph$.
We set $\norme[*]{\delta_{h}^+} \eqq \sup_{\shwh\in\shVkph}
\frac{|\langle \delta_{h}^+,\shwh\rangle|}{\norme[\shVkph]{\shwh}}$ with
$\norme[\shVkph]{\shwh}\eqq \big(\sum_{T\in\Th} \snorme[\shVkpT]{\shwT}^2\big)^{\frac12}$.

\bLem[Consistency]
Assume $u\in H^{t+1}(\Omega)$ with $t\in[1,k+1]$.
There is $C$ such that 
\begin{equation}
\|\delta_h^+\|_{*} \le C\Big( \sum_{T\in\Th} \big(\snorme[\sharp,T]{u-\calJ_T^{k+1}(u)}^2
+ h_T^{-1} \|u-J_T^{k+1}(u)\|_{L^2(\dT)}^2\big)\Big)^{\frac12}.
\end{equation}
\eLem

\bproof
We have $\langle \delta_{h}^+,\shwh\rangle=\term_1-\term_2$ with
\begin{align*}
\term_1&\eqq \sum_{T\in\Th} \big( (f,w_T)_{L^2(T)} - (\GRAD \calJ_T^{k+1}(u),\GRAD\opRec^+(\shwT))_{\bL^2(T)}\big),\\
\term_2&\eqq \sum_{T\in\Th} \big(h_T^{-1}(\opSti^+(\shIkpT(u)),\opSti^+(\shwT))_{L^2(\dTi)}
+  h_T^{-1}(J_T^{k+1}(u),w_T)_{L^2(\dTb)}\big).
\end{align*}
Since $f=-\Delta u$, integration by parts and the definition~\eqref{eq:def_Rec_HHO_curved_ipp} of $\opRec^+(\shwT)$ give
\begin{align*}
\term_1 ={}& \sum_{T\in\Th} \big( (\GRAD u,\GRAD w_T)_{\bL^2(T)} - (\nT\SCAL\GRAD u,w_T-w_{\dT})_{L^2(\dTi)} - (\nT\SCAL\GRAD u,w_T)_{L^2(\dTb)}\big) \\
&- \sum_{T\in\Th} \big( (\GRAD \calJ_T^{k+1}(u),\GRAD w_T)_{\bL^2(T)} 
- ( \nT\SCAL\GRAD \calJ_T^{k+1}(u),w_T-w_{\dT})_{L^2(\dTi)} 
\ifSp \cuthere \fi
- ( \nT\SCAL\GRAD \calJ_T^{k+1}(u),w_T)_{L^2(\dTb)} \big) \\
={}& \sum_{T\in\Th} \big( (\GRAD \xi,\GRAD w_T)_{\bL^2(T)} - (\nT\SCAL\GRAD \xi,w_T-w_{\dT})_{L^2(\dTi)} - (\nT\SCAL\GRAD \xi,w_T)_{L^2(\dTb)}\big),
\end{align*}
with $\xi\eqq u-\calJ_T^{k+1}(u)$, 
where we used that $\sum_{T\in\Th}(\nT\SCAL\GRAD u,w_{\dT})_{L^2(\dTi)}=0$. The term $\term_1$
is now bounded by using the Cauchy--Schwarz inequality.
To bound $\term_2$, we proceed as in the proof of~\eqref{eq:approx_S_p} for the
interior faces
(here, we use again that the face component of $\shIkpT(u)$ is defined using 
$\PikdTis$), and we use $u_{|\dTb}=0$ and the 
Cauchy--Schwarz inequality for the boundary faces.
\eproof

Using the above stability, approximation, 
and consistency results and reasoning as in the proofs of
Lemma~\ref{lem:esterr_HHO} and
Theorem~\ref{th:esterr_HHO} leads to optimally converging
$H^1$-error estimates; see \cite[Thm.~3.10]{BCDE:21}.

\section{Finite element and virtual element viewpoints}
\label{sec:FEM_HHO}

Our goal here is to bridge the HHO method with the finite element
and virtual element viewpoints. For simplicity, we focus on the equal-order HHO method.
Recall that a finite element is defined on a mesh cell $T\in\Th$ as a triple
$(T,P_T,\Sigma_T)$, where $P_T$ is a finite-dimensional space composed of
functions defined on $T$ and $\Sigma_T$ are the degrees of freedom, \ie
a collection of linear forms
on $P_T$ forming a basis of $\calL(P_T;\Real)$. The material of this 
section originates from ideas in \cite{CoDPE:16,CiErL:19}; see
also \cite{Lemaire:21}.

Let $k\ge0$ be the polynomial degree.
Consider the finite-dimensional functional space
\begin{equation} \label{eq:V_virtual_HHO}
\calV^k_T \eqq \big\{v\in H^1(T) \tq \LAP v \in \Pkd(T), \
\nT\SCAL \GRAD v_{|\dT} \in \sVkdT\big\}.
\end{equation}
We observe that $\Pkpd(T)\subset \calV^k_T$,
but there are other functions in $\calV^k_T$, and these functions 
are in general not accessible to direct computation. For this reason, the members
of $\calV^k_T$ are called virtual functions.

\bLem[$\calV^k_T \leftrightarrow \shVkT$]
The linear spaces $\calV^k_T$ and $\shVkT$ are isomorphic. Consequently,
$\dim(\calV^k_T)=\dim(\shVkT)={k+d\choose d}+{k+d-1\choose d-1}\#\FT$.
 \label{lem:iso_virtual}
\eLem

\bProof
Let us set
\begin{align*}
(\shVkT)^\perp&\eqq\{\shvT\in \shVkT\tq (v_T,1)_{L^2(T)}+(v_{\dT},1)_{L^2(\dT)}=0\},\\
(\calV^k_T)^\perp&\eqq \{v\in\calV^k_T\tq (v,1)_{L^2(T)}=0\}.
\end{align*}
To prove the assertion, it suffices to build
an isomorphism $\Phi_T:(\shVkT)^\perp\to (\calV^k_T)^\perp$. For all $\shvT\in (\shVkT)^\perp$,
$\varphi\eqq \Phi_T(\shvT)$ is the unique function in $(\calV^k_T)^\perp$ such that
$-\Delta \varphi = v_T$ in $T$ and $\nT\SCAL\GRAD\varphi=v_{\dT}$ on $\dT$.
This Neumann problem is well-posed since
$(v_T,1)_{L^2(T)}+(v_{\dT},1)_{L^2(\dT)}=0$ and $(\varphi,1)_{L^2(T)}=0$. 
This directly implies that the
map $\Phi_T$ is bijective.
\eProof

We define the virtual reconstruction operator
$\opRecv : \shVkT \to \calV^k_T$ such that for all
$\shvT\in\shVkT$, the function
$\opRecv(\shvT) \in \calV^k_T$ is uniquely defined by the following equations:
\begin{align}
&\psv[T]{\GRAD \opRecv(\shvT)}{\GRAD w} =
- \pss[T]{\svT}{\Delta w} + \pss[\dT]{\svdT}{\nT\SCAL\GRAD w},\label{eq:def_Recv_HHO} \\
&\pss[T]{\opRecv(\shvT)}{1} = \pss[T]{\svT}{1},\label{eq:def_Recv_HHO_mean}
\end{align}
where \eqref{eq:def_Recv_HHO} holds for all $w\in (\calV^k_T)^\perp$.
Notice that $\opRecv(\shvT)$ is well-defined, but it is not explicitly computable
(it can be approximated to a desired accuracy by using, say, a
finite element method on a subgrid of $T$).
Let $\shIvkT:\calV^k_T\to \shVkT$ denote the restriction of the reduction operator
$\shIkT$ to $\calV^k_T$. We slightly abuse the terminology by calling the operator
$\shIvkT$ the degrees of freedom on $\calV^k_T$ (one could define more rigorously the
degrees of freedom by choosing bases of $\Pkd(T)$ and $\PkF(F)$ for all $F\in\FT$).
Let $I_V$ denote the identity operator on a generic space $V$.

\bLem[Finite element] 
We have $\shIvkT\circ \opRecv=I_{\shVkT}$ and $\opRecv\circ \shIvkT=I_{\calV^k_T}$.
Consequently, the triple $(T,\calV^k_T,\shIvkT)$ is a finite element with interpolation
operator $\opRecv\circ \shIkT$.
\label{lem:virtual_HHO}
\eLem

\bProof
We only need to prove the two identities regarding $\shIvkT$ and $\opRecv$.
(Actually, proving just one identity is sufficient since 
$\dim(\calV^k_T)=\dim(\shVkT)$, but we provide two proofs for completeness.)
\\
\textup{(i)} Let $\shvT\in \shVkT$. To prove that $\shIvkT\circ \opRecv=I_{\shVkT}$,
we need to show that
\[\Theta\eqq (\opRecv(\shvT)-v_T,q)_{L^2(T)}+(\opRecv(\shvT)-v_{\dT},r)_{L^2(\dT)}=0,
\quad \forall (q,r)\in \shVkT.
\]
Let us write $(q,r)=(q',r)+(c,0)$ with $c\eqq\frac{1}{\mes{T}}\big((q,1)_{L^2(T)}
+(r,1)_{L^2(\dT)}\big)$ so that $(q',r)\in (\shVkT)^\perp$. Using the isomorphism
from Lemma~\ref{lem:iso_virtual}, we set $\psi\eqq \Phi_T(q',r)\in (\calV^k_T)^\perp$
and observe that
\begin{align*}
\Theta &=(\opRecv(\shvT)-v_T,q')_{L^2(T)}+(\opRecv(\shvT)-v_{\dT},r)_{L^2(\dT)}\\
&=-(\opRecv(\shvT)-v_T,\Delta \psi)_{L^2(T)}+(\opRecv(\shvT)-v_{\dT},\nT\SCAL\GRAD\psi)_{L^2(\dT)}\\
&=\psv[T]{\GRAD \opRecv(\shvT)}{\GRAD \psi} +
 \pss[T]{\svT}{\Delta \psi} - \pss[\dT]{\svdT}{\nT\SCAL\GRAD \psi} = 0,
\end{align*}
where we used~\eqref{eq:def_Recv_HHO_mean} in the first line, the definition of $\Phi_T$ on
the second line, and integration by parts and~\eqref{eq:def_Recv_HHO} on the third line.
\\
\textup{(ii)} Let $v\in \calV^k_T$. The definition~\eqref{eq:def_Recv_HHO} of $\opRecv$ implies that for all $w\in (\calV^k_T)^\perp$, we have
\begin{align*}
(\GRAD \opRecv(\shIvkT(v)),\GRAD w)_{\bL^2(T)} &= -(\Pi_T^k(v),\Delta w)_{L^2(T)}+(\Pi_{\dT}^k(v),\nT\SCAL\GRAD w)_{L^2(\dT)} \\
&=-(v,\Delta w)_{L^2(T)}+(v,\nT\SCAL\GRAD w)_{L^2(\dT)} = (\GRAD v,\GRAD w)_{\bL^2(T)},
\end{align*}
since $\Delta w\in \Pkd(T)$ and $\nT\SCAL\GRAD w_{|\dT}\in \sVkdT$. 
Since $\opRecv(\shIvkT(v))-v\in (\calV^k_T)^\perp$ owing to \eqref{eq:def_Recv_HHO_mean},
and $w$ is arbitrary in $(\calV^k_T)^\perp$, $\opRecv(\shIvkT(v))=v$.
\eProof

\bRem[Right inverse]
Lemma~\ref{lem:virtual_HHO} implies that $\opRecv:\shVkT \to \calV^k_T$ 
is a right inverse of $\shIkT$. Right inverses with other codomains can be 
devised. For instance, the right inverse devised in \cite{ErnZa:20} using
bubble functions maps onto $\mathbb{P}^{k+d+1}_d(T)$ on simplicial meshes
and allows one to build a globally $H^1$-conforming function. The construction
can be extended to general meshes.
\eRem

Let us define the high-order Crouzeix--Raviart-type finite element space
\begin{equation} \label{eq:def_CR_global}
\calV_{h}^k \eqq  \{ v_h\in \Ldeux \tq v_{h|T}\in \calV^k_T, \ \forall T\in\Th,
\ \Pi_F^k(\jump{v_h})=0,\ \forall F\in\calFi\},
\end{equation}
where $\jump{v_h}$ denotes the jump of $v_h$ across the mesh interface $F\in\calFi$,
and let us set $\calV_{h,0}^k\eqq \{v_h\in\calV_h^k\tq \Pi_F^k(v_h)=0,\ \forall F\in\calFb\}$.
The global degrees of freedom of a function $v_h\in \calV_{h,0}^k$ are 
$\shIvkh(v_h)\eqq (\PikcTs(v_h),\PikcFs(v_{h|\calF}))\in \shVkhz$, which is 
meaningful owing to jump condition in~\eqref{eq:def_CR_global}.
A natural way to use the finite element identified in Lemma~\ref{lem:virtual_HHO}
to approximate the model problem~\eqref{eq:weak_diff} is to seek $u_h\in\calV_{h,0}^k$
such that
\begin{equation} \label{eq:weak_virtual}
a_h^{\textsc{v}}(u_h,w_h) = \ell(\Pi_\calT^k(w_h)), \quad \forall w_h\in \calV_{h,0}^k,
\end{equation}
with $a_h^{\textsc{v}}:\calV_{h,0}^k\times \calV_{h,0}^k\to\Real$
such that $a_h^{\textsc{v}}(v_h,w_h)\eqq (\GRAD_\calT v_h,\GRAD_\calT w_h)_{\Ldeuxd}$
(notice that $v_h\mapsto\|\GRAD_\calT v_h\|_{\Ldeuxd}$ defines a norm on $\calV_{h,0}^k$).
The use of the $L^2$-orthogonal projection to evaluate the right-hand side
of~\eqref{eq:weak_virtual} is to stick to the multiscale HHO method proposed in
\cite{CiErL:19}. Using the above reduction and reconstruction operators,
an equivalent reformulation of~\eqref{eq:weak_virtual} is to seek $\shuh\in \shVkhz$
such that $\hat a_h^{\textsc{v}}(\shuh,\shwh) = \ell(\wcT)$ for all $\shwh\in \shVkhz$,
with
$\hat a_h^{\textsc{v}}(\shvh,\shwh)\eqq (\GRAD_\calT\opRecvT(\shvh),\GRAD_\calT\opRecvT(\shwh))_{\Ldeuxd}$ and $\opRecvT(\shvh)_{|T}\eqq \opRecv(\shvT)$ for all $T\in\Th$.

The discrete problem~\eqref{eq:weak_virtual} is
not easily tractable since the computation
of the basis functions in $\calV_{h,0}^k$
is possible only by using some subgrid discretization method
in each mesh cell. This approach is reasonable when dealing with a diffusion
problem characterized by subgrid scales that are not captured by the mesh $\calT$.
Instead, in the absence of multiscale features, it is more efficient to use the
original HHO method presented in Chapter~\ref{chap:diffusion}. To bridge the two methods,
we first notice that an equivalent formulation of the original HHO method is
to seek
$u_h\in\calV_{h,0}^k$
such that $a_h^{\textsc{hho}}(u_h,w_h) = \ell(\Pi_\calT^k(w_h))$ 
for all $w_h\in \calV_{h,0}^k$, 
with $a_h^{\textsc{hho}}:\calV_{h,0}^k\times \calV_{h,0}^k\to\Real$
such that
\begin{equation}
a_h^{\textsc{hho}}(v_h,w_h) \eqq
(\GRAD_\calT \calE_\calT^{k+1}(v_h),\GRAD_\calT \calE_\calT^{k+1}(w_h))_{\Ldeuxd}
+ s_h(\shIvkh(v_h),\shIvkh(w_h)).
\end{equation}
The quantity $a_h^{\textsc{hho}}(v_h,w_h)$ is computable even if one only
knows the global degrees of freedom $\shIvkh(v_h)$ and $\shIvkh(w_h)$ of 
$v_h$ and $w_h$,
without the need to explicitly knowing these functions (notice that $\calE_T^{k+1}=\opRec\circ\shIvkT$ for all $T\in\Th$, where $\opRec$ is the computable reconstruction operator defined in~\eqref{eq:def_Rec_HHO}-\eqref{eq:def_Rec_HHO_mean}).
The role of the stabilization in the original HHO method can then be understood
as a computable way of ensuring that $a_h^{\textsc{hho}}$ remains $H^1$-coercive on
$\calV_{h,0}^k$, in the spirit of the seminal ideas
developed for the virtual element method (see, \eg \cite{BBCMMR:13}).

\bLem[Coercivity on $\calV_{h,0}^k$]
There is $\eta>0$ such that $a_h^{\textsc{hho}}(v_h,v_h)\ge \eta \|\GRAD_\calT v_h\|_{\Ldeuxd}^2$
for all $v_h\in \calV_{h,0}^k$.
\eLem

\bProof
\textup{(i)} Let us first prove the following inverse inequalities: There is $C$ such that
for all $T\in\Th$ and all $v\in\calV^k_T$,
\begin{equation}\label{eq:disc_trace_virtual}
\|\nT\SCAL\GRAD v\|_{L^2(\dT)} + h_T^{\frac12}\|\Delta v\|_{L^2(T)}
\le Ch_T^{-\frac12}\|\GRAD v\|_{\bL^2(T)}.
\end{equation}
Let $F\in\FT$ and set $q\eqq (\nT\SCAL\GRAD v)_{|F} \in \PkF(F)$.
For simplicity, we assume that $T$ is a simplex (otherwise,
for every simplicial subface of $F$, one carves a simplex inside $T$ of
diameter uniformly equivalent to $h_T$). It results from \cite[Lem.~A.3]{ErnVo:20}
that there is $C_{\textsc{div}}$ such that for all $q\in \PkF(F)$, there is
$\btheta_q\in \RT_d^k(T)$ (the Raviart--Thomas finite element space of order $k$ in $T$) satisfying $\nT\SCAL\btheta_q=q$ on $F$, $\DIV\btheta_q=\Delta v\in \Pkd(T)$, and 
\[
\|\btheta_q\|_{\bL^2(T)} \le C_{\textsc{div}} \min_{\ba\in \bV_q} \|\ba\|_{\bL^2(T)},
\]
with $\bV_q\eqq\{\ba\in \bH(\text{\rm div};T)\tq (\nT\SCAL\ba)_{|F}=q, \, \DIV\ba=\Delta v\}$.
Since $\GRAD v\in \bV_q$, invoking the discrete trace inequality~\eqref{eq:disc_trace} for the polynomial $\btheta_q$, we infer that
\begin{equation}
\|\nT\SCAL\GRAD v\|_{L^2(F)} = \|\nT\SCAL\btheta_q\|_{L^2(F)}
\le Ch_T^{-\frac12}  \|\btheta_q\|_{\bL^2(T)} \le CC_{\textsc{div}}h_T^{-\frac12}\|\GRAD v\|_{\bL^2(T)}. \label{eq:trace_virtual}
\end{equation}
Moreover, setting $r\eqq \Delta v\in \Pkd(T)$, an integration by parts gives
$\|\Delta v\|_{L^2(T)}^2 = -(\GRAD v,\GRAD r)_{\bL^2(T)}+(\nT\SCAL\GRAD v,r)_{L^2(\dT)}$.
Invoking the Cauchy--Schwarz inequality, the discrete inverse inequalities~\eqref{eq:disc_inv}-\eqref{eq:disc_trace}, and the bound~\eqref{eq:trace_virtual} readily gives 
$\|\Delta v\|_{L^2(T)}
\le Ch_T^{-1}\|\GRAD v\|_{\bL^2(T)}$. This completes the proof of~\eqref{eq:disc_trace_virtual}.  
\\
\textup{(ii)} Let $v_h\in \calV_{h,0}^k$. Integrating by parts 
in~\eqref{eq:def_Recv_HHO}, invoking the Cauchy--Schwarz inequality, and 
using~\eqref{eq:disc_trace_virtual} shows that $\|\GRAD \opRecv(\shvT)\|_{\bL^2(T)}
\le C|\shvT|_{\shVkT}$ for all $\shvT\in \shVkT$ and all $T\in\Th$
(we only use the bound on the normal derivative in~\eqref{eq:disc_trace_virtual}). 
Owing to Lemma~\ref{lem:virtual_HHO},
we infer that $\|\GRAD v_{h|T}\|_{\bL^2(T)} \le C|\shIvkT(v_{h|T})|_{\shVkT}$. The assertion now
follows by invoking the lower bound in Lemma~\ref{lem:stab_HHO} and 
summing cellwise since $a_h^{\textsc{hho}}(v_h,v_h)
= \sum_{T\in\Th} a_T(\shIvkT(v_{h|T}),\shIvkT(v_{h|T}))$.
\eProof

\chapter{Linear elasticity and hyperelasticity}
\label{Chap:Elas}

In this chapter, we show how to discretize using HHO methods 
linear elasticity and nonlinear hyperelasticity problems. 
In particular, we pay particular attention to the robustness of the discretization
in the quasi-incompressible limit.
For linear elasticity, we reconstruct the strain tensor in the space composed of
symmetric gradients of vector-valued polynomials. 
For nonlinear hyperelasticity, 
we reconstruct the deformation gradient in a full tensor-valued
polynomial space, and not just in a space composed of polynomial gradients.
We also consider a second gradient reconstruction
in an even larger space built using Raviart--Thomas polynomials,
for which no additional stabilization is necessary. Finally, we present
some numerical examples.

\section{Continuum mechanics}
\label{Sec:Elas}

We are interested in finding the static equilibrium configuration of an
elastic continuum body that occupies the domain $\Dom$ in the
reference configuration.
Here, $\Dom \subset \Rd$, $d \in \{2,3\}$, is a 
bounded Lipschitz domain with unit outward normal $\bn$ and
boundary partitioned 
as $\front = \overline{\Bn} \cup \overline{\Bd}$
with two relatively open and disjoint subsets $\Bn$ and $\Bd$.
The body undergoes deformations under the
action of a body force $\loadext:\Dom\to\Rd$, a traction force $\Gn:\Bn\to\Rd$,
and a prescribed displacement $\vuD:\Bd\to\Rd$.
We assume that $\Bd$ has positive measure so as to prevent rigid-body
motions. Due to the deformation, a point $\bx \in \Dom$ in the reference
configuration is mapped to a point $\bx' = \bx + \vu(\bx)$ in the equilibrium configuration, 
where $\vu : \Dom \rightarrow \Rd$ is the displacement field.

\subsection{Infinitesimal deformations and linear elasticity}
\label{sec:model-linear_elas}

Since we are concerned here with infinitesimal deformations, 
a relevant measure of the deformations of the body is 
the linearized strain tensor such that
\begin{equation}
\strain(\vu) \eqq \frac{1}{2} ( \GRAD \vu + \GRAD \vu\tr).
\label{Eq4:def-eps}
\end{equation}
Notice that $\strain(\vu)$ 
takes values in the space $\Msym$ composed of symmetric tensors of order $d$. 
Moreover, in the framework of linear isotropic elasticity, 
the internal stresses in the body are described at any point in $\Dom$
by the stress tensor $\stress$ which depends on the linearized strain tensor $\strain$
at that point. The constitutive stress-strain relation is linear and takes the form
\begin{equation}
\stress(\strain) \eqq 2\mu\strain+\lambda \trace(\strain){\matrice{I}}_{d},
\label{Def_of_pols}
\end{equation}
where $\lambda$ and $\mu$ are material parameters called
Lam\'e coefficients, and ${\matrice{I}}_{d}$ is the
identity tensor in $\Real^{d\CROSS d}$. Notice that $\stress(\strain)$
also takes values in $\Msym$ (the symmetry of $\stress(\strain)$ is actually
a consequence of the balance of angular momentum in infinitesimal deformations).
For simplicity, we assume that $\lambda$ and $\mu$ are constant in $\Dom$.
Owing to thermodynamic stability, we have 
$\mu>0$ and $\lambda+\frac23 \mu >0$.
The coefficient $\kappa\eqq\lambda+\frac23 \mu$,
called bulk modulus, describes the compressibility of the material.
Very large values relative to $\mu$, \ie $\lambda \gg \mu$,
correspond to almost incompressible materials.  

In the above setting, the displacement field 
$\vu : \Dom \rightarrow \Rd$ satisfies the following equations:
\begin{alignat}{2}
-\mdivergence \stress(\strain(\vu))&=\loadext&\quad&\text{in $\Dom$},\\
\vu&=\vuD&\quad&\text{on $\Bd$},\label{Eq4:elas3}\\
\stress(\strain(\vu)) \bn&=\Gn&\quad&\text{on $\Bn$},\label{Eq4:elas4}
\end{alignat}
together with~\eqref{Eq4:def-eps} and \eqref{Def_of_pols}.
Setting $\vecteur{H}^1(\Dom)\eqq H^1(\Dom;\Rd)$, 
the functional space composed of the kinematically admissible displacements and its tangent
space are
\begin{align}
\VD &\eqq \bset \vv \in \vecteur{H}^1(\Dom)\tq \vv_{|\Bd} = \vuD \eset, \\
\Vz &\eqq \bset \vv \in \vecteur{H}^1(\Dom)\tq \vv_{|\Bd} = \vzero \eset.
\end{align}
Recall that the $\vecteur{H}^1$-norm is defined as
$\norme[\vecteur{H}^1(\Dom)]{\vv}\eqq \big(\|\vv\|_{\Ldeuxd}^2
+\ell_\Dom^2\normem[\Dom]{\GRAD\vv}^2\big)^{\frac12}$, where 
the length scale $\ell_\Dom\eqq\diam(\Dom)$ is introduced to be
dimensionally consistent. Since $\mes{\Bd}>0$, Korn's inequality
implies that there is $C_{\textsc{k}}>0$ such that $C_{\textsc{k}}\normem[\Dom]{\GRAD \vv}
\le \normem[\Dom]{\strain(\vv)}$ for all $\vv\in \Vz$
(see, \eg \cite{Horgan:95}, \cite[Thm.~10.2]{MacLean_2000}). Moreover, 
the Poincar\'e--Steklov inequality
applied componentwise shows that there is $C_{\textsc{ps}}>0$ such that 
$C_{\textsc{ps}}\normev[\Dom]{\vv}\le \ell_\Dom \normem[\Dom]{\GRAD \vv}$ for all
$\vv\in\Vz$. Assuming $\loadext\in\Ldeuxd$ and $\Gn\in \bL^2(\Bn)$,
the weak formulation of the linear elasticity problem is as follows:
Seek $\vu \in \VD$ such that
\begin{equation}
a(\vu,\vw) = \ell(\vw)\eqq \psv[\Dom]{\loadext}{\vw} + \psv[\Bn]{\Gn}{\vw},
\quad \forall \vw \in \Vz,
\label{eq:weak-elas}
\end{equation}
with the bilinear form $a:\Hund\times\Hund\to\Real$ such that
\begin{equation}
a(\vv,\vw) \eqq \psm[\Dom]{\stress(\strain(\vv))}{\strain(\vw)} =
2 \mu \psm[\Dom]{\strain(\vv)}{\strain(\vw)} + 
\lambda \pss[\Dom]{\divergence \vv}{\divergence \vw},
\end{equation}
where we used that $\trace(\strain(\vv))=\divergence\vv$.
Simple manipulations show that
\begin{equation}
2\mu\norme[\ell^2]{\strain}^2 + \lambda |\trace(\strain)|^2 \ge
\min(2\mu,3\kappa)\norme[\ell^2]{\strain}^2,
\end{equation}
where $\norme[\ell^2]{\strain}^2 \eqq \strain{:}\strain = \sum_{1\le i,j\le d}|\strain_{ij}|^2$.
Combining this bound with the Korn and Poincar\'e--Steklov inequalities
shows that the bilinear form $a$ is coercive on $\Vz$. Hence, after lifting
the Dirichlet datum, one can show that the model problem~\eqref{eq:weak-elas} is
well-posed by invoking the Lax--Milgram lemma.

Standard convexity arguments show that 
the weak solution $\bu\in\VD$ to~\eqref{eq:weak-elas} is the unique minimizer
in $\VD$ of the energy functional $\engy:\VD\to\Real$ such that
\begin{equation} \label{eq:elast_engy}
\engy(\vv) \eqq \frac12
\Big(2\mu\normem[\Dom]{\strain(\vv)}^2 + \lambda\normes[\Dom]{\divergence \vv}^2\Big)
- \ell(\bv).
\end{equation}
Moreover, the weak formulation~\eqref{eq:weak-elas} expresses the 
principle of virtual work, wherein the test function $\vw$ plays
the role of a virtual displacement.

\bRem[Rigid-body motions]
An important fact in continuum mechanics is that the gradient and strain operators have
different kernels. In fact, $\GRAD \vv=\vzero$ if and only if there is $\ba\in\Real^d$
such that $\vv=\ba$, \ie the displacement field $\vv$ represents a translation. 
Instead, $\strain(\vv)=\vzero$ if and only if $\vv\in \RM$ where
\begin{equation} \label{eq:def_rigid_body}
\RM\eqq \begin{cases}
\bpolP_{3}^0 + \vecteur{x} \CROSS \bpolP_{3}^0&\text{if $d=3$}; \quad \dim(\RM)=6,\\
\bpolP_{2}^0 + \vecteur{x}^\perp \polP_{2}^0&\text{if $d=2$}; \quad \dim(\RM)=3,
\end{cases}
\end{equation}
where $\bpolP_d^0\eqq \polP_d^0(\Real^{d};\Real^d)$, $\bx$ is the position vector in $\Real^d$, 
$\CROSS$ denotes the cross product in $\Real^3$ and $\bx^\perp=(x,y)^\perp=(-y,x)$ if $d=2$.
Notice that for $d\in\{2,3\}$, we have $\bpolP_d^0\subsetneq \RM \subsetneq \bpolP_d^1$.
Fields in $\RM$ are called rigid-body motions (or translation-rotation motions).
\eRem

\subsection{Finite deformations and hyperelasticity}
\label{sec:model-hyperelas}

We are now concerned with finite deformations. 
We adopt the Lagrangian description so that all the differential operators are
taken with respect to the coordinates in the reference configuration. 
The deformations are measured using the deformation gradient
\begin{equation} \label{eq:def_Fdef}
\Fdef(\vu) \eqq \matrice{I}_d +\grad\vu,
\end{equation}
taking values in the set $\Rdd_{+}$ of $d \times d$ matrices with positive determinant. 
In the setting of homogeneous hyperelastic materials, the internal efforts in the body are
described at any point in $\Dom$ by the first Piola--Kirchhoff tensor $\PK$
which depends (nonlinearly) on the deformation gradient $\Fdef$ at that point. 
The constitutive relation between $\PK$ and $\Fdef$ is derived by postulating a
strain energy density $\Psi : \Rdd_{+} \rightarrow \Reel$ and setting 
\begin{equation} \label{eq:def_PK}
\PK(\Fdef) \eqq \partial_{\Fdef} \Psi(\Fdef).
\end{equation} 
We will mainly deal with hyperelastic materials of Neohookean type 
extended to the compressible range such that
\begin{equation}\label{NeoLaw}
\Psi(\Fdef) \eqq \frac{\mu}{2}  ( \Fdef{:}\Fdef - d) - \mu \ln J + \frac{\lambda}{2} (\ln J)^2,
\qquad J\eqq \det(\Fdef),
\end{equation}
where $\mu$ and $\lambda$ are material constants. Since $\partial_{\Fdef} J = J\Fdef\mtr$, 
\eqref{NeoLaw} gives
\begin{equation}
\PK(\Fdef) = \mu ( \Fdef - \Fdef\mtr) +  \lambda  \ln J  \Fdef\mtr.
\end{equation}

In the above setting, the displacement field 
$\vu : \Dom \rightarrow \Rd$ satisfies the following equations:
\begin{alignat}{2} 
-\mdivergence{\PK(\Fdef(\vu))} &= \loadext & \quad &\text{in $\Dom$}, \label{eq_pnba} \\
\vu &= \vuD & \quad  &\text{on $\Bd$}, \label{eq_pnbb}\\
\PK(\Fdef(\vu)) \,  \bn &= \Gn & \quad &\text{on $\Bn$},\label{eq_pnbc}
\end{alignat}
together with \eqref{eq:def_Fdef} and \eqref{eq:def_PK}, and 
where we assumed so-called dead external forces $\loadext$ and $\Gn$ (\ie independent of 
the deformed configuration).
Defining the energy functional $\engy: \VD \rightarrow \Reel$ such that
\begin{equation}\label{eq:def_calE}
\engy(\vv) \eqq \int_{\Dom} \Psi(\Fdef(\vv)) \dV - \ell(\vv),
\end{equation}
with the linear form $\ell$ defined in~\eqref{eq:weak-elas},
the static equilibrium problem~\eqref{eq_pnba}--\eqref{eq_pnbc} 
consists of seeking the stationary points of the energy functional $\engy$ which satisfy the following weak form of the Euler--Lagrange equations:
\begin{equation} \label{eq:Euler_Lag_hypereal}
0 =  D \engy (\vu) [ \vv] = \int_{\Dom} \PK(\Fdef(\vu)) : \grad \vv  \dV - \ell(\vv),
\end{equation}
for all virtual displacements $\vv\in \Vz$. We assume that the strain energy density function $\Psi$ is polyconvex, so that local minimizers of the energy functional exist (see ~\cite{Ball1976}). We refer the reader to the textbooks 
\cite{Bonet1997,Ciarlet1988,Ogden1997} for further insight into the physical modeling.

\section{HHO methods for linear elasticity}
\label{sec:elas_HHO}

The goal of this section is to present and analyze the HHO method to
discretize the linear elasticity problem introduced in Sect.~\ref{sec:model-linear_elas}.
We assume in the whole section that $\lambda\ge0$.

\subsection{Discrete unknowns, reconstruction, and stabilization}
\label{sec:lin_elas_setting}

Let $\Th$ be a mesh of $\Dom$ belonging to a shape-regular mesh sequence 
(see Sect.~\ref{sec:mesh} and~\ref{sec:mesh_reg}). We additionally require that
every mesh cell $T\in\Th$ is star-shaped with respect to every point in a ball with radius
uniformly equivalent to $h_T$; this will allow us to invoke a local Korn inequality. 
We assume that $\Dom$ is a polyhedron so that the mesh covers $\Dom$ exactly. 
Moreover, we assume that every mesh boundary face belongs either to $\Bd$ or to
$\Bn$; the corresponding subsets of $\calFb$ are denoted by $\calFbD$ and $\calFbN$. 
Recall that in HHO methods, the discrete unknowns are polynomials attached
to the mesh cells and the mesh faces. In the context of continuum mechanics,
both unknowns are vector-valued: the cell unknowns approximate the displacement field
in the cell, and the face unknowns approximate its trace on the mesh faces. 
For brevity, we only consider the equal-order setting for the cell and face unknowns.
One important difference with the diffusion model problem is that we now take 
the polynomial degree $k\ge1$. The reason for excluding the case $k=0$ is related to
the necessity to control the rigid-body motions in each mesh cell (for a lowest-order
nonconforming method, see \cite{BDiPG:19}). 

\begin{figure}
    \centering
     \subfloat[$k=1$]{
        \centering
        \includegraphics[scale=0.4]{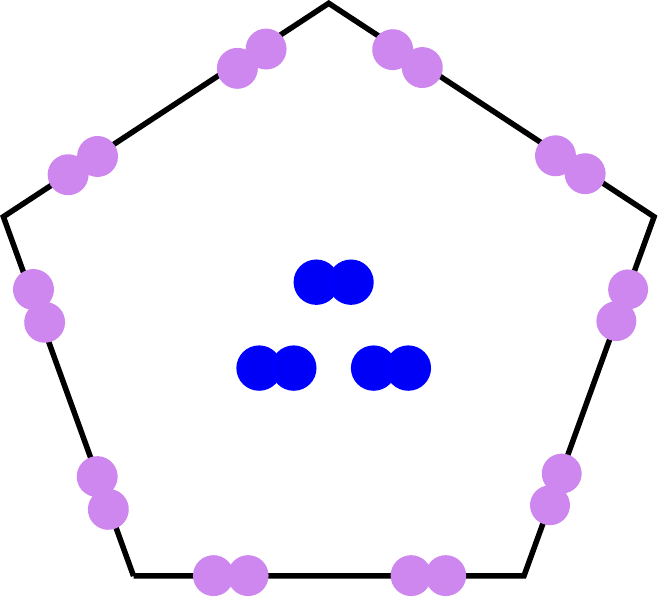}
  }
    ~ 
    \subfloat[$k=2$]{
        \centering
	    \includegraphics[scale=0.4]{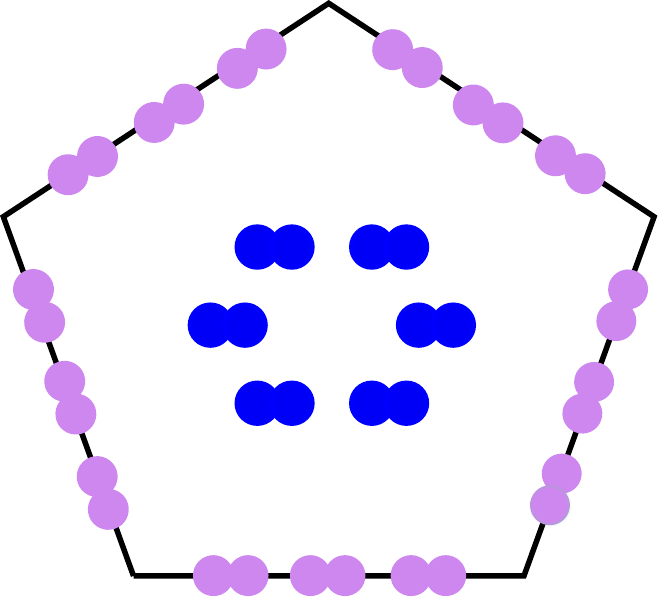}
    }
     ~ 
    \subfloat[$k=3$]{
        \centering
	    \includegraphics[scale=0.43]{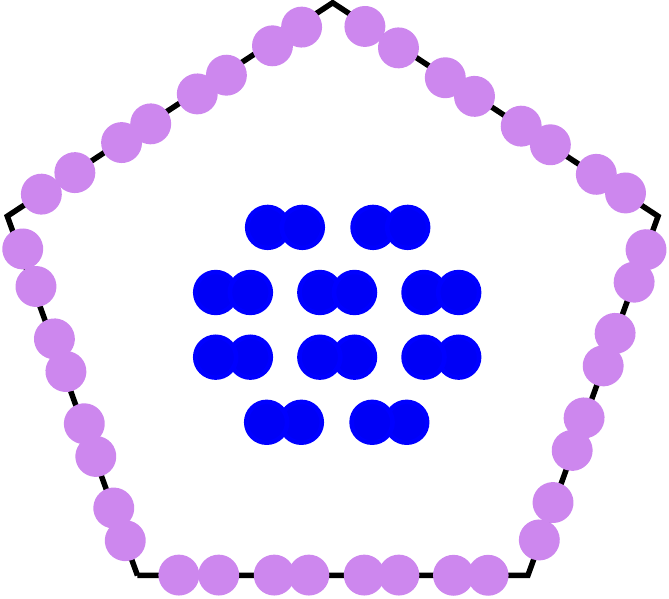}
    }
\caption{Face (violet) and cell (blue) unknowns in $\hVkT$ in a pentagonal cell ($d=2$) for $k\in\{1,2,3\}$ (each dot in a pair represents one basis function associated with one Cartesian component).}
\label{fig_HHO_dofs_elas}
\end{figure}

For every mesh cell $T\in\Th$, we set
\begin{equation}
\hVkT \eqq \vPkd(T) \times \vVkdT, \qquad
\vVkdT \eqq \bigtimes_{F\in\FT} \vPkF(F),
\label{eq:def_VkT}
\end{equation}
with $\vPkd(T)\eqq \Pkd(T;\Real^d)$ and $\vPkF(F)\eqq \PkF(F;\Real^d)$;
see Figure~\ref{fig_HHO_dofs_elas}.
A generic element in $\hVkT$ is denoted by $\hvT\eqq (\vT,\vdT)$.
The HHO space is defined as 
\begin{equation} \label{eq:def_hVkh}
\hVkh \eqq \Vkh \times \VkFh,
\qquad
\Vkh \eqq \bigtimes_{T\in\calT} \vPkd(T),
\qquad
\VkFh \eqq \bigtimes_{F\in\calF} \vPkF(F),
\end{equation}
so that $\dim(\hVkh) = d{k+d\choose d}\#\calT + d{k+d-1\choose d-1}\#\calF$.
A generic element in $\hVkh$ is denoted by $\hvh\eqq (\vTh,\vFh)$ with
$\vTh\eqq (\vv_T)_{T\in\calT}$ and $\vFh\eqq (\vv_F)_{F\in\calF}$, and we 
localize the components of $\hvh$ associated with a mesh cell $T\in\Th$ 
and its faces by using the notation
$\hvT\eqq \big(\vv_T,\vv_{\dT}\eqq (\vv_F)_{F\in\FT}\big)\in \hVkT$.
The local HHO reduction operator  
$\hIkT:\vecteur{H}^1(T)\to \hVkT$ is such that
\begin{equation} \label{eq:def_hIkT}
\hIkT(\vv)\eqq (\PikTv(\vv),\PikdTv(\vv_{|\dT})),
\end{equation}
and its global counterpart $\Ikh : \vecteur{H}^1(\Dom) \to \hVkh$ is 
$\Ikh(\vv)\eqq (\PikcTv(\vv),\PikcFv(\vv_{|\calF}))$,
where the $L^2$-orthogonal projections act componentwise.

The present HHO method for linear elasticity hinges on three local operators: 
(i) a displacement reconstruction operator, (ii) a divergence reconstruction operator,
and (iii) a stabilization operator.  The displacement reconstruction operator
is based on inverting the Gram matrix associated with the strain tensor. For this reason,
we need to specify additionally the rigid-body motions of the reconstructed displacement.
For a polynomial degree $l\ge1$, we set 
\begin{equation}
\bpolP_d^l(T)^\perp \eqq \{ \bq \in \bpolP_d^l(T)\tq
{\PizTv(\bq)=\vzero,\ \PizTm(\gradss\bq)=\mzero}\},
\end{equation}
{where $\gradss\vv\eqq \frac12(\GRAD\vv-\GRAD\vv\tr)$ denotes the skew-symmetric part of the gradient (recall that the $\bL^2$-orthogonal projection $\PizTv$ returns the
mean-value over $T$)}.
It is easy to see that $\bpolP_d^l(T)=\bpolP_d^l(T)^\perp\oplus\RM$, where
$\RM$ is the space of rigid-body motions defined in~\eqref{eq:def_rigid_body}.
The displacement reconstruction operator $\opDep : 
\hVkT \to \bpolP_d^{k+1}(T)$ is defined such that for all $\hvT\in\hVkT$,
$\opDep(\hvT)\in \bpolP_d^{k+1}(T)$ is uniquely specified as follows:
\begin{align}
&\psm[\T]{\strain(\opDep(\hvT))}{\strain(\bq)} =
-\psv[\T]{\vT}{\divergence\strain(\bq)} + \psv[\dT]{\vdT}{\strain(\bq)\nT},
\label{eq:def_Rec_HHO_elas} \\
&{\PizTv(\opDep(\hvT))=\PizTv(\vv_T), \quad
\PizTm(\gradss\opDep(\hvT))=\frac{1}{\mes{T}}\int_{\dT}\vdT\otimes_{\textrm{ss}}\nT\ds,}
\label{eq:def_Rec_HHO_mean_elas}
\end{align}
where \eqref{eq:def_Rec_HHO_elas} holds for all $\bq\in \bpolP_d^{k+1}(T)^\perp$ and
$\ba\otimes_{\textrm{ss}}\bb \eqq \frac12(\ba\otimes\bb-\bb\otimes\ba)$ 
for all $\ba,\bb\in\Real^d$.
Proceeding as in the proof of Lemma~\ref{lem:ell_proj_HHO} and since
$\PizTm(\gradss\opDep(\hIkT(\vv)))=\frac{1}{\mes{T}}\int_{\dT}\PikdTv(\vv_{|\dT})\otimes_{\textrm{ss}}\nT\ds=\PizTm(\gradss\vv)$ (owing to \eqref{eq:def_Rec_HHO_mean_elas},
\eqref{eq:def_hIkT}, and integration by parts), one readily
establishes the following identity.

\bLem[Elliptic projection]
We have $\opDep\circ\hIkT= \bcalE_T^{k+1}$ where
$\bcalE_T^{k+1} : \vecteur{H}^1(T) \to \bpolP_d^{k+1}(T)$ 
is the elliptic strain projection such that
$\psm[\T]{\strain(\bcalE_T^{k+1}(\vv)-\vv)}{\strain(\bq)} = 0$ for all
$\bq\in \bpolP_d^{k+1}(T)^\perp$, $\PizTv(\bcalE_T^{k+1}(\vv)-\vv)=\vzero$, and
$\PizTm(\gradss(\bcalE_T^{k+1}(\vv)-\vv))=\mzero$.
\eLem

The divergence reconstruction operator $\opDiv : \hVkT \to \Pkd(T)$
is defined by solving the following well-posed problem: For all $\hvT\in\hVkT$,
\begin{equation} \label{eq:op_rec_div_HHO}
\pss[\T]{\opDiv(\hvT)}{q} = -\psv[\T]{\vT}{\GRAD q} + \psv[\dT]
{\vdT}{q\nT},
\end{equation}
for all $q\in \Pkd(T)$. In practice, the computation
of $\opDiv(\hvT)$ entails inverting the mass matrix in $\Pkd(T)$. We have the following
important commutation result. 

\bLem[Commuting with divergence]  \label{lem:commut_div_HHO}
The following holds true:
\eq{ \label{eq:commut_div_HHO}
\opDiv(\hIkT(\vv)) = \PikTs(\divergence\vv), \quad
\forall \vv \in \vecteur{H}^1(T).
}
\eLem

\bproof
Let $q\in \Pkd(T)$. Since $\GRAD q\in \vPkd(T)$ and $q\nT \in \vVkdT$, we have
\begin{align*}
\pss[T]{\opDiv(\hIkT(\vv))}{q} &= -\psv[\T]{\PikTv(\vv)}{\GRAD q} +
\psv[\dT]{\PikdTv(\vv_{|\dT})}{q\nT} \\
&= -\psv[T]{\vv}{\GRAD q} + \psv[\dT]{\vv}{q\nT} = \pss[T]{\divergence\vv}{q}.
\end{align*}
Since $\opDiv(\hIkT(\vv))\in \Pkd(T)$ and $q$ is arbitrary in $\Pkd(T)$, 
this proves~\eqref{eq:commut_div_HHO}.
\eproof

Finally, the stabilization operator is inspired from the one 
devised for the Poisson model problem in Sect.~\ref{sec:stab_diff}. Here, 
we define $\opStv : \hVkT \to \vVkdT$ such that for all $\hvT \in \hVkT$, 
\begin{equation}
\opStv(\hvT) \eqq \PikdTv\Big(\vv_{T|\dT}
-\vdT+\big((I-\PikTv)\opDep(\hvT)\big)_{|\dT}\Big), \label{eq:def_SKk_HHO_elas}
\end{equation}
where $I$ is the identity operator. Letting $\vecteur{\delta}_{\dT}\eqq \vv_{T|\dT}-\vdT$,
we observe that $\opStv(\hvT) = \PikdTv\big(\vecteur{\delta}_{\dT}-((I-\PikTv)\opDep(\vzero,
\vecteur{\delta}_{\dT}))_{|\dT}\big)$.

\subsection{Discrete problem, energy minimization, and traction recovery}
\label{sec:disc_pb_elas}

The global bilinear form $a_h:\hVkh\times\hVkh\to\Real$ is assembled cellwise
as in Sect.~\ref{sec:assembling_HHO} by setting
$a_h(\hvh,\hwh) \eqq \sum_{T\in\Th} a_T(\hvT,\hwT)$ where for all $T\in\Th$,
the local bilinear forms $a_T:\hVkT\times\hVkT\to\Real$ are such that 
\begin{align}
a_T(\hvT,\hwT) \eqq {}&
2\mu \psm[T]{\strain(\opDep(\hvT))}{\strain(\opDep(\hwT))} 
+ \lambda \pss[T]{\opDiv(\hvT)}{\opDiv(\hwT)} \nonumber \\
&+ 2\mu h_T^{-1}\psv[\dT]{\opStv(\hvT)}{\opStv(\hwT)}.
\end{align}
Notice that the stabilization term is weighted by the Lam\'e parameter $\mu$.
To account for the Dirichlet boundary condition, we define the subspaces
\begin{align}
\VkFhd &\eqq \{ \vFh \in \VkFh \tq \vF = \PikFv(\vuD), \ \forall F \in \calFbD \}, \label{eq:def_VkFhd}\\
\VkFhz &\eqq \{ \vFh \in \VkFh \tq \vF = \vecteur{0}, \ \forall F \in \calFbD \},\label{eq:def_VkFhz}
\end{align}
as well as $\hVkhd \eqq \Vkh \times \VkFhd$ and $\hVkhz \eqq \Vkh \times \VkFhz$.
The discrete problem is as follows:
\begin{equation}
\left\{
\begin{array}{l}
\text{Find $\huh \in \hVkhd$ such that}\\[2pt]
a_h(\huh,\hwh) = \ell_h(\hwh)\eqq \psv[\Dom]{\loadext}{\vw_\calT}+\psv[\Bn]{\Gn}{\vw_\calF},
\quad \forall \hwh \in \hVkhz.
\end{array}
\right.
\label{eq:weak_HHO_elas}
\end{equation}
Notice that the cell component of the test function is used against the body force 
$\loadext$, whereas the face component is used against the traction force $\Gn$.
We will see in the next section that the bilinear form $a_h$ is coercive on $\hVkhz$
so that the discrete problem~\eqref{eq:weak_HHO_elas} is well-posed. Moreover,
let us define the discrete energy functional $\engy_h:\hVkhd\to\Real$ such that
\begin{equation} \label{eq:elast_engy_HHO}
\engy_h(\hvh) \eqq \frac12
\Big(2\mu\normem[\Dom]{\strain(\opDeph(\hvh))}^2 + \lambda\normes[\Dom]{\opDivh(\hvh)}^2+2\mu s_h(\hvh,\hvh)\Big)
- \ell_h(\hvh),
\end{equation}
where the global reconstructions $\opDeph(\hvh)$ and $\opDivh(\hvh)$ are such that
$\opDeph(\hvh)_{|T}\eqq \opDep(\hvT)$ and $\opDivh(\hvh)_{|T}\eqq\opDiv(\hvT)$ for all
$T\in\Th$, and with the global stabilization bilinear form 
\begin{equation} \label{eq:def_sh_elas}
s_h(\hvh,\hwh)\eqq \sum_{T\in\Th} h_T^{-1}\psv[\dT]{\opStv(\hvT)}{\opStv(\hwT)}.
\end{equation} 
Then, the same arguments as in the proof of Proposition~\ref{prop:HHO_min_diff}
show that $\huh \in \hVkhd$ solves~\eqref{eq:weak_HHO_elas} if and only if $\huh$ 
minimizes $\engy_h$ in $\hVkhd$.

The algebraic realization of~\eqref{eq:weak_HHO_elas} leads to a linear system
with symmetric positive-definite stiffness matrix having the same block-structure
as in~\eqref{eq:HHO_alg}. The right-hand side vector can now have nonzero
face components due to the Neumann boundary condition. In any case, a computationally
effective way to solve the linear system is again to use static condensation:
one eliminates locally all the cell unknowns, solves the global transmission 
problem coupling all the face unknowns, and finally recovers the cell unknowns by local
post-processing. Moreover, the result of Proposition~\ref{prop:trans_HHO} on the 
global transmission problem can be 
readily extended to the setting of linear elasticity provided the right-hand side
of~\eqref{eq:glob_trans_pb_HHO} is modified to include the contribution of
the Neumann boundary condition.

The material of Sect.~\ref{sec:flux_recovery} on flux recovery can be readily 
adapted to the present setting leading to the important notion of equilibrated
tractions defined on the boundary of the mesh cells and at the Neumann boundary faces.
Let $\optStv : \vVkdT\to \vVkdT$ be defined such that
\begin{equation} \label{eq:def_optStv}
\optStv(\mu) \eqq \PikdTv\Big(\bmu-\big((I-\PikTv)\opDep(0,\bmu)\big)_{|\dT}\Big),
\end{equation}
and let 
$\optStv^{*}: \vVkdT\to \vVkdT$ be its adjoint such that
$\psv[\dT]{\optStv^{*}(\blambda)}{\bmu} = \psv[\dT]{\blambda}{\optStv(\bmu)}$
for all $\blambda,\bmu\in \vVkdT$.
Then, for all $\hvh\in \hVkh$, we can define numerical tractions 
$\vecteur{T}_{\dT}(\hvT)\in \vVkdT$ 
at the boundary of every mesh cell $T\in\calT$ by setting
\begin{equation} \label{eq:def_num_traction}
\vecteur{T}_{\dT}(\hvT) \eqq -\stress(\strain(\opDep(\hvT)))_{|\dT} \nT
+ 2\mu h_T^{-1}(\optStv^*\circ \optStv) (\vv_{T|\dT}-\vdT).
\end{equation}
A direct adaptation of the proof of Proposition~\ref{prop:flux_recovery}
establishes the following result.

\bProp[Rewriting with tractions] \label{prop:flux_recovery_elas}
Let $\huh\in \hVkhd$ solve~\eqref{eq:weak_HHO_elas} and 
let the tractions $\vecteur{T}_{\dT}(\huT)\in \vVkdT$
be defined as in~\eqref{eq:def_num_traction} for all $T\in\calT$. The following holds:\\
\textup{(i)} Equilibrium at every mesh interface $F=\partial T_-\cap \partial T_+ \cap H_F\in \calFi$:
\begin{equation} \label{eq:equil_traction}
\vecteur{T}_{\dT_-}(\hat{\vu}_{T_-})_{|_F} + \vecteur{T}_{\dT_+}(\hat{\vu}_{T_+})_{|_F} = \vzero,
\end{equation}
and at every Neumann boundary face $F=\partial T_-\cap \Bn\cap H_F \in \calFbN$:
\begin{equation} \label{eq:equil_Neumann}
\vecteur{T}_{\dT_-}(\hat\vu_{T_-})_{|_F}  + \PikFv(\Gn) = \vzero.
\end{equation}
\textup{(ii)} Balance with the source term in every mesh cell $T\in\calT$: For all $\bq\in \vPkd(T)$,
\begin{equation} \label{eq:balance_traction}
\psm[T]{\stress(\strain(\opDep(\huT)))}{\strain(\bq)} + \psv[\dT]{\vecteur{T}_{\dT}(\huT)}{\bq} = \psv[T]{\loadext}{\bq}.
\end{equation}
\textup{(iii)}
\eqref{eq:equil_traction}-\eqref{eq:equil_Neumann}-\eqref{eq:balance_traction} are an equivalent rewriting of~\eqref{eq:weak_HHO_elas}. 
\eProp
 
\bRem[Literature]
HHO methods for linear elasticity were introduced in \cite{DiPEr:15}.
HDG methods for linear elasticity were developed, among others, in
\cite{SoCoS:09,Nguyen2012a,FuCoS:15}, weak Galerkin methods 
in \cite{WaWWZ:16}, and a hybridizable weakly conforming Galerkin method
in \cite{KrWWW:16}.
\eRem

\subsection{Stability and error analysis}
\label{sec:analysis_elas}

The stability and error analysis relies on the vector-valued version of the
inequalities from Sect.~\ref{sec:basic_tools}. In particular, we need
the multiplicative trace inequality from Lemma~\ref{lem:mtr} and the discrete
inverse inequalities from Lemma~\ref{lem:disc_inv}. Moreover, in addition to the
local Poincar\'e--Steklov inequality from Lemma~\ref{lem:PS}, 
we need the following local Korn inequality (see 
\cite[App.~A.1]{BDiPD:18}): 
There is $C_{\textsc{k}}$ such that for all $T\in\Th$,
\begin{equation} \label{eq:local_Korn}
\normem[\T]{\GRAD \vv} \le C_{\textsc{k}}\normem[\T]{\strain(\vv)}, \quad
\forall \vv\in \vecteur{H}^1(T)^\perp,
\end{equation}
with $\vecteur{H}^1(T)^\perp\eqq \{\vv\in \vecteur{H}^1(T) \tq
{\PizTv(\vv)=\vzero,\ \PizTm(\gradss\vv)=\mzero}\}$.
Combining the local Poincar\'e--Steklov and Korn inequalities yields
$\normev[\T]{\vv} \le C_{\textsc{ps}}C_{\textsc{k}} h_T\normem[\T]{\strain(\vv)}$
for all $\vv\in \vecteur{H}^1(T)^\perp$.

\bLem[Stability]  \label{lem:stab_HHO_elas}
Equip $\hVkT$ with the seminorm
$\snorme[\strain,T]{\hvT}^2 \eqq
\normem[\T]{\strain(\vT)}^2 + h_{T}^{-1}\normev[\dT]{\vT-\vdT}^2$.
Assume $k\ge1$.
There are $0<\alpha\le\omega < +\infty$ such that for all $T\in\Th$ and all $\hvT\in \hVkT$,
\begin{align}
\alpha \snorme[\strain,T]{\hvT}^2 &\le
\normem[\T]{\strain(\opDep(\hvT))}^2
+ h_{T}^{-1}\normev[\dT]{\opStv(\hvT)}^2 \le
\omega \snorme[\strain,T]{\hvT}^2.
\end{align}
\eLem

\bproof The only difference with the proof of Lemma~\ref{lem:stab_HHO}
arises in the proof of~\eqref{eq:bnd_on_trace_rT} when bounding 
$(I-\PikTv)(\opDep(\hvT))$. We first notice that there is (a unique)
$\vr_{\vU}\in \RM$ such that $\PizTv(\opDep(\hvT)-\vr_{\vU})=\vzero$ and
$\PizTm(\gradss(\opDep(\hvT)-\vr_{\vU}))=\mzero$. Since $(I-\PikTv)(\vr_{\vU})=\vzero$
(because $k\ge1$), $\opDep(\hvT)-\vr_{\vU}\in \vecteur{H}^1(T)^\perp$, and
$\strain(\vr_{\vU})=\mzero$, we infer that
\begin{align*}
\normev[T]{(I-\PikTv)(\opDep(\hvT))} 
\ifSp \else \alhere \fi 
=\normev[T]{(I-\PikTv)(\opDep(\hvT)-\vr_{\vU})}
\le \normev[T]{\opDep(\hvT)-\vr_{\vU}}\ifSp \alhere \fi \\
\ifSp \else \alhere \fi 
\le C_{\textsc{ps}}C_{\textsc{k}} h_T\normem[\T]{\strain(\opDep(\hvT)-\vr_{\vU})}
= Ch_T\normem[\T]{\strain(\opDep(\hvT))},\ifSp \alhere \fi
\end{align*}
owing to the
combined Poincar\'e--Steklov and Korn inequalities. All the other arguments in the proof
are the vector-valued version of those invoked for Lemma~\ref{lem:stab_HHO}.
\eproof


\bLem[Coercivity, well-posedness]  \label{lem:coer_HHO_elas}
The map $\hVkhz\ni \hvh \mapsto \norme[\strain,h]{\hvh}\eqq \big( \sum_{T\in\Th} 
\snorme[\strain,T]{\hvT}^2\big)^{\frac12}\in [0,+\infty)$ defines a norm on $\hVkhz$.
Moreover, the discrete bilinear form $a_h$ satisfies the coercivity property
\begin{equation}
a_h(\hvh,\hvh) \ge 2\mu \norme[\strain,h]{\hvh}^2,
\quad \forall \hvh\in\hVkhz, \label{eq:coer_HHO_elas}
\end{equation}
and the discrete problem~\eqref{eq:weak_HHO_elas} is well-posed.
\eLem 

\bproof
The only nontrivial property is the definiteness of the map.
Let $\hvh\in \hVkhz$ be such that
$\norme[\strain,h]{\hvh}=0$. Then, for all
$T\in\Th$, $\vT$ is a rigid displacement whose trace on
$\dT$ is $\vdT$.
Since two rigid displacements that coincide on a mesh face
are identical, we infer that
$\vTh$ is a global rigid displacement, and since $\vv_F=\vzero$ for all $F\in\calFbD$,
we conclude that $\vTh$ and $\vFh$ are zero. Hence, $\norme[\strain,h]{\SCAL}$
defines a norm on $\hVkhz$. The coercivity property~\eqref{eq:coer_HHO_elas}
follows by summing over the mesh cells 
the lower bound from Lemma~\ref{lem:stab_HHO_elas} and recalling that
$\lambda\ge0$ by assumption. 
Finally, the well-posedness of~\eqref{eq:weak_HHO_elas} results from
the Lax--Milgram lemma. 
\eproof

To derive an error estimate, we introduce the consistency
error such that 
\begin{equation}
\langle \delta_h,\hwh\rangle 
\eqq \ell_h(\hat\vw_h)-a_h(\Ikh(\vu),\hwh), \quad \forall \hwh\in\hVkhz,
\end{equation}
and we bound the dual norm $\|\delta_h\|_*\eqq \sup_{\hwh\in\hVkhz}
\frac{|\langle \delta_h,\hwh\rangle|}{\norme[\strain,h]{\hwh}}$.
For all $T\in\Th$ and all $\vv\in \vecteur{H}^{1+r}(T)$, $\phi\in H^{1+r}(T)$, $r > \frac12$,
we define the local (semi)norms 
\begin{equation} \label{eq:sharp_dagger_elas}
\snorme[\sharp,T]{\vv} \eqq \normem[\T]{\strain(\vv)}
+ h_T^{\frac12}\normem[\dT]{\strain(\vv)},
\quad
\norme[\dagger,T]{\phi} \eqq \normes[\T]{\phi}
+ h_T^{\frac12}\normes[\dT]{\phi}.
\end{equation} 
The global counterparts are
$\snorme[\sharp,\calT]{\vv}\eqq \big(\sum_{T\in\Th}
\snorme[\sharp,T]{\vv}^2\big)^{\frac12}$ for all $\vv\in \vecteur{H}^{1+r}(\calT)$ and
$\norme[\dagger,\calT]{\phi}\eqq \big(\sum_{T\in\Th}
\norme[\dagger,T]{\phi}^2\big)^{\frac12}$ for all $\phi\in H^{1+r}(\calT)$.
Let $\bcalE_\calT^{k+1}$ be the global elliptic strain projection such that
$\bcalE_\calT^{k+1}(\vv)_{|T}=\bcalE_T^{k+1}(\vv_{|T})$ for all $T\in\Th$ and all $\vv\in\Hund$.

\bLem[Consistency] \label{lem:consist_HHO_elas}
Assume that the exact solution satisfies 
$\vu\in \vecteur{H}^{1+r}(\Dom)$, $r>\frac12$. There is $C$,
uniform with respect to $\mu$ and $\lambda$, such that
\begin{equation}
\norme[*]{\delta_h}\le C\big( 
\mu \snorme[\sharp,\calT]{\vu-\bcalE_\calT^{k+1}(\vu)} +
\lambda \norme[\dagger,\calT]{\divergence\vu - \PikcTs(\divergence\vu)}\big).
\end{equation}
\eLem

\bproof
Since the exact solution 
satisfies $-\mdivergence \stress(\strain(\vu))=\loadext$ in $\Dom$ and 
$\stress(\strain(\vu))\bn=\Gn$ on $\Bn$, integrating by parts and using that the 
normal component of $\stress(\strain(\vu))$ is single-valued at every mesh interface and that
$\vw_F$ vanishes at every Dirichlet boundary face, we infer that
\begin{align*}
\ell(\hwh) &= \sum_{T\in\Th} \big((\stress(\strain(\vu)),\strain(\wT))_{\bL^2(T)}
- (\stress(\strain(\vu))\nT,\wT)_{\bL^2(\dT)}\big) +\psv[\Bn]{\Gn}{\vw_\calF} \\
&=\sum_{T\in\Th} \big((\stress(\strain(\vu)),\strain(\wT))_{\bL^2(T)}
+ (\stress(\strain(\vu))\nT,\vw_{\dT}-\wT)_{\bL^2(\dT)}\big).
\end{align*}
Similar manipulations to the Poisson model problem (see the proof
of Lemma~\ref{lem:bnd_HHO}) and the commuting property from Lemma~\ref{lem:commut_div_HHO}
give
$\langle \delta_h,\hwh\rangle = 
-\sum_{T\in\Th}(\term_{1,T}+\term_{2,T}+\term_{3,T})$, where
\begin{align*}
\term_{1,T}&\eqq 2\mu \psm[\dT]{\strain(\vu-\bcalE_T^{k+1}(\vu))\nT}{\wT-\vw_{\dT}}, \\
\term_{2,T}&\eqq \lambda \psv[\dT]{(\divergence\vu - \PikTs(\divergence\vu)))\nT}{\wT-\vw_{\dT}}, \\
\term_{3,T}&\eqq 2\mu h_{T}^{-1} \psv[\dT]{ \opStv(\hIkT(\vu))}{\opStv(\hwT)}.
\end{align*}
The first two terms are bounded by using the Cauchy--Schwarz inequality. 
The third term is bounded by proceeding as in the
proof of Lemma~\ref{lem:approx_HHO_St}, except that we additionally invoke
the local Korn inequality~\eqref{eq:local_Korn}. 
\eproof

As in Lemma~\ref{lem:esterr_HHO} for the Poisson model problem, 
stability (Lemma~\ref{lem:coer_HHO_elas}) and consistency
(Lemma~\ref{lem:consist_HHO_elas}) imply the bound 
\begin{equation}
\mu \norme[\strain,h]{\heh}\le C\big( 
\mu \snorme[\sharp,\calT]{\vu-\bcalE_\calT^{k+1}(\vu)} +
\lambda \norme[\dagger,\calT]{\divergence\vu - \PikcTs(\divergence\vu)}\big),
\end{equation}
with the discrete error $\heh\eqq \huh-\Ikh(\vu)$.
Using the triangle inequality and the approximation properties of the 
elliptic strain projection (see \cite[App.~A.2]{BDiPD:18})
leads, as in Theorem~\ref{th:esterr_HHO}, to the following energy-error estimate.
In the spirit of Sect.~\ref{sec:error_anal_HHO},
we use the notation
$|h^\alpha\bphi|_{\bH^\beta(\calT)}\eqq \big( \sum_{T\in\Th}
h_T^{2\alpha}|\bphi|_{\bH^\beta(T)}^2\big)^{\frac12}$ for $\bphi\in \bH^\beta(\calT)$,
and we let $\strain_\calT$ denote the strain operator applied cellwise 
(\ie using the broken gradient). 

\bTheo[Energy-error estimate] 
Let $\vu\in \VD$ be the exact solution and let $\huh\in\hVkhd$ be the HHO solution
solving~\eqref{eq:weak_HHO_elas}.
Assume that $\vu\in \vecteur{H}^{1+r}(\Dom)$, $r>\frac12$.
There is $C$, uniform with respect to $\mu$ and $\lambda$, such that
\begin{equation} \label{eq:en_err_elas1}
\mu \normem[\Dom]{\strain_\calT(\vu-\opDeph(\huh))} 
\le C \big( \mu \snorme[\sharp,\calT]{\vu-\bcalE_\calT^{k+1}(\vu)}
+ \lambda \norme[\dagger,\calT]{\divergence\vu - \PikcTs(\divergence\vu)} \big).
\end{equation}
Moreover, if $\vu\in \vecteur{H}^{t+1}(\calT)$ and $\divergence\vu\in H^{t}(\calT)$
for some $t\in [\frac12,k+1]$, we have
\begin{equation} \label{eq:en_err_elas2}
\mu \normem[\Dom]{\strain_\calT(\vu-\opDeph(\huh))}
\le C\big( \mu \snorme[\vecteur{H}^{t+1}(\calT)]{h^t\vu} +
\lambda \snorme[H^{t}(\calT)]{h^t\divergence\vu}\big).
\end{equation}
\label{th:esterr_HHO_elas}
\eTheo

\bRem[Quasi-incompressible limit]
The remarkable fact about the error estimates~\eqref{eq:en_err_elas1}-\eqref{eq:en_err_elas2} is
that the right-hand side depends on the second Lam\'e parameter $\lambda$ only
through the smoothness of $\divergence\vu$. Furthermore, 
the incompressible limit (\ie the Stokes equations)
can be treated by introducing a pressure
variable attached to the mesh cells \cite{AgBDP:15}. The pressure unknowns, up to the cell mean-value, can be locally eliminated together with the cell velocity unknowns by static condensation. Moreover, as shown in \cite{DPELS:16}, the discretization can be made
pressure-robust.
\eRem

An improved $\bL^2$-error estimate can be derived by adapting the arguments presented in
Sect.~\ref{sec:error_anal_HHO_L2}. We assume the following elliptic regularity property: There are a constant $C_{\mathrm{ell}}$
and a regularity pickup index $s\in(\frac12,1]$ such that for all $\bg\in\Ldeuxd$,
the unique field $\bzeta_{\bg}\in \Vz$ such that $a(\vv,\bzeta_{\bg})=(\vv,\bg)_{\Ldeuxd}$ for all $\vv\in \Vz$ satisfies the regularity estimate 
$\mu\|\bzeta_{\bg}\|_{\bH^{1+s}(\Dom)}+\lambda\|\divergence\bzeta_{\bg}\|_{H^{s}(\Dom)}\le C_{\mathrm{ell}} \ell_\Dom^{2} \|\bg\|_{\Ldeuxd}$. Then, proceeding as in the proof of Lemma~\ref{lem:L2_est_HHO} (see also \cite{DiPEr:15} for the original arguments) leads to the following discrete $\bL^2$-error estimate: There is $C$, uniform with respect to $\mu$ and $\lambda$, such that
\begin{align}
\mu \norme[\Dom]{\be_\calT} \le {}& C\ell_\Dom^{1-s}h^s \Big( 
\mu \snorme[\sharp,\calT]{\vu-\bcalE_\calT^{k+1}(\vu)}
+ \max(\mu,\lambda) \norme[\dagger,\calT]{\divergence\vu - \PikcTs(\divergence\vu)}
\nonumber \\
&+ h\normev[\Dom]{\loadext-\PikTv(\loadext)} + h^{\frac12}\normev[\Bn]{\Gn-\vecteur{\Pi}^k_{\calFbN}(\Gn)}\Big). 
\end{align}

\section{HHO methods for hyperelasticity}
\label{sec:HHO_hyperelasticity}

The goal of this section is to present and analyze two HHO methods to 
discretize the hyperelasticity problem introduced in Sect.~\ref{sec:model-hyperelas}. Following the ideas outlined in Sect.~\ref{sec:variants_rec}, we consider
\textup{(i)} 
a stabilized HHO method reconstructing the deformation gradient in $\Pkd(T;\Rdd)$ 
and \textup{(ii)} 
an unstabilized method reconstructing the deformation gradient in a larger
polynomial space built using Raviart--Thomas polynomials. 
For both methods, the discrete setting is the same as the one considered 
for the linear elasticity problem (see Sect.~\ref{sec:lin_elas_setting}):
the mesh $\Th$ satisfies the assumptions stated therein, 
the discrete unknowns belong to the local space $\hVkT\eqq \vPkd(T) \times \vVkdT$ 
defined in~\eqref{eq:def_VkT} 
for every mesh cell $T\in\Th$ and a polynomial degree $k\ge1$, 
and the global HHO space is the space $\hVkh\eqq \Vkh \times \VkFh$ 
defined in~\eqref{eq:def_hVkh} together with the subspaces 
$\hVkhd \eqq \Vkh \times \VkFhd$ and $\hVkhz \eqq \Vkh \times \VkFhz$ related to the enforcement of the Dirichlet boundary condition (see~\eqref{eq:def_VkFhd}-\eqref{eq:def_VkFhz}). 

\subsection{The stabilized HHO method}
\label{sec:sHHO}

The local gradient reconstruction operator $\opGRec:\hVkT \to \Pkd(T;\Rdd)$ 
is such that for all $\hvT\in\hVkT$,
$\opGRec(\hvT)\in \Pkd(T;\Rdd)$ is uniquely determined by the equations
\begin{equation} \label{eq:def_opGRec_mat}
(\opGRec(\hvT),\bq)_{\bL^2(T)}=-(\vvT,\DIV\bq)_{L^2(T)}+(\vv_\dT,
\bq\nT)_{\bL^2(\dT)}, \quad \forall \bq\in \Pkd(T;\Rdd).
\end{equation}
To compute $\opGRec(\hvT)$, it suffices to invert the mass matrix associated with the scalar-valued polynomial space $\Pkd(T)$ since only the right-hand side changes when computing each entry of the tensor $\opGRec(\hvT)$. Notice in passing that
\begin{equation} \label{eq:tr_G_opDiv}
\trace(\opGRec(\hvT))=\opDiv(\hvT), \quad \forall \hvT\in\hVkT,
\end{equation} 
where $\opDiv$ is the divergence reconstruction operator defined in~\eqref{eq:op_rec_div_HHO} (take in~\eqref{eq:def_opGRec_mat} $\bq\eqq q\matrice{I}_d$ with $q$ arbitrary in $\Pkd(T)$).
Moreover, proceeding as in the proof of Lemma~\ref{lem:grad_rec}(ii), one readily
verifies that 
\begin{equation}
\opGRec(\hIkT(\vv))=\PikTv(\GRAD \vv),
\end{equation}
for all $\vv\in \bH^1(T)$,
where $\PikTv$ is the $\bL^2$-orthogonal projection onto $\Pkd(T;\Rdd)$
and the local reduction operator $\hIkT:\vecteur{H}^1(T)\to \hVkT$ is defined 
in~\eqref{eq:def_hIkT}.
Finally, the stabilization operator $\opStv :
\hVkT \to \vVkdT$ is defined in~\eqref{eq:def_SKk_HHO_elas}
as for linear elasticity (another possibility is to define it as the vector-valued 
version of the one used for the Poisson model problem in \eqref{eq:def_SKk_HHO}).
Adapting the arguments
in the proof of Lemma~\ref{lem:stab_HHO} leads to the following result.

\bLem[Stability]  \label{lem:stab_HHO_hyperelas}
Equip $\hVkT$ with the seminorm
$\snorme[\hVkT]{\hvT}^2 \eqq
\normem[\T]{\GRAD \vvT}^2 + h_{T}^{-1}\normev[\dT]{\vvT-\vdT}^2$.
There are $0<\alpha\le\omega<+\infty$ such that for all $T\in\Th$ and all $\hvT\in \hVkT$,
\begin{align}
\alpha \snorme[\hVkT]{\hvT}^2 &\le
{\normem[\T]{\opGRec(\hvT))}^2}
+ h_{T}^{-1}\normev[\dT]{\opStv(\hvT)}^2 \le
\omega \snorme[\hVkT]{\hvT}^2.
\end{align}
\eLem

For all $\hvT\in \hVkT$, we reconstruct the deformation gradient in every mesh cell $T\in\Th$ as
\begin{equation}
\Fdef_T(\hvT) \eqq \matrice{I}_d + \opGRec(\hvT).
\end{equation}
For all $\hvh\in \hVkh$, we define the global reconstructions 
$\Fdef_\calT(\hvh)$ and $\opGRech(\hvh)$ such that 
\begin{equation} \label{eq:glob_recs_stab}
\Fdef_\calT(\hvh)_{|T} \eqq \Fdef_T(\hvT), \quad 
\opGRech(\hvh)_{|T} \eqq \opGRec(\hvT), \quad \forall T\in\Th,
\end{equation} 
so that $\Fdef_\calT(\hvh) = \matrice{I}_d + \opGRech(\hvh)$.
Recalling the linear form $\ell_h(\hvh)\eqq \psv[\Dom]{\loadext}{\vv_\calT}+\psv[\Bn]{\Gn}{\vv_\calF}$ from~\eqref{eq:weak_HHO_elas} and the stabilization bilinear form $s_h$ from~\eqref{eq:def_sh_elas}, we define the discrete energy functional 
$\engy_h : \hVkhd \rightarrow \Reel$ such that 
(compare with \eqref{eq:def_calE})
\begin{equation}\label{discrete_energy_stab}
\engy_h (\hvh) \eqq \int_{\Dom} \Psi(\Fdef_\calT(\hvh))\dV 
+ \beta_0\mu s_h(\hvh,\hvh) - \ell_h(\hvh),
\end{equation}
with a non-dimensional positive weight $\beta_0>0$.
For linear elasticity, a simple choice (considered in \eqref{eq:elast_engy_HHO})
is $\beta_0=1$; 
for finite deformations of hyperelastic materials, the choice of 
$\beta_0$ is further discussed in Remark~\ref{rem:choice_beta_hyper}.
The discrete problem consists in seeking the stationary points in $\hVkhd$ of the discrete energy functional $\engy_h$: 
Find $\huh\in\hVkhd$ such that (compare with~\eqref{eq:Euler_Lag_hypereal})
\begin{equation}\label{eq_points_critique_stab}
\psm[\Dom]{\PK(\Fdef_\calT(\huh))}{\opGRech(\hwh)} 
+ 2\beta_0\mu s_h(\huh,\hwh) = \ell_h(\hwh), \quad \forall \hwh \in \hVkhz.
\end{equation}

As for the linear elasticity problem, the discrete problem~\eqref{eq_points_critique_stab} can be reformulated in terms of equilibrated tractions. 
Let $\optStv : \vVkdT\to \vVkdT$ be defined in~\eqref{eq:def_optStv}
and let $\optStv^{*}: \vVkdT\to \vVkdT$ be its adjoint.
For all $\hvh\in \hVkh$, we define numerical tractions 
$\vecteur{T}_{\dT}(\hvT)\in \vVkdT$ 
at the boundary of every mesh cell $T\in\calT$ by setting
\begin{equation} \label{eq:def_num_traction_hyper1}
\vecteur{T}_{\dT}(\hvT) \eqq -\PikTv(\PK(\Fdef_T(\hvT)))_{|\dT} \nT
+ 2\beta_0\mu h_T^{-1}(\optStv^*\circ \optStv) (\vv_{T|\dT}-\vdT).
\end{equation}
A direct adaptation of the proof of Proposition~\ref{prop:flux_recovery}
establishes the following result.

\bProp[Rewriting with tractions] \label{prop:flux_recovery_hyper1}
Let $\huh\in \hVkhd$ solve~\eqref{eq_points_critique_stab} and 
let the tractions $\vecteur{T}_{\dT}(\huT)\in \vVkdT$
be defined as in~\eqref{eq:def_num_traction_hyper1} for all $T\in\calT$. 
The following holds:\\
\textup{(i)} Equilibrium at every mesh interface $F=\partial T_-\cap \partial T_+ \cap H_F\in \calFi$ and at every Neumann boundary face $F=\partial T_-\cap \Bn\cap H_F \in \calFbN$: \eqref{eq:equil_traction} and~\eqref{eq:equil_Neumann} hold true.
\\
\textup{(ii)} Balance with the source term in every mesh cell $T\in\calT$: For all $\bq\in \vPkd(T)$,
\begin{equation} \label{eq:balance_traction_hyper1}
\psm[T]{\PK(\Fdef_T(\huT))}{\GRAD \bq} + \psv[\dT]{\vecteur{T}_{\dT}(\huT)}{\bq} = \psv[T]{\loadext}{\bq}.
\end{equation}
\textup{(iii)}
The above identities are an equivalent rewriting of~\eqref{eq_points_critique_stab} that fully characterizes any HHO solution $\huh\in \hVkhd$.
\eProp

\bRem[Literature]
HHO methods for hyperelastic materials undergoing finite deformations were
introduced in \cite{AbErPi:18}, see also \cite{BoDPS:17} for nonlinear
elasticity and small deformations. HDG methods for nonlinear elasticity
were developed in \cite{Soon2008,Nguyen2012a,KaLeC:2015},
discontinuous Galerkin methods in \cite{Eyck2006,Noels2006,Eyck2008},
gradient schemes in \cite{Droniou2015},
virtual element methods in \cite{Chi2017,Wriggers2017},
and a (low-order) hybrid dG method with
conforming traces in \cite{Wulfinghoff2017}.
\eRem

\subsection{The unstabilized HHO method}
\label{sec:uHHO}

In nonlinear elasticity, the use
of stabilization can lead to numerical difficulties since it is not
clear beforehand how large the stabilization parameter ought to be; see 
\cite{Eyck2008,Chi2017} for related discussions. Moreover, 
\cite[Sect.~4]{KaLeC:2015} presents an example where spurious
solutions can appear in an HDG discretization if the stabilization
parameter is not large enough. Motivated by these observations, we
present in this section the unstabilized HHO method devised in 
\cite{AbErPi:18}. We assume for simplicity that the mesh is simplicial.

Let $T\in\Th$ and let $k\ge1$. Let us set
\begin{equation}
\RT_d^k(T;\Rdd)\eqq \Pkd(T;\Rdd) \oplus (\tilde\Poly^{k}_d(T;\Rd) \otimes \bx),
\end{equation} 
where $\tilde\Poly^{k}_d(T;\Rd)$ is the space composed of the restriction to $T$ of $\Rd$-valued $d$-variate homogeneous polynomials of degree $k$.
The local gradient reconstruction operator $\opGRecRT:\hVkT \to \RT_d^k(T;\Rdd)$ 
is such that for all $\hvT\in\hVkT$,
$\opGRecRT(\hvT)\in \RT_d^k(T;\Rdd)$ is uniquely determined by the equations
\begin{equation} \label{eq:def_opGRec_mat_RT}
(\opGRecRT(\hvT),\bq)_{\bL^2(T)}=-(\vvT,\DIV\bq)_{\bL^2(T)}+(\vv_\dT,
\bq\nT)_{\bL^2(\dT)}, \quad \forall \bq\in \RT_d^k(T;\Rdd).
\end{equation}
In practice, the lines of $\opGRecRT(\hvT)$ can be computed separately 
by inverting the mass matrix associated with the space $\RT_d^k(T;\Rd)$. Notice that 
the size of the
linear system resulting from~\eqref{eq:def_opGRec_mat_RT} is larger 
than the one resulting from~\eqref{eq:def_opGRec_mat}; the respective sizes are
$d{k+d\choose d}+{k+d-1\choose d-1}$ vs.~${k+d\choose d}$,
\eg 15 vs.~4 for $d=3$, $k=1$ and 36 vs.~10 for $d=3$, $k=2$.
Adapting the arguments of the proof of Lemma~\ref{lem:grad_rec_RT}(ii)
to the tensor-valued case shows that
\begin{equation}
\opGRecRT(\hIkT(\vv))=\PikTv(\GRAD \vv),
\end{equation}
for all $\vv\in \bH^1(T)$,
where $\PikTv$ is the $\bL^2$-orthogonal projection onto $\Pkd(T;\Rdd)$
and the local reduction operator $\hIkT:\vecteur{H}^1(T)\to \hVkT$ is defined 
in~\eqref{eq:def_hIkT}.
Adapting the arguments of the proof of Lemma~\ref{lem:grad_rec_RT}(iii) for the
lower bound and proceeding as usual for the upper bound leads to the following result.

\bLem[Stability]  \label{lem:stab_HHO_hyperelasRT}
Recall the seminorm $\snorme[\hVkT]{\SCAL}$ from Lemma~\ref{lem:stab_HHO_hyperelas}.
There are $0<\alpha\le\omega<+\infty$ such that
$\alpha \snorme[\hVkT]{\hvT}^2 \le
{\normem[\T]{\opGRec(\hvT))}^2} \le
\omega \snorme[\hVkT]{\hvT}^2$ for all $T\in\Th$ and all $\hvT\in \hVkT$.
\eLem

For all $\hvT\in \hVkT$, we now reconstruct the deformation gradient in every mesh cell $T\in\Th$ as
\begin{equation}
\FdefRT(\hvT) \eqq \matrice{I}_d + \opGRecRT(\hvT).
\end{equation}
For all $\hvh\in \hVkh$, we define the global reconstructions 
$\FdefRTh(\hvh)$ and $\opGRecRTh(\hvh)$ such that 
\begin{equation} \label{eq:glob_recs_unstab}
\FdefRTh(\hvh)_{|T} \eqq \FdefRT(\hvT), \quad
\opGRecRTh(\hvh)_{|T} \eqq \opGRecRT(\hvT), \quad \forall T\in\Th,
\end{equation} 
so that $\FdefRTh(\hvh) = \matrice{I}_d + \opGRecRTh(\hvh)$.
Recalling the linear form $\ell_h(\hvh)\eqq \psv[\Dom]{\loadext}{\vv_\calT}+\psv[\Bn]{\Gn}{\vv_\calF}$ (see~\eqref{eq:weak_HHO_elas}), we define the discrete energy functional 
$\engy_h : \hVkhd \rightarrow \Reel$ such that 
(compare with \eqref{discrete_energy_stab}; we use the same notation for simplicity)
\begin{equation}\label{discrete_energy_unstab}
\engy_h(\hvh) \eqq \int_{\Dom} \Psi(\FdefRTh(\hvh))\dV 
- \ell_h(\hvh).
\end{equation}
The discrete problem consists in seeking the stationary points in $\hVkhd$ of the discrete energy functional $\engy_h$: 
Find $\huh\in\hVkhd$ such that (compare with~\eqref{eq:Euler_Lag_hypereal})
\begin{equation}\label{eq_points_critique_unstab}
\psm[\Dom]{\PK(\FdefRTh(\huh))}{\opGRecRTh(\hwh)} 
= \ell_h(\hwh), \quad \forall \hwh \in \hVkhz.
\end{equation}

The discrete problem~\eqref{eq_points_critique_unstab} can be reformulated in terms of equilibrated tractions. 
For all $\hvh\in \hVkh$, we can define numerical tractions 
$\vecteur{T}_{\dT}(\hvT)\in \vVkdT$ 
at the boundary of every mesh cell $T\in\calT$ by setting
\begin{equation} \label{eq:def_num_traction_hyper2}
\vecteur{T}_{\dT}(\hvT) \eqq - \vecteur{\Pi}_T^{\textsc{rt}}(\PK(\Fdef_T(\hvT)))_{|\dT} \nT,
\end{equation}
where $\vecteur{\Pi}_T^{\textsc{rt}}$ denotes the $\bL^2$-orthogonal projection
onto $\RT_d^k(T;\Rdd)$.
A direct adaptation of the proof of Proposition~\ref{prop:flux_recovery}
establishes the following result.

\bProp[Rewriting with tractions] \label{prop:flux_recovery_hyper2}
Let $\huh\in \hVkhd$ solve~\eqref{eq_points_critique_unstab} and 
let the tractions $\vecteur{T}_{\dT}(\huT)\in \vVkdT$
be defined as in~\eqref{eq:def_num_traction_hyper2} for all $T\in\calT$. 
The following holds:\\
\textup{(i)} Equilibrium at every mesh interface $F=\partial T_-\cap \partial T_+ \cap H_F\in \calFi$ and at every Neumann boundary face $F=\partial T_-\cap \Bn\cap H_F \in \calFbN$: \eqref{eq:equil_traction} and~\eqref{eq:equil_Neumann} hold true.
\\
\textup{(ii)} Balance with the source term in every mesh cell $T\in\calT$: For all $\bq\in \vPkd(T)$,
\begin{equation} \label{eq:balance_traction_hyper2}
\psm[T]{\PK(\Fdef_T(\huT))}{\GRAD \bq} + \psv[\dT]{\vecteur{T}_{\dT}(\huT)}{\bq} = \psv[T]{\loadext}{\bq}.
\end{equation}
\textup{(iii)}
The above identities are an equivalent rewriting of~\eqref{eq_points_critique_unstab}. 
\eProp

\bRem[Divergence]
Notice that the identity~\eqref{eq:tr_G_opDiv} no longer holds if the
gradient is reconstructed using Raviart--Thomas polynomials. Instead, one 
only has $\PikTs(\trace(\opGRec(\hvT)))=\opDiv(\hvT)$ for all $\hvT\in\hVkT$,
indicating that a high-order perturbation may hamper robustness in the 
quasi-incompressible limit. So far, robustness was observed in the numerical
experiments. 
\eRem

\subsection{Nonlinear solver and static condensation}
\label{sec:implement}

The nonlinear problems \eqref{eq_points_critique_stab} and
\eqref{eq_points_critique_unstab} can be solved by using Newton's method. 
This requires evaluating the fourth-order elastic modulus $\moduletangent(\Fdef) \eqq
\partial_{\Fdef \Fdef}^2 \Psi(\Fdef)$. In particular, for Neohookean materials
(see~\eqref{NeoLaw}), we have
\begin{align}
\moduletangent(\Fdef) =& \mu  ( \matrice{I}  \, \overline{\otimes} \, \matrice{I} + \Fdef\mtr \underline{\otimes} \, \Fdef^{-1})  -\lambda   \ln J  \Fdef\mtr \underline{\otimes} \, \Fdef^{-1}  + \lambda  \Fdef\mtr\otimes \Fdef\mtr,
\end{align}
with $(\ba \otimes \bb)_{ijkl} \eqq \ba_{ij}
\bb_{kl}$, $(\ba \, \underline{\otimes} \, \bb)_{ijkl}
\eqq \ba_{il} \bb_{jk}$, and $(\ba \,
\overline{\otimes} \, \bb)_{ijkl} \eqq \ba_{ik} \bb_{jl}$, for all $1\le i,j,k,l\le d$. 
Let $i\ge0$ be the index of the Newton's iteration. Given an initial discrete
displacement $\huh^0\in\hVkhd$, one computes at each Newton's
iteration the incremental displacement $\hduh^i\in \hVkhz$ and
updates the discrete displacement as $\huh^{i+1} = \huh^i+\hduh^i$. 
The linear system of equations to be solved is
\begin{align}\label{eq_stiffness_matrix}
& \psm[\Dom]{\moduletangent(\Fdef_\calT^* (\huh^i)): \opGRech^*(\hduh^i)}{\opGRech^*(\hwh)} +  2\beta_0\mu s_h(\hduh^i,\hwh) = -R_h^i( \hwh),
\end{align}
for all $\hwh\in\hVkhz$, with the residual term
\begin{align}
 R_h^i(\hwh) & \eqq \psm[\Dom]{\PK(\Fdef_\calT^*(\huh^i))}{\opGRech^*(\hwh)} + 2\beta_0\mu s_h(\huh^i,\hwh) - \ell_h(\hwh),
\end{align}
where $\beta_0>0$, $\Fdef_\calT^*=\Fdef_\calT$, and $\opGRech^*=\opGRech$ (see~\eqref{eq:glob_recs_stab}) in the stabilized case and
$\beta_0=0$, $\Fdef_\calT^*=\FdefRTh$, and $\opGRech^*=\opGRecRTh$ (see~\eqref{eq:glob_recs_unstab}) in the unstabilized case.
We notice that in both cases the cell unknowns can be eliminated locally by using static condensation at each Newton's iteration~\eqref{eq_stiffness_matrix}.

\begin{remark}[Choice of $\beta_0$] 
To our knowledge, there is no general theory on the choice of $\beta_0$ in the case of finite deformations of hyperelastic materials. Following ideas developed in \cite{Eyck2008, Eyck2008a,BdVLM:15}, one can consider to take (possibly in an adaptive fashion) the largest eigenvalue (in absolute value) of the elastic modulus $\moduletangent$. This choice introduces additional nonlinearities to be handled by Newton's method, and may require some relaxation. 
Another possibility discussed in \cite{Chi2017} for virtual element methods is based
on the trace of the Hessian of the isochoric part of the strain-energy
density. Such an approach bears similarities with the classic
selective integration for FEM, and for Neohookean materials, 
this choice implies to take $\beta_0 = 1$. Finally,
we mention that too large values of the stabilization parameter $\beta_0$ 
can deteriorate the condition number of the stiffness matrix and can cause 
numerical instabilities in Newton's method. \label{rem:choice_beta_hyper}%
\end{remark}

\section{Numerical examples}\label{sec_numconv}
We present two examples that are close to industrial simulations: one for linear elasticity, a perforated strip subjected to uniaxial extension, and one for hyperelasticity, the pinching of a pipe.
The material parameters are $\mu = 23.3~\MPa$ and  $\lambda = 11650~\MPa$, which correspond to a Young modulus $E=70~\MPa$ and a Poisson ratio $\nu=0.499$. Both simulations are in the quasi-incompressible regime to show the robustess of HHO methods.

We first consider a strip of width $2L=200~\mm$ and height $2H=360~\mm$. The strip is perforated in its middle by a circular hole of radius $R=50~\mm$, and is subjected to a uniaxial extension $\delta=5~\mm$ at its top and bottom ends. We consider the linear elasticity model. For symmetry reasons, only one quarter of the strip is discretized. The Euclidean norm of the displacement field and the trace of stress tensor are plotted in Fig.~\ref{fig::strip} for $k=1$ on a mesh composed of 536 cells with hanging nodes. There is no sign of volumetric locking, thereby comfirming the robustness of HHO methods in the quasi-incompressible limit.
\begin{figure}[htbp]
    \centering
    \subfloat[Euclidean norm of displacement ($\mm$)]{
        \centering
        \includegraphics[scale=0.201, trim= -20 17 70 0, clip=true]{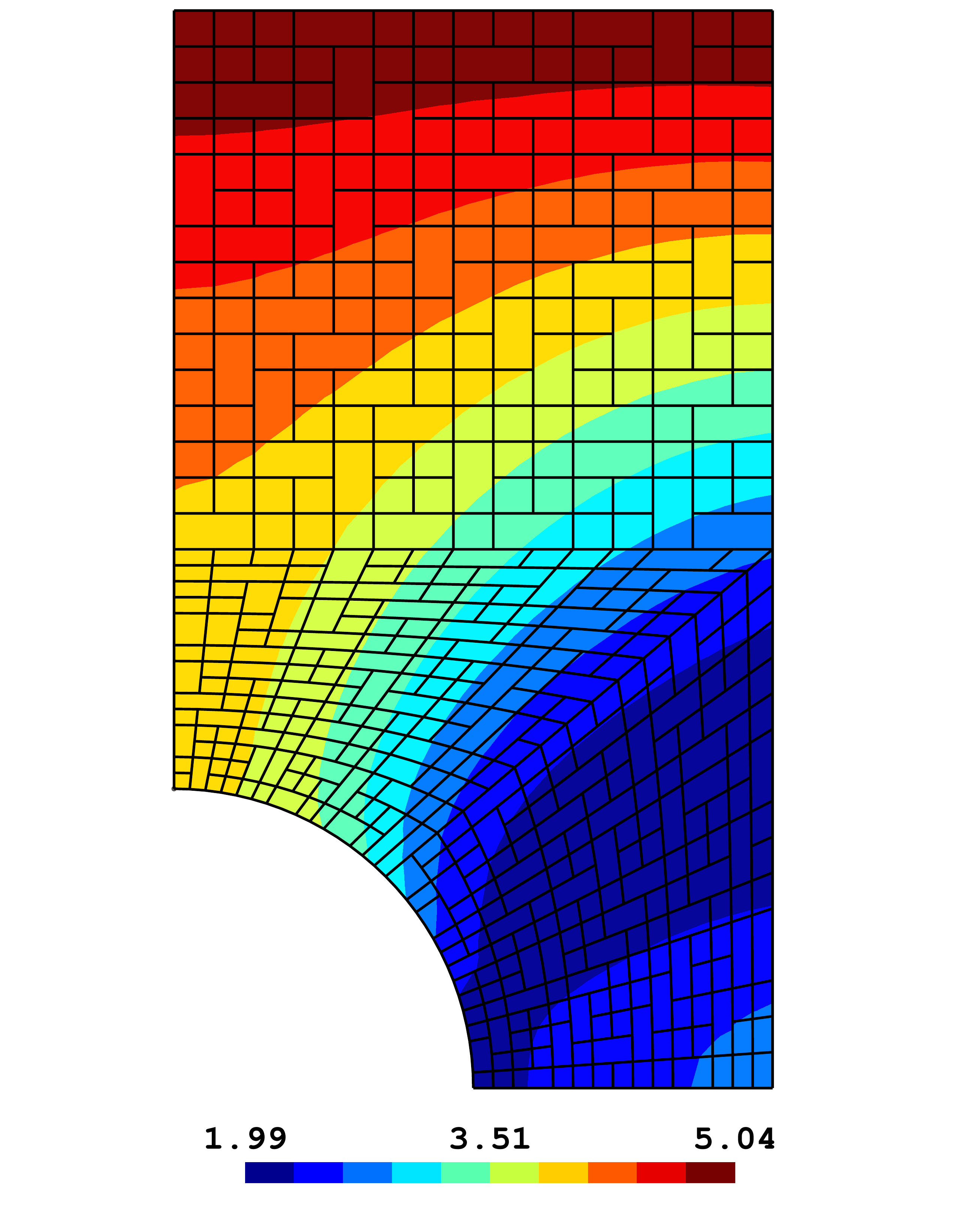}
  }
    ~ 
    \subfloat[Trace of stress tensor ($\MPa$)]{
        \centering
       \includegraphics[scale=0.225, trim= 00 20 120 30, clip=true]{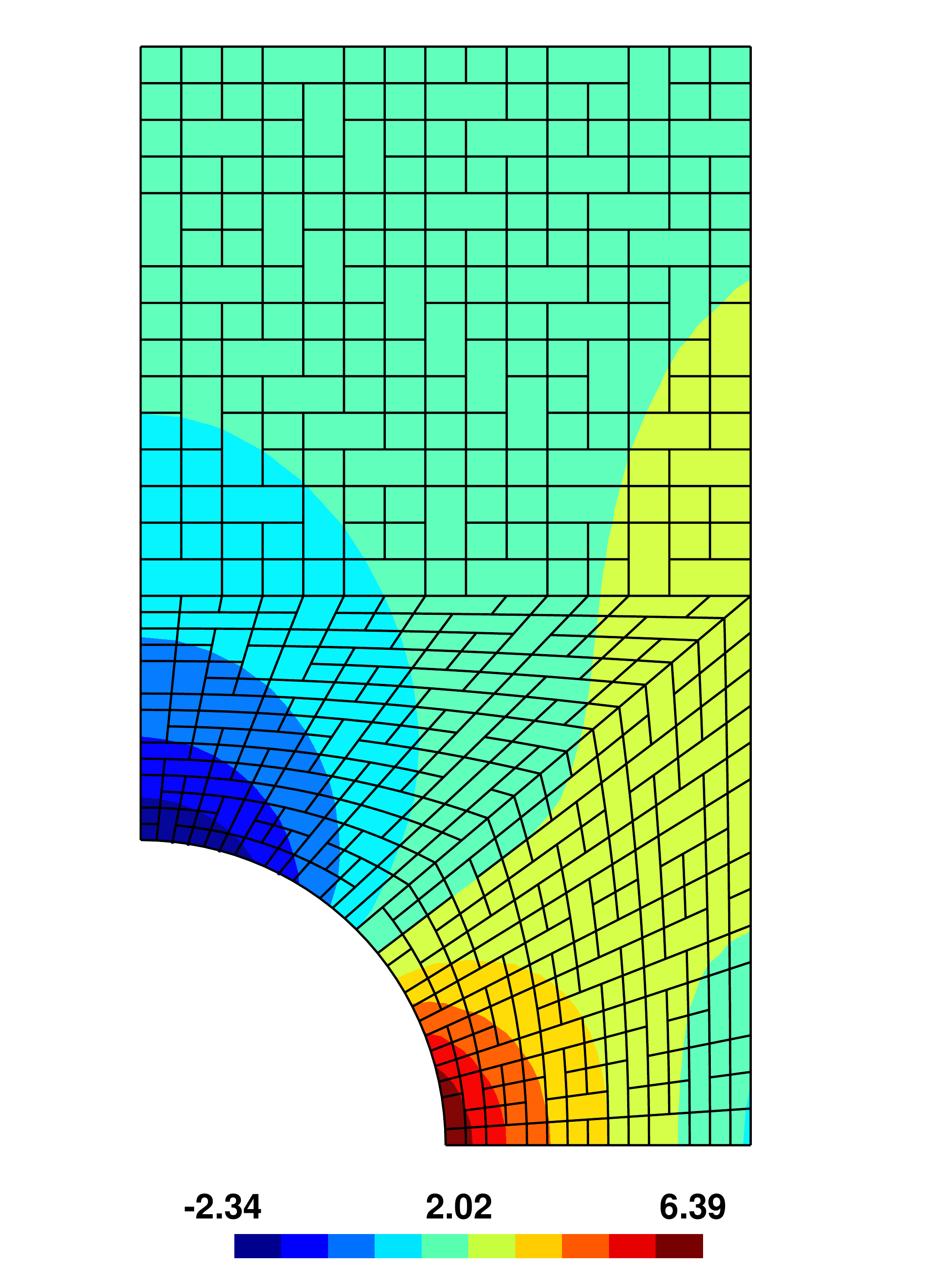}
    }
    \caption{Perforated strip: Euclidean norm of the displacement field (left) and trace of stress tensor (right); polynomial degree $k=1$; the mesh with hanging nodes is treated as a polygonal mesh and is composed of 536 cells (quadrilaterals, pentagons and hexagons).}
    \label{fig::strip}
\end{figure}

%
The second numerical example is the pinching of a pipe due to external forces.  The pipe has an outer radius of $24~\mm$ and an inner radius of $23~\mm$ (the thickness is equal to $1~\mm$) and a length equal to $100~\mm$. One end is clamped, the other end and the inner surface are free, and the outer surface is subjected to a compression force of 0.01~MPa, oriented downwards in the upper half and upwards in the lower half of the outer surface.  Since the geometry as well as the boundary conditions are symmetric, it is sufficient to model one half of the pipe in finite deformations. The mesh is composed of 40{,}500 tetrahedra. The von Mises stress is plotted in Fig.~\ref{fig::strip} for $k=1$ and different gradient reconstructions (using full polynomials or Raviart--Thomas polynomials with $k=1$) on the deformed configuration. Both simulations do not present any sign of volumetric locking.
\begin{figure}[htbp]
    \centering
    \subfloat[$\opGRech^*=\opGRech$ with stabilization]{
        \centering
        \includegraphics[scale=0.25, trim= 25 30  875 60, clip=true]{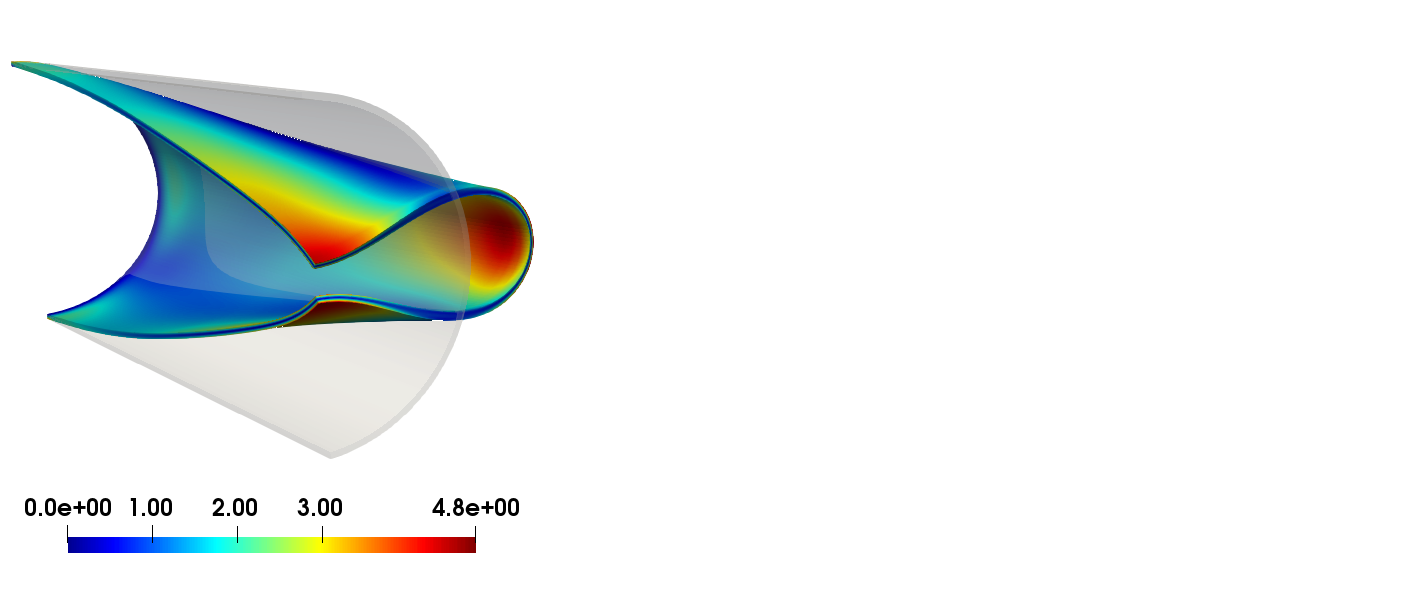}
  }  
    \subfloat[$\opGRech^*=\opGRecRTh$ without stabilization]{
        \centering
       \includegraphics[scale=0.25, trim= 25 30  875 60, clip=true]{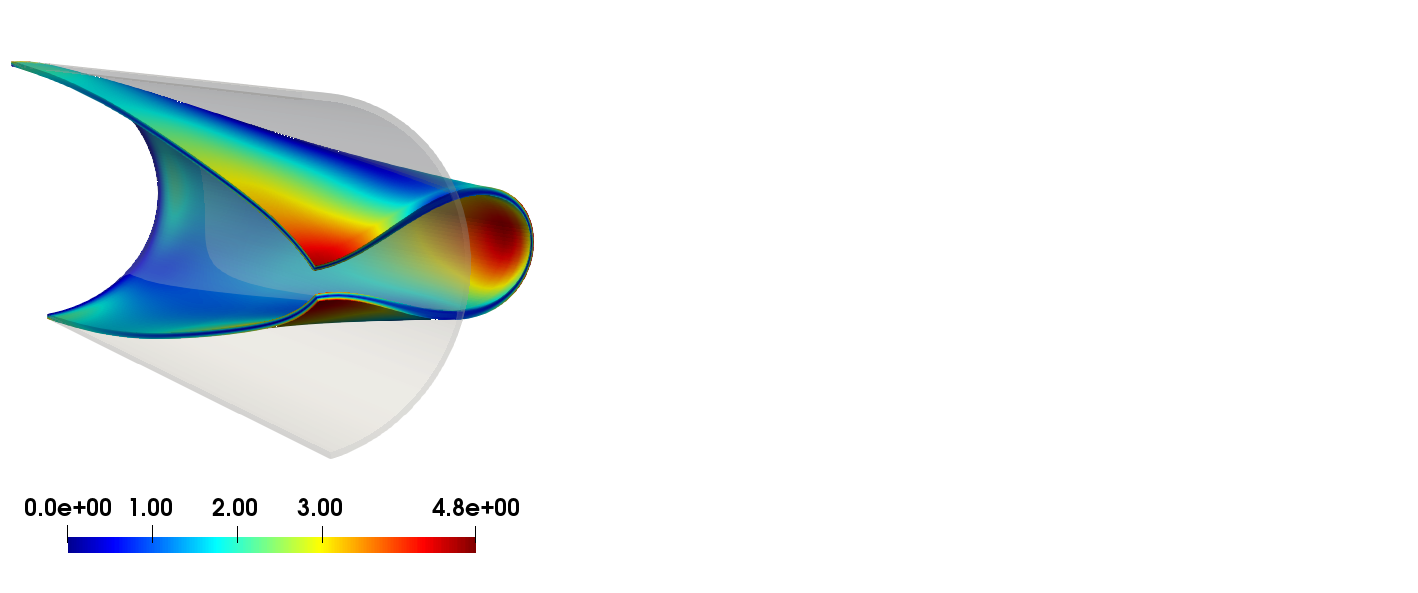}
    }
    \caption{Pinching of a pipe: von Mises stress (in $\MPa$) on the deformed configuration for $k=1$ and different gradient reconstructions. }
    \label{fig::pipe}
\end{figure}

\chapter{Elastodynamics}
\label{chap:wave}

The goal of this chapter is to show how the HHO method can be used
for the space semi-discretization of the elastic wave equation. 
For simplicity, we restrict the scope to 
media undergoing infinitesimal deformations and governed by
a linear stress-strain constitutive relation. We consider first 
the second-order formulation in time and then the mixed formulation
leading to a first-order formulation in time. The time discretization
is realized, respectively, by means of Newmark schemes and diagonally-implicit or explicit
Runge--Kutta schemes. Interestingly, considering the mixed-order HHO
method is instrumental to devise explicit Runge--Kutta schemes.
HHO methods for acoustic and elastic wave propagation
were developed in \cite{BDE:21,BDES:21}, and HDG methods for these
problems were studied in \cite{NgPeC:11,StNPC:16,SCNPC:17,CFHJSS:18}.

\section{Second-order formulation in time}
\label{sec:wave_2nd}

In this section, we consider the second-order formulation in time of the 
elastic wave equation. We refer the reader to Sect.~\ref{sec:model-linear_elas}
for a description of the linear elasticity model in the case of static problems. 
Let $J\eqq (0,T_{\textrm{f}})$ with final time $T_{\textrm{f}}>0$
be the time interval, and let $\Dom$ be an open, bounded, connected, Lipschitz
subset of $\Real^d$ in space dimension $d\ge2$.
The elastic wave equation reads as follows:
\begin{equation}
\rho \partial_{tt}\bu -
\GRAD\SCAL \bsigma(\strain(\vu)) = \loadext \quad\text{in $J\times\Dom$},
\label{eq:elasto_2nd}
\end{equation}
where $\rho$ is the material density,
$\loadext$ the body force, and $\bu$ the displacement field.
The stress tensor $\bsigma(\strain(\vu))$ depends
on the displacement field by means of the linearized strain tensor $\strain\eqq
\strain(\bu)
\eqq \frac12(\GRAD \bu + \GRAD\bu\tr)$ as follows:
\begin{equation} \label{eq:def_sigma_elasto}
\bsigma(\strain) \eqq \polA \strain \eqq 
2\mu \strain + \lambda \trace(\strain) {\matrice{I}}_{d},
\end{equation}
where $\polA$ is the fourth-order stiffness tensor,
$\mu$ and $\lambda$ are the Lam\'e parameters,
and $\matrice{I}_d$ the identity tensor. 
We assume that the coefficients $\rho$, $\mu$, and $\lambda$ are piecewise
constant on a partition of $\Omega$ into a finite collection of polyhedral
subdomains, that $\rho$ and $\mu$ take 
positive values, and that $\lambda$ takes nonnegative values.
The wave equation~\eqref{eq:elasto_2nd} describes the propagation of different
types of elastic waves in the medium. In particular,  
the speeds of P- and S-waves are $c_{\mathrm{P}} \eqq 
\sqrt{\frac{\lambda+2\mu}{\rho}}$ and $c_{\mathrm{S}} \eqq 
\sqrt{\frac{\mu}{\rho}}$.

The wave equation~\eqref{eq:elasto_2nd} is subjected to the initial and boundary
conditions 
\begin{equation} \label{eq:IC_BC_wave_2nd}
\bu(0) =\bu_{0}, \ \partial_t\bu(0) = \bv_{0}\; \text{in $\Dom$}, \qquad
\bu=\vzero\;\text{on $J\times\front$},
\end{equation}
where the homogeneous Dirichlet boundary condition is chosen for simplicity.
Assuming that 
$\loadext\in C^0(\oJ;\bL^2(\Omega))$ with $\oJ\eqq[0,T_{\textrm{f}}]$, $\bu_0\in \bH^1_0(\Omega)$,
and $\bv_0\in \bH^1_0(\Omega)$, and that
$\bu\in C^0(\oJ;\bH^1_0(\Omega))\cap H^2(J;\bL^2(\Omega))$, we have for a.e.~$t\in J$,
\begin{equation} \label{eq:weak_elasto}
(\partial_{tt}\vu(t),\vw)_{\bL^2(\rho;\Omega)} + a(\vu(t),\vw) = (\loadext(t),\vw)_{\bL^2(\Omega)},
\quad \forall \vw\in \bH^1_0(\Omega),
\end{equation}
with the bilinear form 
\begin{equation}
a(\bv,\bw)\eqq (\stress(\strain(\bv)),\strain(\bw))_{\bL^2(\Omega)}
= (\strain(\bv),\strain(\bw))_{\bL^2(2\mu;\Omega)} 
+ (\GRAD\SCAL\bv,\GRAD\SCAL\bw)_{L^2(\lambda;\Omega)}.
\end{equation}
Here, for a weight function $\phi\in L^\infty(\Omega)$ taking nonnegative 
values, we used the 
notation $\|v\|_{L^2(\phi;\Omega)} \eqq \|\phi^{\frac12}v\|_{L^2(\Omega)}$ for all
$v\in L^2(\Omega)$, and a similar notation for vector-valued fields
in $\bL^2(\Omega)$ (this defines a norm if $\phi$ is
uniformly bounded from below away from zero).

An important property of the elastic wave equation is energy balance. 
The time-dependent energy associated with the weak solution is defined for all $t\in \oJ$ as
\begin{equation} \label{eq:def_engy_wave}
\engy(t) \eqq \frac12 \|\partial_t\vu(t)\|_{\bL^2(\rho;\Omega)}^2 + 
\frac12 \|\strain(\bu(t))\|_{\bL^2(2\mu;\Omega)}^2
+ \frac12 \|\GRAD\SCAL\bu(t)\|_{L^2(\lambda;\Omega)}^2.
\end{equation}
Assuming $\bu\in C^1(\oJ;\bH^1_0(\Omega))$ and
testing~\eqref{eq:weak_elasto} against $\bw\eqq \partial_t\vu(t)$ gives $\frac{d}{dt}\engy(t)=(\loadext(t),\partial_t\vu(t))_{\bL^2(\Omega)}$ for all $t\in J$. Integrating in time over $(0,t)$ leads to the energy balance equation 
\begin{equation} \label{energy_bal_wave}
\engy(t)=\engy(0) + \int_0^t (\loadext(s),\partial_t\vu(s))_{\bL^2(\Omega)}\ds,
\end{equation}
where $\engy(0)$ can be evaluated from the 
initial condition~\eqref{eq:IC_BC_wave_2nd}. In the absence of body forces,
\eqref{energy_bal_wave} implies energy conservation, \ie $\engy(t)=\engy(0)$
for all $t\in\oJ$.

\subsection{HHO space semi-discretization}

We consider the discrete setting described in Sect.~\ref{sec:lin_elas_setting}.
In particular, $\Th$ is a mesh of $\Dom$ belonging to a shape-regular mesh sequence,
and we assume that $\Dom$ is a polyhedron so that the mesh covers $\Dom$ exactly. 
Moreover, we assume that the mesh is compatible with the above partition of $\Dom$
regarding the material properties, so that the parameters $\rho$, $\mu$,
and $\lambda$ are piecewise constant on the mesh.

To allow for a bit more generality (this will be handy when 
studying the first-order formulation in time in Sect.~\ref{sec:wave_1st}),
we consider either the equal-order case or the mixed-order case for the
cell and the face unknowns in the HHO method. 
Letting $k\ge1$ denote the degree of the face unknowns, 
the cell unknowns can have degree $k'\eqq k$ (equal-order) or $k'\eqq k+1$
(mixed-order). Only the equal-order case
was considered in Sect.~\ref{sec:lin_elas_setting} for the static problem, 
and we refer the reader to
Sect.~\ref{sec:variants_cell} for a study of the mixed-order HHO method
applied to the Poisson model problem.

We use a unified notation to cover both cases, and for simplicity we use only the
superscript $k$ in the HHO spaces composed of polynomial pairs.  
For every mesh cell $T\in\Th$, we set
\begin{equation}
\hVkT \eqq \vPkprd(T) \times \vVkdT, \qquad
\vVkdT \eqq \bigtimes_{F\in\FT} \vPkF(F),
\label{eq:def_VkT_wave}
\end{equation}
with $\vPkprd(T)\eqq \Pkprd(T;\Real^d)$ and $\vPkF(F)\eqq \PkF(F;\Real^d)$;
see Figure~\ref{fig_HHO_dofs_wave}.
A generic element in $\hVkT$ is denoted by $\hvT\eqq (\vT,\vdT)$.
The HHO space is then defined as
\begin{equation} \label{eq:def_hVkh_wave}
\hVkh \eqq \Vkprh \times \VkFh,
\qquad
\Vkprh \eqq \bigtimes_{T\in\calT} \vPkprd(T),
\qquad
\VkFh \eqq \bigtimes_{F\in\calF} \vPkF(F).
\end{equation}
A generic element in $\hVkh$ is denoted by $\hvh\eqq (\vTh,\vFh)$ with
$\vTh\eqq (\vv_T)_{T\in\calT}$ and $\vFh\eqq (\vv_F)_{F\in\calF}$, and we 
localize the components of $\hvh$ associated with a mesh cell $T\in\Th$ 
and its faces by using the notation
$\hvT\eqq \big(\vv_T,\vv_{\dT}\eqq (\vv_F)_{F\in\FT}\big)\in \hVkT$.
The local HHO reduction operator  
$\hIkT:\vecteur{H}^1(T)\to \hVkT$ is such that
$\hIkT(\vv)\eqq (\vecteur{\Pi}^{k'}_T(\vv),\PikdTv(\vv_{|\dT}))$,
and its global counterpart $\Ikh : \vecteur{H}^1(\Dom) \to \hVkh$ is such that
$\Ikh(\vv)\eqq (\vecteur{\Pi}_\calT^{k'}(\vv),\PikcFv(\vv_{|\calF}))$.

\begin{figure}
    \centering
     \subfloat[equal-order]{
        \centering
        \includegraphics[scale=0.45]{Figures/Fig_inconnues_HHO_penta_k1.pdf}
  }
    \qquad 
    \subfloat[mixed-order]{
        \centering
	    \includegraphics[scale=0.45]{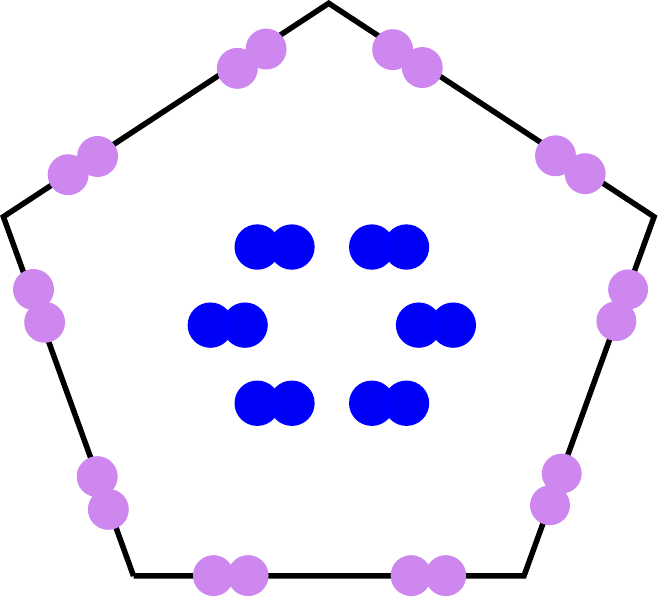}
    }
\caption{Face (violet) and cell (blue) unknowns in $\hVkT$ in a pentagonal cell ($d=2$) for $k=1$; left: equal-order case, right: mixed-order case (each dot in a pair represents one basis function associated with one Cartesian component of the displacement).}
\label{fig_HHO_dofs_wave}
\end{figure}

To reconstruct a strain tensor in every mesh cell $T\in\Th$, we can consider
the displacement reconstruction operator $\opDep : 
\hVkT \to \bpolP_d^{k+1}(T)$ defined in~\eqref{eq:def_Rec_HHO_elas}-\eqref{eq:def_Rec_HHO_mean_elas}. However, to highlight another possibility which is relevant in the case
of nonlinear materials (see, \eg Sect.~\ref{sec:HHO_hyperelasticity}), we consider here
a (symmetric-valued) strain reconstruction operator 
$\opStr:\hVkT \to \polP_d^{k}(T;\Msym)$  
such that for all $\hvT\in \hVkT$,
\begin{equation} \label{eq:def_strain_rec_wave}
(\opStr(\hvT),\bq)_{\bL^2(T)} = -(\vT,\GRAD\SCAL\bq)_{\bL^2(T)} + 
(\vdT,\bq\bn_T)_{\bL^2(\dT)}, 
\end{equation}
for all $\bq \in \polP_d^{k}(T;\Msym)$.
Notice that $\opStr(\hvT)$ can be evaluated componentwise by inverting the mass matrix
associated with a chosen basis of the scalar-valued polynomial space
$\Pkd(T)$. 
Recalling~\eqref{eq:def_SKk_HHO_elas} (equal-order case) and \eqref{eq:stab_HDG}
(mixed-order case),
the local stabilization operator is defined as follows:
\begin{alignat}{2}
&\opStv(\hvT) \eqq \PikdTv\Big(\vv_{T|\dT}
-\vdT+\big((I-\PikTv)\opDep(\hvT)\big)_{|\dT}\Big),&\quad&\text{if $k'=k$},\label{eq:def_stab_equal_wave}\\
&\opStv(\hvT) \eqq \PikdTv(\vv_{T|\dT})
-\vdT,&\quad&\text{if $k'=k+1$},
\end{alignat}
so that in the equal-order case, the displacement reconstruction operator
$\opDep$ needs to be evaluated as well.
Adapting the proof of Lemma~\ref{lem:stab_HHO_elas}, one can show that 
there are $0<\alpha\le\omega<+\infty$ such that for all $T\in\Th$ and all $\hvT\in \hVkT$,
\begin{align}\label{eq:stab_prop_wave}
\alpha \snorme[\strain,T]{\hvT}^2 &\le
\normem[\T]{\opStr(\hvT)}^2
+ h_{T}^{-1}\normev[\dT]{\opStv(\hvT)}^2 \le
\omega \snorme[\strain,T]{\hvT}^2,
\end{align}
recalling the seminorm $\snorme[\strain,T]{\hvT}^2 \eqq
\normem[\T]{\strain(\vT)}^2 + h_{T}^{-1}\normev[\dT]{\vT-\vdT}^2$.

We define the global 
discrete bilinear form $a_h:\hVkh\times \hVkh\to\Real$ such that
$a_h(\hvh,\hwh)\eqq \sum_{T\in\calT_h} a_T(\hvT,\hwT)$ with
$a_T:\hVkT\times \hVkT\to \Real$ such that 
\begin{align}
a_T(\hvT,\hwT) \eqq {}& (\opStr(\hvT),\opStr(\hwT))_{\bL^2(2\mu;T)}
+ (\opDiv(\hvT),\opDiv(\hwT))_{L^2(\lambda;T)} \nonumber \\
&+ \tau_{\dT} (\opStv(\hvT),\opStv(\hwT))_{\bL^2(\dT)},
\end{align}
with the weight $\tau_{\dT}\eqq 2\mu_{|T} h_T^{-1}$,  
and where $\opDiv(\SCAL) \eqq \trace(\opStr(\SCAL))$
coincides with the divergence reconstruction operator 
defined in~\eqref{eq:op_rec_div_HHO}. 
We define the global strain reconstruction operator 
$\opStrh:\hVkh\to \bW_\calT\eqq \Pkd(\Th;\Msym)$ such that 
\begin{equation} \label{eq:def_global_strain}
\opStrh(\hvh)_{|T} \eqq \opStr(\hvT), \quad  \forall T\in\Th,\ \forall \hvh\in \hVkh,
\end{equation} 
and the global divergence reconstruction operator $\opDivh:\hVkh\to \Pkd(\Th;\Real)$
such that $\opDivh(\hvh)\eqq \trace(\opStrh(\hvh))$.
We also define the global stabilization bilinear form
$s_h$ on $\hVkh\times \hVkh$ such that 
\begin{equation}
s_h(\hvh,\hwh)\eqq \sum_{T\in\calT_h}\tau_{\dT} 
(\opStv(\hvT),\opStv(\hwT))_{L^2(\dT)}.
\end{equation} 
Letting $\VkFhz \eqq \{ \vFh \in \VkFh \tq \vF = \vecteur{0}, \ \forall F \in \calFb \}$
and $\hVkhz \eqq \Vkh \times \VkFhz$,
the space semi-discrete HHO scheme for the elastic wave equation is as follows:
Seek $\huh\eqq (\vu_\calT,\vu_\calF)
\in C^2(\oJ;\hVkhz)$ such that for all $t\in \oJ$,
\begin{equation} \label{eq:HHO_2nd_semi_elasto}
(\partial_{tt}\vu_\calT(t),\vw_\calT)_{\bL^2(\rho;\Omega)}
+ a_h(\huh(t),\hwh) = (\loadext(t),\vw_\calT)_{\bL^2(\Omega)},
\end{equation}
for all $\hwh\eqq (\vw_\calT,\vw_\calF)\in \hVkhz$. Notice that the acceleration
term only involves the cell components; the same remark applies to the body force
(as in the static case). 
Consistently with this observation, the initial conditions 
for~\eqref{eq:HHO_2nd_semi_elasto} only concern $\bu_\calT$ and read
$\vu_\calT(0) = \vecteur{\Pi}_\calT^{k'}(\vu_{0})$,
$\partial_t\vu_\calT(0) = \vecteur{\Pi}_\calT^{k'}(\vv_{0})$,
whereas the boundary condition is encoded in the fact that $\huh(t)\in
\hVkhz$ for all $t\in\oJ$. Notice that 
$\vu_\calF(0)\in \VkFhz$ is uniquely
determined by the equations $a_h((\vu_\calT(0),\vu_\calF(0)),(\vzero,\vw_\calF))=0$
for all $\vw_\calF \in \VkFhz$ 
with $\vu_\calT(0)$ specified by the initial condition.

The time-dependent energy associated with the space semi-discrete HHO problem \eqref{eq:HHO_2nd_semi_elasto} is defined for all $t\in \oJ$ as (compare with~\eqref{eq:def_engy_wave})
\begin{align} 
\engy_h(t) \eqq {}&\frac12 \|\partial_t\vu_\calT(t)\|_{\bL^2(\rho;\Omega)}^2 + 
\frac12 \|\opStrh(\huh(t))\|_{\bL^2(2\mu;\Omega)}^2
+ \frac12 \|\opDivh(\huh(t))\|_{L^2(\lambda;\Omega)}^2 \nonumber \\
&+ \frac12 s_h(\huh(t),\huh(t)).\label{eq:def_engy_wave_HHO_2nd}
\end{align}
Then, proceeding as in the continuous case, one shows that
\begin{equation} \label{energy_bal_wave_HHO_2nd}
\engy_h(t)=\engy_h(0) + \int_0^t (\loadext(s),\partial_t\vu_\calT(s))_{\bL^2(\Omega)}\ds,
\end{equation}
so that, in the absence of body forces,
\eqref{energy_bal_wave_HHO_2nd} implies again energy conservation, \ie $\engy_h(t)=\engy_h(0)$ for all $t\in\oJ$.

Let $N_\calT^{k'}\eqq \dim(\Vkprh)$ and $N_{\calF,0}^{k}\eqq \dim(\VkFhz)$.
Let $(\sfU_\calT(t),\sfU_\calF(t))\in \Real^{N_\calT^{k'}\times N_{\calF,0}^k}$ 
be the component vectors of the space semi-discrete solution 
$\huh(t)\in \hVkhz$ once bases 
$\{\bvarphi_i\}_{1\le i\le N_\calT^{k'}}$
and $\{\bpsi_j\}_{1\le j\le N_{\calF,0}^k}$ for $\Vkprh$ and $\VkFhz$, respectively, 
have been chosen. 
Let $\sfF_\calT(t)\in \Real^{N_\calT^{k'}}$ have components
given by $\sfF_i(t)\eqq (\bef(t),\bvarphi_i)_{\bL^2(\Omega)}$ for all $1\le i\le N_\calT^{k'}$ and all 
$t\in \oJ$. The algebraic realization of~\eqref{eq:HHO_2nd_semi_elasto} 
is as follows: For all $t\in \oJ$,
\begin{equation}
\begin{bmatrix}
\sfM_{\calT\calT} & 0\\
0 & 0
\end{bmatrix} \begin{bmatrix}
\partial_{tt}\sfU_\calT(t)\\
\bullet
\end{bmatrix} + \begin{bmatrix}
\sfA_{\mathcal{T}\mathcal{T}} & \sfA_{\mathcal{T}\mathcal{F}}\\
\sfA_{\mathcal{F}\mathcal{T}} & \sfA_{\mathcal{F}\mathcal{F}}
\end{bmatrix} \begin{bmatrix}
\sfU_\calT(t)\\
\sfU_\calF(t)
\end{bmatrix} = \begin{bmatrix}
\sfF_\calT(t)\\
0
\end{bmatrix},
\label{eq:HHO_semi_2nd_alg}
\end{equation}
with the mass matrix $\sfM_{\calT\calT}$ associated with the inner product
in $L^2(\rho;\Omega)$ and the cell basis functions,
whereas the symmetric positive-definite stiffness matrix with blocks
$\sfA_{\mathcal{T}\mathcal{T}}$, $\sfA_{\mathcal{T}\mathcal{F}}$,
$\sfA_{\mathcal{F}\mathcal{T}}$, $\sfA_{\mathcal{F}\mathcal{F}}$ is
associated with the bilinear form $a_h$ and the cell and face
basis functions. The bullet stands for $\partial_{tt}\sfU_\calF(t)$ which is irrelevant 
owing to the structure of the mass matrix.
The matrices $\sfM_{\calT\calT}$
and $\sfA_{\calT\calT}$ are block-diagonal, but this is not the case for the matrix
$\sfA_{\calF\calF}$ since the components $\vu_\calF(t)$ attached to the faces 
belonging to the same cell are coupled together through the strain reconstruction operator
(and the stabilization operator in the equal-order case).

\bRem[Error analysis]
The error analysis for the space semi-discrete problem is performed in \cite[Thm.~3.1\&3.2]{BDES:21} for the acoustic wave equation and can be extended to the elastic wave equation. Following the seminal ideas from \cite{Dupont:73,Baker:76}, the key idea is to exploit the approximation properties of the discrete solution operator in the static case (see Sect.~\ref{sec:analysis_elas} for linear elasticity) and use the stability properties of the wave equation in time. For brevity, we only mention that the energy-error $\|\partial_t\vu-\partial_t\vu_\calT\|_{L^\infty(J;\bL^2(\rho;\Dom))} + \|\strain(\vu)-\opStrh(\huh)\|_{L^\infty(J;\bL^2(2\mu;\Dom))}$ decays as $\calO(h^{k+1})$ if $\vu\in C^1(\oJ;\vecteur{H}^{k+2}(\Dom))$, and assuming full elliptic regularity pickup ($s=1$), $\|\vecteur{\Pi}^{k'}_\calT(\vu)-\vu_\calT\|_{L^\infty(J;\bL^2(\mu;\Dom))}$ decays as $\calO(h^{k+2})$ if additionally $\loadext\in C^1(\oJ;\vecteur{H}^k(\Dom))$.
\eRem

\subsection{Time discretization}

Let $(t^n)_{0\le n\le N}$ be the discrete time nodes with $t^0\eqq0$ and $t^N\eqq T_{\textrm{f}}$.
For simplicity, we consider a fixed time step $\Delta t\eqq \frac{T_{\textrm{f}}}{N}$.
A classical time discretization of~\eqref{eq:HHO_2nd_semi_elasto} relies on the
Newmark scheme with parameters $\beta$ and $\gamma$, which is 
second-order accurate in time, implicit if $\beta>0$, unconditionally stable 
if $\frac12\le \gamma \le 2\beta$ (the classical choice is 
$\gamma=\frac12$ and $\beta=\frac14$) and conditionally stable if $\frac12\le 
\gamma$ and $2\beta<\gamma$. 
We detail the implementation for $\beta>0$. The
Newmark scheme considers a
displacement, a velocity, and an acceleration at each time node, 
which are all hybrid unknowns, say
$\huh^n,\hvh^n,\hxh^n\in \hVkhz$. The scheme is initialized
by setting $\huh^0\eqq \Ikh(\bu_0)$, $\hvh^0\eqq \Ikh(\bv_0)$, and the 
initial acceleration $\hxh^0\in \hVkhz$ is defined by
solving $(\bx_\calT^0,\vw_\calT)_{\bL^2(\rho;\Omega)} + a_h(\huh^0,(\vw_\calT,\vzero))
=(\loadext(0),\vw_\calT)_{\bL^2(\Omega)}$ for all $\vw_\calT\in \Vkprh$ and
$a_h(\hxh^0,(\vzero,\vw_\calF))=0$ for all $\vw_\calF\in\VkFhz$.
Then, given $\huh^n,\hvh^n,\hxh^n$ from the previous time-step
or the initial condition, the HHO-Newmark scheme proceeds as follows: 
For all $n\in\{0,\ldots,N-1\}$,
\begin{enumerate}
\item Predictor step: 
$\huh^{*n}\eqq \huh^n+ \Delta t \hvh^n + \tfrac12\Delta t^2(1-2\beta)\hxh^n$,
$\hvh^{*n}\eqq \hvh^n+\Delta t(1-\gamma)\hxh^n$.
\item Linear solve to find the acceleration 
$\hxh^{n+1}\in \hVkhz$ such that for all $\hwh\in \hVkhz$,
\begin{equation} \label{eq:Newmark_b}
(\bx_\calT^{n+1},\vw_\calT)_{\bL^2(\rho;\Omega)} + \beta \Delta t^2 a_h(\hxh^{n+1},\hwh) = (\loadext(t^{n+1}),\vw_\calT)_{\bL^2(\Omega)} - a_h(\huh^{*n},\hwh).
\end{equation}
\item Corrector step: 
$\huh^{n+1}\eqq \huh^{*n} + \beta \Delta t^2 \hxh^{n+1}$,
$\hvh^{n+1}\eqq \hvh^{*n} + \gamma \Delta t \hxh^{n+1}$.
\end{enumerate}
The algebraic realization 
of \eqref{eq:Newmark_b} amounts to finding $(\sfX_\calT^{n+1},\sfX_\calF^{n+1})
\in \Real^{N_\calT^{k'}\times N_{\calF,0}^k}$ such that 
\begin{equation}
\left( \begin{bmatrix}
\sfM_{\calT\calT} & 0\\
0 & 0
\end{bmatrix} + \beta \Delta t^2 \begin{bmatrix}
\sfA_{\mathcal{T}\mathcal{T}} & \sfA_{\mathcal{T}\mathcal{F}}\\
\sfA_{\mathcal{F}\mathcal{T}} & \sfA_{\mathcal{F}\mathcal{F}}
\end{bmatrix}\right) \begin{bmatrix}
\sfX_\calT^{n+1}\\
\sfX_\calF^{n+1}
\end{bmatrix} = \begin{bmatrix}
\sfG_\calT^{n+1}\\
\sfG_\calF^{n+1}
\end{bmatrix},
\label{eq:Newmark_alg}
\end{equation}
with $\sfG_\calT^{n+1}\eqq \sfF_\calT^{n+1} - (\sfA_{\calT\calT}\sfU^{*n}_\calT+
\sfA_{\calT\calF}\sfU^{*n}_\calF)$, $\sfG_\calF^{n+1}\eqq - (\sfA_{\calF\calT}\sfU^{*n}_\calT+
\sfA_{\calF\calF}\sfU^{*n}_\calF)$, and $(\sfU^{*n}_\calT,\sfU^{*n}_\calF)$ are the
components of the predicted displacement $\huh^{*n}$. Since the matrix
$\sfM_{\calT\calT}+\beta \Delta t^2 \sfA_{\mathcal{T}\mathcal{T}}$ is block-diagonal,
static condensation can be applied to~\eqref{eq:Newmark_alg}, \ie
the cell acceleration unknowns can be eliminated locally, 
leading to a global transmission problem coupling only
the face acceleration unknowns.

An important property of the HHO-Newmark scheme is energy balance. 
For all $n\in\{0,\ldots,N\}$, we define the discrete energy
\begin{align} 
\engy^n \eqq {}& \frac12 \|\bv_\calT^n\|_{\bL^2(\rho;\Omega)}^2
+ \frac12\|\opStrh(\huh^n)\|_{\bL^2(2\mu;\Omega)}^2 
+ \frac12 \|\opDivh(\huh^n)\|_{L^2(\lambda;\Omega)}^2 \nonumber \\
&+\frac12 s_h(\huh^n,\huh^n) + \delta \Delta t^2 \|\bx_\calT^n\|_{L^2(\rho;\Omega)}^2,
\label{eq:energy_Newmark}
\end{align}
with $\delta\eqq \tfrac14(2\beta-\gamma)$, \ie $\delta=0$ for the standard choice
$\beta=\frac14$, $\gamma=\frac12$. 
Using standard manipulations for Newmark schemes (see
\cite[Lemma 3.3]{BDE:21}), one can show that
$\engy^n = \engy^1 + \sum_{m=1}^{n-1} \frac12(\loadext(t^{m+1})+\loadext(t^m),
\bu_\calT^{m+1}-\bu_\calT^{m})_{\bL^2(\Omega)}$,
so that $\engy^n$ is conserved in the absence of body forces. 

\section{First-order formulation in time}
\label{sec:wave_1st}

The first-order formulation of the elastic wave equation is obtained
by introducing the velocity field $\bv\eqq \partial_t \bu$ 
and the stress tensor $\stress\eqq \stress(\strain(\vu))$ as independent unknowns. 
Taking the time derivative of \eqref{eq:def_sigma_elasto} and exchanging the order
of derivatives leads to 
\begin{equation}
\polA^{-1}\partial_t \stress -\strain(\bv)  = \mzero,\qquad
\rho\partial_t \bv- \GRAD\SCAL\stress =\loadext,
\label{eq:elasto_1st}
\end{equation}
with $\polA^{-1}\btau=\frac{1}{2\mu}(\btau-\frac{\lambda}{2\mu+\lambda d}\trace(\btau)\matrice{I}_d)$,
together with the initial conditions $\stress(0) = \polA\strain(\vu_0)$, 
$\vv(0) = \vv_{0}$ in $\Dom$, 
and the boundary condition $\vv =\vzero$ on $J\times\Gamma$.
Assuming that $\vv\in H^1(J;\bL^2(\Omega)) \cap L^2(J;\bH^1_0(\Omega))$
and $\stress\in H^1(J;L^2(\Omega;\Msym))$, we obtain
\begin{equation} \label{eq:weak_1st_elasto}\left\{\begin{aligned}
(\partial_t\stress(t),\btau)_{\bL^2(\polA^{-1};\Omega)} - (\strain(\vv(t)),\btau)_{\bL^2(\Omega)}
&= 0, \\
(\partial_t\vv(t),\vw)_{L^2(\rho;\Omega)} + (\stress(t),\strain(\vw))_{\bL^2(\Omega)}
&= (\loadext(t),\vw)_{\bL^2(\Omega)},
\end{aligned}\right.\end{equation}
for all $(\btau,\vw)\in L^2(\Omega;\Msym)\times \bH^1_0(\Omega)$ and a.e.~$t\in J$.

\subsection{HHO space semi-discretization}

The idea is to approximate $\stress$ by a cellwise unknown $\stress_\calT\in C^1(\oJ;\bW_\calT)$ 
and $\vv$ by a hybrid unknown $\hvh\in C^1(\oJ;\hVkhz)$ 
(recall that $\bW_\calT\eqq \Pkd(\Th;\Msym)$). 
The space semi-discrete problem then reads as follows: For all $t\in\oJ$,
\begin{equation}\label{eq:HHO_1st_semi_elasto}\left\{\begin{aligned}
&(\partial_t \stress_\calT(t),\btau_\calT)_{\bL^2(\polA^{-1};\Omega)} - (\opStrh(\hvh(t)),\btau_\calT)_{\bL^2(\Omega)} =0,\\
&(\partial_t\vv_\calT(t),\vw_\calT)_{\bL^2(\rho;\Omega)} + (\stress_\calT(t),\opStrh(\hwh))_{\bL^2(\Omega)} + \tilde s_h(\hvh(t),\hwh) = (\loadext(t),\vw_\calT)_{\bL^2(\Omega)},
\end{aligned}\right.
\end{equation}
for all $(\btau_\calT,\hwh)\in \bW_\calT\times \hVkhz$, where the global strain
reconstruction operator $\opStrh$ is defined in~\eqref{eq:def_global_strain}.  
The stabilization bilinear form is
$\tilde s_h(\hvh,\hwh) \eqq \sum_{T\in\calT_h} \tilde\tau_{\dT} 
(\bS_{\dT}(\hvT),\bS_{\dT}(\hwT))_{\bL^2(\dT)}$ with
parameter 
$\tilde\tau_{\dT}\eqq \rho c_{\mathrm{S}}\frac{\ell_\Omega}{h_T}$
(that is, $\tilde\tau_{\dT}=\calO(h_T^{-1})$) or
$\tilde\tau_{\dT}\eqq \rho c_{\mathrm{S}}$
(that is, $\tilde\tau_{\dT}=\calO(1)$) where $c_{\mathrm{S}}=\sqrt{\mu/\rho}$.
The initial conditions for~\eqref{eq:HHO_1st_semi_elasto} 
are $\stress_\calT(0) = \polA \opStrh(\Ikh(\vu_0))$
and $\vv_\calT(0) = \Pi_\calT^{k'}(\vv_0)$, whereas 
the boundary condition is encoded in the fact that $\hvh(t)\in\hVkhz$ 
for all $t\in\oJ$.

The space semi-discrete schemes~\eqref{eq:HHO_2nd_semi_elasto} and
\eqref{eq:HHO_1st_semi_elasto} are not equivalent. Indeed, assume that 
$\huh$ solves~\eqref{eq:HHO_2nd_semi_elasto}, $(\stress_\calT,\hvh)$ solves
\eqref{eq:HHO_1st_semi_elasto}, and set $\hzh(t)\eqq \huh(0)+\int_0^t \hvh(s)\ds$.
Then, observing that the first equation in~\eqref{eq:HHO_1st_semi_elasto}
implies that $\partial_t \stress_\calT(t)=\polA\opStrh(\hvh(t))$, using the initial condition
for $\stress_\calT$ and the linearity of $\opStrh$ 
gives $\stress_\calT(t)=\polA \opStrh(\huh(0))+
\int_0^t\polA\opStrh(\hvh(s))\ds=\polA \opStrh(\huh(t))$ for all $t\in\oJ$.
Substituting into the second equation in \eqref{eq:HHO_1st_semi_elasto} and since
$\partial_t\vv_\calT=\partial_{tt}\vecteur{z}_\calT$, we infer that
for all $\hwh\in \hVkhz$,
\begin{multline*}
(\partial_{tt}\vecteur{z}_\calT(t),\vw_\calT)_{\bL^2(\rho;\Omega)}
+ (\opStrh(\hzh(t)),\opStrh(\hwh))_{\bL^2(\polA;\Dom)} + \tilde s_h(\partial_t\hzh(t),\hwh) \\ = (\loadext(t),\vw_\calT)_{\bL^2(\Omega)},
\end{multline*}
which differs from~\eqref{eq:HHO_2nd_semi_elasto} in the form
of the stabilization term. This difference in structure between the two formulations
has an impact on energy conservation. Indeed, defining the discrete energy for all $t\in\oJ$ as
\begin{equation}
\engy^{\ast}_h(t) \eqq \frac12 \|\vv_\calT(t)\|_{\bL^2(\rho;\Omega)}^2 + 
\frac12 \|\stress_\calT(t)\|_{\bL^2(\polA^{-1};\Omega)}^2,
\end{equation}
testing \eqref{eq:HHO_1st_semi_elasto} with $(\btau_\calT,\hwh)\eqq(\stress_\calT(t),\hvh(t))$ for
all $t\in J$ and integrating over time leads to
\begin{equation} \label{energy_bal_wave_HHO_1st}
\engy_h^{\ast}(t) + \int_0^t
\tilde s_h(\hvh(s),\hvh(s))\ds =\engy_h^{\ast}(0) + \int_0^t (\loadext(s),\vv_\calT(s))_{\bL^2(\Omega)}\ds.
\end{equation}
Comparing with~\eqref{energy_bal_wave_HHO_2nd}, we see that in the second-order formulation, the stabilization is included in the definition of the discrete energy and an exact energy balance is obtained, whereas in the first-order formulation, the discrete energy is independent of the stabilization, but the latter plays a dissipative role in the energy balance.  

\bRem[Link with HDG]
Recalling the material in Sect.~\ref{sec:HHO_HDG}, the space semi-discrete problem \eqref{eq:HHO_1st_semi_elasto} can be rewritten as an HDG formulation for the first-order wave equation. HDG methods typically consider the stabilization parameter $\tilde\tau_{\dT}=\calO(1)$
\cite{NgPeC:11}; see also \cite[Tab.~4]{StNPC:16} for a numerical study.
\eRem

Let $M_\calT^k\eqq \dim(\bW_\calT)=\frac{d(d+1)}{2}N_\calT^k$ 
and $\{\bzeta_i\}_{1\le i\le M_\calT^k}$ be 
the chosen basis for $\bW_\calT$. 
Let $\sfZ_\calT(t) \in \Real^{M_{\calT}^{k}}$ and
$(\sfV_\calT(t),\sfV_\calF(t))\in \Real^{N_\calT^{k'}\times N_{\calF,0}^k}$ 
be the component vectors of $\stress_\calT(t)\in \bW_\calT$ and $\hvh(t)\in \hVkhz$, respectively. Let $\sfM_{\calT\calT}^{\bsigma}$ be the mass matrix associated with the inner product in $\bL^2(\polA^{-1};\Omega)$ and the basis functions $\{\bzeta_i\}_{1\le i\le M_\calT^k}$, 
and recall that $\sfM_{\calT\calT}$ is the mass matrix associated with the inner product in $L^2(\rho;\Omega)$ and the basis functions $\{\bvarphi_i\}_{1\le i\le N_\calT^{k'}}$.
Let $\sfS_{\calT\calT}$, $\sfS_{\calT\calF}$, $\sfS_{\calF\calT}$, $\sfS_{\calF\calF}$ be the four blocks composing the matrix representing the stabilization bilinear form $\tilde s_h$. 
Let $\sfK_{\calT}\in \Real^{M_\calT^k\times N_\calT^{k'}}$ and 
$\sfK_{\calF}\in \Real^{M_\calT^k\times N_{\calF,0}^{k}}$ be the (rectangular) matrices 
representing the strain reconstruction operator $\opStrh$.
The algebraic realization of~\eqref{eq:HHO_1st_semi_elasto} is as follows: For all $t\in \oJ$, 
\begin{equation} \label{eq:HHO_1st_semi_alg_elasto}
\begin{bmatrix}
\sfM_{\calT\calT}^{\bsigma} & 0 & 0\\
0 & \sfM_{\calT\calT} & 0\\
0 & 0 & 0
\end{bmatrix} \begin{bmatrix}
\partial_t\sfZ_\calT(t)\\
\partial_t\sfV_\calT(t)\\
\bullet
\end{bmatrix} + \begin{bmatrix}
0 & -\sfK_{\mathcal{T}} & -\sfK_{\mathcal{F}}\\
\sfK_{\mathcal{T}}\tr & \sfS_{\mathcal{TT}} & \sfS_{\mathcal{TF}}\\
\sfK_{\mathcal{F}}\tr & \sfS_{\mathcal{FT}} & \sfS_{\mathcal{FF}}
\end{bmatrix} \begin{bmatrix}
\sfZ_\calT(t)\\
\sfV_\calT(t)\\
\sfV_\calF(t)
\end{bmatrix} = \begin{bmatrix}
0\\
\sfF_\calT(t)\\
0
\end{bmatrix},
\end{equation}
where the bullet stands for $\partial_t\sfV_\calF(t)$ which is irrelevant owing to
the structure of the mass matrix. Notice that the third equation 
in~\eqref{eq:HHO_1st_semi_alg_elasto} implies that 
\begin{equation}
\sfS_{\mathcal{FF}} \sfV_\calF(t) = -(\sfK_\calF\tr \sfZ_\calT(t)
+ \sfS_{\calF\calT}\sfV_\calT(t)),
\end{equation}
and that the submatrix $\sfS_{\mathcal{FF}}$ is 
symmetric positive-definite. A crucial observation is that
this submatrix is additionally block-diagonal in the 
mixed-order case, but this property is lost in the equal-order case owing to the presence
of the displacement reconstruction operator in the stabilization (see~\eqref{eq:def_stab_equal_wave}).

\bRem[Error analysis]
The error analysis for the space semi-discrete problem~\eqref{eq:HHO_1st_semi_elasto} is performed in \cite[Thm.~4.3]{BDES:21} for $\tilde\tau_{\dT}=\calO(h_T^{-1})$ and the acoustic wave equation (it can be extended to the elastic wave equation). In particular, the energy-error $\|\vv-\vv_\calT\|_{L^\infty(J; \bL^2(\rho;\Dom))} + \|\stress-\stress_\calT\|_{L^\infty(J;\bL^2(\polA^{-1};\Dom))}$ decays as $\calO(h^{k+1})$ if $(\stress,\vv)\in C^1(\oJ;\matrice{H}^{k+1}(\Dom)\times\vecteur{H}^{k+2}(\Dom))$. 
We refer the reader to \cite{CoQue:14} for the error analysis in the HDG setting with $\tilde\tau_{\dT}=\calO(1)$, including an improved $L^\infty(J;\Ldeuxd)$-estimate on a post-processed displacement field decaying at rate $\calO(h^{k+2})$.
\eRem

\subsection{Time discretization}

The space semi-discrete problem~\eqref{eq:HHO_1st_semi_elasto} can be discretized in time by means of a Runge--Kutta (RK) scheme. RK schemes are defined by a set of coefficients, $\{a_{ij}\}_{1\le i,j\le s}$, $\{b_i\}_{1\le i\le s}$, $\{c_i\}_{1\le i\le s}$, where $s\ge1$ is the number of stages. We consider diagonally implicit RK schemes (DIRK), where $a_{ij}=0$ if $j> i$, and explicit RK schemes (ERK), where additionally $a_{ii}=0$. The implementation of DIRK and ERK schemes is slightly different owing to the treatment of the face unknowns. 

Let us start with DIRK schemes. For all $n\in\{1,\ldots,N\}$, given
$(\sfZ_\calT^{n-1},\sfV_\calT^{n-1})$ from the previous time-step or the
initial condition and letting $\sfF_\calT^{n-1+c_j}\eqq \sfF_\calT(t_{n-1}+c_j\Delta t)$
for all $1\le j\le s$, one proceeds as follows:
\begin{enumerate}
\item Solve sequentially for all $1\le i\le s$,
\begin{multline} \label{eq:HHO-DIRK}
\begin{bmatrix}
\sfM_{\calT\calT}^{\bsigma} & 0 & 0\\
0 & \sfM_{\calT\calT} & 0\\
0 & 0 & 0
\end{bmatrix} \begin{bmatrix}
\sfZ_\calT^{n,i}\\
\sfV_\calT^{n,i}\\
\bullet
\end{bmatrix} = \begin{bmatrix}
\sfM_{\calT\calT}^{\bsigma} & 0 & 0\\
0 & \sfM_{\calT\calT} & 0\\
0 & 0 & 0
\end{bmatrix} \begin{bmatrix}
\sfZ_\calT^{n-1}\\
\sfV_\calT^{n-1}\\
\bullet
\end{bmatrix}
\\
+ \Delta t \sum_{j=1}^{i} a_{ij} \left( 
\begin{bmatrix}
0\\
\sfF_\calT^{n-1+c_j}\\
0
\end{bmatrix} - \begin{bmatrix}
0 & -\sfK_{\mathcal{T}} & -\sfK_{\mathcal{F}}\\
\sfK_{\mathcal{T}}\tr & \sfS_{\mathcal{TT}} & \sfS_{\mathcal{TF}}\\
\sfK_{\mathcal{F}}\tr & \sfS_{\mathcal{FT}} & \sfS_{\mathcal{FF}}
\end{bmatrix} \begin{bmatrix}
\sfZ_\calT^{n,j}\\
\sfV_\calT^{n,j}\\
\sfV_\calF^{n,j}
\end{bmatrix}
\right).
\end{multline}
This is a linear system for the triple $(\sfZ_\calT^{n,i},
\sfV_\calT^{n,i},\sfV_\calF^{n,i})$ (which appears on both the left- and right-hand sides). 
The upper $2\times 2$ submatrix associated with the cell unknowns
$(\sfZ_\calT^{n,i},\sfV_\calT^{n,i})$ being block-diagonal, 
static condensation can be efficiently performed in \eqref{eq:HHO-DIRK} 
leading to a global transmission problem coupling only the components of $\sfV_\calF^{n,i}$.
\item Finally set
\begin{multline}
\begin{bmatrix}
\sfM_{\calT\calT}^{\bsigma} & 0\\
0 & \sfM_{\calT\calT}
\end{bmatrix} \begin{bmatrix}
\sfZ_\calT^{n}\\
\sfV_\calT^{n}
\end{bmatrix} \eqq \begin{bmatrix}
\sfM_{\calT\calT}^{\bsigma} & 0\\
0 & \sfM_{\calT\calT} 
\end{bmatrix} \begin{bmatrix}
\sfZ_\calT^{n-1}\\
\sfV_\calT^{n-1}
\end{bmatrix}
\\
+ \Delta t \sum_{j=1}^{s} b_{j} \left( 
\begin{bmatrix}
0\\
\sfF_\calT^{n-1+c_j}
\end{bmatrix} - \begin{bmatrix}
0 & -\sfK_{\mathcal{T}} & -\sfK_{\mathcal{F}}\\
\sfK_{\mathcal{T}}\tr & \sfS_{\mathcal{TT}} & \sfS_{\mathcal{TF}}
\end{bmatrix} \begin{bmatrix}
\sfZ_\calT^{n,j}\\
\sfV_\calT^{n,j}\\
\sfV_\calF^{n,j}
\end{bmatrix}
\right). \label{eq:final_update_DIRK}
\end{multline}
\end{enumerate}
For ERK schemes, instead, one proceeds as follows:
\begin{enumerate}
\item Set $(\sfZ_\calT^{n,1},\sfV_\calT^{n,1})\eqq (\sfZ_\calT^{n-1},\sfV_\calT^{n-1})$
and solve $\sfS_{\mathcal{FF}} \sfV_\calF^{n,1} = -(\sfK_\calF\tr \sfZ_\calT^{n,1}
+ \sfS_{\calF\calT}\sfV_\calT^{n,1})$.
\item If $s\ge2$, solve sequentially for all $2\le i\le s$,
\begin{multline}
\begin{bmatrix}
\sfM_{\calT\calT}^{\bsigma} & 0 \\
0 & \sfM_{\calT\calT} 
\end{bmatrix} \begin{bmatrix}
\sfZ_\calT^{n,i}\\
\sfV_\calT^{n,i}
\end{bmatrix} = \begin{bmatrix}
\sfM_{\calT\calT}^{\bsigma} & 0 \\
0 & \sfM_{\calT\calT}
\end{bmatrix} \begin{bmatrix}
\sfZ_\calT^{n-1}\\
\sfV_\calT^{n-1}
\end{bmatrix}
\\
+ \Delta t \sum_{j=1}^{i-1} a_{ij} \left( 
\begin{bmatrix}
0\\
\sfF_\calT^{n-1+c_j}
\end{bmatrix} - \begin{bmatrix}
0 & -\sfK_{\mathcal{T}} & -\sfK_{\mathcal{F}}\\
\sfK_{\mathcal{T}}\tr & \sfS_{\mathcal{TT}} & \sfS_{\mathcal{TF}}
\end{bmatrix} \begin{bmatrix}
\sfZ_\calT^{n,j}\\
\sfV_\calT^{n,j}\\
\sfV_\calF^{n,j}
\end{bmatrix}
\right),
\end{multline}
and $\sfS_{\mathcal{FF}} \sfV_\calF^{n,i} = -(\sfK_\calF\tr \sfZ_\calT^{n,i}
+ \sfS_{\calF\calT}\sfV_\calT^{n,i})$.
\item Finally update the cell unknowns as in~\eqref{eq:final_update_DIRK}. 
\end{enumerate}
We emphasize that the ERK scheme is effective only in the mixed-order
case since the submatrix $\sfS_{\calF\calF}$ is then block-diagonal.
The HHO-ERK scheme is subjected for its stability to a CFL condition
on the time-step. The choice $\tilde\tau_{\dT}=\calO(1)$ is
recommended for the stabilization parameter since it leads to a CFL
condition scaling linearly with the mesh size (the scaling is quadratic
for $\tilde\tau_{\dT}=\calO(h_T^{-1})$). For the HHO-DIRK scheme,
both choices for the stabilization parameter are viable, and 
numerical experiments indicate that the choice $\tilde\tau_{\dT}=\calO(h_T^{-1})$
leads to more accurate solutions, with an $O(h^{k+2})$ decay rate
for the $L^\infty(J;\Ldeuxd)$-norm.

\section{Numerical example}
\label{sec:wave_num}

To illustrate the HHO methods described in the previous sections, we consider
the propagation of an elastic wave in a 
two-dimensional heterogeneous domain $\Omega$ such that
$\overline\Omega=\overline\Omega_{1}\cup\overline\Omega_{2}$ with $\Omega_{1}\eqq 
(-\frac32,\frac32)\times(-\frac32,0)$ and $\Omega_{2}\eqq 
(-\frac32,\frac32)\times(0,\frac32)$. The material properties are
$\rho_1\eqq 1$, $\mu_1=\lambda_1\eqq1$ in $\Omega_1$
and $\rho_2\eqq 1$, $\mu_2=\lambda_2\eqq9$  in $\Omega_2$, so that 
$c_{\mathrm{S},2}\eqq 3c_{\mathrm{S},1}$, $c_{\mathrm{P},2}\eqq
3c_{\mathrm{P},1}$.
The simulation time is $T_{\textrm{f}}\eqq 1$,
and homogeneous Dirichlet boundary conditions are enforced.
The body force is $\loadext\eqq\vzero$, and the initial conditions are $\vu_0\eqq\vzero$
together with
\begin{equation}
\vv_0(x,y) \eqq \theta \exp\big( -\pi^2 \tfrac{r^2}{\lambda^2}\big) 
(\bx-\bx_c),
\end{equation}
with $\theta\eqq 10^{-2}$ [s$^{-1}$], 
$\lambda\eqq \frac{v_{\mathrm{P},2}}{f_c}$ [m] with $f_c\eqq 10$ [s$^{-1}$],
$r^2\eqq \|\bx-\bx_c\|_{\ell^2}^2$, $\bx_c\eqq (0,\frac23)$. The initial
condition corresponds to a Ricker wave centered at the point 
$\bx_c\in \Omega_2$. The wave first propagates in $\Omega_2$, then is partially
transmitted to $\Omega_1$ and later it is also reflected at the boundary of $\Omega$. 

Numerical results are obtained using the Newmark scheme (with $\beta=\frac14$, 
$\gamma=\frac12$), a three-stage singly diagonally implicit RK of order 4 (in short,
SDIRK(3,4)), and a four-stage explicit RK scheme of order 4 (in short, ERK(4)).
The Butcher tableaux for the RK schemes are, respectively, 
\begin{equation}
\begin{array}{c|ccc}
\gamma & \gamma & 0 & 0\\
\frac{1}{2} & \frac{1}{2} - \gamma & \gamma & 0\\
1-\gamma  & 2\gamma & 1-4 \gamma & \gamma\\
\hline
&\delta&1-2\delta&\delta
\end{array}
\qquad\qquad
\begin{array}{c|cccc}
0&0&0&0&0\\
\frac12&\frac12&0&0&0\\[2pt]
\frac12&0&\frac12&0&0\\[2pt]
1&0&0&1&0\\
\hline
&\frac16&\frac13&\frac13&\frac16
\end{array}
\end{equation}
with $\gamma\eqq \frac{1}{\sqrt{3}}\cos\left(\frac{\pi}{18}\right)+\frac{1}{2}$, $\delta\eqq \frac{1}{6\left(2\gamma-1\right)^{2}}$.
We consider a quadrangular mesh of size $h\eqq 2^{-6}$ and a time-step $\Delta t\eqq 0.1\times 2^{-6}$.
Figure~\ref{fig:2d_profiles_wave} reports the velocity profiles over the computational domain at the four simulation times $t\in \{\frac18,\frac14,\frac12,1\}$. These profiles are obtained using the SDIRK(3,4) scheme ($k'=k$, $\tilde\tau_{\dT}=\calO(1)$). We observe the various reflections of the elastic waves at the interface and at the domain boundary. 

\begin{figure}[htb]
\begin{centering}
\includegraphics[width=0.24\columnwidth]{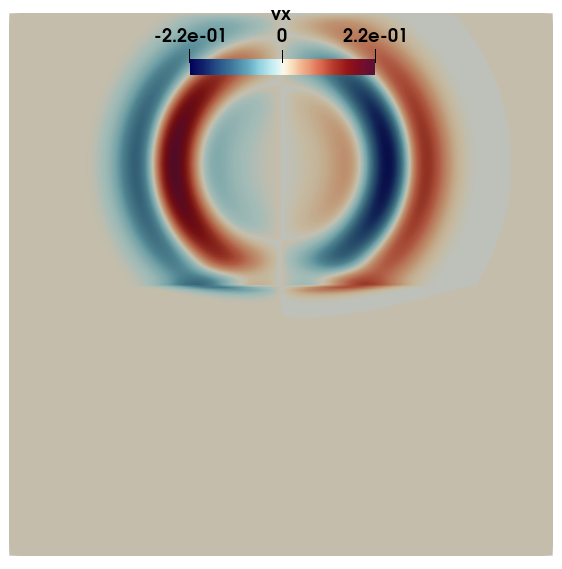}%
\includegraphics[width=0.24\columnwidth]{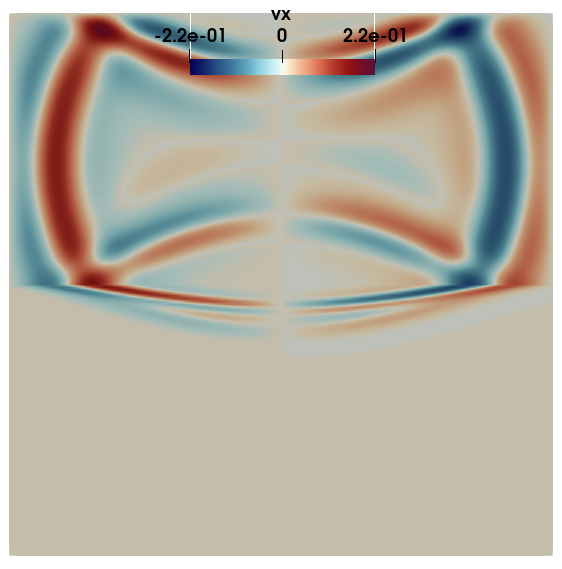}%
\includegraphics[width=0.24\columnwidth]{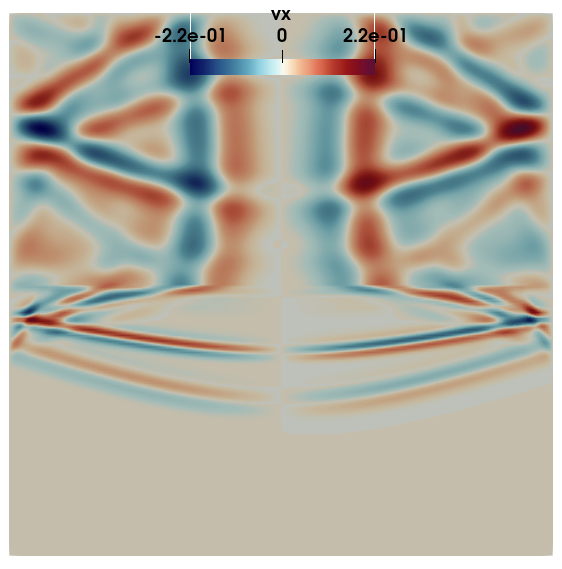}%
\includegraphics[width=0.24\columnwidth]{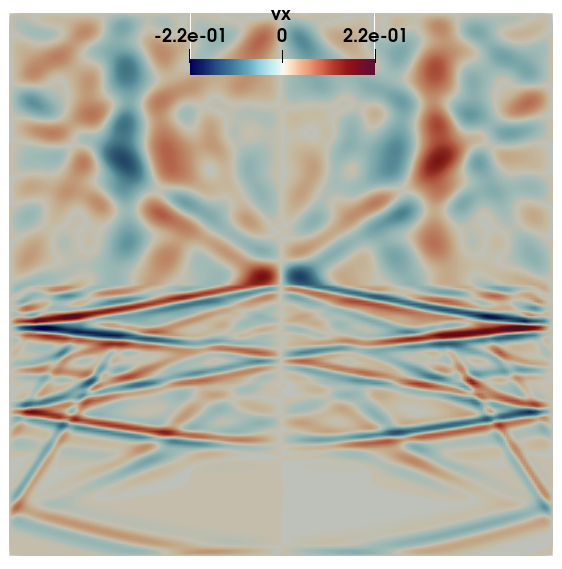}%
\par\end{centering}
\begin{centering}
\includegraphics[width=0.24\columnwidth]{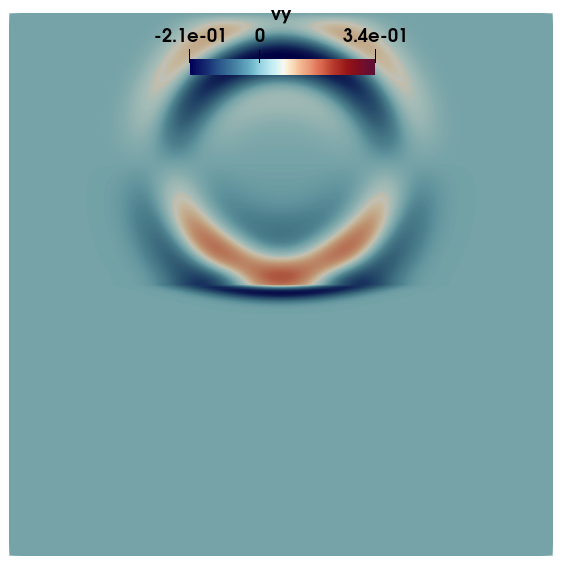}%
\includegraphics[width=0.24\columnwidth]{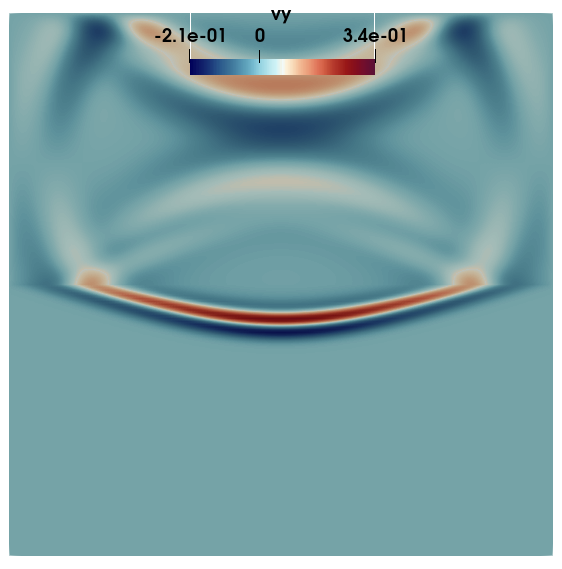}%
\includegraphics[width=0.24\columnwidth]{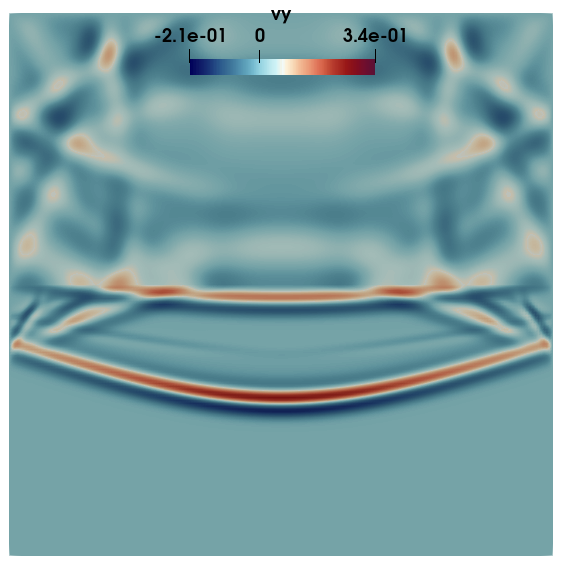}%
\includegraphics[width=0.24\columnwidth]{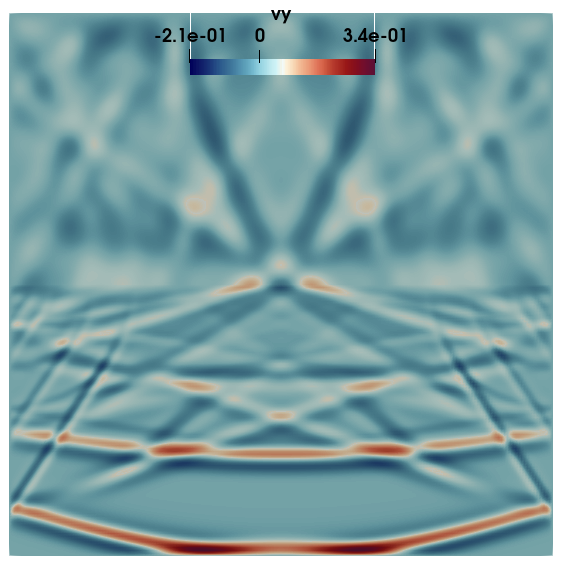}%
\par\end{centering}
\caption{Velocity profiles at the times $t\in \{\frac18,\frac14,\frac12,1\}$ (from left to right). Upper row: $v_x$; bottom row: $v_y$. SDIRK(3,4) scheme, $k'=k$, $\tilde\tau_{\dT}=\calO(1)$.}
\label{fig:2d_profiles_wave}
\end{figure}

These results can be 
compared against semi-analytical solutions obtained using the \texttt{gar6more2d} software.\footnote{\texttt{https://gforge.inria.fr/projects/gar6more2d/}} 
The semi-analytical solution is
based on a reformulation of the problem with zero initial conditions and a Dirac source
term with a time delay of $0.15$ [s] (this value is tuned to match
the choice of the parameter $\theta$, see
\cite{Boillot:14}). The
comparisons are made by tracking the velocity at two
sensors, one located in $\Omega_1$ at the point $S_1\eqq(\frac13,-\frac13)$
and one located in $\Omega_2$ at the point $S_2\eqq(\frac13,\frac13)$. Since 
the semi-analytical solution
assumes propagation in two half-spaces, the comparison with the simulations remains
meaningful until the reflected waves at the boundary reach one
of the sensors (this happens around the times $0.7$
for $S_1$ and $0.45$ for $S_2$). Figure~\ref{fig:sdirk} reports the results
for the cell velocity component $v_x$ with $\Delta t\eqq 0.1\times 2^{-8}$ for the second-order Newmark scheme,
$\Delta t\eqq 0.1\times 2^{-6}$ for the SDIRK(3,4) scheme, and 
$\Delta t\eqq 0.1\times 2^{-9}$ for the ERK(4) scheme (owing to the stability 
condition). Equal-order is used for the Newmark and SDIRK schemes, and mixed-order for the ERK scheme.
For both RK schemes, the stabilization parameter is $\tilde\tau_{\dT}=\calO(1)$.
We observe that increasing the polynomial degree in the HHO discretization is beneficial
for all the time-stepping schemes, and that the predictions overlap with the 
semi-analytical solution for $k=3$. For the RK schemes, the predictions are already quite
accurate for $k=2$, but this is not the case for the Newmark scheme. As expected, the 
profiles at the sensor $S_1$ are more difficult to capture due to the transmission of the
incoming wave across the interface separating the two media. 

\begin{figure}[htb]
\begin{centering}
\includegraphics[width=0.45\columnwidth]{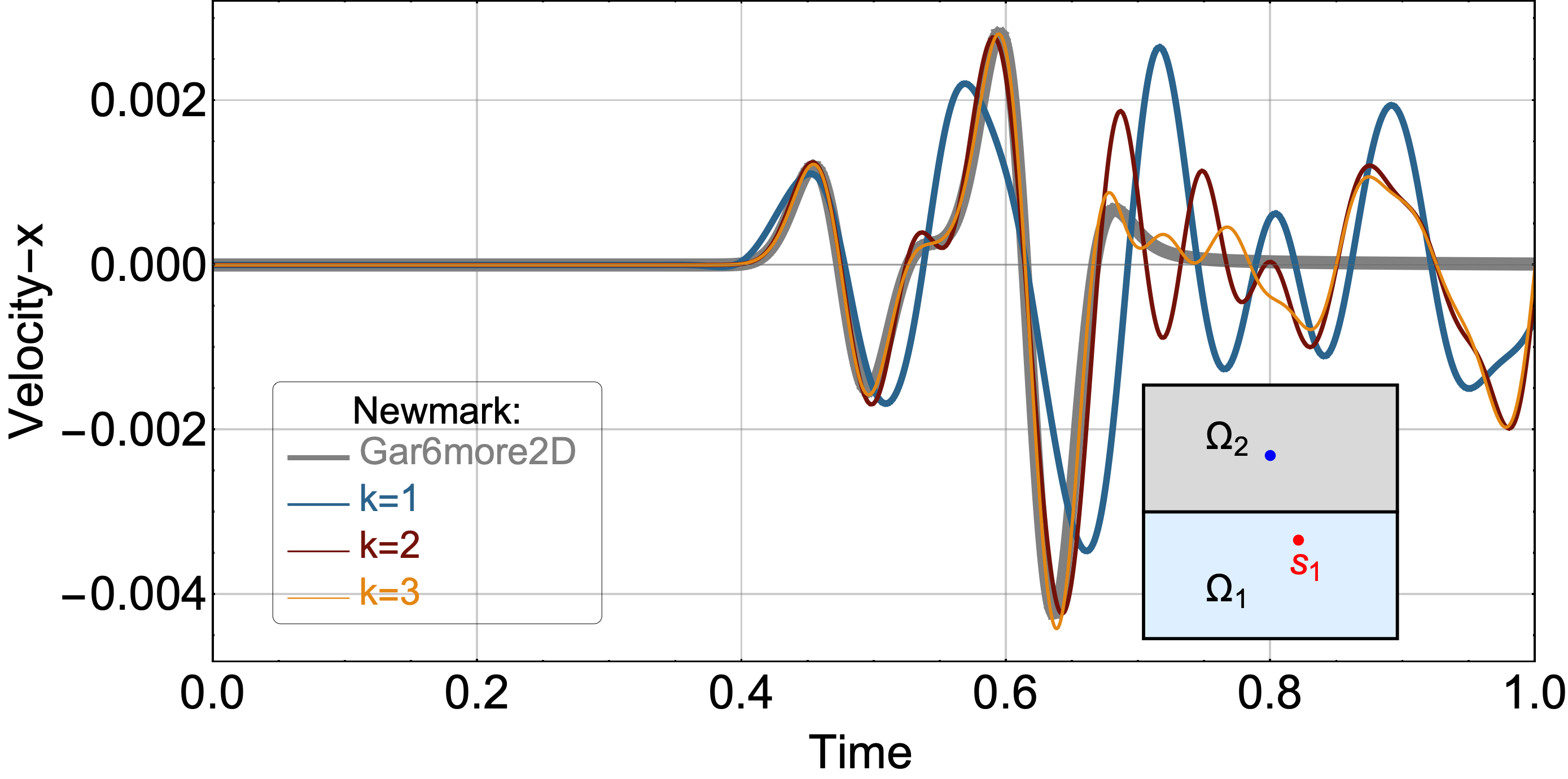}
\includegraphics[width=0.45\columnwidth]{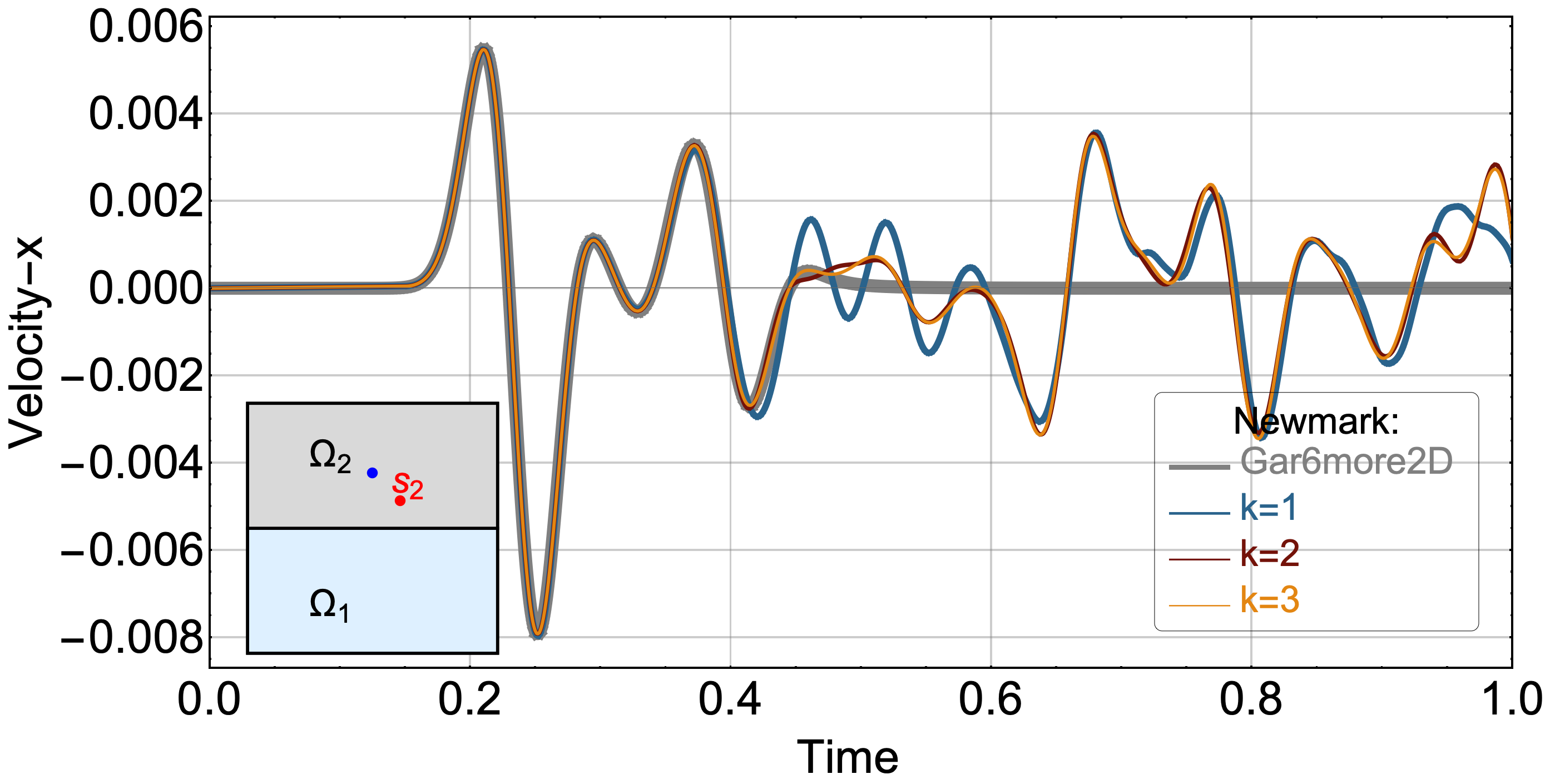}
\includegraphics[width=0.45\columnwidth]{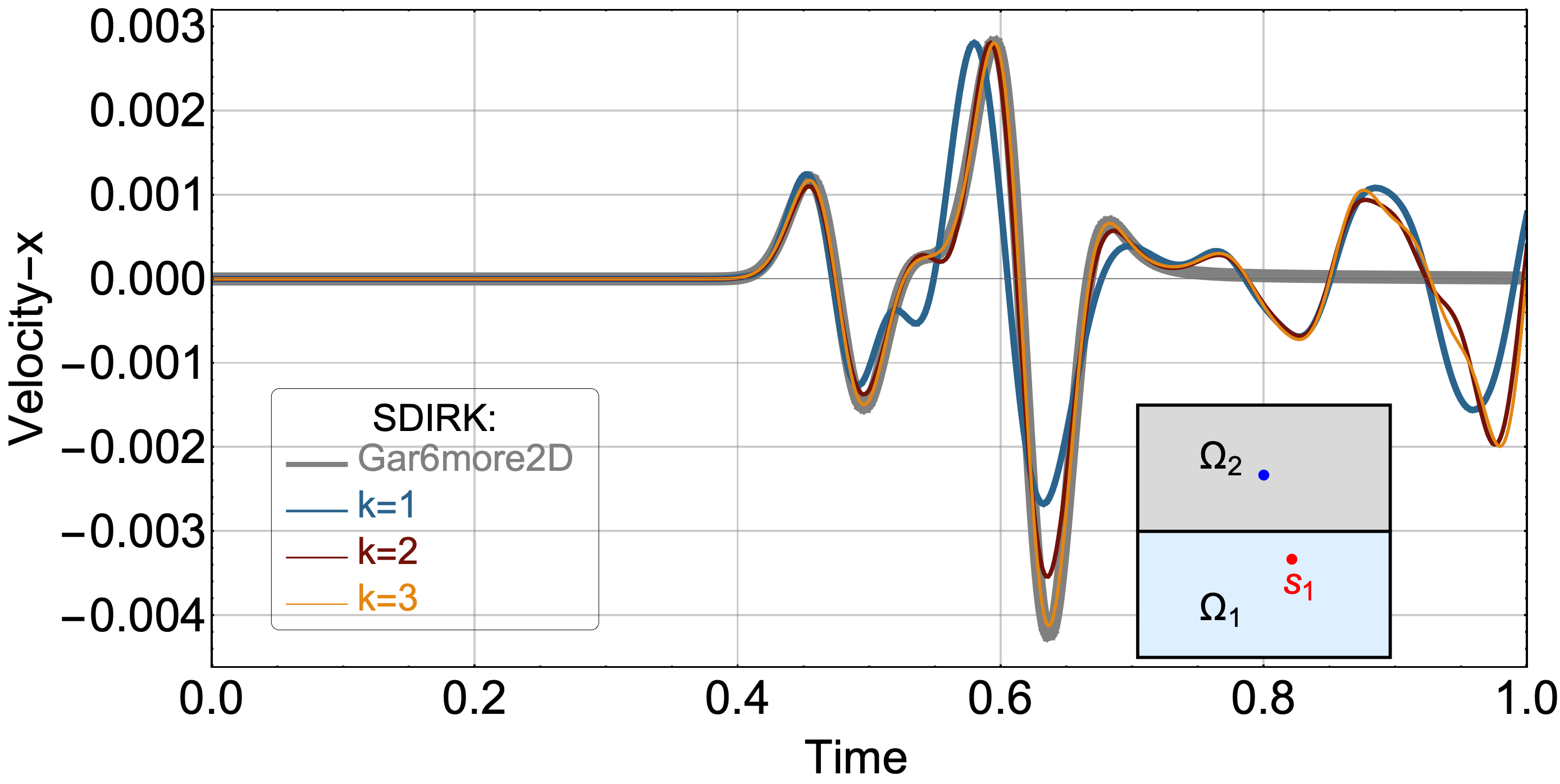}
\includegraphics[width=0.45\columnwidth]{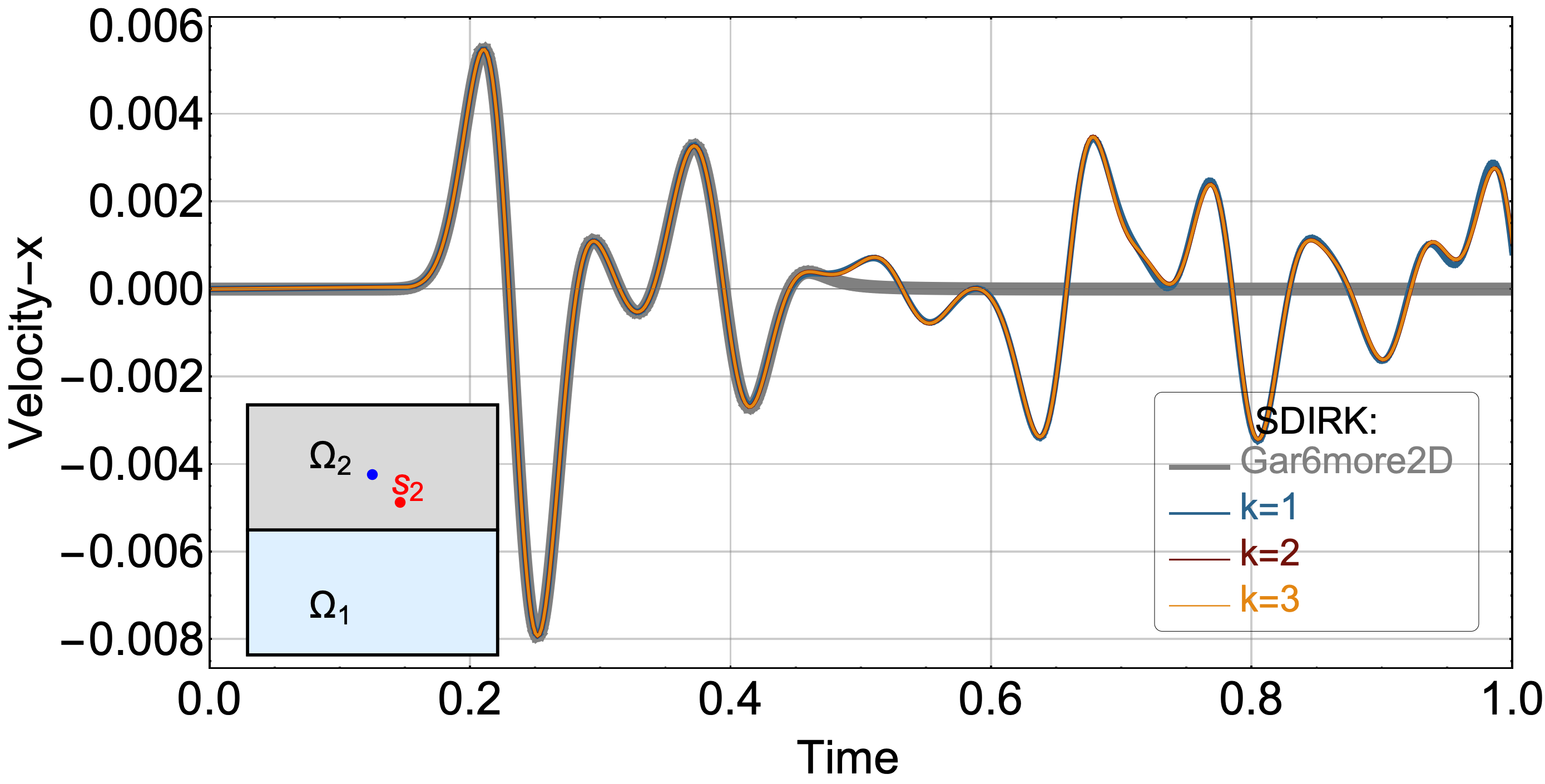}
\includegraphics[width=0.45\columnwidth]{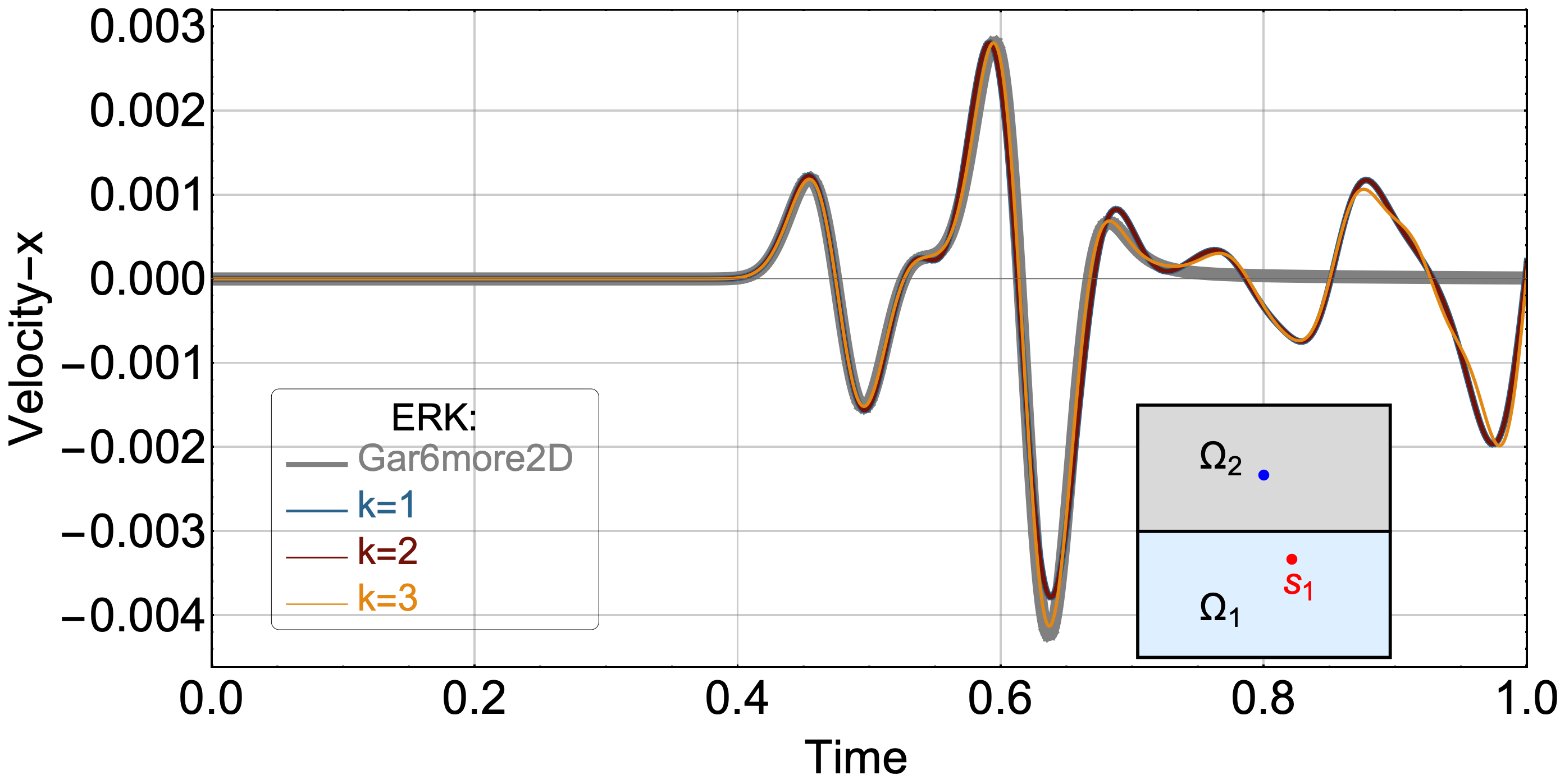}
\includegraphics[width=0.45\columnwidth]{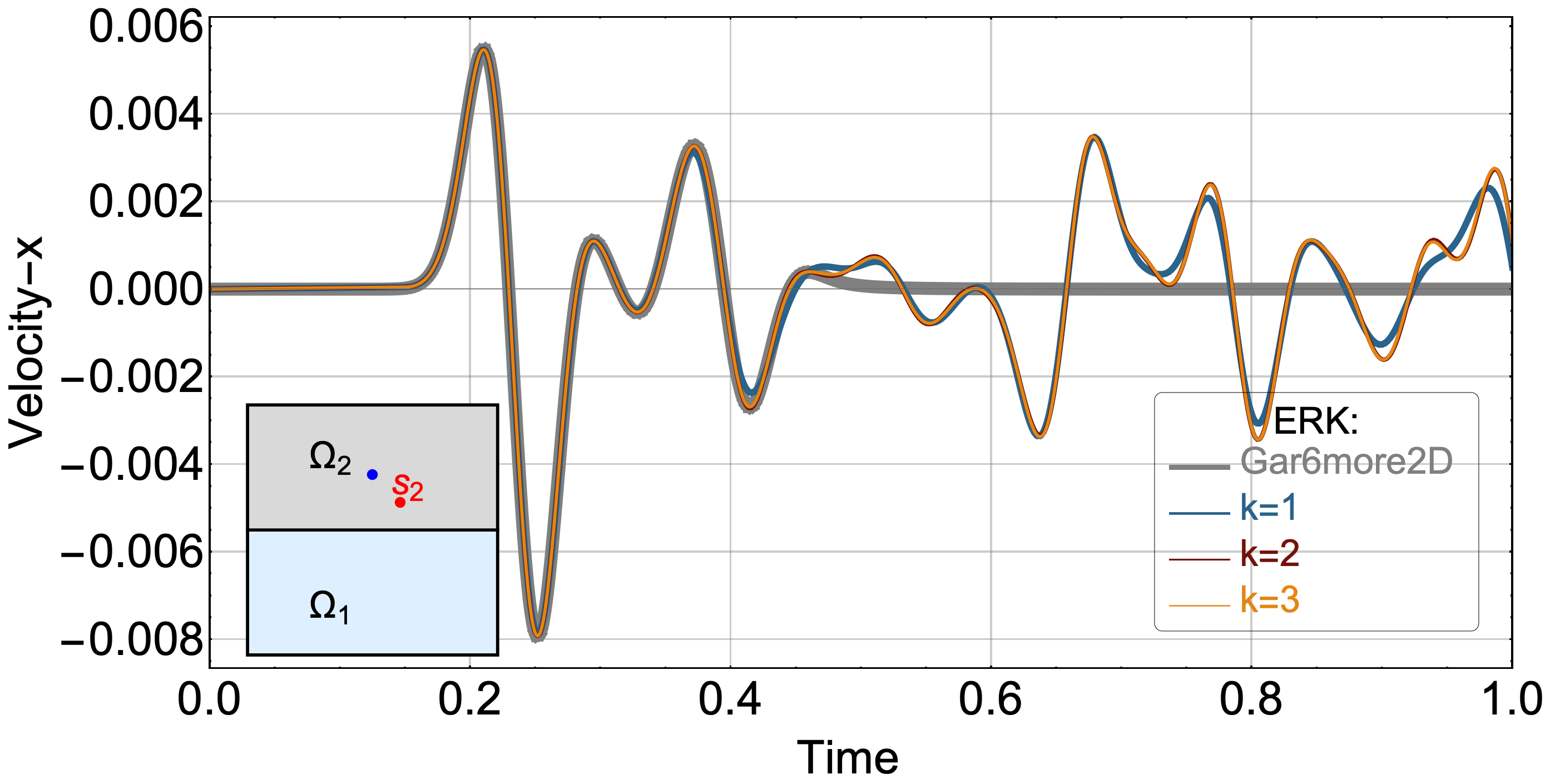}
\par\end{centering}
\caption{Velocity component $v_x$ at the sensor $S_1$ (left column) and at the sensor $S_2$ (right column). The polynomial degree is $k\in\{1,2,3\}$. Upper row: Newmark ($\Delta t\eqq 0.1\times 2^{-8}$), middle row: SDIRK(3,4) ($\Delta t\eqq 0.1\times 2^{-6}$), bottom row: ERK(4) ($\Delta t\eqq 0.1\times 2^{-9}$).}
\label{fig:sdirk}
\end{figure}


\chapter{Contact and friction}\label{chap::contact}

In this chapter, we show how the HHO method can be used to discretize a linear elasticity
problem with nonlinear boundary conditions resulting from contact and friction. 
The main idea is to use a boundary penalty technique to enforce these conditions. 
This approach leads, under some assumptions, to a discrete semilinear form enjoying 
a monotonicity property. The error analysis reveals that the degree of the face unknowns
on the contact/friction boundary has to be raised to $(k+1)$ to ensure optimal estimates.

\section{Model problem}
\label{sec: sig_model_problems}

As in Sect.~\ref{sec:model-linear_elas}, we consider an elastic body occupying the bounded Lipschitz 
domain $\Omega\subset\Rd$, $d\in\{2,3\}$, in the reference configuration. 
The boundary $\front$ is now partitioned into three disjoint subsets: the Dirichlet
boundary $\Bd$, the Neumann boundary $\Bn$, and the contact/friction boundary $\Bc$,
with $\meas(\Bd) > 0$ (to prevent rigid-body motions) and $\meas(\Bc) > 0$. 
The body undergoes infinitesimal deformations due to volume forces $\loadext \in L^2(\Omega; \Rd)$
and surface loads $\Gn \in L^2(\Bn; \Rd)$, and it is clamped on
$\Bd$ (for simplicity). Recall that the linearized strain tensor
associated with a displacement field $\vv:\Omega \rightarrow \Reel^d$ is
$\strain(\vv) := \frac12 (\grad \vv + \grad \vv\tr)\in \Msym$.
Assuming a linear elastic behaviour, 
the Cauchy stress tensor resulting from a strain tensor $\strain$
is given by
\begin{equation} \label{eq:stress_strain}
\stress(\strain)= 2\mu \strain + \lambda \trace (\strain) \matrice{I}{}_d \in \Msym,
\end{equation}
where $\mu$ and $\lambda$ are the Lam\'e coefficients of the material satisfying
$\mu > 0$ and $3\lambda +2\mu> 0$, and $\matrice{I}{}_d$ is the identity tensor of order $d$.

Let $\NO$ be the unit outward normal vector to $\Omega$. On the boundary, we consider the
following decompositions into normal and tangential components of a displacement field $\vv$ and
a stress tensor $\stress$:
\begin{equation}
\vv = \vn \NO + \vvt \quad \hbox{ and } \qquad 
\stressn \eqq \stress\NO = \stressnn\NO + \stresst,
\end{equation}
where $\vn \eqq \vv \SCAL \NO$ and $\stressnn \eqq \stressn \SCAL \NO$
(so that $\vvt\SCAL\NO=0$ and $\stresst\SCAL\NO=0$).
The model problem consists in finding the displacement field 
$\vu: \Omega \rightarrow \Reel^d$ such that, using the shorthand notation $\stress(\vu)\eqq \stress(\strain(\vu))$,
\begin{align} 
\label{eq}
&-\mdivergence{\stress(\vu)} = \loadext\; \text{in $\Omega$},\\ 
\label{eq:Dir_Neu_contact}
&\vu = \vzero \; \text{on $\Bd$} \;\wedge\;
\stressn(\vu) = \Gn \; \text{on $\Bn$},\\
\label{t.cont}
&\un  \le 0\;\wedge\;
\stressnn(\vu) \le 0 \;\wedge\;
\stressnn(\vu)\, \un = 0 \quad \text{on $\Bc$},\\
\label{tresca}
&\abs{ \stresst(\vu) } \le \Sc\ \text{if $\vut = \vzero$}\;\wedge\;
\stresst(\vu) =  -\Sc \frac{\vut}{\abs{\vut}}\ \text{if $\abs{\vut}>0$}\quad \text{on $\Bc$},
\end{align}
where \eqref{t.cont} are called unilateral contact conditions 
and \eqref{tresca} Tresca friction conditions. 
The first condition in \eqref{t.cont} expresses non-interpenetration, whereas the last condition, called complementarity condition, means that either there is contact ($\un = 0$) or there is no normal force ($\stressnn(\vu)=0$).
In~\eqref{tresca}, 
$\Sc\ge0$ is a given threshold parameter (more generally, $s$ can be a nonnegative function on $\Bc$),
and $\abs{\SCAL}$ stands for the 
Euclidean norm in $\Reel^{d}$ (or the absolute value depending on the context).
The conditions in \eqref{tresca} mean that sliding cannot occur as long as the magnitude of the tangential stress $\abs{ \stresst(\vu) }$ is lower than the threshold $s$. When the threshold is reached, sliding can happen, in a direction opposite to $\stresst (\vu)$ (see, e.g., \cite[Chapter 10]{kikuchi-oden-88}). The case of frictionless contact is recovered by setting $\Sc:=0$ in \eqref{tresca}.

Let us briefly discuss some variants of the above model. On the one hand,
bilateral contact with Tresca friction can be considered by keeping \eqref{tresca},
whereas \eqref{t.cont} is substituted by the condition
\begin{equation} \label{bilateral}
\un  = 0. 
\end{equation}
This setting is relevant to model persistent contact.
In the case of unilateral contact, 
nonzero tangential stress ($\abs{ \stresst(\vu) } $ $>$ $0$) can occur in regions with no-adhesion ($\un < 0$), which is not expected physically.
The setting of bilateral contact 
prevents such situations. Indeed, since $\un=0$, there are no 
regions with no-adhesion.
On the other hand, substituting \eqref{tresca} by 
\begin{equation}
\abs{ \stresst(\vu) } \le \ \Fc \abs{ \stressnn(\vu) }\ \text{if $\vut = \vzero$} \;\wedge\;
\stresst(\vu) =  -\Fc \abs{ \stressnn(\vu) } \frac{\vut}{\abs{\vut}}\ \text{if $\abs{\vut}>0$},
\label{coulomb}
\end{equation}
where $\Fc \geq 0$ is a given friction coefficient, leads to static Coulomb friction. The condition \eqref{coulomb} is an adaptation of the quasi-static (or dynamic) Coulomb's law, in which the tangential velocity $\dot\vut$ plays the same role as the displacement $\vut$. 
In the rest of this chapter, we focus on 
the Tresca friction model. This choice is motivated more by mathematical simplicity than physical reasons. Moreover, the Tresca friction model can be useful when Coulomb friction is approximated iteratively.

Recalling the notation $\Hund\eqq H^1(\Dom;\Rd)$, we introduce the Hilbert space $\VD$ and the convex 
cone $\vecteur{K}$ such that
\begin{equation*}
\VD := \left\{ \vv \in \Hund \:|\: 
\vv_{|\Bd}=\vzero\right\},
\quad
\vecteur{K} := \left\{ \vv \in \VD \:|\: 
\vn \leq 0 \hbox{ on } \Bc \right\},
\end{equation*}
i.e., the Dirichlet condition on $\Bd$ is explicitly enforced in the space
$\VD$
and the non-interpenetration condition on $\Bc$ is explicitly enforced
in the cone $\vecteur{K}$. We define the following bilinear form and the following
linear and nonlinear forms:
\begin{align}
a(\vv, \vw) {}&:=  \psm[\Omega]{\stress(\strain(\vv))}{\strain(\vw)} 
= 2\mu \psm[\Omega]{\strain(\vv)}{\strain(\vw)} + \lambda \pss[\Omega]{\divergence{\vv}}{\divergence{\vw}}, \label{eq:def_bilin_a}\\
\ell(\vw) &:=  \psv[\Omega]{\loadext}{\vw} 
+ \psv[\Bn]{\Gn}{\vw},  \qquad
j(\vw) :=  \int_{\Bc} s\abs{\vwt}\ds, %
\end{align}
for all $\vv,\vw\in \VD$. The weak formulation of \eqref{eq}--\eqref{tresca} leads to the following variational inequality:
\begin{equation} \label{weak}
\left\{
\begin{array}{l}
\text{Find $\vu \in \vecteur{K}$ such that}\\
a(\vu,\vw - \vu) + j(\vw) - j(\vu) \geq \ell(\vw - \vu), \quad \forall \vw \in \vecteur{K}.
\end{array}
\right.
\end{equation}
This problem admits a unique solution; see, e.g., \cite[Theorem~10.2]{kikuchi-oden-88}.
Moreover, this solution is the unique minimizer in $\vecteur{K}$ of the energy functional
$\engy:\VD\to\Reel$ such that
\begin{equation} \label{energy-tresca}
\engy(\vv) \eqq \frac12 a(\vv,\vv) + j(\vv) - \ell(\vv).
\end{equation}
 
\section{HHO-Nitsche method} \label{sec:hho_face_elas}

The HHO-Nitsche method presented in this section to approximate the model problem \eqref{weak} is inspired by the FEM-Nitsche method devised in \cite{ChoHi:13,Chouly:14}. 
Therefore, we first start with a brief description of the ideas underlying this latter method.

\subsection{FEM-Nitsche method}
\label{sec:FEM-Nitsche}

The two keys ideas in the FEM-Nitsche method are on the one hand a reformulation 
due to \cite{CurnierA:88}
of the conditions~\eqref{t.cont}-\eqref{tresca} as nonlinear equations and
on the other hand the use of a consistent boundary-penalty method inspired by Nitsche
\cite{nit71} to enforce these conditions in the discrete problem. 

For all $x\in\Reel$,
let $\projRm{x} \eqq \min(x,0)$ denote its projection onto  
$\Reel^-:=(-\infty,0]$, and for all $\vecteur{x} \in \Reel^{d}$,
let $[\vecteur{x}]_{\alpha} \eqq \vecteur{x}$ if $\abs{\vecteur{x}} \leq \alpha$ and 
$[\vecteur{x}]_{\alpha} \eqq \alpha \frac{\vecteur{x}}{\abs{\vecteur{x}}}$ 
if $\abs{\vecteur{x}}>\alpha$ denote its projection onto 
the closed ball $\mathbb{B}(\vzero,\alpha)$ centered at $\vzero$ and of radius $\alpha>0$.
Let $\Upsn$ and $\Upst$ be positive functions on $\Bc$. 
Then, as pointed out in \cite{CurnierA:88} (see also \cite{Chouly:14}),
the conditions \eqref{t.cont}-\eqref{tresca} are equivalent to the following statements:
\begin{alignat}{2}
\stressnn(\vu) &= \projRm{\taun(\vu)},&\qquad \taun(\vu)&\eqq\stressnn (\vu) - \Upsn \un, \label{t.contnitsche}\\
\stresst (\vu) &= \pTr{\taut(\vu)},&\qquad \taut(\vu)&\eqq \stresst (\vu) - \Upst \vut. \label{trescanitsche} 
\end{alignat}

Let $\Th$ be a simplicial mesh of $\Dom$. 
We assume that $\Dom$ is a polyhedron so that the mesh covers $\Dom$ exactly, and 
that every mesh boundary face belongs either to $\Bd$, $\Bn$, or
$\Bc$. The corresponding subsets of $\calFb$ are denoted by $\calFbD$, $\calFbN$, and $\calFbC$.
Let $\ThC\subset \Th$ be the collection of the mesh cells having
at least one boundary face on $\Bc$ and set $\dTC:=\dT\cap\Bc$ for all $T\in\ThC$.
In what follows, we need the following discrete trace inequality which is a slight variant of Lemma~\ref{lem:disc_inv} specialized to $\ThC$: There is $C_{\rm dt}$ such that for all $T\in\ThC$ and all $\vecteur{q}\in \polP^k_d(T;\Real^q)$, $q\in\{1,d\}$,
\begin{equation}
\|\vecteur{q}\|_{\bL^2(\dTC)} \le C_{\rm dt}h_T^{-\frac12}\|\vecteur{q}\|_{\bL^2(T)}.\label{eq:disc_trace_C}
\end{equation}

For the time being, we consider an $H^1$-conforming finite element subspace $\VhD\subset \VD$. Then, as shown in \cite{ChoHi:13,Chouly:14}, the FEM-Nitsche method leads to the discrete semilinear form $a^{\textsc{fem}}_h:\VhD\times\VhD\to\Real$ such that $a^{\textsc{fem}}_h(\SCAL;\SCAL)\eqq a(\SCAL,\SCAL)+n^{\textsc{fem}}_h(\SCAL;\SCAL)$ with
\begin{align}
&n^{\textsc{fem}}_h(\vv_h;\vw_h) \eqq  \label{eq:FEM-Nitsche}\\
&- \theta(\Upsn^{-1}\stressnn(\vv_h),\stressnn(\vw_h))_{L^2(\Bc)} + (\Upsn^{-1}\projRm{\taun(\vv_h)},(\taun+(\theta-1)\stressnn)(\vw_h))_{L^2(\Bc)} \nonumber \\
&-\theta(\Upst^{-1}\stresst(\vv_h),\stresst(\vw_h))_{\bL^2(\Bc)} + (\Upst^{-1}\pTr{\taut(\vv_h)},(\taut+(\theta-1)\stresst)(\vw_h))_{\bL^2(\Bc)}, \nonumber
\end{align}
$\stress(\vv_h)\eqq \stress(\strain(\vv_h))$, $\stress(\vw_h)\eqq \stress(\strain(\vw_h))$,
$\taun(\vv_h)$ and $\taut(\vv_h)$ defined as in~\eqref{t.contnitsche}-\eqref{trescanitsche},
and $\theta\in\{1,0,-1\}$ is a symmetry parameter. Choosing $\theta:=1$ 
leads to a symmetric formulation with a variational structure, 
choosing $\theta:=0$ is interesting to simplify the implementation by avoiding some 
terms in the formulation, and choosing $\theta:=-1$ allows one to improve on the stability
of the method by exploiting its skew-symmetry (see~\eqref{eq:cond_ell} where the lower bound vanishes for $\theta=-1$).

The discrete semilinear form $a^{\textsc{fem}}_h$ enjoys two key properties: (conditional) monotonicity and consistency. On the one hand, monotonicity holds true under a minimal condition on the penalty parameters. We assume that $\Upsn$ and $\Upst$ are piecewise constant on $\Bc$ with $\UpsnF\eqq \gamman h_{T_-}^{-1}$ and $\UpstF\eqq \gammat h_{T_-}^{-1}$ with positive parameters $\gamman$ and $\gammat$, for all $F\eqq \partial T_-\cap \Bc \in\calFbC$. Then, assuming that
\begin{equation}\label{eq:cond_ell}
\min(\varrho^{-1}\gamman,2\gammat) \ge 3(\theta+1)^2 C_{\rm dt}^2 \mu,
\end{equation}
with $\varrho:=\max(1,\frac{\lambda}{2\mu})$ and $C_{\rm dt}$ from~\eqref{eq:disc_trace_C},
we have (see Lemma~\ref{lem:monoton_HHO} for the arguments of the proof) 
\begin{equation}
a^{\textsc{fem}}_h(\vv_h;\bdelta_h)-a^{\textsc{fem}}_h(\vw_h;\bdelta_h)
\ge \frac13 a(\bdelta_h,\bdelta_h),
\end{equation}
for all $\vv_h,\vw_h\in \VhD$ with $\bdelta_h\eqq \vv_h-\vw_h$.
Concerning consistency, the key observation is that assuming that the exact solution satisfies $\vu\in \bH^{1+r}(\Dom)$, $r>\frac12$, we have $a^{\textsc{fem}}_h(\vu;\vw_h)=\ell(\vw_h)$ for all $\vw_h\in \VhD$. Indeed, integration by parts gives 
$a(\vu,\vw_h)-\ell(\vw_h)=(\stressn(\vu),\vw_h)_{\bL^2(\Bc)}$, while \eqref{t.contnitsche}-\eqref{trescanitsche} imply that
\begin{equation*}
n^{\textsc{fem}}_h(\vu;\vw_h)=(\stressnn(\vu),w_{h,n})_{L^2(\Bc)}+(\stresst(\vu),\vw_{h,t})_{\bL^2(\Bc)}
=(\stressn(\vu),\vw_h)_{\bL^2(\Bc)},
\end{equation*}
since $(\taun-\stressnn)(\vw_h)=\Upsn w_{h,n}$, $(\taut-\stresst)(\vw_h)=\Upst\vw_{h,t}$,
$\vw_h=w_{h,n}\NO+\vw_{h,t}$, and $\stressn(\vu)\SCAL\vw_h=\stressnn(\vu)w_{h,n}
+\stresst(\vu)\SCAL\vw_{h,t}$.

\subsection{Discrete setting for HHO-Nitsche} \label{sec: discrete_set}

The discrete setting for the HHO-Nitsche method is the same as for 
the linear elasticity problem
in Sect.~\ref{sec:lin_elas_setting}. 
As for FEM-Nitsche, 
we assume that every mesh boundary face belongs either to $\Bd$, $\Bn$, or 
$\Bc$, and the corresponding subsets of $\calFb$ are again denoted by $\calFbD$, 
$\calFbN$ and $\calFbC$. The HHO-Nitsche method uses the same key ideas as FEM-Nitsche:
the nonlinear reformulation \eqref{t.contnitsche}-\eqref{trescanitsche} of the contact and
friction conditions, and the weak enforcement of these nonlinear conditions by means of
a consistent boundary-penalty method inspired by Nitsche and originally developed in the 
context of HHO methods in \cite{CaChE:20}.

Our starting point is the equal-order HHO method devised for the linear elasticity problem,
where the discrete unknowns are polynomials of degree at most $k\ge1$ attached to
the mesh cells and to the mesh faces. One modification is that the 
degree of the face unknowns is raised to $(k+1)$ on the boundary faces in $\calFbC$.
This choice is motivated by the fact that these face unknowns are used 
to evaluate the quantities
$\taun$ and $\taut$ in Nitsche's formulation, so that the error estimate depends on how well these unknowns approximate the trace of the exact solution $\vu$ on $\Bc$. 
At the same time, this choice increases 
only marginally the computational costs.
For every mesh cell $T\in \Th$, let
$\FT$ be the collection of the mesh faces that are subsets of $\dT$,
which we partition as $\FT=\FTc\cup\FTnc$ with $\FTc \eqq \FT \cap \Fhbc$
(the subset $\FTc$ is empty for all $T\not\in\ThC$).
Then, the local HHO discrete space is
\begin{equation}
\hVkT \eqq \vPkd(T) \times \vVkkpdT, \quad
\vVkkpdT \eqq \bigtimes_{F\in\FTnc} \vPkF(F) \times  \bigtimes_{F\in\FTc} \vPkpF(F).
\end{equation}
A generic element in $\hVkT$ is denoted by $\hvT \eqq (\vT,\vdT)$. The discrete unknowns are illustrated in Fig.~\ref{fig_HHO_dofs_cont}.
\begin{figure}[htbp]
    \centering
    \subfloat[Pentagonal cell with no contact face ($\FTnc = \FT)$]{
        \centering 
       \includegraphics[scale=0.52]{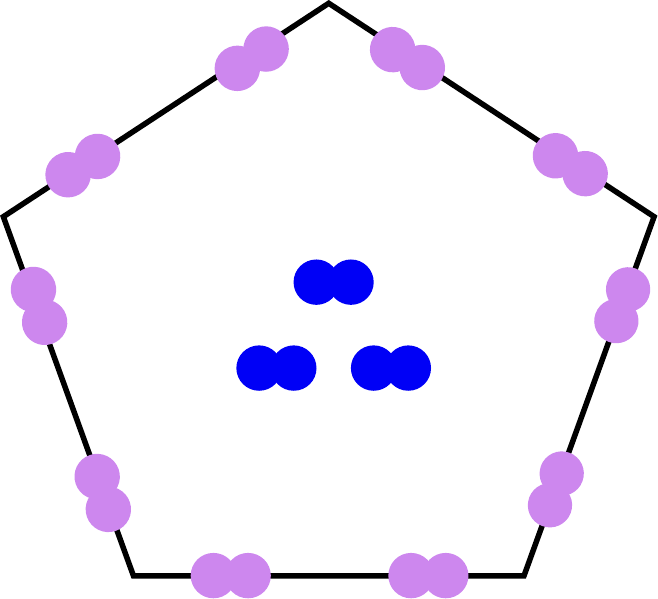} 
  }
    \qquad 
    \subfloat[Pentagonal cell with one contact face in red ($\FTnc \varsubsetneq \FT$)]{
        \centering
	    \includegraphics[scale=0.52]{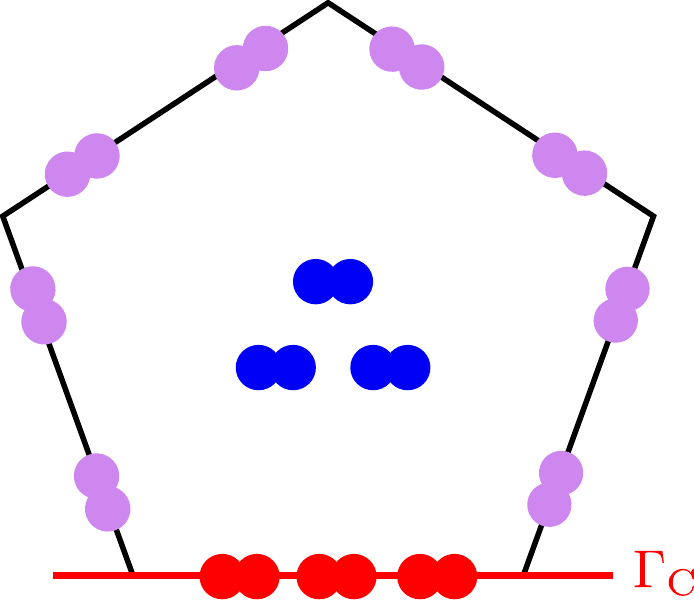}
    }
\caption{Face (red or violet) and cell (blue) unknowns in $\hVkT$ for $k=1$ and $d=2$ (each dot represents a basis function). }
\label{fig_HHO_dofs_cont}
\end{figure}

We consider as in Chap.~\ref{chap:wave} the local strain reconstruction operator $\opStr : \hVkT \rightarrow \Pkd(T;\Msym)$ such that for all $\hvT \in \hVkT$, 
\begin{equation}\label{eq_reconstruction_grad_tres}
\psm[T]{\opStr(\hvT)}{\matrice{q}} = {-\psm[T]{\vT}{\mdivergence \matrice{q}} + \psv[\dT]{\vdT}{\matrice{q} \nT}},
\end{equation}
for all $\matrice{q} \in \Pkd(T;\Msym)$.
The local discrete divergence operator $\opDiv: \hVkT \rightarrow \Pkd(T)$ is simply defined by taking the trace of the reconstructed strain tensor, i.e., for all $\hvT \in \hVkT$, we set
$\opDiv (\hvT ) \eqq \trace (\opStr(\hvT))$.
The local stabilization operator $\opStv: \hVkT \rightarrow  \vVkkpdT$ is readily adapted from the one considered for linear elasticity by setting for all $\hvT \in \hVkT$,
\begin{equation} \label{eq:stab_tres}
\opStv(\hvT) \eqq \PikkpdTv\Big(\vv_{T|\dT}
-\vdT+\big((I-\PikTv)\opDep(\hvT)\big)_{|\dT}\Big),
\end{equation}
where $\PikkpdTv$ is the $L^2$-orthogonal projections onto $ \vVkkpdT$ and the displacement reconstruction operator $\opDep : \hVkT \to \bpolP_d^{k+1}(T)$ is defined in \eqref{eq:def_Rec_HHO_elas}.
Using the above operators leads to the following local bilinear form defined on $\hVkT \times  \hVkT$:
\begin{align}\label{eq:a_T_hho}
a_T(\hvT,\hwT) \eqq {}& 2 \mu \psm[T]{ \opStr(\hvT)}{ \opStr(\hwT)} 
+ \lambda \pss[T]{ \opDiv(\hvT) }{ \opDiv(\hwT) } \nonumber \\
&+ 2 \mu \WT \psv[\dT]{\opStv(\hvT)}{\opStv(\hwT)}.
\end{align}

For simplicity, we employ the Nitsche technique only on the subset $\Bc$ where the nonlinear frictional contact conditions are enforced, whereas we resort to a strong enforcement of the homogeneous Dirichlet condition on the subset $\Bd$.
The global discrete spaces for the HHO-Nitsche method are
\begin{align}
\hVkh &\eqq  \bigtimes_{T\in\calT} \vPkd(T) \times \bigtimes_{F\in\calFi \cup \calFbN\cup\calFbD} \vPkF(F)  \times \bigtimes_{F\in\calFbC} \vPkpF(F). \label{eq:def_Ukh_tres} \\
\hVkhd &\eqq \{ \hvh \in \hVkh \: \vert \: \vF = \vecteur{0}, \; \forall F \in \Fhbd \},
\end{align}
leading to the notation $\hvh :=\big((\vT)_{T\in \Th}, (\vF)_{F\in\Fh}\big)$ for a generic element $\hvh\in \hVkh$. For all $T\in\Th$, we denote by $\hvT:=(\vT,\vdT\eqq (\vF)_{F\in\FT})\in\hVkT$ the local components of $\hvh$ attached to the mesh cell $T$ and the faces composing $\dT$, and for any mesh face $F\in\Fh$, we denote by $\vF$ the component of $\hvh$ attached to the face $F$. The global discrete bilinear form related to the linear elasticity part of the problem is, as usual, assembled cellwise by setting $a_h(\hvh,\hwh)\eqq \sum_{T\in\Th} a_T(\hvT,\hwT)$, and it remains to extend to the HHO setting the Nitsche-like semilinear form $n_h^{\textsc{fem}}$ defined in~\eqref{eq:FEM-Nitsche}. To this purpose, we set for all $T\in\ThC$ and all $\hwT\in \hVkT$,
\begin{equation} \label{eq:disc_stress_contact}
\stress(\hwT) \eqq 2\mu \opStr(\hwT)+ \lambda \opDiv( \hwT) \matrice{I}_d \in \Pkd(T; \Msym),
\end{equation}
with the decomposition 
$\stress(\hwT)\nT \eqq \stressnn(\hwT)\nT + \stresst(\hwT)$. 
Inspired by \eqref{t.contnitsche}-\eqref{trescanitsche}, we also introduce 
the linear operators $\taun:\hVkT\to \PkpF(\FTc)$ and $\taut:\hVkT\to \bpolP^{k+1}_{d-1}(\FTc)$ such that (notice the use of the face component on the right-hand side)
\begin{align}
\taun(\hwT) &\eqq \stressnn(\hwT) - \Upsn w_{\dT,n} , \qquad
\taut(\hwT) \eqq \stresst(\hwT) - \Upst \vecteur{w}_{\dT,t},
\end{align}
together with the decomposition $\vecteur{w}_{\dT}\eqq w_{\dT,n}\nT 
+ \vecteur{w}_{\dT,t}$ for the face component. We then set
$n_h^{\textsc{hho}}(\hvh,\hwh) \eqq \sum_{T\in\ThC} n^{\textsc{hho}}_T(\hvT,\hwT)$ 
for all $\hvh,\hwh\in\hVkh$ with
\begin{align}
& n_T^{\textsc{hho}}(\hvT,\hwT) \eqq \\ 
&-\theta (\Upsn^{-1}\stressnn(\hvT),\stressnn(\hwT))_{L^2(\dTC)} 
+ (\Upsn^{-1}\projRm{\taun(\hvT)},(\taun+(\theta-1)\stressnn)(\hwT))_{L^2(\dTC)} \nonumber \\
&-\theta (\Upst^{-1}\stresst(\hvT),\stresst(\hwT))_{\bL^2(\dTC)} 
+ (\Upst^{-1}\pTr{\taut(\hvT)})((\taut+(\theta-1)\stresst)(\hwT))_{\bL^2(\dTC)}, \nonumber
\end{align}
where $\theta\in\{1,0,-1\}$ is again the symmetry parameter. 
This leads to the following discrete HHO-Nitsche problem: 
\begin{equation}\label{eq:DP_N_Tresca} 
\left\lbrace 
\begin{alignedat}{2}
&\text{Find }\huh\in \hVkhd \text{ such that}\\	
&a_h^{\textsc{hho}}(\huh;\hwh) =  \ell_h(\hwh) \quad\forall \hwh\in \hVkhd,
\end{alignedat}
\right. 
\end{equation}
with $a^{\textsc{hho}}_h(\SCAL;\SCAL)\eqq a_h(\SCAL,\SCAL)+n^{\textsc{hho}}_h(\SCAL;\SCAL)$
and the linear form on the right-hand side is defined as
$\ell_h(\hwh)\eqq \psv[\Dom]{\loadext}{\vw_\calT}+\psv[\Bn]{\Gn}{\vw_\calF}$.

\bRem[Literature]
The above HHO-Nitsche method for contact and friction problems is devised and analyzed
in \cite{ChErPi:20}. This is, to our knowledge, so far 
the only discretization method supporting
polyhedral meshes that benefits from the same features as the FEM-Nitsche method devised
in \cite{ChoHi:13,Chouly:14}, namely optimal error estimates without additional assumptions
on the contact/friction set (see also \cite{ChHiR:15} for the analysis of FEM-Nitsche). 
Notice also
that \cite{ChErPi:20} tracks the dependency of the penalty parameters
and error estimates on the Lam\'e parameters $\mu$ and $\lambda$. Other polyhedral
discretization methods for contact/friction problems, that however do not hinge on Nitsche's approach,
include virtual element \cite{WrRuR:16,WangWei:18}, weak Galerkin \cite{GuGuZ:18}, and hybridizable discontinuous Galerkin \cite{ZhWuX:19} methods.
\eRem

\subsection{Stability and error analysis}
\label{sec:analysis}

In this section we outline the stability and error analysis for the above HHO-Nitsche method,
and we refer the reader to \cite{ChErPi:20} for more details. 

\bLem[Monotonicity, well-posedness] 
Assume that $\Upsn$ and $\Upst$ are piecewise constant on $\Bc$ with $\UpsnF\eqq \gamman h_{T_-}^{-1}$ and $\UpstF\eqq \gammat h_{T_-}^{-1}$ with positive parameters $\gamman$ and $\gammat$, for all $F\eqq \partial T_-\cap \Bc \in\calFbC$. Then, assuming that the minimality condition~\eqref{eq:cond_ell} on $\gamman$ and $\gammat$ holds true, we have
\begin{equation} \label{eq:monoton_HHO}
a^{\textsc{hho}}_h(\hvh;\hbdelta_h)-a^{\textsc{hho}}_h(\hwh;\hbdelta_h)
\ge \frac13 a_h(\hbdelta_h,\hbdelta_h),
\end{equation}
for all $\hvh,\hwh\in \hVkhd$ with $\hbdelta_h\eqq \hvh-\hwh$. Moreover,
the discrete problem~\eqref{eq:DP_N_Tresca} is well-posed.
\label{lem:monoton_HHO}
\eLem

\bproof
\textup{(i)} 
We have $n^{\textsc{hho}}_h(\hvh;\hbdelta_h)-n^{\textsc{hho}}_h(\hwh;\hbdelta_h)
= -\sum_{T\in\ThC}(A_{T,n}+A_{T,t})$ with
\begin{align*}
\frac{\gamman}{h_T}A_{T,n}&\eqq\theta\|\stressnn(\hbdelta_T)\|_{L^2(\dTC)}^2
-(\delta\taun,\taun(\hbdelta_T))_{L^2(\dTC)}
-(\theta-1)(\delta\taun,\stressnn(\hbdelta_T))_{L^2(\dTC)},\\
\frac{\gammat}{h_T}A_{T,t}&\eqq \theta\|\stresst(\hbdelta_T)\|_{\bL^2(\dTC)}^2
-(\delta\taut,\taut(\hbdelta_T))_{\bL^2(\dTC)}
-(\theta-1)(\delta\taut,\stresst(\hbdelta_T))_{\bL^2(\dTC)},
\end{align*}
with $\delta\taun\eqq \projRm{\taun(\hvT)}-\projRm{\taun(\hwT)}$
and $\delta\taut\eqq \pTr{\taut(\hvT)}-\pTr{\taut(\hwT)}$.
Using that $(\projRm{x} - \projRm{y})(x-y) \geq (\projRm{x} - \projRm{y})^2 \geq 0$ 
for all $x,y \in \Reel$, Young's inequality and the identity $\theta+\frac14(\theta-1)^2=\frac14(\theta+1)^2$ shows that
\[
A_{T,n} \le \frac14(\theta+1)^2\|\stressnn(\hbdelta_T)\|_{L^2(\dTC)}^2
\le \frac14(\theta+1)^2\frac{C_{\rm dt}^2}{\gamman} \|\stressnn(\hbdelta_T)\|_{L^2(T)}^2,
\]
where the last bound follows from the discrete trace inequality~\eqref{eq:disc_trace_C}.
Using the definition~\eqref{eq:disc_stress_contact} of the discrete stress, the triangle and Young's inequalities gives
\[
A_{T,n} \le (\theta+1)^2\frac{C_{\rm dt}^2}{\gamman} \mu\varrho \times
\big( 2\mu\|\opStr(\hbdelta_T)\|_{\bL^2(T)}^2+\lambda\|\opDiv(\hbdelta_T)\|_{L^2(T)}^2\big),
\]
recalling that $\varrho:=\max(1,\frac{\lambda}{2\mu})$. Using similar arguments, and in particular that $(\pTr{\vx} - \pTr{\vy})\SCAL (\vx-\vy) \geq | \pTr{\vx} - \pTr{\vy} |^2 \geq 0$ for all $\vx,\vy \in \Rd$, shows that 
\[
A_{T,t} \le (\theta+1)^2\frac{C_{\rm dt}^2}{\gammat} \mu \times
2\mu\|\opStr(\hbdelta_T)\|_{\bL^2(T)}^2.
\]
Putting these bounds together and using the condition~\eqref{eq:cond_ell} on the penalty parameters $\gamman$ and $\gammat$ proves that
\begin{align*}
n^{\textsc{hho}}_h(\hvh;\hbdelta_h)-n^{\textsc{hho}}_h(\hwh;\hbdelta_h) &\ge 
-\frac23 \sum_{T\in\ThC} \big( 2\mu\|\opStr(\hbdelta_T)\|_{\bL^2(T)}^2+\lambda\|\opDiv(\hbdelta_T)\|_{L^2(T)}^2 \big),
\end{align*}
so that $n^{\textsc{hho}}_h(\hvh;\hbdelta_h)-n^{\textsc{hho}}_h(\hwh;\hbdelta_h)\ge -\frac23 a_h(\hbdelta_h,\hbdelta_h)$.
This proves~\eqref{eq:monoton_HHO} since 
$a^{\textsc{hho}}_h(\hvh;\hbdelta_h)-a^{\textsc{hho}}_h(\hwh;\hbdelta_h)
=a_h(\hbdelta_h,\hbdelta_h)+n^{\textsc{hho}}_h(\hvh;\hbdelta_h)-n^{\textsc{hho}}_h(\hwh;\hbdelta_h)$.
\\
\textup{(ii)} Recalling~\eqref{eq:stab_prop_wave} shows that $a_h$ is coercive on $\hVkhd$ with respect to the norm $\snorme[\strain,h]{\hvh}^2\eqq \sum_{T\in\Th}\snorme[\strain,T]{\hvT}^2$ with $\snorme[\strain,T]{\hvT}^2 \eqq
\normem[\T]{\strain(\vT)}^2 + h_{T}^{-1}\normev[\dT]{\vT-\vdT}^2$. Therefore, combining the monotonicity property~\eqref{eq:monoton_HHO} with the arguments from \cite[Corollary 15, p.~126]{brezis-68} (see 
also \cite{Chouly:14}) proves that \eqref{eq:DP_N_Tresca} is well-posed.
\eproof

Let us finally state without proof an $\bH^1$-error estimate. Referring to \cite{ChErPi:20}
for more details, we observe that the bound on the consistency error combines the arguments from the proof of Lemma~\ref{lem:consist_HHO_elas} (for linear elasticity) and the arguments at the end of Sect.~\ref{sec:FEM-Nitsche} (for FEM-Nitsche). Let $\opStrh$ and $\opDivh$ be the global reconstruction operators such that $\opStrh(\hvh)_{|T}\eqq\opStr(\hvT)$ and $\opDivh(\hvh)_{|T}\eqq\opDiv(\hvT)$ for all $T\in\Th$ and all $\hvh\in \hVkhd$. Let $\matrice{\Pi}^l_\calT$, $l\in\{k,k+1\}$, denote the global $\bL^2$-orthogonal projection onto the corresponding piecewise polynomial space. 

\bTheo[$\bH^1$-error estimate] 
Assume that the penalty parameters satisfy the tighter condition 
\ifSp
$\min(\varrho^{-1}\gamman,2\gammat) \ge 3\big((\theta+1)^2 
+ (4+(\theta-1)^2)\big) C_{\rm dt}^2 \mu$.
\else
\[\min(\varrho^{-1}\gamman,2\gammat) \ge 3\big((\theta+1)^2 
+ (4+(\theta-1)^2)\big) C_{\rm dt}^2 \mu.\]
\fi
Let $\huh$ be the discrete solution of~\eqref{eq:DP_N_Tresca} with local
components $\huT$ for all $T\in\Th$.
Assume that the exact solution satisfies $\vu\in \bH^{1+r}(\Omega)$, $r>\frac12$. 
There is $C$, uniform with respect to $\mu$ and $\lambda$, such that
\begin{equation*}
2\mu\normem[\Dom]{\strain(\vu)-\opStrh(\huh)}^2 
+ \lambda\normes[\Dom]{\divergence{\vu}-\opDivh(\huh)}^2
\le C\big(\Psi^{\textsc{el}}(\vu)+\Psi^{\textsc{co}}(\vu)+\Psi^{\textsc{fr}}(\vu)\big),
\end{equation*}
with $\Psi^{\textsc{el}}(\vu)\eqq 2\mu|\strain(\vu)-\matrice{\Pi}^k_\calT(\strain(\vu))|_{\sharp,\calT}^2 + 2\mu\|\strain_\calT(\vu-\matrice{\Pi}_\calT^{k+1}(\vu))\|_{\bL^2(\Dom)}^2 + \lambda\varrho\|\DIV\vu-\Pi^k_\calT(\DIV \vu)\|_{\dagger,\calT}^2$ with the (semi)norms $|\SCAL|_{\sharp,\calT}$ and $\|\SCAL\|_{\dagger,\calT}$ defined in~\eqref{eq:sharp_dagger_elas}, and
\begin{align}
\Psi^{\textsc{co}}(\vu) &\eqq \sum_{T\in\ThC} \Big(\frac{h_T}{\gamman}
\|\stressnn(\vu)-\stressnn(\hIkT(\vu))\|_{L^2(\dTC)}^2
+ \frac{\gamman}{h_T} \|\delta u_{T,n}\|_{L^2(\dTC)}^2\Big), \label{eq:def_Psi_co}\\
\Psi^{\textsc{fr}}(\vu) &\eqq \sum_{T\in\ThC} \Big(\frac{h_T}{\gammat}
\|\stresst(\vu)-\stresst(\hIkT(\vu))\|_{\bL^2(\dTC)}^2
+ \frac{\gammat}{h_T} \|\delta \vu_{T,t}\|_{\bL^2(\dTC)}^2\Big),\label{eq:def_Psi_fr}
\end{align}
where the local reduction operator is defined such that
$\hIkT(\vu)\eqq (\PikTv(\vu),\PikkpdTv(\vu_{|\dT}))$,
and $(\delta u_{T,n},\delta \vu_{T,t})$ are the normal and tangential components of 
$\vu_{|\dT}-\PikkpdTv(\vu_{|\dT})$.
\label{th:Signorini_error}
\eTheo

An error estimate on the satisfaction of the contact/friction conditions is also
given in \cite[Thm.~12]{ChErPi:20}. Moreover, provided the exact solution satisfies
$\vu\in \bH^{1+r}(\calT)$ and $\divergence{\vu}\in H^{r}(\calT)$ with $r\in(\frac12,k+1]$,
Theorem~\ref{th:Signorini_error} implies that the $\bH^1$-error decays optimally with rate $\calO(h^{k+1})$. Notice however that in general, when there is a transition between contact and no-contact, the best expected regularity exponent is $r=\frac52-\varepsilon$, $\varepsilon > 0$, so that the maximal convergence rate is $\mathcal{O}(h^{\frac32-\varepsilon})$ 
and is reached for $k=1$. Finally, we notice that using face polynomials of degree $(k+1)$ on the faces in $\calFbC$ is crucial to estimate optimally the rightmost terms in~\eqref{eq:def_Psi_co}-\eqref{eq:def_Psi_fr}.

\bRem[Quasi-incompressible limit] 
In this situation, the factor $\varrho$ can be very large. 
The minimality condition~\eqref{eq:cond_ell} is robust with respect to the quasi-incompressible limit in the two following situations:
\textup{(i)} for the 
skew-symmetric variant $\theta=-1$, since the penalty parameters $\gamman$ and $\gammat$ need only to be positive real numbers 
(instead, for $\theta\in\{0,1\}$, this property is lost for
$\gamman$ which needs to scale as $\mu\varrho$); 
\textup{(ii)} for bilateral contact and any value of $\theta$,
since only the parameter $\gammat$ is used and its value
remains independent of $\varrho$. In contrast, the error estimate from Theorem~\ref{th:Signorini_error} is affected by large values of $\varrho$. The numerical experiments reported in~\cite{ChErPi:20} do not indicate, however, any sign of lack of robustness. 
\eRem

\section{Numerical example}\label{ssec::gv}

We consider a prototype for an industrial application that 
simulates the installation of a notched plug in a rigid pipe. The mesh is composed of 21{,}200 hexahedra and 510 prisms (for symmetry reasons, only one quarter of the pipe is discretized). The notched plug has a length of $56~\mm$ and an outer radius of $8~\mm$. The pipe is supposed to be rigid and has an inner radius of $8.77~\mm$ (there is an initial gap of $0.77~\mm$ between the plug and the pipe). The contact zone with Tresca's friction ($\Sc:=3{,}000~\MPa$) is between the rigid pipe and the ten notches of the plug. In the actual industrial setting, an indenter imposes a displacement to the upper surface of the plug. To simplify, sufficiently large vertical and horizontal forces are applied to the upper surface of the plug to impose a contact between the pipe and the notches.  The material parameters for the plug are $\mu:=80,769~\MPa$ and $\lambda:=121,154~\MPa$ (which correspond to a Young modulus $E=210,000~\MPa$ and a Poisson ratio $\nu=0.3$). The simulation is performed using $k\eqq1$, the symmetric variant $\theta:=1$, and the penalty parameters $\gamman=\gammat:=2\mu$). The discrete nonlinear problem \eqref{eq:DP_N_Tresca} is solved by a generalized Newton's method as in \cite{CurnierA:88}.
The von Mises stress is plotted in Fig.~\ref{fig:gv_VM} on the deformed configuration (a zoom on the contact zone is shown). 
We remark that there is contact between the notches and the pipe. Finally, the normal stress $\stressnn$ is visualized in Fig.~\ref{fig:gv_pres_cont} on the inferior surface of the plug.  We remark that 
all the notches are in contact except the first three (from left to right) and the last one (where $\stressnn=0$), and that a transition between contact and non-contact is located at the fourth notch.
Moreover, the maximal value of the normal stress is reached at the extremity of the notches.


\begin{figure}
   \centering
   \includegraphics[scale=0.14, trim =  50 20 20 390, clip=true]{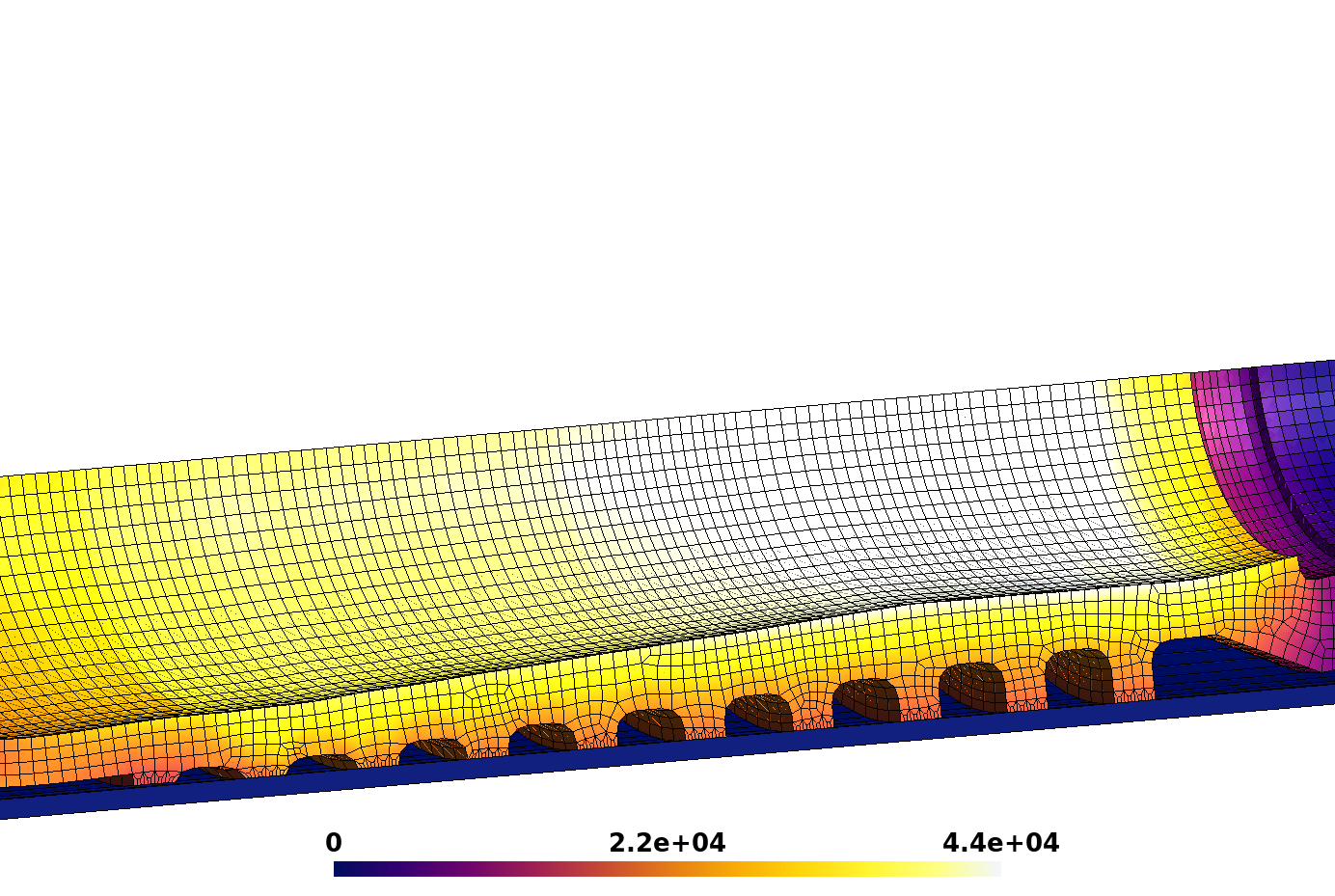} 
\caption{Notch plug (zoom on contact zone): von Mises stress (\MPa) on the deformed configuration.}
\label{fig:gv_VM}
\end{figure}

\begin{figure}
   \centering
   \includegraphics[scale=0.15, trim = 0 20 0 370, clip=true]{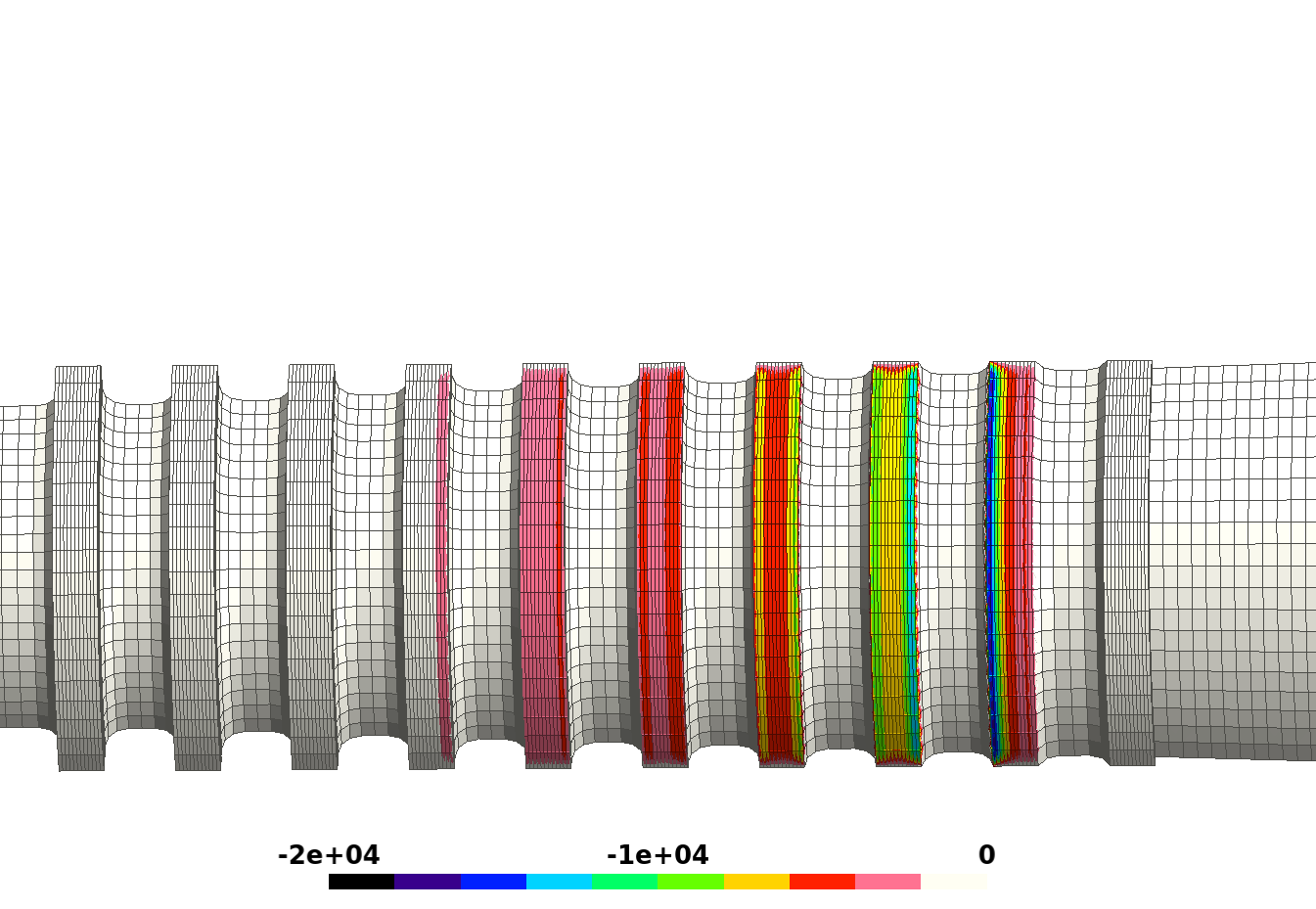} 
\caption{Notch plug: normal stress $\stressnn$ (\MPa) on the contact zone.}
\label{fig:gv_pres_cont}
\end{figure}

\chapter{Plasticity}\label{chap::plashpp}

Modeling plasticity problems is particularly relevant in nonlinear solid mechanics since plasticity can have a major influence on the behavior of a mechanical structure. One difficulty is that the plastic deformations are generally assumed to be incompressible, leading to volume-locking problems if (low-order) $H^1$-conforming finite elements are used. Mixed methods avoid these problems, but need additional globally coupled unknowns to enforce the incompressibility of the plastic deformations. 
Discontinuous Galerkin methods also avoid locking problems, but generally require to perform the integration of the behavior law at quadrature nodes located on the mesh faces, and not only in the mesh cells. In contrast, HHO methods are free of volume locking, only handle primal unknowns, and integrate the behavior law only at quadrature nodes in the mesh cells.

\section{Plasticity model}\label{sec::model}

Contrary to the elastic and hyperlastic models, the elastoplastic model is based on the assumption that the deformations are no longer reversible. We place ourselves within the framework of generalized standard materials \cite{Halphen1975,Lemaitre1994}. Moreover, the plasticity model is assumed to be strain-hardening (or perfect) and rate-independent, i.e., the speed of the deformations has no influence on the solution. For this reason, only the incremental plasticity problem with a pseudo-time is considered.

\subsection{Kinematics and additive decomposition}
We consider an elastoplastic material body that  occupies the domain $\Dom$ in the reference configuration.  Here, $\Dom \subset \Rd$, $d \in \{2,3\}$, is a bounded connected Lipschitz domain with unit outward normal $\bn$ and
boundary partitioned as $\front = \overline{\Bn} \cup \overline{\Bd}$
with two relatively open and disjoint subsets $\Bn$ and $\Bd$.  Due to the deformation, a point $\bx \in \Dom$ is mapped to a point $\bx'(t) = \bx + \vu(t,\bx)$ in the equilibrium configuration, where $\vu : J\times \Dom \rightarrow \Rd$ is the displacement field and $J$ is the pseudo-time interval. The deformation gradient $\Fdef(\vu) = \matrice{I}_d +\grad \vu$ takes values in $\Rdd_{+}$, which is the set of $\Rdd$-matrices with positive determinant. 

The regimes of infinitesimal and finite deformations are condidered here. For infinitesimal deformations, we consider the linearized strain tensor (see Sect.~\ref{sec:model-linear_elas})
\begin{equation}\label{eq:def_strain_inf_plas} 
\strain(\vu)\eqq \frac12(\grad\vu+\grad\vu\tr).
\end{equation}
For finite deformations, we adopt the logarithmic strain framework \cite{Miehe2002} leading to the following strain tensor:
\begin{equation}\label{eq:def_strain_finite_plas} 
\Elog(\vu) := \frac{1}{2} \ln\big(\Fdef(\vu)\tr \Fdef(\vu)\big) \qqe \mathcal{L}(\Fdef(\vu)),
\end{equation}
with the transformation $\mathcal{L}:\Rdd_{+}\to \Msym$ such that
$\mathcal{L}(\matrice{t})\eqq \frac12\ln(\matrice{t}\tr\matrice{t})$. Evaluating 
$\mathcal{L}(\matrice{t})$ requires to perform an eigenvalue decomposition of 
$\matrice{t}\tr\matrice{t}$.

Both strain tensors defined in~\eqref{eq:def_strain_inf_plas}-\eqref{eq:def_strain_finite_plas} are symmetric, and we notice that for infinitesimal deformations where 
$\| \grad\vu \|_{\ell^2} \ll 1$ with  $\| \grad\vu \|_{\ell^2}\eqq (\grad\vu:\grad\vu)^{\frac12}$, we have $\Elog(\vu)\approx\strain(\vu)$. 
To avoid the proliferation of cases, we work in this chapter
with the tensor $\Elog(\vu)$, keeping in mind that everything can be adapted to 
infinitesimal deformations.

\subsection{Helmholtz free energy and yield function}\label{ss:helmhotz}

In the framework of generalized standard materials, the material state is described locally
by the strain tensor $\Elog\in\Msym$ (we drop the dependency on $\vu$), the plastic strain tensor $\pElog\in\Msym$ which is trace-free, 
and a finite collection of internal variables $\IV:= ( \alpha_1, \ldots, \alpha_m) \in \Reel^m$.
The elastic strain tensor is then defined as follows:
\begin{equation}
\eElog \eqq \Elog - \pElog \in \Msym, \qquad \trace( \pElog )= 0.
\end{equation}
The Helmholtz free energy $\Psi: \Msym \times \Reel^m \rightarrow \Reel$ acts on a generic pair $( \matrice{e}, \ba)$ representing the elastic strain tensor and the internal variables. We assume that this function satisfies the following hypothesis.
\bHypo[Helmholtz free energy] 
$\Psi$ can be decomposed additively into an elastic and a plastic part as follows:\label{hypo:free_energy}
\begin{equation}
 \Psi(\matrice{e}, \ba) \eqq \frac{1}{2}\,  \matrice{e} : \elasticmodule:  \matrice{e} + \Psi^{\textrm{p}}( \ba)
 \end{equation}
where $\Psi^{\textrm{p}}:\Real^m\to\Real$ is convex (and strongly convex for strain-hardening plasticity), and the elastic modulus is $\elasticmodule \eqq 2\mu \tenseur{4}{I}^{\textrm{s}} + \lambda \matrice{I} \otimes \matrice{I}$, with $\mu>0$, $3\lambda + 2 \mu > 0$, $(\tenseur{4}{I}^{\textrm{s}})_{ij,kl}\eqq  \frac{1}{2}(\delta_{ik}\delta_{jl} + \delta_{il}\delta_{jk})$, and $(\matrice{I} \otimes \matrice{I})_{ij,kl}=\delta_{ij}\delta_{kl}$ for all $1\le i,j,k,l\le d$. The elastic modulus $\elasticmodule$ is isotropic, constant, and positive definite with $\matrice{e} : \elasticmodule : \matrice{e} = 2 \mu \, \matrice{e}: \matrice{e} +\lambda\trace(\matrice{e})^2$ for all $\matrice{e} \in \Msym$. \hfill$\square$
\label{hyp_free_ener}
\eHypo

Owing to the second principle of thermodynamics, the (logarithmic) stress tensor $\stressT \in \Msym$ and the internal forces $\IF \in \Reel^{m}$ are derived from $\Psi$ as follows:
\begin{align}
\stressT(\matrice{e}) \eqq \partial_{\matrice{e}} \Psi(\matrice{e}) = \elasticmodule : \matrice{e}, \qquad 
\IF(\ba) \eqq \partial_{\ba} \Psi^{\textrm{p}}(\ba). \label{eq_IF}
\end{align}
(Notice that $\stressT(\strain(\vu))=\elasticmodule : \strain(\vu)$ coincides with the usual stress tensor in the case of infinitesimal deformations and no plasticity.) 

The criterion to determine whether the deformations are plastic hinges on a scalar yield function $\Phi: \Msym \times \Reel^m \rightarrow \Reel$, which is a continuous and convex function of the stress tensor $\stressT$ and the internal forces $\IF$. The convex set of admissible states (or plasticity admissible domain) is
\begin{equation}
\mathcal{A} := \left\lbrace (\stressT, \IF) \in \Msym \times \Reel^m \: | \: \Phi(\stressT, \IF) \leq 0 \right\rbrace.
\end{equation}
This set is partitioned into the elastic domain $\mathcal{A}^{\textrm{e}}\eqq \left\lbrace (\stressT, \IF) \in\mathcal{A} \: | \: \Phi(\stressT, \IF) < 0 \right\rbrace=\interior(\mathcal{A})$ and the yield surface $\partial \mathcal{A}\eqq \left\lbrace (\stressT, \IF) \in \mathcal{A} \: | \: \Phi(\stressT, \IF) = 0 \right\rbrace$.
\bHypo[Yield function]
The yield function satisfies the following properties: \textup{(i)} $\Phi$ is piecewise analytical; \textup{(ii)} the point $(\matrice{0}, \vecteur{0})$  lies in the elastic domain, i.e., $\Phi(\matrice{0}, \vecteur{0}) < 0$; \textup{(iii)} $\Phi$ is differentiable at all points on the yield surface $\partial\mathcal{A}$. \hfill$\square$
\label{hyp_yield}
\eHypo

\bExp[Nonlinear isotropic hardening with von Mises yield criterion]
The internal variable is
$\IV\eqq p$, where $p\ge0$ is the equivalent plastic strain. The plastic part of the 
free energy is $\Psi^{\textrm{p}}(p) \eqq \sigma_{y,0}p +   \frac{H}{2} p^2 + (\sigma_{y,\infty} -\sigma_{y,0})(p - \frac{1-e^{-\delta p}}{\delta})$, where $H \geq 0$ is the isotropic hardening modulus, $ \sigma_{y,0} >0$, resp. $ \sigma_{y,\infty} \geq 0$, is the initial, resp. infinite, yield stress and $\delta \geq 0$ is the saturation parameter. The internal force is $q \eqq \sigma_{y,0} + Hp + (\sigma_{y,\infty} -\sigma_{y,0})(1 - e^{-\delta p})$. The perfect plasticity model is retrieved by taking $H\eqq 0$ and $\sigma_{y,\infty} \eqq \sigma_{y,0}$. Finally, the $J_2$-plasticity model with a von Mises criterion uses the yield function $\Phi(\stressT,q)\eqq\sqrt{3/2}\|\dev(\stressT)\|_{\ell^2}-q$, where  $\dev(\matrice{t})\eqq\matrice{t}-\frac{1}{d}\trace(\matrice{t})\matrice{I}_d$, $\|\matrice{t}\|_{\ell^2}\eqq (\matrice{t}:\matrice{t})^{\frac12}$ for any tensor $\matrice{t}\in\Rdd$.\label{exp:strain_vonMises}
\eExp

\subsection{Plasticity problem in incremental form}\label{ss:mech_model}

We are interested in finding the quasi-static evolution in the pseudo-time interval $J\eqq [0,T_{\textrm{f}}]$, $T_{\textrm{f}} > 0$, of the elastoplastic material body. We focus on the incremental form of the problem so that $J$ is discretized into $N$ subintervals defined by the discrete pseudo-time nodes $t^0\eqq 0 < t^1 < \ldots < t^N \eqq T_{\textrm{f}}$. The evolution occurs, for all $1\le n\le N$, under the action of a body force $\loadext^n:\Dom \rightarrow \Rd $, a traction force $\Gn^n: \Bn \rightarrow \Rd$ on the Neumann boundary $\Bn$, and a prescribed displacement $\vuD^n: \Bd \rightarrow \Rd$ on the Dirichlet boundary $\Bd$ ($\Bd$ has positive measure to prevent rigid-body motions). Recalling that $ \bH^1(\Dom) := H^1(\Dom;\Rd)$,
we denote by $\VD^n$, resp. $\Vz$, the set of all kinematically admissible displacements which satisfy the Dirichlet conditions, resp. homogeneous Dirichlet conditions on $\Bd$:
\begin{align}
 \VD^n = \left\lbrace \vv \in \bH^1(\Dom) \: | \: \vv_{|\Bd} = \vuD^n \right\rbrace,  \quad \Vz = \left\lbrace \vv \in \bH^1(\Dom) \: | \: \vv_{|\Bd} = \vecteur{0}  \right\rbrace.
\end{align}

It is customary to regroup the plastic strain tensor and the internal variables into the so-called generalized internal variables so that
\begin{equation} \label{eq:def_space_X}
\IVG \eqq ( \pElog, \IV )\in \IVGS \eqq 
\left\lbrace (\matrice{p}, \IV) \in  \Msym \times \Reel^{m}  \: | \: \trace(\matrice{p})= 0 \right\rbrace.
\end{equation}
The incremental plasticity problem proceeds as follows: For all $1 \leq n \leq N$, given $\vu^{n-1} \in \VD^{n-1}$ and $\IVG^{n-1}\eqq(\pElognm,\IV^{n-1}) \in  L^2(\Dom; \IVGS)$ from the previous pseudo-time step or the initial condition, find $\vu^n \in \VD^n$ and $\IVG^n\eqq(\pElogn,\IV^n) \in  L^2(\Dom; \IVGS)$ such that
\begin{align}
\label{weak_equil}
&\psm[\Dom]{ \PK^n }{ \GRAD \vw } = \ell^n(\vw)\eqq \psv[\Dom]{\loadext^n}{\vw} +  \psv[\Bn]{\Tn^n}{\vw}, \quad  \forall \vw\in \Vz,
\\
\label{weak_plas}
&( \IVG^n, \PK^n) \eqq  \texttt{PLASTICITY}(\IVG^{n-1}, \Fdef^{n-1}, \Fdef^n)\quad \text{pointwise in $\Dom$},
\end{align}
where $\Fdef^m\eqq \Fdef(\vu^m)$, $m\in\{n-1,n\}$. Letting $\Elog^m\eqq \mathcal{L}(\Fdef^m)$,
the procedure $\texttt{PLASTICITY}$ finds $\IVG^n$ and the Lagrange multiplier $\Lambda^n$ solving the following constrained nonlinear problem:
\begin{align}
&\pElogn-\pElognm=\Lambda^n \partial_{\stressT}\Phi(\stressT^n,\IF^n), \quad
\IV^n-\IV^{n-1}=-\Lambda^n \partial_{\IF}\Phi(\stressT^n,\IF^n), \label{eq:nonlin_opt_plas1}\\
&\Lambda^n\ge0, \quad \Phi(\stressT^n,\IF^n)\le0,\quad \Lambda^n\Phi(\stressT^n,\IF^n)=0,\label{eq:nonlin_opt_plas2}
\end{align}
where $\stressT^n\eqq \partial_{\matrice{e}} \Psi(\Elogn-\pElogn)=\elasticmodule:(\Elogn-\pElogn)$ and $\IF^n\eqq \partial_{\ba} \Psi^{\textrm{p}}(\IV^n)$. The first Piola--Kirchhoff stress tensor is then defined as $\PK^n \eqq \stressT^n : \partial_{\matrice{t}} \mathcal{L}(\Elogn)$, noting that for infinitesimal deformations, $\PK^n\approx \elasticmodule:(\Elogn-\pElogn)$.
One example of procedure for solving \eqref{eq:nonlin_opt_plas1}-\eqref{eq:nonlin_opt_plas2} is the standard radial return mapping \cite{Simo1992b, Simo1998}. For strain-hardening plasticity and infinitesimal deformations, the weak formulation \eqref{weak_equil}-\eqref{weak_plas} is well-posed, see \cite[Sect.~6.4]{Han2013}. For perfect plasticity, under additional hypotheses on the loads, the existence of a solution with bounded infinitesimal deformation is studied in \cite{Maso2006}. 

The incremental problem \eqref{weak_equil}-\eqref{weak_plas} can be reformulated as an incremental variational inequality by introducing a dissipative function \cite{Miehe2002,Djoko2007a}.
Given $(\vu^{n-1},\IVG^{n-1})\in \VD^{n-1}\times L^2(\Dom; \IVGS)$, 
we define the energy functional $\engy^n:\VD^n\times L^2(\Dom; \IVGS)\to \Real$ such that
\begin{equation}
\engy^n(\vv,\btheta) \eqq \int_\Dom \Psi^n(\Fdef(\vv),\btheta)\dx - \ell^n(\vv),
\end{equation}
with the incremental pseudo-energy density $\Psi^n:\Rdd_+\times\IVGS\to\Real$ such that
\begin{equation}
\Psi^n(\Fdef,\btheta) \eqq \Psi(\eElog,\IV)-\Psi(\eElognm,\IV^{n-1})+D_{\IVG^{n-1}}(\btheta),
\end{equation}
where $\btheta\eqq(\pElog,\IV)$, $\eElog\eqq \Elog-\pElog$ with $\Elog\eqq \mathcal{L}(\Fdef)$,
and with the incremental dissipation function $D_{\IVG^{n-1}}(\btheta)\eqq \sup_{(\stressT,\IF)\in\mathcal{A}} \big((\stressT:(\pElog-\pElognm)-\IF\cdot(\IV-\IV^{n-1})\big)$ ($D$ is convex and positively homogeneous of degree one). Then, a pair $(\vu^n,\IVG^n) \in \VD^n\times L^2(\Dom; \IVGS)$ solving~\eqref{weak_equil}-\eqref{weak_plas} satisfies the Euler--Lagrange equations of the minimization problem $\min_{(\vv,\btheta)\in \VD^n\times L^2(\Dom; \IVGS)} \engy^n(\vv,\btheta)$.

\section{HHO discretizations}
\label{sec:HHO_plas}

In this section, we present HHO methods to solve nonlinear plasticity problems. 

\subsection{Discrete unknowns}

Let $\Th$ be a mesh of $\Dom$ belonging to a shape-regular mesh sequence 
(see Sect.~\ref{sec:mesh} and~\ref{sec:mesh_reg}).
We assume that $\Dom$ is a polyhedron so that the mesh covers $\Dom$ exactly. 
Moreover, we assume that every mesh boundary face belongs either to $\Bd$ or to
$\Bn$. The corresponding subsets of $\calFb$ are denoted by $\calFbD$ and $\calFbN$. 
Recall that in HHO methods, the discrete unknowns are polynomials attached
to the mesh cells and the mesh faces. In the context of continuum mechanics,
both unknowns are vector-valued: the cell unknowns approximate the displacement field
in the cell, and the face unknowns approximate its trace on the mesh faces; see Figure~\ref{fig_HHO_dofs_elas}. 

For simplicity, we consider only the equal order-case where $k\ge1$ is the polynomial degree of both face and cell unknowns.
For every mesh cell $T\in\Th$, we set
\begin{equation}
\hVkT \eqq \vPkd(T) \times \vVkdT, \qquad
\vVkdT \eqq \bigtimes_{F\in\FT} \vPkF(F),
\end{equation}
with $\vPkd(T)\eqq \Pkd(T;\Real^d)$ and $\vPkF(F)\eqq \PkF(F;\Real^d)$.
A generic element in $\hVkT$ is denoted by $ \hvT\eqq (\vT,\vdT)$.
The HHO space is then defined as follows:
\begin{equation} 
\hVkh \eqq \Vkh \times \VkFh,
\qquad
\Vkh \eqq \bigtimes_{T\in\calT} \vPkd(T),
\qquad
\VkFh \eqq \bigtimes_{F\in\calF} \vPkF(F).
\end{equation}
A generic element in $\hVkh$ is denoted by $\hvh\eqq (\vTh,\vFh)$ with
$\vTh\eqq (\vv_T)_{T\in\calT}$ and $\vFh\eqq (\vv_F)_{F\in\calF}$, and we 
localize the components of $\hvh$ associated with a mesh cell $T\in\Th$ 
and its faces by using the notation
$\hvT\eqq \big(\vv_T,\vv_{\dT}\eqq (\vv_F)_{F\in\FT}\big)\in \hVkT$.
The Dirichlet boundary condition on the displacement field is enforced explicitly on the discrete unknowns attached to the mesh boundary faces in $\Fhbd$. Letting $\PikFv$ denote the $L^2$-orthogonal projection onto $\vPkF(F)$, we set
\begin{align}
\hVknhd &\eqq \left\lbrace \hvh \in \hVkh \: \vert \: \vF = \PikFv(\vu_{\textrm{D}}^n), \; \forall F \in \Fhbd   \right \rbrace, \\
\hVkhz &\eqq \left\lbrace \hvh \in \hVkh \: \vert \: \vF = \vecteur{0}, \; \forall F \in \Fhbd   \right \rbrace.
\end{align}

The discrete generalized internal variables are computed locally at the quadrature points of every mesh cell. We introduce the quadrature points $\Qp_T= (\Qp_{T,j})_{1 \leq j \leq m_{\textsc{q}}}$ and the weights $\vecteur{\Wp}_T= (\Wp_{T,j})_{1 \leq j \leq m_{\textsc{q}}}$, with $\Qp_{T,j} \in T$ and $\Wp_{T,j} \in \Reel$ for all $1 \leq j \leq m_{\textsc{q}}$ and all $T \in \Th$. We denote by $k_{\textsc{q}}$ the order of the quadrature. Then, the discrete generalized internal variables are sought in the space
\begin{equation}
\IVGTh := \bigcross_{T \in\Th}  \big(\underbrace{\IVGS\times \cdots \times \IVGS}_{\text{$m_{\textsc{q}}$ times}}\big),
\end{equation}
that is, for all $T\in\Th$, the generalized internal variables attached to $T$ form a vector $\IVG_T$ whose components are (a bit abusively) denoted by $(\IVG_T(\Qp_{T,j}))_{1\le j\le m_{\textsc{q}}}$ with $\IVG_T(\Qp_{T,j})\in \IVGS$ for all $1\le j\le m_{\textsc{q}}$.
In what follows, we use the following notation:
\begin{equation}\label{eq_prods}
\psmd[T]{\matrice{p}}{\matrice{q}} \eqq \sum_{j=1}^{m_{\textsc{q}}} \Wp_{T,j} \, \matrice{p}(\Qp_{T,j}) : \matrice{q}(\Qp_{T,j}),
\end{equation}
where, according to the context, the arguments can be either a continuous, tensor-valued function defined on $T$ or a vector in $(\Rdd)^{m_{\textsc{q}}}$. The global counterpart $\psmd[\Dom]{\matrice{p}}{\matrice{q}}$ is obtained by summing \eqref{eq_prods} over the mesh cells.

\subsection{Discrete plasticity problem in incremental form}

Recall the local gradient reconstruction $\opGRec:\hVkT \to \polP^{k}(T;\Rdd)$ defined in \eqref{eq:def_opGRec_mat} and the deformation gradient operator such that $\opFRec(\hvT) 
\eqq \matrice{I}_d + \opGRec(\hvT)$ for all $T\in\Th$. The global counterparts of these operators, 
which are defined in every mesh cell as above, are tensor-valued piecewise polynomials 
in $\polP^{k}(\calT;\Rdd)$ denoted by $\opGRech$ and $\opFRech$.
The global stabilization bilinear form  $s_h :
\hVkh \times \hVkh \to \Reel$ is defined in \eqref{eq:def_sh_elas} 
as for the linear elasticity problem, and we consider a positive weight $\beta_0>0$ 
(the choice $\beta_0=1$ was made for linear elasticity in Sect.~\ref{sec:disc_pb_elas}).

The discrete plasticity problem in incremental form proceeds as follows:
For all $ 1 \leq n \leq N$, given $\huh^{n-1} \in \hVknmhd$ and 
$\IVGh^{n-1}\eqq (\pElognm_\calT,\IV^{n-1}_\calT)\in \IVGTh$ from the previous pseudo-time step or the initial condition, find $\huhn \in \hVknhd$ and $\IVGh^n\eqq (\pElogn_\calT,\IV^{n}_\calT) \in \IVGTh$ such that
\begin{align}\label{discrete_problem_ptv_fp}
&\psmd[\Dom]{\PK_\calT^{n}}{\opGRech(\hwh)}
+ 2\beta_0\mu  s_h(\huh^n,\hwh) = \ell^n(\hwh),
\quad \forall \hwh \in \hVkhz, \\
\label{discrete_problem_plast_fp}
&(\IVG^n_{T,j},\PK^n_{T,j}) \eqq  \texttt{PLASTICITY}(\IVG^{n-1}_{T,j}, \Fdef^{n-1}_{T,j}, \Fdef^n_{T,j}), \; \forall T\in\Th,\forall j\in\{1,\ldots,m_{\textsc{q}}\},
\end{align}
where
$\IVG^m_{T,j}\eqq \IVG^m_T(\Qp_{T,j})$, $\Fdef^m_{T,j}\eqq \Fdef^m_\calT(\Qp_{T,j})$
with $m\in\{n-1,n\}$, and  $\PK^n_\calT(\Qp_{T,j})\eqq \PK^n_{T,j}$ for all $T\in\Th$ and all $j\in\{1,\ldots,m_{\textsc{q}}\}$.
Notice that the same procedure \texttt{PLASTICITY} is used as in the continuous setting.

\bRem[Litterature]
HHO methods for plasticity were developed in \cite{AbErPi:19,AbErPi:19a}.
Discontinuous Galerkin methods have been developed in \cite{Hansbo2010, Liu2010,Liu2013},
and virtual element methods in \cite{BdVLM:15, Artioli2017,Wriggers2017a}.
\eRem

\subsection{Nonlinear solver}
\label{sec:newton}

The nonlinear problem \eqref{discrete_problem_ptv_fp}-\eqref{discrete_problem_plast_fp} can be solved by using Newton's method. This requires evaluating the consistent (nominal) elastoplastic tangent modulus $\depmodule$ at every Gauss point in every mesh cell. The evaluation of $\depmodule$ can be included within the procedure \texttt{PLASTICITY}. To this purpose, we rewrite \eqref{weak_plas} as
\begin{equation}
( \IVG^n, \PK^n,\depmodule^n) \eqq  \texttt{PLASTICITY}(\IVG^{n-1}, \Fdef^{n-1}, \Fdef^n).
\end{equation}
Referring to the constrained nonlinear problem~\eqref{eq:nonlin_opt_plas1}-\eqref{eq:nonlin_opt_plas2} and recalling that $\elasticmodule$ denotes the (state-independent) elastic modulus (see \eqref{hypo:free_energy}), one first computes the infinitesimal elastoplastic tangent modulus $\depmoduleinf^n$ such that
\begin{equation}
\depmoduleinf^{n} \eqq \elasticmodule - \frac{(\elasticmodule:\partial_{\stressT}\Phi)\otimes (\elasticmodule:\partial_{\stressT}\Phi)}{\partial_{\stressT}\Phi:\elasticmodule:\partial_{\stressT}\Phi + \partial_{\IF}\Phi:\partial^2_{\ba\ba}\Psi^{\textrm{p}}:\partial_{\IF}\Phi},
\end{equation}
with the partial derivatives of $\Phi$ evaluated at $(\stressT^n,\IF^n)$ and the second derivative of $\Psi^{\textrm{p}}$ evaluated at $\IV^n$. Then, one sets
\begin{equation}
\depmodule^n \eqq (\partial_{\matrice{t}}\mathcal{L})\tr : \depmoduleinf^n : \partial_{\matrice{t}}\mathcal{L} + \stressT^n : \partial^2_{\matrice{t}\matrice{t}}\mathcal{L},
\end{equation}
where the partial derivatives of $\mathcal{L}$ are evaluated at $\Elog^n$.

Let $i\ge0$ be the index of the Newton's iteration and recall that $\huh^{n-1} \in \hVknmhd$ and 
$\IVGh^{n-1}\in \IVGTh$ are given from the previous pseudo-time step or the initial condition. The Newton's method is initialized 
by setting $\huh^{n,0}\eqq \huh^{n-1}$ (up to the update of the Dirichlet condition) 
and $\IVGh^{n,0}\eqq \IVGh^{n-1}$. Then, for all $i\ge0$, given $\huh^{n,i}\in \hVknhd$, one computes at each Newton's iteration 
the incremental displacement $\hduh^{n,i} \in \hVkhz$ such that
\begin{align}\label{eq_stiffness_matrix2}
&\psmd[\Dom]{\depmodulecalT^{n,i}{:} \opGRech(\hduh^{n,i})}{\opGRech(\hwh)} + 2\beta_0\mu s_h(\hduh^{n,i},\hwh) = -R_h^{n,i}(\hwh), \\
\label{eq:plasticity_Newton}
&(\IVG^{n,i}_{T,j},\PK^{n,i}_{T,j},\depmoduleTj^{n,i}) \eqq  \texttt{PLASTICITY}(\IVG^{n-1}_{T,j}, \Fdef^{n-1}_{T,j}, \Fdef^{n,i}_{T,j}),
\end{align}
where \eqref{eq_stiffness_matrix2} holds for all $\hwh \in \hVkhz$ with the residual term
\begin{equation}
R_h^{n,i}(\hwh) \eqq \psmd[\Dom]{\PK_\calT^{n,i}}{\opGRech(\hwh)} 
+ 2\beta_0\mu  s_h(\huh^{n,i},\hwh) - \ell^n(\hwh),
\end{equation}
and where \eqref{eq:plasticity_Newton} holds for all $T\in\Th$ and all $j\in\{1,\ldots,m_{\textsc{q}}\}$, with $\Fdef^{n-1}_{T,j}$, $\Fdef^{n,i}_{T,j}$ evaluated from $\huh^{n-1}$, $\huh^{n,i}$, respectively, and $\depmodulecalT^{n,i}(\Qp_{T,j})\eqq \depmoduleTj^{n,i}$.
At the end of each Newton's iteration, one updates the discrete displacement as 
$\huh^{n,i+1} = \huh^{n,i}+\hduh^{n,i}$. The discrete generalized internal variables do not need to
be updated at the end of the iteration, but only once Newton's method has converged.

For strain-hardening plasticity, the consistent elastoplastic tangent modulus is symmetric positive-definite. The following result gives some sufficient conditions for the linear system~\eqref{eq_stiffness_matrix2} to be coercive.

\begin{theorem}[Coercivity]\label{th::coer}
Assume the following: \textup{(i)} $k_{\textsc{q}}\ge 2k$ and all the quadrature weights are positive; \textup{(ii)} $\beta_0>0$; \textup{(iii)} the plastic model is strain-hardening. Let $\theta>0$ be the smallest eigenvalue of the fourth-order symmetric positive-definite tensors $(2\mu)^{-1}\depmoduleTj^{n,i}$ for all $T\in\Th$ and all
$j\in\{1,\ldots,m_{\textsc{q}}\}$. Then, the linear system~\eqref{eq_stiffness_matrix2} in each Newton's iteration is coercive, i.e., there is $\alpha>0$, independent of $h$, such that for all $\hvh \in \hVkhz$,
\begin{equation}\label{eq_stable}
\psmd[\Dom]{\depmodulecalT^{n,i}{:} \opGRech(\hvh)}{\opGRech(\hvh)}+  2\beta_0\mu s_h(\hvh,\hvh) \geq \alpha \min(\beta_0, \theta) 2\mu\norme[\hVkhz]{\hvh}^2,
\end{equation}
where $\norme[\hVkhz]{\hvh}^2\eqq \sum_{T\in\Th} \snorme[\hVkT]{\hvT}^2$ with $\snorme[\hVkT]{\hvT}^2\eqq
\normem[\T]{\GRAD \vvT}^2 + h_{T}^{-1}\normev[\dT]{\vvT-\vdT}^2$.
\end{theorem}
\bProof
Since the material is strain-hardening, we have $\theta >0$. Let $\hvh\in \hVkhz$.
Since $\opGRec (\hvT) \in \Pkd(T, \Rdd)$ for all $T\in\Th$,  
since all the quadrature weights are positive, and
$k_{\textsc{q}}\ge 2k$, we infer that
\begin{align*}
& \psmd[\Dom]{\depmodulecalT^{n,i}{:} \opGRech(\hvh)}{\opGRech(\hvh)}+  2\beta_0\mu s_h(\hvh,\hvh)\\
& \geq  \sum_{T\in\Th}   \sum_{1\le j\le m_{\textsc{q}}} 2\mu\theta\Wp_{T,j} \|\opGRec(\hvT)(\Qp_{T,j}) \|_{\ell^2}^2  + 2\beta_0\mu s_h(\hvh,\hvh) \\
&\geq 2\mu \min(\theta, \beta_0) \sum_{T\in\Th}  \Big(  \normem[T]{\opGRec (\hvT)}^2 + h_T^{-1}\normev[\dT]{\SdTk(\hvT)}^2 \Big).
\end{align*}
We conclude by using the stability result from Lemma~\ref{lem:stab_HHO_hyperelas}.
\eProof

\bRem[Choice of $\beta_0$]
Theorem~\ref{th::coer} indicates that the smallest eigenvalue $\theta$ is a natural target for the value of the weight parameter $\beta_0$ in the stabilization. A numerical study on
the influence of $\beta_0$ is presented in \cite[Sec.~5.3]{AbErPi:19a}. Another possibility
considered for virtual element methods in \cite{Wriggers2017a} is a piecewise constant stabilization parameter depending on the shape of the cell, the value of $\theta$, and a minimal user-defined value when $\theta \leq 0$.
\eRem

%
%

%
\section{Numerical examples}\label{sec::sec_numexp}
The goal of this section is to illustrate the above HHO method on two industrial applications where finite plasticity is present: a torsion of a square-section bar and an hydraulic pump under internal forces. For both examples, we use the nonlinear isotropic hardening model described in Example~\ref{exp:strain_vonMises}.

\subsection{Torsion of a square-section bar}\label{ss::torsion_fp}

This first example allows one to test the robustness of HHO methods under large torsion. The bar has a square-section of length $L\eqq 1~\mm$ and of height $H\eqq 5~\mm$ along the $z$-direction. The bottom end is clamped and the top end is subjected to a rotation of angle $\Theta=360^{\circ}$ around its center along the $z$-direction. The following material parameters ared used: Young modulus $E \eqq206.9~\GPa$, Poisson ratio $\nu\eqq0.29$, hardening parameter $H \eqq 129.2~\MPa$, initial yield stress $\sigma_{y,0} \eqq 450~\MPa$, infinite yield stress $\sigma_{y,\infty} \eqq 715~\MPa$, and saturation parameter $\delta \eqq 16.93$. The equivalent plastic strain $p$ is plotted at the quadrature points in Fig.~\ref{fig::torsion_p} for $k=2$. There is no sign of localization of the plastic deformations even for large rotations and large plastic deformations (around $50\%$). Moreover, the trace of the Cauchy stress tensor $\stress$ is plotted at the quadrature points on the final configuration in Fig.~\ref{fig::torsion_trace} for $k=2$. As expected, there is no sign of volume locking (no oscillation of the trace of the stress tensor, except at both ends which are fully constrained by Dirichlet conditions, so that stress concentrations are present).

\begin{figure}[htbp]
\centering
        \includegraphics[scale=0.44, trim=0 30 0 0]{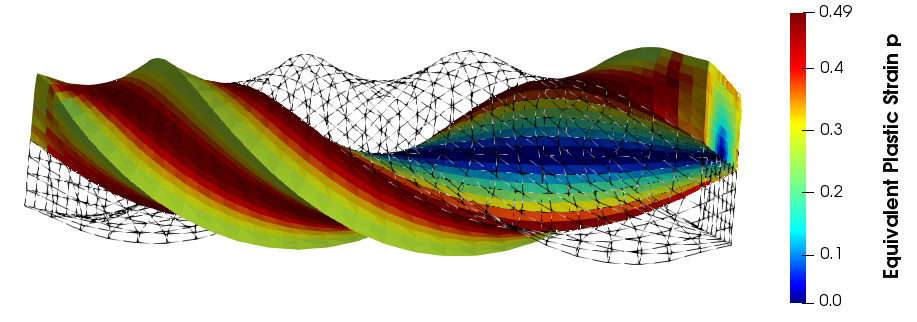}
    \caption{Torsion of a square-section bar: equivalent plastic strain $p$ at the quadrature points for $k=2$ and a rotation of angle $\Theta=360^{\circ}$.}
    \label{fig::torsion_p}
\end{figure}
\begin{figure}[htbp]
\centering
\ifSp
        \includegraphics[scale=0.5, trim= 0 0 0 90]{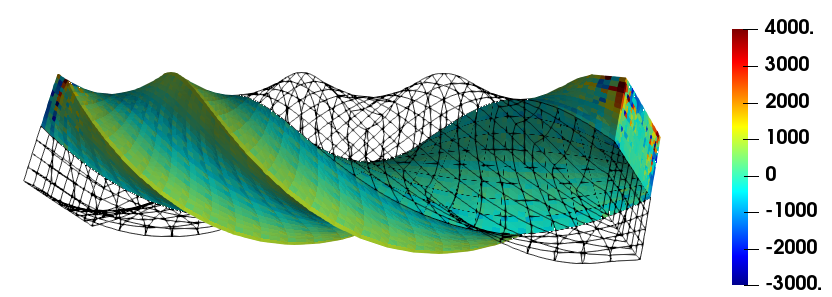}
\else
        \includegraphics[scale=0.5, trim= 0 0 0 0]{Figures/Fig_torsion_hho222_trace.png}
\fi
    \caption{Torsion of a square-section bar: trace of the Cauchy stress tensor $\stress$ (in $\MPa$) at the quadrature points for $k=2$ and a rotation of angle $\Theta=360^{\circ}$.}
    \label{fig::torsion_trace}
\end{figure}

\subsection{Hydraulic pump under internal forces}
This test case based on an industrial problem focuses on the deformation of an hydraulic pump and two of its pipes under the influence of a pressurized fluid. Since the study is restricted to the structural part of the problem, the force applied by the fluid on the walls of the pump and its pipes is replaced by an equivalent internal force. This surface force corresponds to a pressure of $14~\MPa$ in the reference configuration. Moreover, the bottom of the pump is clamped and the other surfaces are free. The description of the geometry and the mesh is given on the \CA~web site\footnote{Test PERF009: https://www.code-aster.org/V2/doc/default/fr/man\_v/v1/v1.01.262.pdf}.  Strain-hardening plasticity with a von Mises yield criterion is considered with the following material parameters: Young modulus $E \eqq 200~\GPa$, Poisson ratio $\nu\eqq 0.3$, hardening parameter $H \eqq 200~\MPa$, initial and infinite yield stresses $\sigma_{y,0} = \sigma_{y,\infty} \eqq  500~\MPa$, and saturation parameter $\delta \eqq 0$. The mesh is composed of 23,837 tetrahedra and 41,218 triangular faces.
The discrete global problem to solve has around 500,00 dofs for $k=1$. The Euclidean norm of the displacement and the equivalent plastic strain $p$ are plotted in Fig.~\ref{fig:pump_fp} on the deformed configuration. Note that the upper left part of the pump has the largest displacement.  Moreover, we remark that the plastic deformations are mainly present in the pipes and, in particular, at the junction between the pump and its pipes with nearly 97\% of equivalent plastic strain $p$. 
\begin{figure}[htbp]
    \centering
    \subfloat[Euclidean norm of the displacement]{
        \centering
        \includegraphics[scale=0.25, trim= 345 30 475 30, clip=true]{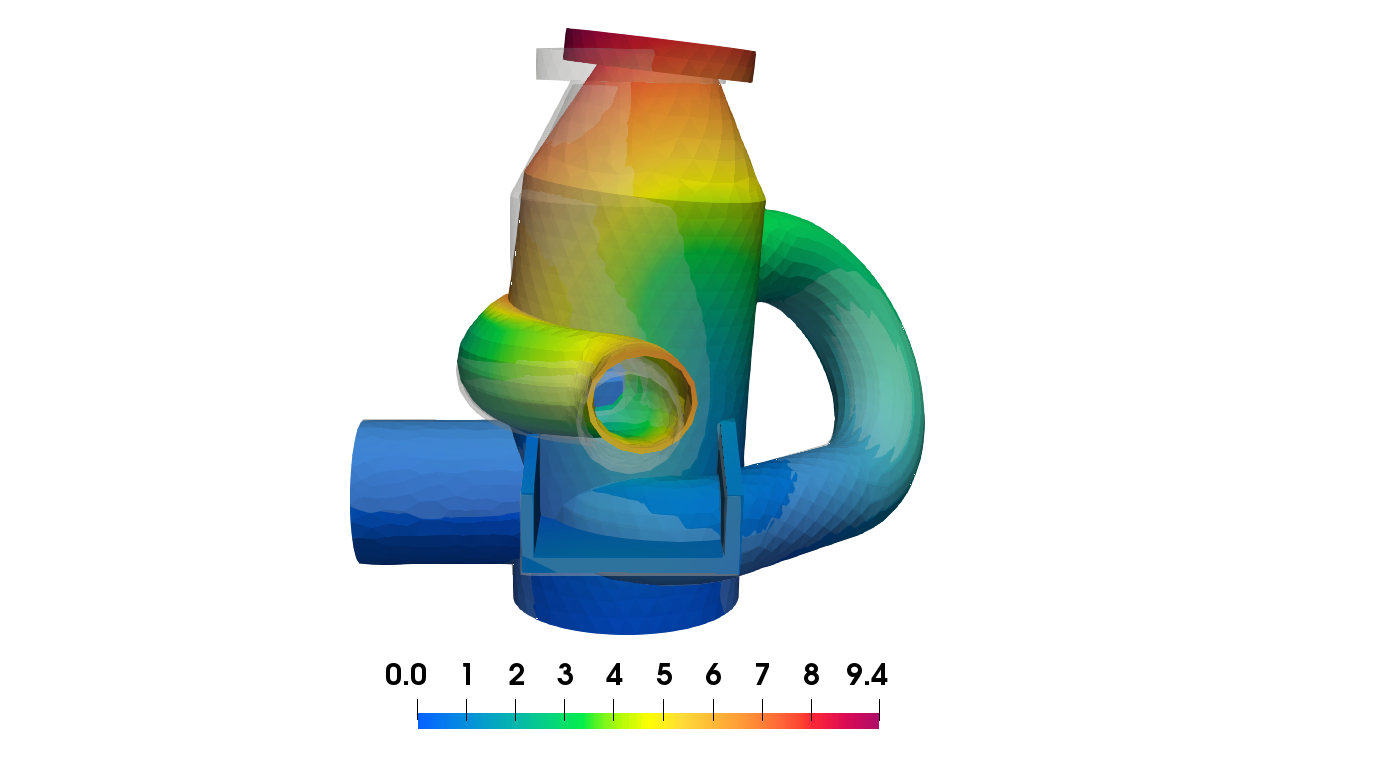}
  }
    ~ 
    \subfloat[Equivalent plastic strain $p$]{
        \centering
	    \includegraphics[scale=0.25,  trim= 340 20 480 0, clip=true]{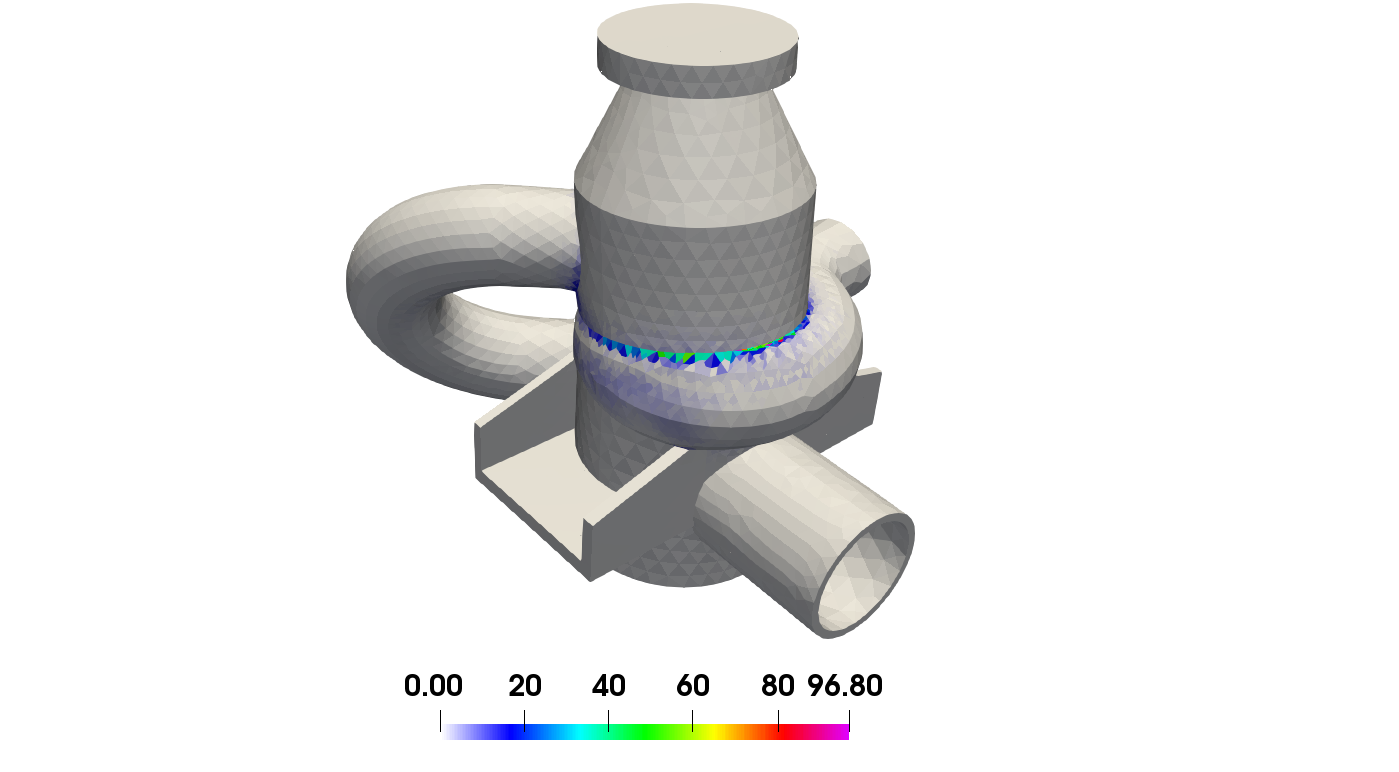}
    }
   \caption{Pump under internal forces:   (a) Euclidean norm of the displacement (in $\mm$) for $k=1$ on the deformed configuration with transparent reference configuration  (b) Equivelent plastic strain $p$ (in \%) on the deformed configuration.}
   \label{fig:pump_fp}
\end{figure}

\newcommand{\poly}[2]{\Poly^{#1}_{#2}}
\newcommand{\polysz}[2]{N^{#1}_{#2}}
\newcommand{\hhosz}[2]{\hat{N}^{#1}_{#2}}
\newcommand{\hhomixsz}[1]{\hat{N}^{k',k}_{#1}}
\newcommand*{\bigO}{\mathcal{O}}
\newcommand*{\shVlkT}{\hat{V}_{\T}^{\ell,k}}
\newcommand*{\shVlkappaT}{\hat{V}_{\T}^{\ell,\kappa}}
\newcommand*{\sVlT}{V_{\T}^{\ell}}
\newcommand*{\sVkappadT}{V_{\partial\T}^{\kappa}}
\newcommand*{\shIlkT}{I^{\ell,k}_T}
\newcommand{\hhodof}[1]{\mathsf{#1}}
\newcommand{\hhodofs}[1]{\bm{\mathsf{#1}}}
\newcommand{\stdvec}[1]{\bm{#1}}
\newcommand{\point}[1]{\bm{#1}}
\newcommand{\proj}[1]{\Pi^{#1}}
\newcommand{\eval}{\diamond}

\newcommand{\Mesh}{\calT}
\newcommand{\Meshes}{\calI}

\newcommand{\Cell}{{T}}
\newcommand{\Cells}{ \mathcal{T} }

\newcommand{\Face}{{F}}
\newcommand{\Faces}{ \mathcal{F} }

\newcommand{\Element}{{K}}

\newcommand{\F}{{\Face}}
\newcommand{\M}{{\Mesh}}
\newcommand*{\zeroF}{\bm{0}_{\Face}}

\newcommand*{\Rmn}{\mathbb{R}^{m \times n}}

\newcommand*{\Rbas}{\varrho}
\newcommand*{\Tbas}{\phi}
\newcommand*{\Fbas}{\phi}
\newcommand*{\Thhobas}{\mu}
\newcommand*{\TFhhobas}{\eta}

\newcommand*{\Rbasv}{\bm{\Rbas}}
\newcommand*{\Tbasv}{\bm{\Tbas}}
\newcommand*{\Fbasv}{\bm{\Fbas}}
\newcommand*{\Thhobasv}{\bm{\Thhobas}}
\newcommand*{\TFhhobasv}{\bm{\TFhhobas}}


\newcommand*{\Qset}{Q}
\newcommand*{\Qsetr}{\hat{Q}}
\newcommand*{\Qw}{\omega}
\newcommand*{\Qx}{\bm{x}}
\newcommand*{\Qwr}{\hat{\omega}}
\newcommand*{\Qxr}{\hat{\bm{x}}}
\newcommand*{\Qindex}{q}
\newcommand*{\Qwq}{\Qw_{\Qindex}}
\newcommand*{\Qxq}{\Qx_{\Qindex}}
\newcommand*{\Qindexo}{q_1}
\newcommand*{\Qwqo}{\Qw_{\Qindexo}}
\newcommand*{\Qxqo}{\Qx_{\Qindexo}}
\newcommand*{\Qwqr}{\Qwr_{\Qindex}}
\newcommand*{\Qxqr}{\Qxr_{\Qindex}}
\newcommand*{\Qdim}{|\Qset|}
\newcommand*{\matr}[1]{\mathsf{#1}}
\newcommand*{\stabhdg}{{S_{\partial T}^+}}
\newcommand*{\stabhdgF}{{S_F^+}}
\newcommand*{\stabhho}{{S_{\partial T}}}
\newcommand*{\stabhhoF}{{S_F}}
\newcommand*{\Vmix}{\hat{V}^{k',k}_T}
\newcommand*{\kcellmix}{k'}
\newcommand*{\kface}{k}
\newcommand*{\polyTm}{\poly{\kcellmix}{d}}
\newcommand*{\polyF}{\poly{\kface}{d-1}}
\newcommand*{\FacesT}{\calF_T}
\newcommand*{\PiTm}{{\Pi}^{\kcellmix}_T}

\chapter{Implementation aspects}

In this chapter, we outline the steps needed to bring the abstract formulation of the HHO method to an actual implementation. For simplicity, we focus on the Poisson model problem (see Chapter \ref{chap:diffusion}). We show how the local HHO operators (reconstruction and stabilization) are translated into matrices that can be used in the actual computation, and we give some criteria to test the implementation. Then we discuss the assembly of the discrete problem and the handling of the boundary conditions. We conclude with a brief overview on computational costs. Along the chapter, we provide some snippets of Matlab\textsuperscript{\tiny\textregistered}/Octave code to show a possible implementation (in 1D) of the critical parts.\footnote{The full source is available at \texttt{https://github.com/wareHHOuse/demoHHO}.} A 3D/polyhedral code called \texttt{DiSk++} fully supporting HHO and discontinuous Galerkin methods is downloadable at the address \texttt{https://github.com/wareHHOuse/diskpp}.\footnote{HHO methods are also implemented in the industrial software \texttt{code\_aster} \cite{codeaster} 
and the academic codes \texttt{SpaFEDTe} and \texttt{HArD{:}{:}Core} available on \texttt{github}.} We also refer the reader to \cite{CiDPE:18} for a description of the implementation of HHO methods using generic programming.

\section{Polynomial spaces} \label{sec:polyspaces}

The HHO method employs polynomials attached to the mesh cells and to the mesh faces. 
These polynomials are represented by their components in chosen polynomial bases. 
The evaluation of the cell basis functions can be done directly in the physical element by manipulating $d$-variate polynomials where $d\ge1$ is the space dimension. Instead, the evaluation of the face basis functions is done by means of affine geometric mappings that transform the $d$-dimensional points composing a face to a $(d-1)$-dimensional reference system associated with the face so that one manipulates $(d-1)$-variate polynomials; see \eqref{eqn:face_affine_mapping}. 

Let us consider first the cell basis functions. 
Let $k\ge0$ be the polynomial degree and  
recall that $\Pkd$ is composed of the $d$-variate polynomials of total degree at most $k$ with $\dim(\Pkd)={k+d\choose d}\qqe N^k_d$. Let $T\in\calT$ be a mesh cell and let $\{\phi_{T,i}\}_{1\le i\le N^k_d}$ be a basis of $\Pkd(T)$. Then, any polynomial $p\in \Pkd(T)$ can be decomposed in this basis as
\begin{align}
    p(\bm{x}) = \sum_{1\le i\le N^k_d} p_i \phi_{T,i}(\bm{x}), \label{eqn:poly-lincomb}
\end{align}
where the coefficients $p_i\in\Real$ are the components of $p$ in the chosen basis. 
These coefficients are the actual information that is stored and manipulated during the computations.
A simple and useful example of basis functions are the scaled monomials. Let 
$\bx_{T} = (x_{T,i})_{1\le i\le d} \in \Rd$ denote the barycenter of $T$ and $h_{T}$ its diameter. Recall that for a multi-index $\alpha\in\mathbb{N}^d$, $|\alpha|\eqq \sum_{1\le i\le d}\alpha_i$ denotes its length. Then, for all $\alpha\in\mathbb{N}^d$ with $|\alpha|\le k$, we set
\begin{align}
    \mu_{T,\alpha}(\mathbf{x}) \eqq \prod_{1\le i\le d}\left(\frac{2(x_i - x_{T,i})}{h_T}\right)^{\alpha_i}, \label{eqn:scaled_monomial}
\end{align}
leading to the basis $\{\mu_{T,\alpha}\}_{\alpha\in\mathbb{N}^d,|\alpha|\le k}$ of $\Pkd(T)$.
The two-dimensional scaled monomial basis is depicted 
in Figure~\ref{fig:basis2d} (up to degree 2 and with $h_T=2$ rather than $2\sqrt{2}$). The code in Listing~\ref{code:basis} implements (\ref{eqn:scaled_monomial}): the function evaluates the basis up to degree \texttt{max\_k} and its derivatives in the element with center \texttt{x\_bar} and size \texttt{h}. It returns two vectors containing the values of the basis functions and their derivatives at the point \texttt{x}.
\begin{listing}
\begin{minted}[mathescape, linenos, framesep=2mm]{matlab}
% Evaluate scalar monomial basis
function [phi, dphi] = basis(x, x_bar, h, max_k)
    k       = (0:max_k)';
    x_tilde = 2*(x-x_bar)/h;
    
    phi         = x_tilde .^ k;
    dphi        = zeros(max_k+1,1);
    dphi(2:end) = (2*k(2:end)/h).*(x_tilde.^k(1:end-1));
end
\end{minted}
    \caption{Possible implementation of a function evaluating the scaled monomial scalar basis and its derivatives in a 1D cell.}
    \label{code:basis}
\end{listing}

The face basis functions can be constructed in an analogous way by working on $\Real^{d-1}$ if $d\ge2$ and using the affine geometric mapping $\bT_F:\Real^{d-1}\to H_F$, where $H_F$ is the affine hyperplane in $\Real^d$ supporting $F$. In particular, scaled monomials can be built by using the point $\bT_F^{-1}(\bx_F)$, where $\bx_F$ is the barycenter of $F$.


%
\begin{figure}[htbp]
    \centering
    \includegraphics[width=0.75\textwidth]{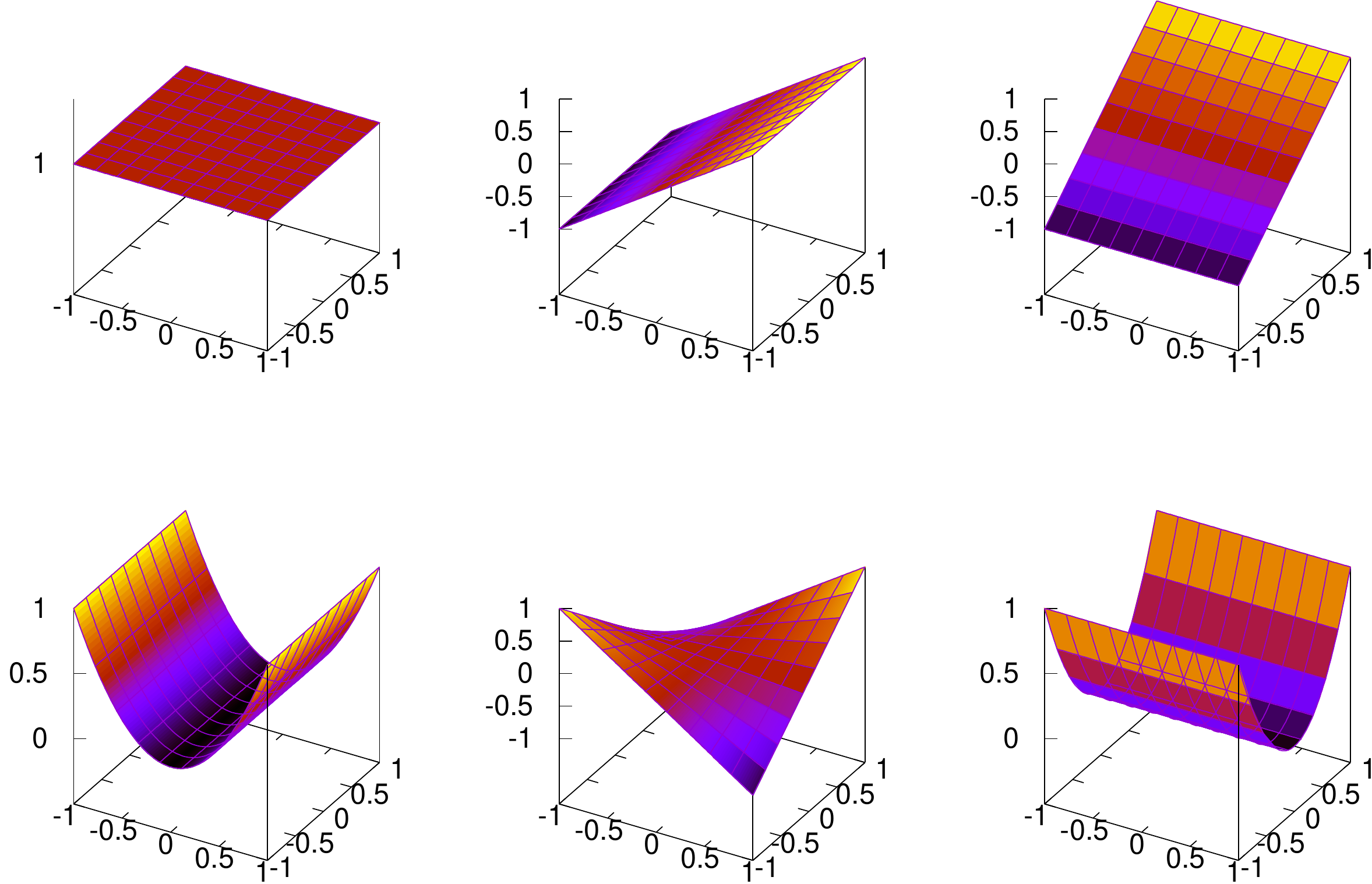}
    \caption{The set of functions of the scaled monomial basis of order 2 in the 2D element $(-1,1)^2$.}
    \label{fig:basis2d}
\end{figure}

\bRem[High order]
The choice of the basis functions is particularly important when working with high-order
polynomials, and its effects can be seen typically for $k\ge3$ and beyond (see, \eg \cite[Sect.~3.1]{HesHa:08} and \cite[Sect.~6.3.5\&Rmk.~7.14]{ErnGu:21a} for general discussions).
It can be beneficial to work with $L^2$-orthogonal bases. 
Such bases are easily devised
for $d=1$ using Legendre polynomials, and for $d\ge2$ if the cells are rectangular cuboids. If other shapes are used, an orthogonalization procedure can be considered, although it can be expensive.
One should bear in mind that the scaled monomial basis suffers from ill-conditioning for high polynomial degrees. 
\eRem

\bRem[Vector-valued case]
In continuum mechanics, HHO methods hinge on vector- and tensor-valued polynomials. 
Bases for
such polynomial spaces can be readily defined as tensor-products of a scalar polynomial basis and 
the Cartesian basis of $\Rd$ or $\Rdd$. For example, if we apply this procedure to $\poly{1}{2}$ with the basis $\{ 1, x, y \}$, we obtain the following vector-valued basis:
\[
    \begin{bmatrix} 1 \\ 0 \end{bmatrix},
    \begin{bmatrix} 0 \\ 1 \end{bmatrix},
    \begin{bmatrix} x \\ 0 \end{bmatrix},
    \begin{bmatrix} 0 \\ x \end{bmatrix},
    \begin{bmatrix} y \\ 0 \end{bmatrix},
    \begin{bmatrix} 0 \\ y \end{bmatrix}.
\]
The same procedure can be readily extended to tensor-valued polynomials.
\eRem

\section{Algebraic representation of the HHO space}

Let $T\in\calT$ be a mesh cell and let $\FacesT$ be the collection of its faces.
Let $k\ge0$ be the degree of the face polynomials. To allow for some generality, we let $k'\in\{k,k+1\}$ be the degree of the cell polynomials (the value $k'=k-1$ can also be considered for $k\ge1$). The local HHO space is
\begin{equation}
    \Vmix \eqq \polyTm(T) \times \left\{ \bigtimes_{F \in \FacesT} \polyF(\Face) \right\}. \label{eq:mixed_space_faces}
\end{equation}
The members of $\Vmix$ are of the form $\shvT\eqq (v_T, v_{F_1}, \ldots, v_{F_n})$, where $v_T\in\polyTm(T)$
and $v_{F_i}\in \polyF(\Face_i)$ for all $1\le i\le n\eqq\#\FacesT$. Notice that for $d=1$, the mesh faces coincide with the mesh vertices, so that the unknown associated with each face is a constant (see Sect.~\ref{sec:1D}); in this case, the degree of the cell unknowns is denoted by $k$. 
Having chosen bases for the above polynomial spaces, we collect all the coefficients in an array of size $\hhosz{k',k}{d}\eqq N^{k'}_d+nN^k_{d-1}$ structured as follows (see Figure~\ref{fig:dof-positioning}):
\begin{align}
    \hhodofs{v}_T := \left[ \hhodof{v}_{\T,1}, \ldots, \hhodof{v}_{\T,\polysz{\kcellmix}{d}} | \hhodof{v}_{\F_1,1}, \ldots, \hhodof{v}_{\F_1,\polysz{k}{d-1}} | \ldots | \hhodof{v}_{\F_n,1}, \ldots, \hhodof{v}_{\F_n,\polysz{k}{d-1}} \right]\tr, \label{eqn:dofsarray}
\end{align}
so that $\hhodofs{v}_T\in \hhodofs{V}^{k',k}_T\eqq \Real^{N^{k'}_d} \times (\Real^{N^k_{d-1}})^n$.
These coefficients are called \emph{degrees of freedom} (DoFs).
The structure of the array in~\eqref{eqn:dofsarray} will guide us in the understanding of the structure of the matrices realizing the HHO operators. 

\begin{figure}
    \centering
    \includegraphics[width=0.85\textwidth]{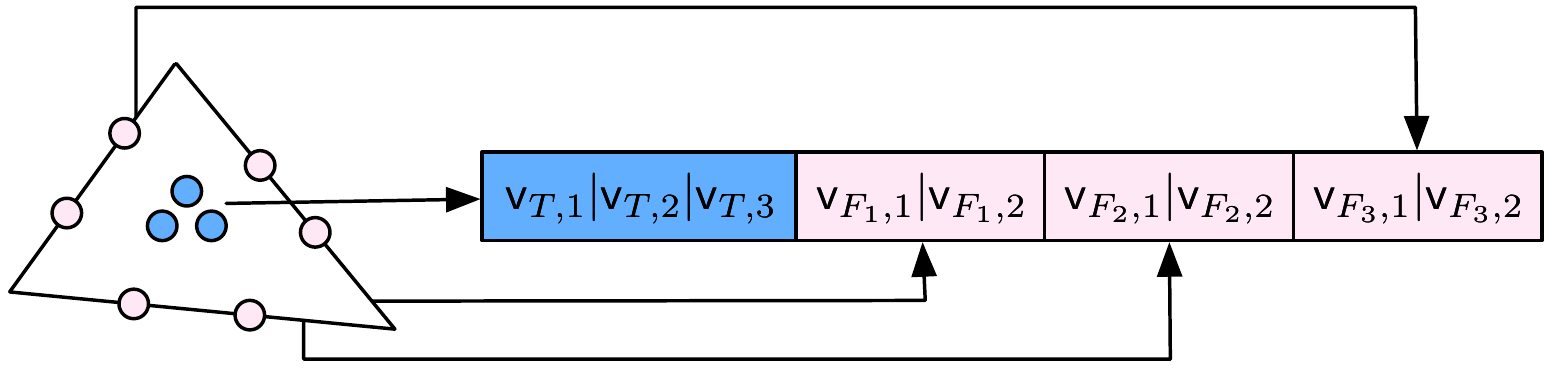}
    \caption{Formation of the local vector of DoFs in the case of a triangle with $\kcellmix = \kface \eqq 1$.}
    \label{fig:dof-positioning}
\end{figure}


\bRem[p-refinement]
The setting can be generalized to account for different polynomial orders on each face $F_i \in \FacesT$. This way, it becomes possible to use neighboring elements with different polynomial orders, opening the way to local $p$-refinement. The only required modification in the implementation is that the size of the sub-arrays in (\ref{eqn:dofsarray}) needs to account for the different polynomial degrees.
\eRem


\section{$L^2$-orthogonal projections}

$L^2$-orthogonal projections allow one to approximate functions belonging to a certain functional space with functions in a finite-dimensional polynomial space. Let us use a common notation $K \in \{T, F\}$ to denote a generic mesh cell or mesh face, with $d(T)\eqq d$ and $d(F)\eqq d-1$. 
Given a function $v \in L^2(K)$, its projection $\proj{k}_K(v)$ on $\poly{k}{d(K)}(K)$ is such that $(\proj{k}_K(v) - v, w)_{L^2(K)} = 0$
for all $w \in \poly{k}{d(K)}(K)$. For notational convenience, let $p_K \eqq \proj{k}_K(v)$. To compute $p_K$, we set up the problem
\begin{align}
\int_K \bigg(\sum_{i\in\mathcal{N}^k_K} w_i \phi_{K,i}(\point{x}) \sum_{j\in\mathcal{N}^k_K} p_{K,j} \phi_{K,j}(\point{x})\bigg) \dx = \int_K \bigg(v(\point{x}) \sum_{i\in\mathcal{N}^k_K} w_i \phi_{K,i}(\point{x})\bigg) \dx, \label{eqn:projection-expanded}
\end{align}
where $\mathcal{N}_K^k\eqq\{1,\ldots,\polysz{k}{d(K)}\}$ and 
the functions $\{\phi_{K,i}\}_{i\in\mathcal{N}^k_K}$ are a set of basis functions
attached to the geometric object $K$.
By defining similarly the coefficient column vectors $\bm{p}_K = \{p_{K,i}\}_{i \in \mathcal{N}^k_K}$ and $\bm{w}$, and the basis function column vector $\bm{\phi}_K(\bm{x}) = \{\phi_{K,i}(\bm{x})\}_{i \in \mathcal{N}^k_K}$, the expression (\ref{eqn:projection-expanded}) can be rewritten in matrix form as
\begin{align}
\stdvec{w}\tr \left(\int_K \stdvec{\phi}_K(\point{x})\stdvec{\phi}_K(\point{x})\tr\dx\right) \stdvec{p}_K = \stdvec{w}\tr \int_K v(\point{x})\stdvec{\phi}_K(\point{x})\dx. \label{eqn:proj-system}
\end{align}
Since $p_K-v$ has to be orthogonal to all the test functions $w$, $\stdvec{p}_K$ is found by setting up and solving the linear system of $\polysz{k}{d(K)}$ equations and $\polysz{k}{d(K)}$ unknowns
\begin{align}
\matr{M}_K\stdvec{p}_K = \int_K v(\point{x})\stdvec{\phi}_K(\point{x})\dx, \label{eqn:proj-problem}
\end{align}
with the mass matrix $\matr{M}_K := \int_K \stdvec{\phi}_K(\point{x})\stdvec{\phi}_K(\point{x})\tr\dx$ (by construction, $\matr{M}_K$ is symmetric positive-definite). One efficient way of solving the linear system~\eqref{eqn:proj-problem} is to compute the Cholesky decomposition of $\matr{M}_K$.

\subsection{Quadratures}

Integrals appearing in (\ref{eqn:proj-system}) are computed numerically using \emph{quadrature rules}. A quadrature rule allows one to approximate integrals over the geometric element $\Element$ as a weighted sum of evaluations of the integrand function $f$ at certain points in $K$:
\begin{align}
  \int_{\Element} f(\point{x})\dx \approx \sum_{\Qindex=1}^{|\Qset|} \Qwq f(\Qxq), \label{eqn:quadrature}
 \end{align}
where $\Qset$ is a set composed of $|\Qset|$ pairs $(\Qxq,\Qwq)$; for each pair, the first element is named \emph{quadrature point}, and the second element is named \emph{weight}. Quadratures are available for simplices, quadrilaterals, and hexahedra.
These quadratures are conceived in a reference cell and mapped to the physical cell by an affine geometric mapping. Quadratures allow exact integration of polynomials up to a certain degree called the \emph{quadrature order}. Integration on geometric objects having a more complex shape can be done by 
triangulating the geometric object and then employing a simplicial quadrature. 
Extensive literature about quadratures exists. Apart from the classical Gauss quadrature points \cite{AbrSt:72}, we mention \cite{grundmol1978,Dunavant1985,Keast1986} for quadratures on simplices and \cite{SomVi:07,SomVi:09,SudWa:13,ChLaS:15} for quadratures on polygons and polyhedra based on various techniques that avoid the need to invoke a sub-triangulation.

By using the tools just introduced, the linear system (\ref{eqn:proj-problem}) is set up numerically as
\begin{align}
\left(\sum_{\Qindex=1}^{\Qdim} \Qwq \stdvec{\phi}_K(\Qxq)\stdvec{\phi}_K(\Qxq)\tr\right)\stdvec{p}_K = \sum_{\Qindex=1}^{\Qdim} \Qwq v(\Qxq) \stdvec{\phi}_K(\Qxq), \label{eqn:massmatrix_using_quad}
\end{align}
where $\Qset$ needs to have a sufficient order to integrate exactly the product of the basis functions. For instance, if $\bm{\phi}_K$ is the vector of basis functions of $\poly{k}{d(K)}(K)$, the quadrature needs to have the sufficient number of points 
to integrate exactly polynomials of degree $2k$. 

\subsection{Reduction operator}
Let $T\in\Th$ be a mesh cell.
The local HHO reduction operator can be rewritten in expanded form as
\begin{equation}
\shIkpkT(v)\eqq (\PiTm(v),\proj{k}_{F_1}(v_{|F_1}), \ldots, \proj{k}_{F_n}(v_{|F_n})) \in \Vmix, \quad \forall v\in H^1(T).
\end{equation}
The reduction is thus the collection of the projections of $v$ on the cell $T$ and on its $n$ faces $F_1, \ldots, F_n$. At the algebraic level, this translates into obtaining the coefficients of $(n+1)$ polynomials by solving $(n+1)$ problems of the form (\ref{eqn:proj-problem}). More precisely, let $\Tbasv_T$ be the vector of cell-based basis functions on the mesh cell $T$ and let $\Fbasv_{F_i}$ be the vector of face-based basis functions on the $i$-th face of $T$. Moreover, let $\matr{M}_T$ and $\matr{M}_{F_i}$ be the corresponding mass matrices. The algebraic version of applying the reduction operator $\shIkpkT$ to a function $v\in H^1(T)$ amounts to finding the array vector $\bm{\mathsf{I}}^{k',k}_T(v)\in \hhodofs{V}^{k',k}_T$ solving the following block-diagonal system:
\begin{align}
\begin{bmatrix}
\matr{M}_T & & & \\
& \matr{M}_{F_1} & & \\
& & \ddots & \\
& & & \matr{M}_{F_n} \\
\end{bmatrix}\bm{\mathsf{I}}^{k',k}_T(v) =
\begin{bmatrix}
\int_{T}v(x)\Tbasv_T(x)\dx\\
\int_{F_1}v(x)\Fbasv_{F_1}(x)\dx\\
\vdots \\
\int_{F_n}v(x)\Fbasv_{F_n}(x)\dx\\
\end{bmatrix}.
\end{align}
Even though it is not used in the actual HHO computations, the computation of $\bm{\mathsf{I}}^{k',k}_T$ is essential to verify the correctness of the implementation of the reconstruction and stabilization operators detailed in the next section.

A possible implementation of the local reduction operator in 1D is shown in Listing~\ref{code:reduction}. The function \texttt{hho\_reduction()} takes the parameters \texttt{pd}, \texttt{elem} and \texttt{fun}, which are respectively a structure containing the computation parameters (in particular the polynomial degree and the cell diameter, which are taken here uniform on the whole mesh), the current element index, and the function to reduce.
At line 7, we ask for a quadrature, obtaining the points, the weights and the size in the variables \texttt{qps}, \texttt{qws}, and \texttt{nn}, respectively. We then proceed with the \texttt{for} loop (line 10) building the mass matrix and the right-hand side; this loop corresponds to the summations in (\ref{eqn:massmatrix_using_quad}). The projection on the cell is finally computed at line 16 (in 1D, we just need to evaluate the function at the faces): compare the structure of the returned vector \texttt{I} with \eqref{eqn:dofsarray}.

\begin{listing}[htp]
\begin{minted}[linenos, framesep=2mm, texcomments]{matlab}
% The HHO reduction operator
function I = hho_reduction(pd, elem, fun)
    % pd.h: cell diameter, uniform for all elements
    % pd.K: polynomial degree, equal for all elements
    x_bar = cell_center(pd, elem);
    [xF1, xF2] = face_centers(pd, elem);
    [qps, qws, nn] = integrate(2*pd.K, pd.h, elem);
    MM = zeros(pd.K+1, pd.K+1);
    rhs = zeros(pd.K+1, 1);
    for ii = 1:nn % This loop is the counterpart of (8.9)
        [phi, ~] = basis(qps(ii), x_bar, pd.h, pd.K);
        MM = MM + qws(ii) * (phi * phi');           % Mass matrix
        rhs = rhs + qws(ii) * phi * fun(qps(ii));   % Right-hand side
    end
    I = zeros(pd.K+3, 1);
    I(1:pd.K+1) = MM\rhs;   % Project on the cell
    I(pd.K+2) = fun(xF1);    % Project on faces: in 1D we just need
    I(pd.K+3) = fun(xF2);    %   to evaluate the function at the faces
end
\end{minted}
    \caption{Possible implementation of the reduction operator in 1D.}
    \label{code:reduction}
\end{listing}

\bRem[Verifying the implementation]
Let us consider a sequence of successively refined meshes $\meshfam\eqq (\calT_i)_{i\in\mathbb{N}}$ and let $h_i$ denote the maximum diameter of the cells composing $\calT_i$.
For each geometric object $\Element \in \{\Cell,\Face\}$ of $\Mesh_i \in \meshfam$, the projection on $\poly{k}{d(K)}(K)$ of a function $v \in H^1(\Dom)$ is computed by solving the problem (\ref{eqn:proj-problem}), obtaining a vector of DoFs $\bm{p}_K$. Such a vector is subsequently used to compute the global quantity
\begin{align*}
    \epsilon_i := \left( \sum_K \int_K \big(v - \proj{k}_K(v)\big)^2 \dx \right)^{1/2} = \left( \sum_K \sum_{\Qindex=1}^{|\Qset_K|} \Qwq \big(v(\Qxq) - \bm{\phi}_K(\Qxq)\tr\bm{p}_K\big)^2 \right)^{1/2},
\end{align*}
where $\Qset_K$ is a quadrature of sufficient order on $K$ and $\bm{\phi}_K$ is the vector of basis functions attached to $K$. The quantity $\epsilon_i$ has to decay, for increasing $i$, with rate $\bigO(h_i^{k+1})$ if the summation is over the mesh cells, whereas it has to decay with a rate of $\bigO(h_i^{k+1/2})$ if the summation is over the mesh faces (see Lemma \ref{lem:poly_approx}).
\eRem


\section{Algebraic realization of the local HHO operators}

Recalling Sect.~\ref{sec:heart}, the reconstruction and stabilization operators lie at the heart of HHO methods. Both operators are locally defined in every mesh cell $T\in\calT$ and map from the local HHO space $\Vmix$ to some polynomial space: the reconstruction operator maps 
to $\poly{k+1}{d}(T)$, and the stabilization operator restricted to each face $F\in\FacesT$ maps 
to $\poly{k}{d-1}(F)$. Since at the discrete level the elements of $\Vmix$ translate to vectors of the form \eqref{eqn:dofsarray}, both operators are represented by matrices that multiply a vector $\hhodofs{v} \in \hhodofs{V}^{k',k}_T$ to yield a vector representing either an element of $\poly{k+1}{d}(T)$ or of $\poly{k}{d-1}(F)$. This means that on a mesh cell with $n$ faces, both matrices have $\hhosz{k',k}{d} = \polysz{\kcellmix}{d} + n\polysz{\kface}{d-1}$ columns, which in turn form $n+1$ horizontally-juxtaposed blocks. We call $T$-block the first and leftmost block, whereas the remaining blocks are called $F_i$-blocks (see Figure~\ref{fig:opmatrix_structure}).

\begin{figure}[htbp]
    \centering
    \includegraphics[width=0.6\textwidth]{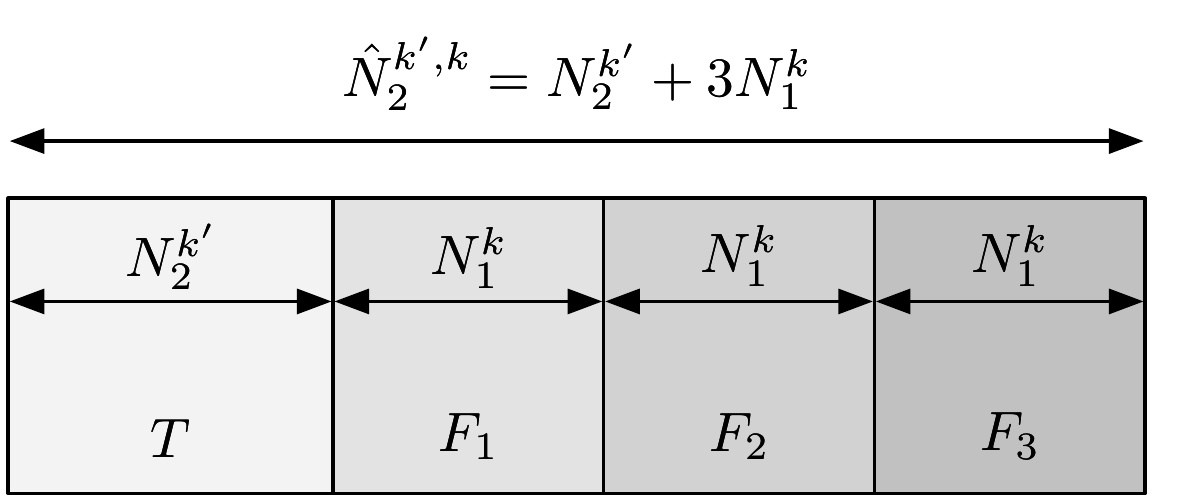}
    \caption{Block structure of an HHO operator matrix on a triangle ($d=2$, $n=3$ faces). In particular, there is one $T$-block and three $F_i$-blocks, giving an horizontal size of $\hhosz{k',k}{2} = \polysz{\kcellmix}{2} + 3\polysz{\kface}{1}$. The vertical size depends on the actual operator, as explained in the text.}
    \label{fig:opmatrix_structure}
\end{figure}

\subsection{Local reconstruction operator}
\label{sec:impl_rec}
Let $T\in\Th$. The local reconstruction operator satisfies (\ref{eq:def_Rec_HHO_bis}), where we expand here the boundary term as a summation on the faces of the mesh cell as follows:
\begin{equation} \label{eq:def_Rec_HHO_bis_copy}
\psv[T]{\GRAD \opRec(\shvT)}{\GRAD q} = \psv[T]{\GRAD \svT}{\GRAD q} - \sum_{F \in \FT}\pss[\F]{\svT-\svF}{\nT\SCAL\GRAD q},
\end{equation}
for all $q\in \Pkpd(T)^\perp\eqq \{q\in \Pkpd(T)\tq (q,1)_{L^2(T)}=0\}$. Moreover,
we have $\pss[T]{\opRec(\shvT)}{1} = \pss[T]{\svT}{1}$ (see \eqref{eq:def_Rec_HHO_mean}).
It is however not necessary to work with the polynomial space $\Pkpd(T)^\perp$, 
and one can consider any subspace $\poly{k+1}{*d}(T)$ leading to a
direct sum $\Pkpd(T) = \mathbb{P}_d^0(T)\oplus \poly{k+1}{*d}(T)$ (notice that
$\dim(\poly{k+1}{*d}(T))=\polysz{k+1}{*d}\eqq N^{k+1}_d-1$). One possibility
is to consider basis functions of $\Pkpd(T)$ such that the first basis function is constant,
and let the remaining basis functions span $\poly{k+1}{*d}(T)$.
Let $\bm{\varrho}(\point{x})$ be the vector of basis functions of $\poly{k+1}{*d}(T)$. 
Using a quadrature $\Qset_T$ of order at least $2k$, the left-hand side of \eqref{eq:def_Rec_HHO_bis_copy} is a plain stiffness matrix such that
\begin{align}
    \matr{K_*} \eqq \sum_{\Qindex=1}^{|\Qset_T|} \Qwq \GRAD\bm{\varrho}(\Qxq) \cdot \GRAD\bm{\varrho}(\Qxq)\tr,
\end{align}
where $\GRAD$ is applied componentwise to $\bm{\varrho}(\point{x})$ and the dot product only to the gradients. Notice that this computation results in a standard stiffness matrix, where the column and the row corresponding to the constant basis function have been dropped.


We next build the right-hand side of~\eqref{eq:def_Rec_HHO_bis_copy} in multiple steps. For simplicity, we assume that we are building the operator for a triangular element, so that $n \eqq 3$ and $d \eqq 2$. Let $\Tbasv_T(\bm{x})$ be the column vector of cell-based basis functions attached to $\Cell$, $\Fbasv_{F_i}(\bm{x})$ the column vector of face-based basis functions attached to the face $F_i$ (recall that these basis functions are computed using a geometric mapping from $\Real^{d-1}$ to the hyperplane supporting $F_i$) and $\bm{0}_F$ a zero column vector of size $\polysz{k}{d-1}$.
We start with $\psv[T]{\GRAD \svT}{\GRAD q}$, where $\svT \in \poly{\kcellmix}{d}(T)$ and $q \in \poly{k+1}{*d}(T)$. In order to evaluate the cell-based part of a DoFs vector of the form (\ref{eqn:dofsarray}), we form a column vector of basis functions $\bm{\mu}(\bm{x}) := [\Tbasv_T(\bm{x})\,|\,\zeroF\,|\,\zeroF\,|\,\zeroF\,]$, where $|$ denotes the vertical concatenation of column vectors. Then, we form the matrix
\begin{align}
    \matr{T} \eqq \sum_{\Qindex=1}^{|\Qset_T|} \Qwq \GRAD\bm{\varrho}(\Qxq) \cdot \GRAD\bm{\mu}(\Qxq)\tr.
\end{align}
This computation yields a matrix where only the $T$-block has nonzero values. Its effect can be intuitively understood by looking separately at the roles of $\GRAD\bm{\mu}(\Qxq)\tr$ and $\GRAD\bm{\varrho}(\Qxq)$ when $\matr{T}$ multiplies a vector $\hhodofs{v} \in \hhodofs{V}_T^{k', k}$. For each quadrature point, $\bm{\mu}$ evaluates the gradients of the cell-based part of $\hhodofs{v}$, whereas $\bm{\varrho}$ tests the value of the polynomial with the gradients of the basis functions of the reconstruction space. In practice, this returns the right-hand side of a projection-like problem where the gradients of $\poly{k+1}{*d}(T)$ are used as test functions.

We continue with the contributions from $\pss[F]{\svT-\svF}{\nT\SCAL\GRAD q}$ for the $n=3$ faces of $T$. For example, in order to compute the contribution due to the face $F_1$, we consider the vector of basis functions $\bm{\eta}_1(\bm{x}) := [\Tbasv_T(\bm{x})\,|\,-\Fbasv_{F_1}(\bm{x})\,|\,\zeroF\,|\,\zeroF\,]$. The contribution to the right-hand side is then computed as
\begin{equation}
    \matr{F}_1 \eqq \sum_{\Qindexo=1}^{|\Qset_{F_1}|} \Qwqo \bm{n}\SCAL\GRAD\bm{\varrho}(\Qxqo) \bm{\eta}_1(\Qxqo)\tr.
\end{equation}
Indeed, taking an array of the form (\ref{eqn:dofsarray}) representing a member of $\Vmix$ and computing the dot-product with $\bm{\eta}_1(\Qxqo)$ corresponds to obtaining the value of the difference of the cell-based and $F_1$-based polynomials at the point $\Qxqo\in F_1$. Notice also that the matrix $\matr{F}_1$ contains nonzero elements only in the $T$-block and in the $F_1$-block. The matrices $\matr{F}_2$ and $\matr{F}_3$ are computed in a similar fashion by taking $\bm{\eta}_2(\bm{x}) := [\Tbasv_T(\bm{x})\,|\,\zeroF\,|\,-\Fbasv_{F_2}(\bm{x})\,|\,\zeroF]$ and $\bm{\eta}_3(\bm{x}) := [\Tbasv_T(\bm{x})\,|\,\zeroF\,|\,\zeroF\,|\,-\Fbasv_{F_3}(\bm{x})]$, respectively. If cells with more than three faces are used, the procedure is easily generalized by computing the remaining $\matr{F}_i$ matrices.


We finally compute the algebraic realization $\matr{R}$ of $\opRec$ (up to the mean-value constraint) by inverting the matrix $\matr{K}_*$ and setting
\begin{align}
\matr{R} := \matr{K}_{*}^{-1}\matr{H} \quad \text{with} \quad \matr{H}\eqq \matr{T} - \sum_{i=1}^3 \matr{F}_i, 
\label{eqn:discoperR}
\end{align}
and the mean-value constraint can be satisfied by adding a suitable contribution from the constant basis function (and increasing by one the size of the vector $\matr{R}$).
Once we have computed $\matr{R}$, we can readily obtain the matrix representing the stiffness term in \eqref{eq:local_aT} as
\begin{align}
    \matr{A} := \matr{R}\tr\matr{K_*}\matr{R} = \matr{H}\tr\matr{R}. \label{eqn:lapl_lc}
\end{align}
Take a moment to analyze the roles of the matrices composing $\matr{A}$. $\matr{K_*}$ is a plain stiffness matrix on $\poly{k+1}{*d}(T)$ and, as such, it operates on polynomials in $\poly{k+1}{*d}(T)$ to compute a standard local stiffness term. In HHO however, DoFs live in the space $\hhodofs{V}^{k',k}_T$: the reconstruction matrix $\matr{R}$ ``translates'' HHO DoFs to the higher-order space $\poly{k+1}{*d}(T)$, on which $\matr{K_*}$ can operate. 
Listing \ref{code:reconstruction} shows a possible realization of the computation of $\matr{R}$ in 1D. At lines 9-13, the stiffness matrix of $\poly{k+1}{d}(T)$ is computed using a quadrature of order $2k$. It is subsequently trimmed to obtain $\matr{K}_{*}$ (line 16) and $\matr{T}$ (line 18). Starting from line 24, the boundary terms are computed. Finally, the reconstruction operator and the matrix $\matr{A}$ are obtained at lines 33 and 34, respectively. An illustration of the action of the reconstruction operator is shown in Figure~\ref{fig:reconstruction-demo}.

\begin{listing}
\begin{minted}[mathescape, linenos, framesep=2mm]{matlab}
% The HHO reconstruction operator
function [A, R] = hho_reconstruction(pd, elem)
    x_bar = cell_center(pd, elem);
    [xF1, xF2] = face_centers(pd, elem);

    stiff_mat = zeros(pd.K+2, pd.K+2);
    gr_rhs = zeros(pd.K+1, pd.K+3);
    
    [qps, qws, nn] = integrate(2*pd.K, pd.h, elem); 
    for ii = 1:nn
        [~, dphi] = basis(qps(ii), x_bar, pd.h, pd.K+1);
        stiff_mat = stiff_mat + qws(ii) * (dphi * dphi');
    end
    
    % Set up local Neumann problem
    gr_lhs = stiff_mat(2:end, 2:end); % Left-hand side
    % Right-hand side, cell part
    gr_rhs(:,1:pd.K+1) = stiff_mat(2:end,1:pd.K+1); % $\psv[T]{\GRAD \svT}{\GRAD q}$
    
    [phiF1, dphiF1] = basis(xF1, x_bar, pd.h, pd.K+1);
    [phiF2, dphiF2] = basis(xF2, x_bar, pd.h, pd.K+1);
    
    % Right-hand side, boundary part
    gr_rhs(1:end, 1:pd.K+1) = gr_rhs(1:end, 1:pd.K+1) + ...
        dphiF1(2:end)*phiF1(1:pd.K+1)'; % $\pss[F_1]{\svT}{\nT\SCAL\GRAD q}$
    
    gr_rhs(1:end, 1:pd.K+1) = gr_rhs(1:end, 1:pd.K+1) - ...
        dphiF2(2:end)*phiF2(1:pd.K+1)'; % $\pss[F_2]{\svT}{\nT\SCAL\GRAD q}$
    
    gr_rhs(1:end, pd.K+2) = - dphiF1(2:end); % $\pss[F_1]{\svF}{\nT\SCAL\GRAD q}$
    gr_rhs(1:end, pd.K+3) = + dphiF2(2:end); % $\pss[F_2]{\svF}{\nT\SCAL\GRAD q}$
    
    R = gr_lhs\gr_rhs;   % Solve problem (up to a constant)
    A = gr_rhs'*R;       % Compute $\psv[T]{\GRAD\opRec(\cdot)}{\GRAD\opRec(\cdot)}$
end
\end{minted}
    \caption{Possible implementation of the reconstruction operator in 1D.}
    \label{code:reconstruction}
\end{listing}

\begin{figure}[htbp]
    \centering
    \includegraphics[width=0.49\textwidth]{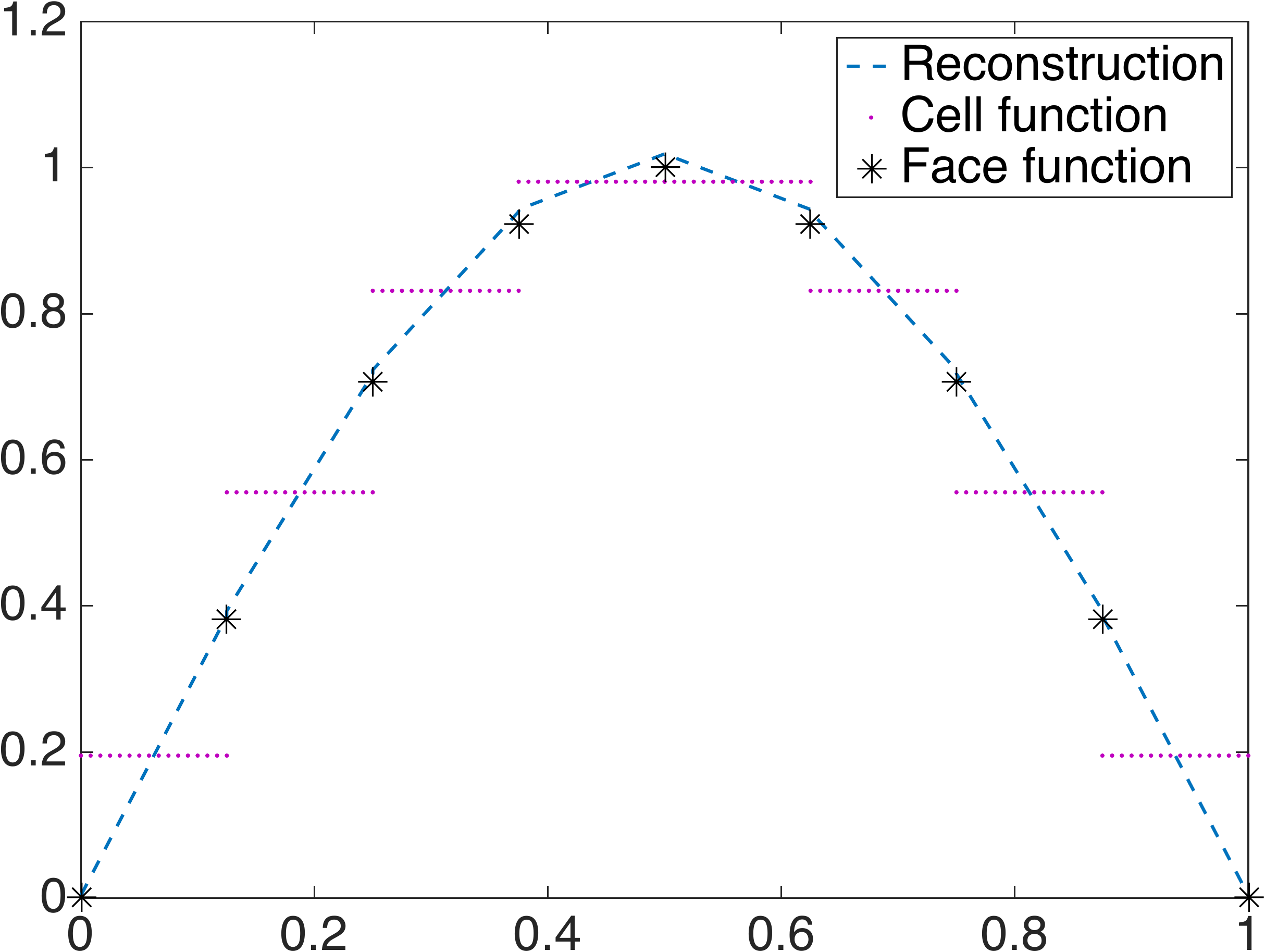}%
    \hfill%
    \includegraphics[width=0.49\textwidth]{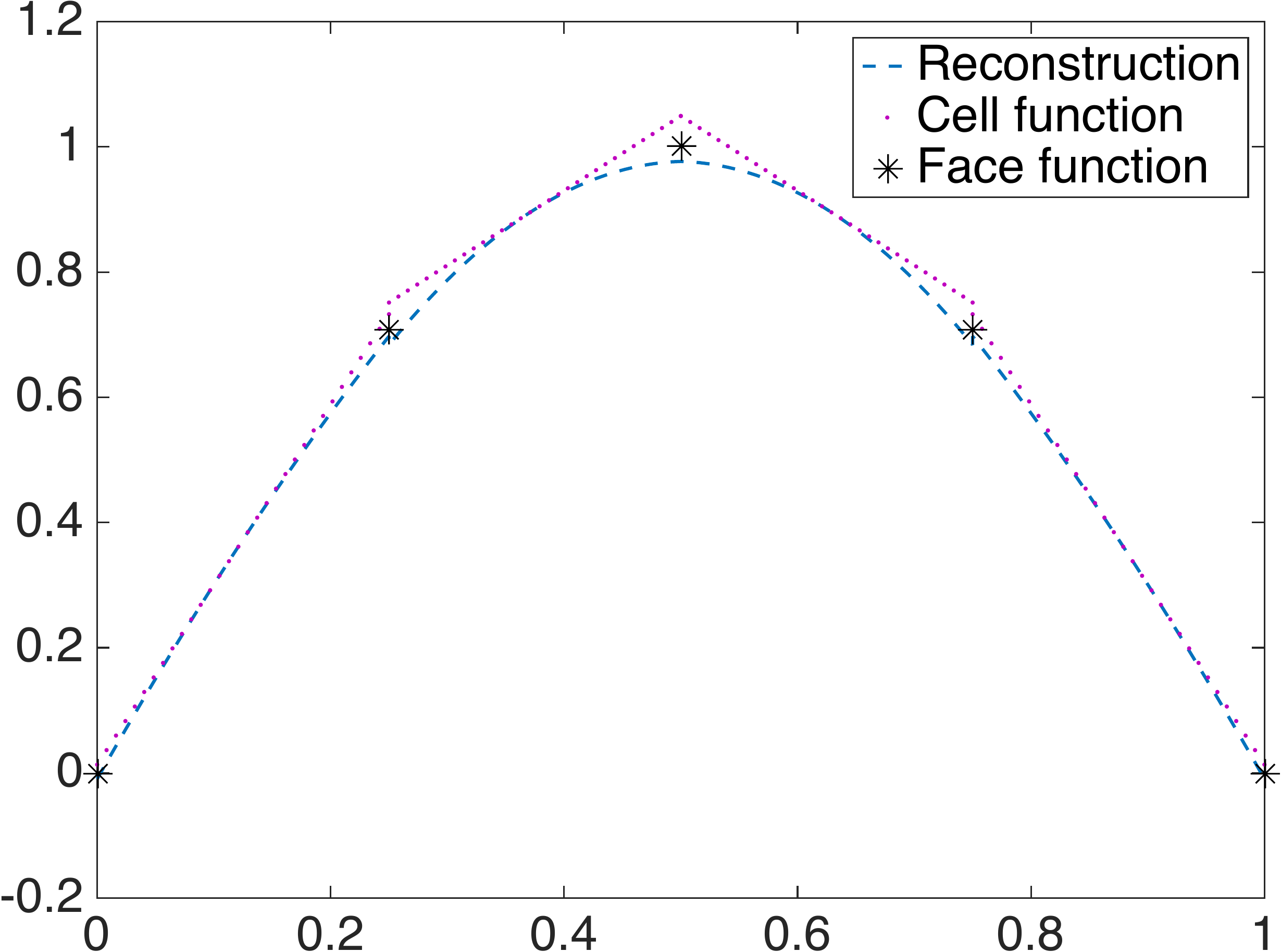}%

    \caption{Illustration of the action of the reconstruction operator $\matr{R}$ in 1D on the DoFs resulting from the computation of $\shIkh(\sin(\pi x))$. On the left panel, the operator acts on cell polynomials of degree 0 (dotted line) and face values (stars) to reconstruct a piecewise polynomial of degree 1 (dashed line). On the right panel, starting from cell polynomials of degree 1 (dotted line) and face values (stars), a piecewise polynomial of degree 2 (dashed line) is reconstructed. Recall that in the 1D case, there is only one DoF per face.}
    \label{fig:reconstruction-demo}
\end{figure}


\bRem[Verifying the implementation]
Given a sequence of successively refined meshes $\meshfam=(\calT_i)_{i\in\polN}$ and a target function $v\in \Hun$, the vector $\bm{\mathsf{I}}_T^{k',k}(v) = [\hhodofs{v}_T | \hhodofs{v}_{F_1} | \ldots | \hhodofs{v}_{F_n}]\tr$ is computed for every mesh cell $T\in\calT_i$ and all $i\in\polN$. Subsequently, we compute the matrix-vector product $\hhodofs{v}^{*} = \matr{R}\bm{\mathsf{I}}_T^{k',k}(v)$, where $\hhodofs{v}^{*}$ are the components of the polynomial $\opRec(\shIkpkT(v)) \in \poly{k+1}{*d}(T)$. The average in $T$ of the reconstructed function is fixed by forming the vector $\hhodofs{v} = [v_{\phi} | \hhodofs{v}^{*}]$, which collects the components of the reconstruction of $v$ in $\poly{k+1}{d}(T)$ and where $v_{\phi}$ is a constant ensuring the condition \eqref{eq:def_Rec_HHO_mean}. We finally compute the $L^2$-error between the reconstruction of $v$ and $v$ itself, which should decay with rate $\bigO(h_i^{k+2})$.
\eRem

\subsection{The stabilization operator}
The computation of the HHO stabilization is a relatively involved task, and for this reason, it is discussed in two steps. In the first step, the Lehrenfeld--Sch\"oberl (LS) stabilization is considered (recall that this stabilization is sufficient when working with mixed-order HHO methods, i.e., $k'=k+1$). In the second step, the equal-order HHO stabilization is discussed as an extension of the LS stabilization. 

\subsubsection{Step 1:  Lehrenfeld--Sch\"oberl stabilization}
The idea behind the LS stabilization is to penalize just the difference between the polynomial attached to a face $\Face_i$ and the trace on $\Face_i$ of the polynomial attached to $\Cell$. This is accomplished by using the operator $Z_F : \Vmix \rightarrow \poly{k}{d-1}(F)$ defined as
\begin{align}
    Z_F(\shvT) := \proj{k}_F(v_T) - v_F, \label{eqn:Zstab_oper}
\end{align}
which is used to build the bilinear form $z_T: \Vmix \times \Vmix \rightarrow \mathbb{R}$ such that
\begin{align}
    z_T(\shvT, \shwT) := \sum_{F \in \FT} h_T^{-1} (Z_F(\shvT), Z_F(\shwT))_{L^2(F)}. \label{eqn:Zstab}
\end{align}

The operator $Z_F$ actually subtracts two polynomials, and at the algebraic level, this is done by subtracting their DoFs. This is accomplished by a matrix $\matr{Z}_i$
of size $\polysz{k}{d-1} \times \hhosz{k}{d}$ constructed as follows.
A matrix $\matr{I}_i$ of size $\polysz{k}{d-1} \times \hhosz{k}{d}$ is first formed by placing a diagonal of ones in correspondence to the $F_i$-block (see Figure~\ref{fig:opmatrix_structure}) of $\matr{I}_i$ (note that the $F_i$-blocks of the stabilization operator are all square of size $\polysz{k}{d-1} \times \polysz{k}{d-1}$). The matrix $\matr{I}_i$ can be thought as a selection matrix such that when left-multiplying a vector of the form (\ref{eqn:dofsarray}), it yields the subvector containing only the DoFs $\hhodof{v}_{F_i,1}, \ldots, \hhodof{v}_{F_i, \polysz{k}{d-1}}$.
The second step consists in computing the DoFs of the polynomial which represents the restriction on $\Face_i$ of the polynomial attached to $\Cell$, and this is done by means of a projection. We first construct the \emph{trace matrix} of size $\polysz{k}{d-1} \times \hhosz{k}{d}$ such that
\begin{align}
    \matr{T}_i \eqq \sum_{\Qindex=1}^{|\Qset_{F_i}|} \Qwq \Fbasv_{F_i}(\Qxq) \Thhobasv(\Qxq)\tr,
\end{align}
whose role is explained as follows: for each quadrature point $\Qxq \in Q_{F_i}$, if a vector of the form (\ref{eqn:dofsarray}) is left-multiplied by $\Thhobasv(\Qxq)\tr$, the operation yields the value of the cell-based polynomial at the point $\Qxq$ (which lies on $F_i$). The subsequent multiplication by $\bm{\phi}_{F_i}(\Qxq)$ then tests the cell-based polynomial with the basis functions of $\poly{k}{d-1}(F_i)$, effectively forming a right-hand side suitable for a projection problem like (\ref{eqn:proj-system}). The left-hand side of the projection problem is the mass matrix of the face $F_i$ of size $\polysz{k}{d-1} \times \polysz{k}{d-1}$ such that
\begin{align}
    \matr{M}_i \eqq \sum_{\Qindex=1}^{|\Qset_{F_i}|} \Qwq \bm{\phi}_{F_i}(\Qxq) \bm{\phi}_{F_i}(\Qxq)\tr,
\end{align}
with which we form the additional matrix $\matr{M}_i^{-1}\matr{T}_i$. This last matrix, applied to a vector of HHO DoFs, yields the sought restriction. Using the matrices just computed, we finally obtain the discrete counterpart of (\ref{eqn:Zstab_oper}) as
\begin{align}
    \matr{Z}_i := \matr{M}_i^{-1}\matr{T}_i - \matr{I}_i, \label{eqn:stab_basic}
\end{align}
which, if applied to a vector $\hhodofs{v} \in \hhodofs{V}^{k',k}_T$, yields the difference between the polynomial on $\Face_i$ and the polynomial on $\Cell$ projected on the face $\Face_i$. This allows us to compute the algebraic counterpart of (\ref{eqn:Zstab}) as
\begin{align}
\matr{Z} := \sum_{i=1}^{n} h_{T}^{-1} \matr{Z}_i\tr \matr{M}_i \matr{Z}_i. \label{eqn:basic_stab_lc}
\end{align}

\subsubsection{Step 2: Equal-order stabilization}
To obtain the equal-order HHO stabilization where $k'=k$, we need to enhance \eqref{eqn:Zstab_oper} by introducing a penalty on the high-order contribution due to the reconstruction. We consider \eqref{eq:def_SKk_HHO}, which we rewrite here by specifying the face $F\in\FacesT$, leading to the operator $S_F : \shVkT \rightarrow \poly{k}{d-1}(F)$ defined as
\begin{align}
S_F(\shvT) := Z_F(\shvT) + \proj{k}_F\big(\opRec(\shvT)-\proj{k}_T \opRec(\shvT)\big).
\end{align}
This operator is used to build the bilinear form $s_T: \shVkT \times \shVkT \rightarrow \mathbb{R}$ such that
\begin{align}
    s_T(\shvT, \shwT) := \sum_{F \in \FT} h_T^{-1} (S_F(\shvT), S_F(\shwT))_{L^2(F)}. \label{eqn:full_stab}
\end{align}
We start by translating in matrix form the term $\proj{k}_F\opRec(\shvT)$. First, we compute
\begin{align}
    \matr{T}_i' \eqq \sum_{q=1}^{|Q_{F_i}|} \Qw_q \Fbasv_{F_i}(\bm{x}_q)\Rbasv(\bm{x}_q)\tr,
\end{align}
which has size $\polysz{k}{d-1} \times \polysz{k+1}{*d}$, to subsequently construct the term
\begin{align}
    \matr{M}_i^{-1}\matr{T}_i'\matr{R}, \label{eqn:stab_R}
\end{align}
where $\matr{R}$ is the reconstruction defined in (\ref{eqn:discoperR}) (notice that it is not necessary to take into account the mean-value correction in this construction). The matrix we just built can be understood by reading it backwards as follows: by applying $\matr{R}$ to an object in $\hhodofs{V}^k_T$, we get its reconstruction in $\poly{k+1}{*d}(T)$. Subsequently, when the trace matrix $\matr{T}_i'$ is applied to the DoFs of the reconstructed polynomial, its columns evaluate the DoFs of the reconstructed function on $\Face_i$, whereas the rows test it with the basis functions of $\poly{k}{d-1}(F_i)$. The final multiplication by $\matr{M}_i^{-1}$ yields the DoFs of the reconstructed polynomial restricted to $\Face_i$.

We proceed similarly to translate the term $\proj{k}_F \proj{k}_T \opRec(\shvT)$ in matrix form. This requires the introduction of the cell mass matrix $\matr{M}$ and the matrix
\begin{align}
    \matr{Q} \eqq \sum_{q=1}^{|Q_T|} \Qw_q \Tbasv_T(\bm{x}_q)\Rbasv(\bm{x}_q)\tr,
\end{align}
which has size $\polysz{k}{d} \times \polysz{k+1}{*d}$. We construct the expression
\begin{align}
    \matr{M}_i^{-1}\matr{\widetilde{T}}_i\matr{M}^{-1}\matr{Q}\matr{R}, \label{eqn:stab_projR}
\end{align}
where $\matr{\widetilde{T}}_i$ is the matrix $\matr{T}_i$ restricted to its first $\polysz{k}{d}$ columns. Again, this last expression is better understood by reading it backwards, and keeping in mind the role of the rows and the columns of each matrix: $\matr{Q}$ evaluates the DoFs of the reconstruction and tests it with the basis functions of $\poly{k}{d}(T)$, whereas the multiplication by $\matr{M}^{-1}$ yields the DoFs corresponding to the result of the projection $\proj{k}_T$. The DoFs of the projection on the face are finally obtained by applying $\matr{M}_i^{-1}\matr{\widetilde{T}}_i$. 

Putting everything together, the matrix form of the equal-order HHO stabilization is computed by combining (\ref{eqn:stab_basic}), (\ref{eqn:stab_R}), and (\ref{eqn:stab_projR}) as follows:
\begin{align}
    \matr{S}_i := \matr{Z}_i + \matr{M}_i^{-1}\matr{T}_i'\matr{R} - \matr{M}_i^{-1}\matr{\widetilde{T}}_i\matr{M}^{-1}\matr{Q}\matr{R}. \label{eqn:stab_full}
\end{align}
It is now possible to build the discrete counterpart of (\ref{eqn:full_stab}) as
\begin{align}
\matr{S} := \sum_{i=1}^{n} h_{T}^{-1} \matr{S}_i\tr \matr{M}_i \matr{S}_i. \label{eqn:full_stab_lc}
\end{align}
We propose in Listing~\ref{code:stabilization} a practical implementation of the equal-order stabilization operator in 1D. On lines 7-8, the matrices \texttt{M} and \texttt{Q} are cut from an order $(k+1)$ mass matrix (\texttt{mass\_mat}); an optimized construction would use an order $k$ basis for the rows and a quadrature of order $2k+1$.

\begin{listing}
\begin{minted}[mathescape, linenos, framesep=2mm]{matlab}
function S = hho_stabilization(pd, elem, R)
    x_bar = cell_center(pd, elem); 
    [xF1, xF2] = face_centers(pd, elem);
    mass_mat = make_mass_matrix(pd, elem, pd.K+1); 
    
    % Compute the term tmp1 = $u_T - \proj{k}_T\opRec(\shuT)$
    M = mass_mat(1:pd.K+1,1:pd.K+1); 
    Q = mass_mat(1:pd.K+1,2:pd.K+2); 
    tmp1 = - M\(Q*R);
    tmp1(1:pd.K+1, 1:pd.K+1) = tmp1(1:pd.K+1, 1:pd.K+1) + eye(pd.K+1);
    
    [phiF1, ~] = basis(xF1, x_bar, pd.h, pd.K+1); 
    Mi = 1;
    Ti = phiF1(2:end)';
    Ti_tilde = phiF1(1:pd.K+1)'; 
    tmp2 = Mi \ (Ti*R);             % tmp2 = $\proj{k}_F\opRec(\shuT)$
    tmp2(pd.K+2) = tmp2(pd.K+2)-1;  % tmp2 = $\proj{k}_F\opRec(\shuT) - u_F$
    tmp3 = Mi \ (Ti_tilde * tmp1);  % tmp3 = $\proj{k}_F(u_T - \proj{k}_T\opRec(\shuT))$
    Si = tmp2 + tmp3;               % Si = $\proj{k}_F\opRec(\shuT) - u_F + \proj{k}_F(u_T - \proj{k}_T\opRec(\shuT))$
    S = Si' * Mi * Si / pd.h;       % Accumulate on S
    
    [phiF2, ~] = basis(xF2, x_bar, pd.h, pd.K+1);
    Mi = 1;
    Ti = phiF2(2:end)';
    Ti_tilde = phiF2(1:pd.K+1)'; 
    tmp2 = Mi \ (Ti*R);             % tmp2 = $\proj{k}_F(\opRec(\shuT))$
    tmp2(pd.K+3) = tmp2(pd.K+3)-1;  % tmp2 = $\proj{k}_F\opRec(\shuT) - u_F$ 
    tmp3 = Mi \ (Ti_tilde * tmp1);  % tmp3 = $\proj{k}_F(u_T - \proj{k}_T\opRec(\shuT))$
    Si = tmp2 + tmp3;               % Si = $\proj{k}_F\opRec(\shuT) - u_F + \proj{k}_F(u_T - \proj{k}_T\opRec(\shuT))$
    S = S + Si' * Mi * Si / pd.h;   % Accumulate on S
end
\end{minted}
    \caption{Possible implementation of the equal-order stabilization operator in 1D.}
    \label{code:stabilization}
\end{listing}

\bRem[Verifying the implementation]
The correctness of the implementation of the stabilization operator is verified as before by taking a sequence of successively refined meshes $\meshfam=(\calT_i)_{i\in\polN}$ and a target function $v \in H^1(\Dom)$. For every mesh cell $\T\in\calT_i$ and all $i\in\polN$, the local vector of DoFs $\bm{\mathsf{I}}^{k',k}_T(v) = [\hhodofs{v}_T | \hhodofs{v}_{F_1} | \ldots | \hhodofs{v}_{F_n}]\tr$ is computed. This vector is then used to compute the quantity
$\epsilon_i := (\sum_{T \in \calT_i} \hhodofs{v}\tr \, \matr{S} \, \hhodofs{v})^{\frac12}$
which should converge to zero with decay rate $\bigO(h_i^{k+1})$. The same result
should be obtained for the LS stabilization.
\eRem

\section{Assembly and boundary conditions}

Using either the mixed-order or the equal-order HHO method, the local contributions in every mesh cell $T\in\Th$ are computed as $\matr{L}_T \eqq \matr{A}_T + \matr{Z}_T$ (using \eqref{eqn:lapl_lc} and \eqref{eqn:basic_stab_lc}) or $\matr{L}_T \eqq \matr{A}_T + \matr{S}_T$ (using \eqref{eqn:lapl_lc} and \eqref{eqn:full_stab_lc}), respectively. Here, we added a subscript referring to the mesh cell $T$ for more clarity. The resulting local matrix $\matr{L}_T$ is statically condensed (see Sect.~\ref{sec:static_HHO}), leading to the condensed matrix $\matr{L}_T^{\textsc{c}}$ of size $(nN^k_{d-1})\times (nN^k_{d-1})$ and the condensed right-hand side $\matr{b}_T^{\textsc{c}}$ of size $(nN^k_{d-1})\times 1$ (recall that $n\eqq \#\FacesT$ is the number of faces of $T$).

The assembly of the global problem requires a local-to-global correspondence array denoted by $\mathcal{G}_T : \{1,\ldots,n\} \rightarrow \{1,\ldots,\#\Faces\}$ for all $T\in\Th$, between the local enumeration of the faces of $T$ and their global enumeration as mesh faces. This array is usually provided by the mesh generator. 
%
In the first stage of the assembly process, one does not bother about boundary conditions (this amounts to assemble a problem with pure Neumann boundary conditions). 
The global matrix $\matr{L}^{\textsc{g}}$ is composed of $\#\Faces \times \#\Faces$ blocks of size $\polysz{k}{d-1} \times \polysz{k}{d-1}$ and the global right-hand side $\matr{b}^{\textsc{g}}$ is composed of $\#\Faces$ blocks of size $\polysz{k}{d-1} \times 1$. Then the local contributions are assembled as follows: For all $T\in\Th$,
\begin{align}
\matr{L}^{\textsc{g}}_{\mathcal{G}_T(i), \mathcal{G}_T(j)} \hookleftarrow \matr{L}_{T;i,j}^{\textsc{c}} \quad \text{and} \quad \matr{b}^{\textsc{g}}_{\mathcal{G}_T(i)} \hookleftarrow \matr{b}_{T;i}^{\textsc{c}}, \quad \forall i,j \in \{1,\ldots, n\},
\end{align}
where we denote by $a \hookleftarrow b$ the operation of accumulating the value $b$ on $a$, i.e. the statement \texttt{a = a + b} of the commonly used imperative programming languages. In other words, the local block $(i,j)$ is summed to the global block in position $(\mathcal{G}_T(i), \mathcal{G}_T(j))$.

It remains to apply the boundary conditions. As discussed in \cite{DiPEL:16}, HHO methods can handle all the classical boundary conditions for the Poisson model problem (see also \cite{BusMu:20}). To fix the ideas, let us assume that the boundary is partitioned 
as $\front = \overline{\Bn} \cup \overline{\Bd}$ leading to the following model problem:
\begin{equation} \label{sec:impl_BC}
-\Delta u=f\;\text{in $\Dom$},\qquad
u=\suD \;\text{on $\Bd$}, \qquad \bm{n}\SCAL\GRAD u=\sGn\;\text{on $\Bn$}.
\end{equation}
We assume that every mesh boundary face belongs either to $\Bd$ or to
$\Bn$; the corresponding subsets of $\calFb$ are denoted by $\calFbD$ and $\calFbN$.
Let us consider an idealized 1D situation with
a simple mesh containing only four faces (vertices), i.e., 
$\calF\eqq \{F_1,F_2,F_3,F_4\}$ with $\calFi\eqq\{F_2,F_3\}$ and $\calFb\eqq\{F_1,F_4\}$.
Then, assuming that only Neumann boundary conditions are enforced (i.e., $\front=\Bn$, $\Bd=\emptyset$ in \eqref{sec:impl_BC}), the global problem takes the form 
\begin{align}
\begin{bmatrix}
\matr{L}_{11} & \matr{L}_{12} &               &               \\
\matr{L}_{21} & \matr{L}_{22} & \matr{L}_{23} &               \\
              & \matr{L}_{32} & \matr{L}_{33} & \matr{L}_{34} \\
              &               & \matr{L}_{43} & \matr{L}_{44} \\
\end{bmatrix}
\begin{bmatrix}
\hhodofs{u}_1 \\
\hhodofs{u}_2 \\
\hhodofs{u}_3 \\
\hhodofs{u}_4 \\
\end{bmatrix}
=
\begin{bmatrix}
\stdvec{b}_1 \\
\stdvec{b}_2 \\
\stdvec{b}_3 \\
\stdvec{b}_4 \\
\end{bmatrix}. \label{eqn:sys_example_dirichlet}
\end{align}
(Notice that in this 1D case, all the entries are actually scalars.)
Assume now that the Neumann boundary condition is applied only on $F_1$ and that the Dirichlet condition is applied on $F_4$. Then we have $\hhodofs{u}_4=\matr{M}_{F_4}^{-1}\stdvec{d}_4$ with $\stdvec{d}_4 \eqq \int_{F_4} \suD(\bx)\Fbasv_{F_4}(\bx)\ds$. Eliminating $\hhodofs{u}_4$ from the first three rows of~\eqref{eqn:sys_example_dirichlet} gives the reduced system
\begin{align}
\begin{bmatrix}
\matr{L}_{11} & \matr{L}_{12} &               \\
\matr{L}_{21} & \matr{L}_{22} & \matr{L}_{23} \\
              & \matr{L}_{32} & \matr{L}_{33} \\
\end{bmatrix}
\begin{bmatrix}
\hhodofs{u}_1 \\
\hhodofs{u}_2 \\
\hhodofs{u}_3 \\
\end{bmatrix}
=
\begin{bmatrix}
\stdvec{b}_1\\
\stdvec{b}_2\\
\stdvec{b}'_3\\
\end{bmatrix},
\end{align}
where $\stdvec{b}'_3 := \stdvec{b}_3 - \matr{L}_{34}\matr{M}_{F_4}^{-1}\stdvec{d}_4$. This process can be conveniently done on the fly during the assembly, but a new mapping $\calG^{\circ}$ has to be used. Such a mapping is computed like $\calG$, but removing the Dirichlet faces. Once the solution of the reduced system is found, the full solution is recovered by plugging $\hhodofs{d}_4$ after $\hhodofs{u}_3$ in the solution vector. An alternative approach is to introduce a Lagrange multiplier to enforce the Dirichlet condition:
\begin{align}
\begin{bmatrix}
\matr{L}_{11} & \matr{L}_{12} &               &                &                \\
\matr{L}_{21} & \matr{L}_{22} & \matr{L}_{23} &                &                \\
              & \matr{L}_{32} & \matr{L}_{33} & \matr{L}_{34}  &                \\
              &               & \matr{L}_{43} & \matr{L}_{44}  & \matr{M}_{F_4} \\
              &               &               & \matr{M}_{F_4} &                \\
\end{bmatrix}
\begin{bmatrix}
\hhodofs{u}_1 \\
\hhodofs{u}_2 \\
\hhodofs{u}_3 \\
\hhodofs{u}_4 \\
\hhodofs{\lambda}_5 \\
\end{bmatrix}
=
\begin{bmatrix}
\stdvec{b}_1 \\
\stdvec{b}_2 \\
\stdvec{b}_3 \\
\stdvec{b}_4 \\
\stdvec{d}_4 \\
\end{bmatrix}. \label{eqn:sys-lagrange}
\end{align}
This second technique leads to a slightly larger system having a saddle-point structure, but it could be easier to implement in a first HHO code.

%

\bRem[Neumann boundary conditions]
The Neumann boundary condition in the model problem~\eqref{sec:impl_BC} leads to a modification of the linear form on the right-hand side of the discrete problem, which reads $\ell(\shwh)\eqq (f,w_\calT)_{L^2(\Dom)}+(\sGn,w_\calF)_{L^2(\Bn)}$ (see Sect.~\ref{sec:disc_pb_elas} for the linear elasticity problem). At the implementation level, the Neumann condition reduces to a contribution on the right-hand side of the linear system positioned according to the Neumann face unknowns. Such a contribution is computed as
$\stdvec{g}_i = \sum_{\Qindex=1}^{|Q_{F_i}|} \Qwq \sGn(\Qxq) \stdvec{\phi}_{F_i}(\Qxq)$,
where $\stdvec{\phi}_{F_i}$ is the vector of basis functions of the globally-numbered $i$-th face. Those contributions are then added to the $i$-th block of the right-hand side. 
\eRem

\section{Remarks on the computational cost of HHO methods}
The computational costs in HHO methods are of two kinds: local costs associated with the assembly and global costs associated with the solution of the global linear system. We focus as before on the Poisson model problem.

The local costs include the computation of the operators and the static condensation, and they differ in the mixed-order and equal-order methods. In the mixed-order method, the computation of the reconstruction and the static condensation are slightly more expensive compared to the equal-order case, essentially because of the increased number of cell-based DoFs. The costs of the stabilization, however, differ substantially between the two variants of the method. This fact can be deduced by comparing the structure of (\ref{eqn:stab_basic}) and (\ref{eqn:stab_full}). The mixed-order stabilization requires $n$ inversions of the face mass matrices $\matr{M}_i$, for a cost of $n \cdot \calO( (\polysz{k}{d-1})^3 ) \approx \calO(k^{3d-3})$, together with the construction of the trace matrices $\matr{T}_i$, for a cost of $n \cdot \calO( (k+1)^{d-1} \cdot \polysz{k+1}{d} \cdot \polysz{k}{d-1}) \approx \calO(k^{3d-2})$. Instead, the equal-order stabilization requires the inversion of the cell mass matrix and other operations which are at least cubic in the size of the cell basis.
\begin{figure}[htbp]
  \centering
  \includegraphics[width=.49\linewidth]{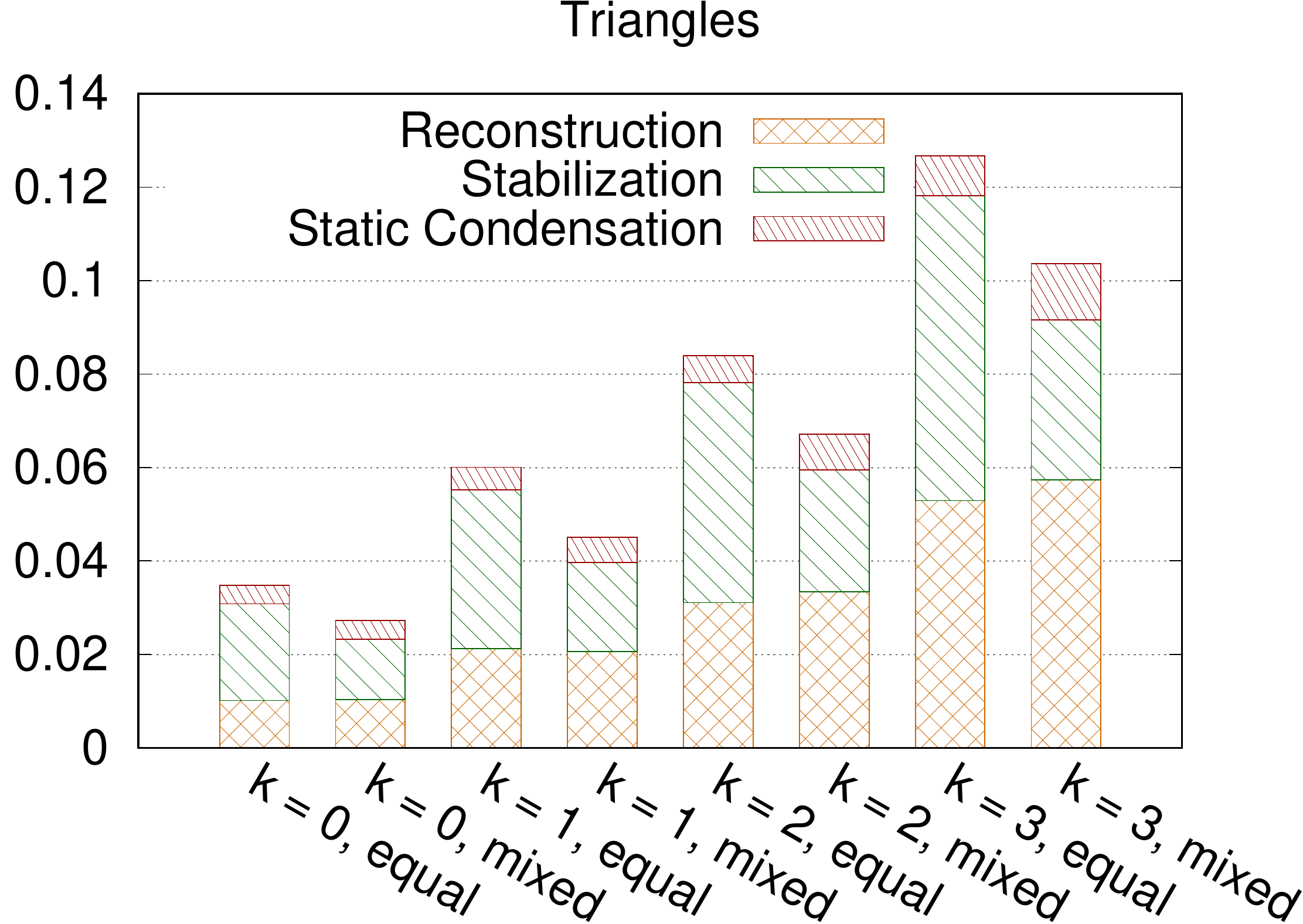}%
  \hfill%
  \includegraphics[width=.49\linewidth]{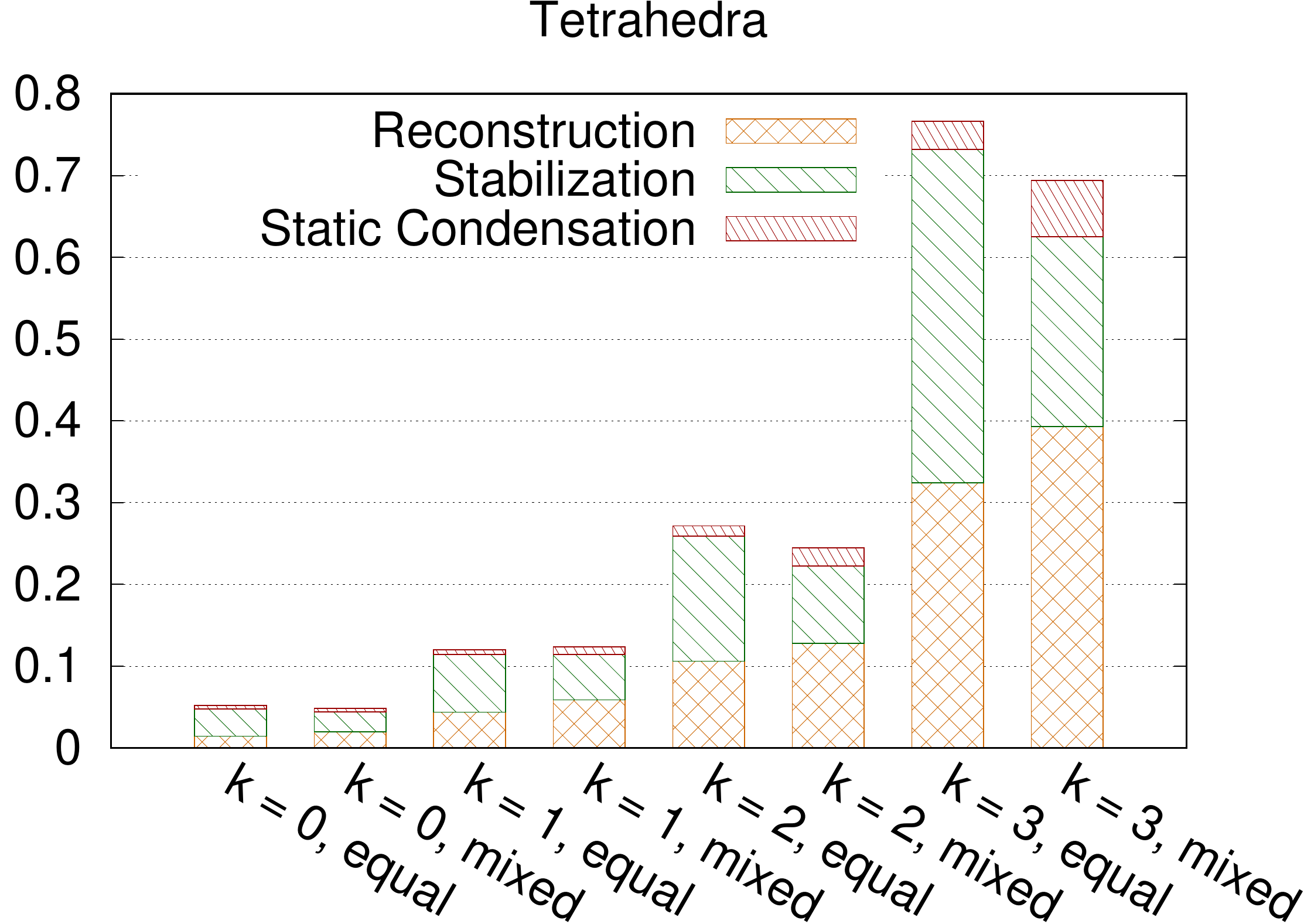}%
  \vspace{2mm}
  \includegraphics[width=.49\linewidth]{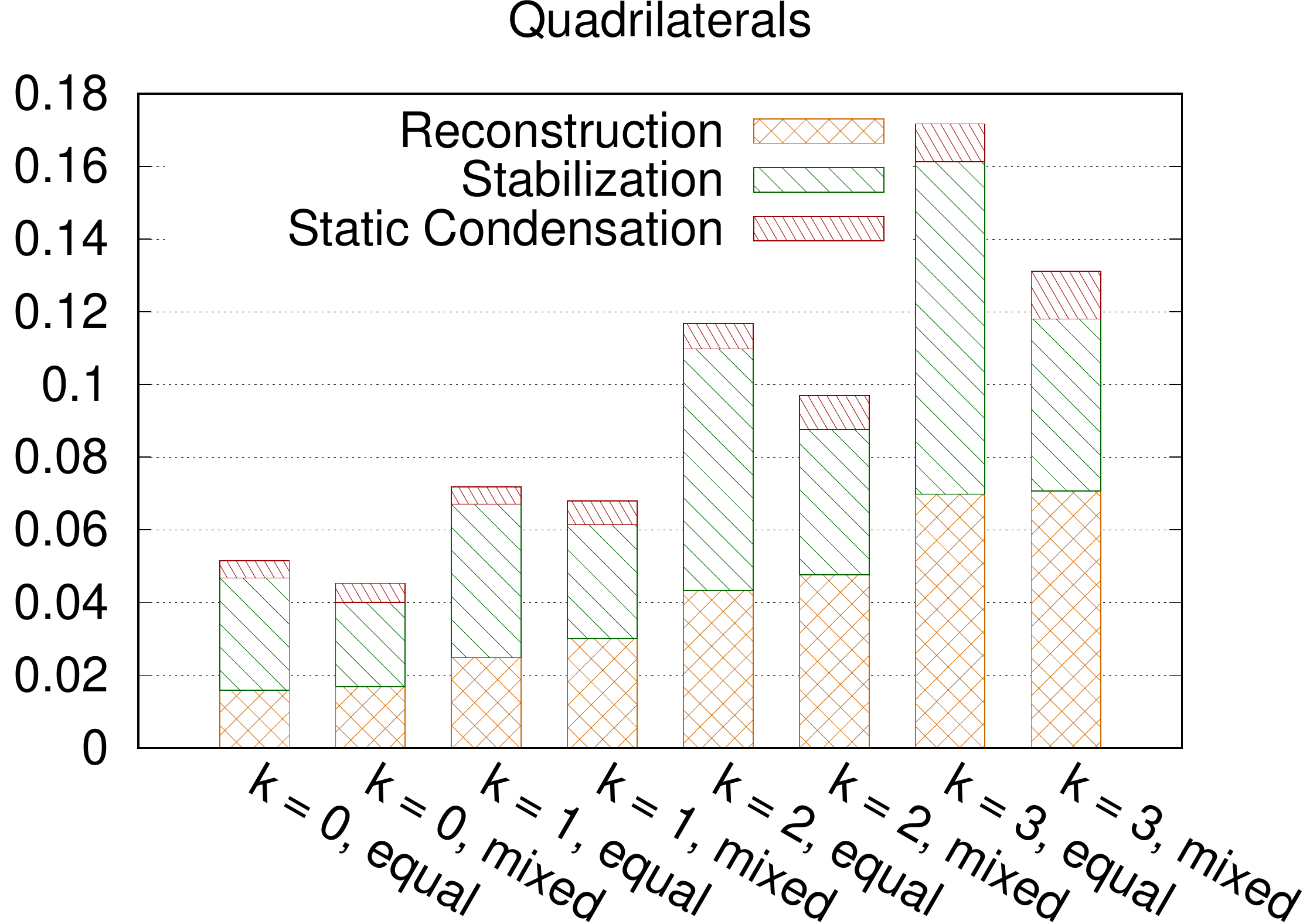}%
  \hfill%
  \includegraphics[width=.49\linewidth]{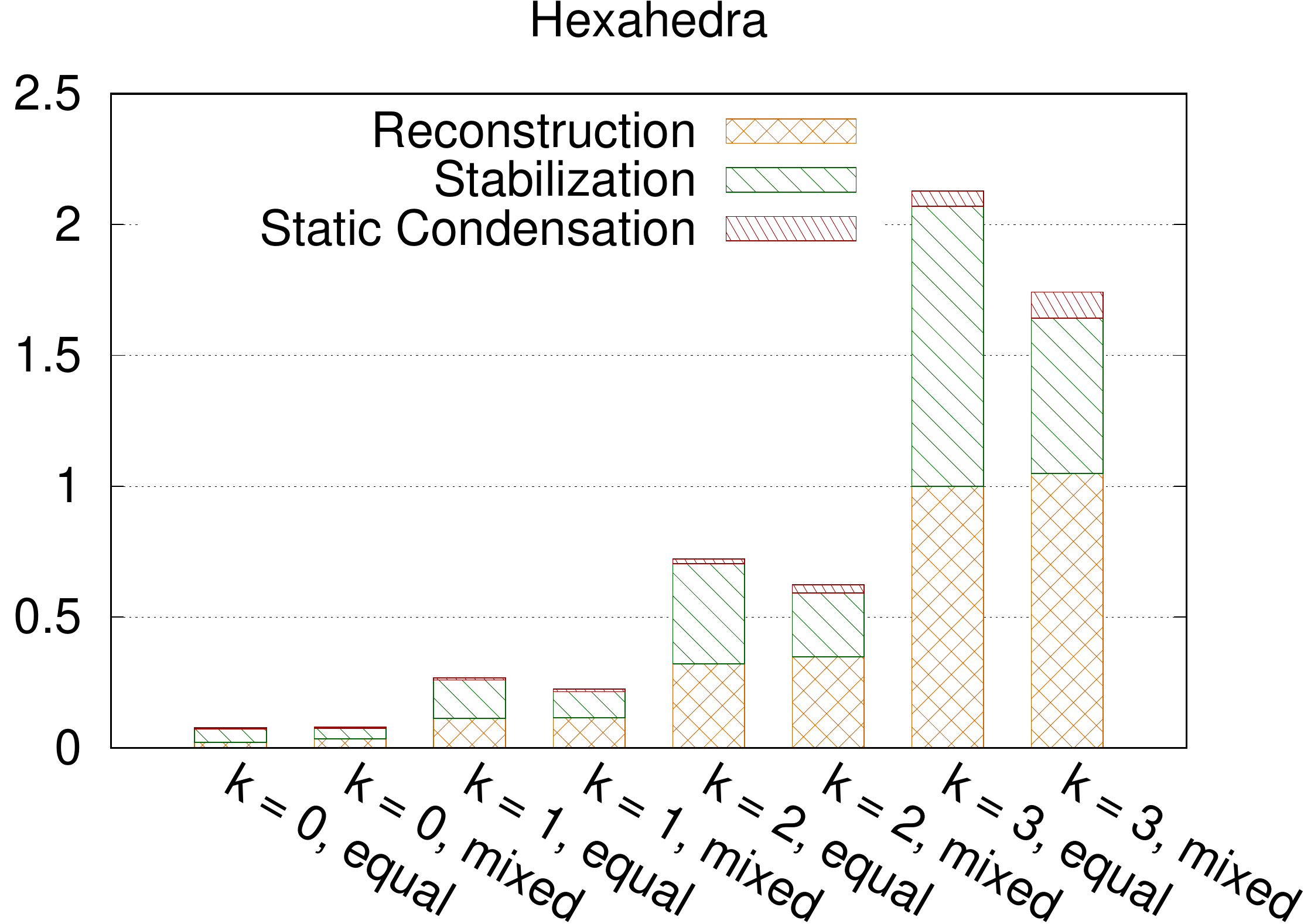}%
  \caption{Comparison of average computational times (in milliseconds) for the construction of the HHO operators on a single mesh cell, including static condensation. In the mixed-order HHO method, even if the cost of reconstruction and static condensation is slightly increased, the lower cost of stabilization results in a reduction of total computational time.}
  \label{fig:hho_oper_cost}
\end{figure}
To illustrate this fact, we performed some computational experiments on common element types, namely triangles and quadrangles in 2D, and tetrahedra and hexahedra in 3D (see Figure~\ref{fig:hho_oper_cost}). In all cases, we observe that when using the mixed-order HHO method, even if one pays a bit more in reconstruction and static condensation, one pays a lot less in stabilization. This turns in an overall reduction of the cost of the computation of the local contributions.

Concerning the global costs, we illustrate the differences between the equal-order HHO method and the well-established symmetric interior-penalty discontinuous Galerkin (SIP-DG) method (see \cite{ArBCM:01} or \cite[Sect.~4.2]{DiPEr:12}) for the Poisson model problem posed in the unit cube $(0,1)^3$. For HHO, we use polynomials of one degree less than in SIP-DG, so that both methods deliver the same error decay rates. We ran the experiments on 3D meshes of tetrahedra (3{,}072 elements) and hexahedra (4{,}096 elements). The global linear systems were solved using the PARDISO linear solver from the Intel MKL library. Memory usage was estimated via the \texttt{getrusage()} system call. The results reported in Tables~\ref{tbl:cost_tetra}-\ref{tbl:cost_hexa} indicate that the HHO discretization is more favorable in terms of linear solver operations and memory usage.

\begin{table}
    \centering
    \begin{tabular}{c||c|c|c|c||c|c|c|c}
      & \multicolumn{4}{c||}{HHO($k$,$k$)} & \multicolumn{4}{c}{SIP-DG($k+1$)} \\
    \hline
    $k$ & $L^2$-error & DoFs  & Mflops    & Memory & $L^2$-error  & DoFs      & Mflops & Memory   \\
    \hline
    0 & 1.73e-2     & 5760  & 38        & 39 MB  & 2.14e-2      & 12288     & 787    & 85 MB    \\
    1 & 1.06e-3     & 17280 & 1006      & 106 MB & 4.61e-4      & 30720     & 11429  & 319 MB   \\
    2 & 9.05e-5     & 34560 & 8723      & 292 MB & 2.14e-5      & 61440     & 92799  & 1108 MB  \\
    3 & 6.45e-6     & 57600 & 40389     & 719 MB & 1.04e-6      & 107520    & 497245 & 3215 MB  \\
    \end{tabular}
    \caption{Comparative cost assessment between HHO and SIP-DG on a tetrahedral mesh composed of 3{,}072 elements.}
    \label{tbl:cost_tetra}
\end{table}
\begin{table}
    \centering
    \begin{tabular}{c||c|c|c|c||c|c|c|c}
      & \multicolumn{4}{c||}{HHO($k,k$)} & \multicolumn{4}{c}{SIP-DG($k+1$)} \\
    \hline
    $k$ & $L^2$-error & DoFs  & Mflops    & Memory & $L^2$-error  & DoFs      & Mflops & Memory   \\
    \hline
    0 & 9.07e-3     & 11520  & 310       & 64 MB   & 6.03e-2      & 16384     & 6677    & 168 MB    \\
    1 & 3.04e-4     & 34560  & 9671      & 293 MB  & 1.72e-4      & 40960     & 104199  & 765 MB   \\
    2 & 1.73e-5     & 69120  & 58977     & 884 MB  & 1.29e-6      & 81920     & 845545  & 2844 MB  \\
    3 & 7.13e-7     & 115200 & 349664    & 2412 MB & 5.24e-8      & 143360    & 4592328 & 8490 MB  \\
    \end{tabular}
    \caption{Comparative cost assessment between HHO and SIP-DG on a hexahedral mesh composed of 4{,}096 (16x16x16) elements.}
    \label{tbl:cost_hexa}
\end{table}

\backmatter

%
%

\ifSp \else
\fancypagestyle{plain}{
    \renewcommand{\headrulewidth}{0pt}%
    \fancyhf{}%
}
\addcontentsline{toc}{chapter}{Bibliography}
\fancyhead[LE]{}
\fancyhead[RE]{Bibliography}
\fancyhead[LO]{Bibliography}
\fancyhead[RO]{}
\fi
\bibliographystyle{plain}
\bibliography{Bibliographie}

\begin{thebibliography}{100}

\bibitem{AbErPi:18}
M.~Abbas, A.~Ern, and N.~Pignet.
\newblock {H}ybrid {H}igh-{O}rder methods for finite deformations of
  hyperelastic materials.
\newblock {\em Comput. Mech.}, 62(4):909--928, 2018.

\bibitem{AbErPi:19}
M.~Abbas, A.~Ern, and N.~Pignet.
\newblock {A} {H}ybrid {H}igh-{O}rder method for incremental associative
  plasticity with small deformations.
\newblock {\em Comput. Methods Appl. Mech. Engrg.}, 346:891--912, 2019.

\bibitem{AbErPi:19a}
M.~Abbas, A.~Ern, and N.~Pignet.
\newblock A {H}ybrid {H}igh-{O}rder method for finite elastoplastic
  deformations within a logarithmic strain framework.
\newblock {\em Internat. J. Numer. Methods Engrg.}, 120(3):303--327, 2019.

\bibitem{AbrSt:72}
M.~Abramowitz and I.~A. Stegun.
\newblock {\em Handbook of Mathematical Functions with Formulas, Graphs, and
  Mathematical Tables}.
\newblock Dover, New York, NY, 1972.

\bibitem{AdaFo:03}
R.~A. Adams and J.~J.~F. Fournier.
\newblock {\em Sobolev spaces}, volume 140 of {\em Pure and Applied Mathematics
  (Amsterdam)}.
\newblock Elsevier/Academic Press, Amsterdam, second edition, 2003.

\bibitem{AgBDP:15}
J.~Aghili, S.~Boyaval, and D.~A. Di~Pietro.
\newblock Hybridization of mixed high-order methods on general meshes and
  application to the {S}tokes equations.
\newblock {\em Comput. Methods Appl. Math.}, 15(2):111--134, 2015.

\bibitem{AndDr:18}
D.~Anderson and J.~Droniou.
\newblock An arbitrary-order scheme on generic meshes for miscible
  displacements in porous media.
\newblock {\em SIAM J. Sci. Comput.}, 40(4):B1020--B1054, 2018.

\bibitem{ArBCM:01}
D.~N. Arnold, F.~Brezzi, B.~Cockburn, and L.~D. Marini.
\newblock Unified analysis of discontinuous {G}alerkin methods for elliptic
  problems.
\newblock {\em SIAM J. Numer. Anal.}, 39(5):1749--1779, 2001/02.

\bibitem{Artioli2017}
E.~Artioli, L.~{Beir\~ao da Veiga}, C.~Lovadina, and E.~Sacco.
\newblock Arbitrary order 2{D} virtual elements for polygonal meshes: part
  {II}, inelastic problem.
\newblock {\em Comput. Mech.}, 60(4):643--657, 2017.

\bibitem{AyLiM:16}
B.~Ayuso~de Dios, K.~Lipnikov, and G.~Manzini.
\newblock The nonconforming virtual element method.
\newblock {\em ESAIM Math. Model. Numer. Anal.}, 50(3):879--904, 2016.

\bibitem{Baker:76}
G.~A. Baker.
\newblock Error estimates for finite element methods for second order
  hyperbolic equations.
\newblock {\em SIAM J. Numer. Anal.}, 13(4):564--576, 1976.

\bibitem{Ball1976}
J.~M. Ball.
\newblock Convexity conditions and existence theorems in nonlinear elasticity.
\newblock {\em Arch. Rational Mech. Anal.}, 63(4):337--403, 1976/77.

\bibitem{Beben:03}
M.~Bebendorf.
\newblock A note on the {P}oincar\'e inequality for convex domains.
\newblock {\em Z. Anal. Anwendungen}, 22(4):751--756, 2003.

\bibitem{BdVLM:15}
L.~{Beir\~ao da Veiga}, C.~Lovadina, and D.~Mora.
\newblock A virtual element method for elastic and inelastic problems on
  polytope meshes.
\newblock {\em Comput. Methods Appl. Mech. Engrg.}, 295:327--346, 2015.

\bibitem{BBCMMR:13}
L.~Beir{\~a}o~da Veiga, F.~Brezzi, A.~Cangiani, G.~Manzini, L.~D. Marini, and
  A.~Russo.
\newblock Basic principles of virtual element methods.
\newblock {\em M3AS Math. Models Methods Appl. Sci.}, 199(23):199--214, 2013.

\bibitem{BoBDP:16}
D.~Boffi, M.~Botti, and D.~A. Di~Pietro.
\newblock A nonconforming high-order method for the {B}iot problem on general
  meshes.
\newblock {\em SIAM J. Sci. Comput.}, 38(3):A1508--A1537, 2016.

\bibitem{BoBrF:13}
D.~Boffi, F.~Brezzi, and M.~Fortin.
\newblock {\em Mixed finite element methods and applications}, volume~44 of
  {\em Springer Series in Computational Mathematics}.
\newblock Springer, Heidelberg, 2013.

\bibitem{Boillot:14}
L.~Boillot.
\newblock {\em {Contributions to the mathematical modeling and to the parallel
  algorithmic for the optimization of an elastic wave propagator in anisotropic
  media}}.
\newblock PhD thesis, {Universit{\'e} de Pau et des Pays de l'Adour, France},
  2014.

\bibitem{BDPGK:18}
F.~Bonaldi, D.~A. Di~Pietro, G.~Geymonat, and F.~Krasucki.
\newblock A {H}ybrid {H}igh-{O}rder method for {K}irchhoff-{L}ove plate bending
  problems.
\newblock {\em ESAIM Math. Model. Numer. Anal.}, 52(2):393--421, 2018.

\bibitem{Bonet1997}
J.~Bonet and R.~D Wood.
\newblock {\em Nonlinear continuum mechanics for finite element analysis}.
\newblock Cambridge university press, Cambridge, 1997.

\bibitem{BoDiP:18}
L.~Botti and D.~A. Di~Pietro.
\newblock Assessment of {H}ybrid {H}igh-{O}rder methods on curved meshes and
  comparison with discontinuous {G}alerkin methods.
\newblock {\em J. Comput. Phys.}, 370:58--84, 2018.

\bibitem{BDiPD:18}
L.~Botti, D.~A. Di~Pietro, and J.~Droniou.
\newblock A {H}ybrid {H}igh-{O}rder discretisation of the {B}rinkman problem
  robust in the {D}arcy and {S}tokes limits.
\newblock {\em Comput. Methods Appl. Mech. Engrg.}, 341:278--310, 2018.

\bibitem{BoDiPD:19}
L.~Botti, D.~A. Di~Pietro, and J.~Droniou.
\newblock A {H}ybrid {H}igh-{O}rder method for the incompressible
  {N}avier-{S}tokes equations based on {T}emam's device.
\newblock {\em J. Comput. Phys.}, 376:786--816, 2019.

\bibitem{BoCDH:21}
M.~Botti, D.~Casta{\~n}{\'o}n~Quiroz, D.~A. Di~Pietro, and A.~Harnist.
\newblock A hybrid high-order method for creeping flows of non-{N}ewtonian
  fluids.
\newblock {\em ESAIM Math. Model. Numer. Anal.}, 55(5):2045--2073, 2021.

\bibitem{BDiPG:19}
M.~Botti, D.~A. Di~Pietro, and A.~Guglielmana.
\newblock A low-order nonconforming method for linear elasticity on general
  meshes.
\newblock {\em Comput. Methods Appl. Mech. Engrg.}, 354:96--118, 2019.

\bibitem{BoDPS:17}
M.~Botti, D.~A. Di~Pietro, and P.~Sochala.
\newblock A {H}ybrid {H}igh-{O}rder method for nonlinear elasticity.
\newblock {\em SIAM J. Numer. Anal.}, 55(6):2687--2717, 2017.

\bibitem{BoDPS:20}
M.~Botti, D.~A. Di~Pietro, and P.~Sochala.
\newblock A hybrid high-order discretization method for nonlinear
  poroelasticity.
\newblock {\em Comput. Methods Appl. Math.}, 20(2):227--249, 2020.

\bibitem{Brenn:03}
S.~C. Brenner.
\newblock {P}oincar\'e-{F}riedrichs inequalities for piecewise {$H^1$}
  functions.
\newblock {\em SIAM J. Numer. Anal.}, 41(1):306--324, 2003.

\bibitem{brezis-68}
H.~Brezis.
\newblock \'{E}quations et in\'equations non lin\'eaires dans les espaces
  vectoriels en dualit\'e.
\newblock {\em Ann. Inst. Fourier (Grenoble)}, 18(fasc. 1):115--175, 1968.

\bibitem{Brezis:11}
H.~Brezis.
\newblock {\em Functional analysis, {S}obolev spaces and partial differential
  equations}.
\newblock Universitext. Springer, New York, NY, 2011.

\bibitem{BCDE:21}
E.~Burman, M.~Cicuttin, G.~Delay, and A.~Ern.
\newblock An unfitted hybrid high-order method with cell agglomeration for
  elliptic interface problems.
\newblock {\em SIAM J. Sci. Comput.}, 43(2):A859--A882, 2021.

\bibitem{BuDeE:21}
E.~Burman, G.~Delay, and A.~Ern.
\newblock An unfitted hybrid high-order method for the {S}tokes interface
  problem.
\newblock {\em IMA J. Numer. Anal.}, 2021.
\newblock hal-02280426.

\bibitem{BDE:21}
E.~Burman, O.~Duran, and A.~Ern.
\newblock Hybrid high-order methods for the acoustic wave equation in the time
  domain.
\newblock {\em Commun. Appl. Math. Comput.}, 2021.
\newblock hal-02922702.

\bibitem{BDES:21}
E.~Burman, O.~Duran, A.~Ern, and M.~Steins.
\newblock Convergence analysis of hybrid high-order methods for the wave
  equation.
\newblock {\em J. Sci. Comput.}, 87(3):Paper No. 91, 30, 2021.

\bibitem{BurEr:18}
E.~Burman and A.~Ern.
\newblock An unfitted hybrid high-order method for elliptic interface problems.
\newblock {\em SIAM J. Numer. Anal.}, 56(3):1525--1546, 2018.

\bibitem{BurEr:19}
E.~Burman and A.~Ern.
\newblock A cut cell hybrid high-order method for elliptic problems with curved
  boundaries.
\newblock In {\em Numerical mathematics and advanced applications---{ENUMATH}
  2017}, volume 126 of {\em Lect. Notes Comput. Sci. Eng.}, pages 173--181.
  Springer, Cham, 2019.

\bibitem{BusMu:20}
R.~Bustinza and J.~Munguia-La-Cotera.
\newblock A hybrid high-order formulation for a {N}eumann problem on polytopal
  meshes.
\newblock {\em Numer. Methods Partial Differential Equations}, 36(3):524--551,
  2020.

\bibitem{CaCDE:19}
V.~Calo, M.~Cicuttin, Q.~Deng, and A.~Ern.
\newblock Spectral approximation of elliptic operators by the hybrid high-order
  method.
\newblock {\em Math. Comp.}, 88(318):1559--1586, 2019.

\bibitem{CaDoG:21}
A.~Cangiani, Z.~Dong, and E.~H. Georgoulis.
\newblock {$hp$}-version discontinuous {G}alerkin methods on essentially
  arbitrarily-shaped elements.
\newblock \emph{Math. Comp.}, published online, arXiv preprint 1906.01715,
  2021.

\bibitem{CaDGH:17}
A.~Cangiani, Z.~Dong, E.~H. Georgoulis, and P.~Houston.
\newblock {\em {$hp$}-version discontinuous {G}alerkin methods on polygonal and
  polyhedral meshes}.
\newblock SpringerBriefs in Mathematics. Springer, Cham, 2017.

\bibitem{CaErP:21}
C.~Carstensen, A.~Ern, and S.~Puttkammer.
\newblock Guaranteed lower bounds on eigenvalues of elliptic operators with a
  hybrid high-order method.
\newblock {\em Numer. Math.}, 149(2):273--304, 2021.

\bibitem{CarFu:00}
C.~Carstensen and S.~A. Funken.
\newblock Constants in {C}l\'ement-interpolation error and residual based a
  posteriori error estimates in finite element methods.
\newblock {\em East-West J. Numer. Math.}, 8(3):153--175, 2000.

\bibitem{CaBCE:18}
K.~L. Cascavita, J.~Bleyer, X.~Chateau, and A.~Ern.
\newblock Hybrid discretization methods with adaptive yield surface detection
  for {B}ingham pipe flows.
\newblock {\em J. Sci. Comput.}, 77(3):1424--1443, 2018.

\bibitem{CaChE:20}
K.~L. Cascavita, F.~Chouly, and A.~Ern.
\newblock Hybrid high-order discretizations combined with {N}itsche's method
  for {D}irichlet and {S}ignorini boundary conditions.
\newblock {\em IMA J. Numer. Anal.}, 40(4):2189--2226, 2020.

\bibitem{CaDiP:20}
D.~Casta{\~n}{\'o}n~Quiroz and D.~A. Di~Pietro.
\newblock A {Hybrid High-Order} method for the incompressible {Navier--Stokes}
  problem robust for large irrotational body forces.
\newblock {\em Comput. Math. Appl.}, 79(8):2655--2677, 2020.

\bibitem{ChELV:21}
T.~Chaumont-Frelet, A.~Ern, S.~Lemaire, and F.~Valentin.
\newblock Bridging the multiscale hybrid-mixed and multiscale hybrid high-order
  methods.
\newblock hal-03235525, 2021.

\bibitem{ChDPF:18}
F.~Chave, D.~A. Di~Pietro, and L.~Formaggia.
\newblock A {Hybrid High-Order} method for {Darcy} flows in fractured porous
  media.
\newblock {\em SIAM J. Sci. Comput.}, 40(2):1063--1094, 2018.

\bibitem{ChDPL:20}
F.~Chave, D.~A. Di~Pietro, and S.~Lemaire.
\newblock A discrete {Weber} inequality on three-dimensional hybrid spaces with
  application to the {HHO} approximation of magnetostatics.
\newblock hal-02892526, 2020.

\bibitem{ChDMP:16}
F.~Chave, D.~A. Di~Pietro, F.~Marche, and F.~Pigeonneau.
\newblock A {H}ybrid {H}igh-{O}rder method for the {C}ahn-{H}illiard problem in
  mixed form.
\newblock {\em SIAM J. Numer. Anal.}, 54(3):1873--1898, 2016.

\bibitem{Chi2017}
H.~Chi, L.~{Beir\~ao da Veiga}, and G.~H. Paulino.
\newblock Some basic formulations of the virtual element method ({VEM}) for
  finite deformations.
\newblock {\em Comput. Methods Appl. Mech. Engrg.}, 318:148--192, 2017.

\bibitem{ChLaS:15}
E.~B. Chin, J.~B. Lasserre, and N.~Sukumar.
\newblock Numerical integration of homogeneous functions on convex and
  nonconvex polygons and polyhedra.
\newblock {\em Comput. Mech.}, 56(6):967--981, 2015.

\bibitem{Chouly:14}
F.~Chouly.
\newblock An adaptation of {N}itsche's method to the {T}resca friction problem.
\newblock {\em J. Math. Anal. Appl.}, 411:329--339, 2014.

\bibitem{ChErPi:20}
F.~Chouly, A.~Ern, and N.~Pignet.
\newblock A hybrid high-order discretization combined with {N}itsche's method
  for contact and {T}resca friction in small strain elasticity.
\newblock {\em SIAM J. Sci. Comput.}, 42(4):A2300--A2324, 2020.

\bibitem{ChoHi:13}
F.~Chouly and P.~Hild.
\newblock A {N}itsche-based method for unilateral contact problems: numerical
  analysis.
\newblock {\em SIAM J. Numer. Anal.}, 51(2):1295--1307, 2013.

\bibitem{ChHiR:15}
F.~Chouly, P.~Hild, and Y.~Renard.
\newblock Symmetric and non-symmetric variants of {N}itsche's method for
  contact problems in elasticity: theory and numerical experiments.
\newblock {\em Math. Comp.}, 84(293):1089--1112, 2015.

\bibitem{Ciarlet1988}
P.~G. Ciarlet.
\newblock {\em Mathematical elasticity. {V}ol. {I}}, volume~20 of {\em Studies
  in Mathematics and its Applications}.
\newblock North-Holland Publishing Co., Amsterdam, 1988.

\bibitem{Ciarlet1978}
P.~G. Ciarlet.
\newblock {\em The finite element method for elliptic problems}, volume~40 of
  {\em Classics in Applied Mathematics}.
\newblock Society for Industrial and Applied Mathematics (SIAM), Philadelphia,
  PA, 2002.
\newblock Reprint of the 1978 original [North-Holland, Amsterdam].

\bibitem{CiDPE:18}
M.~Cicuttin, D.~A. Di~Pietro, and A.~Ern.
\newblock {Implementation of Discontinuous Skeletal methods on
  arbitrary-dimensional, polytopal meshes using generic programming}.
\newblock {\em J. Comput. Appl. Math.}, 344:852--874, 2018.

\bibitem{CiErG:20}
M.~Cicuttin, A.~Ern, and T.~Gudi.
\newblock Hybrid high-order methods for the elliptic obstacle problem.
\newblock {\em J. Sci. Comput.}, 83(1):Paper No. 8, 18, 2020.

\bibitem{CiErL:19}
M.~Cicuttin, A.~Ern, and S.~Lemaire.
\newblock A {Hybrid High-Order} method for highly oscillatory elliptic
  problems.
\newblock {\em Comput. Methods Appl. Math.}, 19(4):723--748, 2019.

\bibitem{Cockburn:16}
B.~Cockburn.
\newblock Static condensation, hybridization, and the devising of the {HDG}
  methods.
\newblock In G.~R. Barrenechea, F.~Brezzi, A.~Cangiani, and E.~H. Georgoulis,
  editors, {\em Building Bridges: Connections and Challenges in Modern
  Approaches to Numerical Partial Differential Equations}, volume 114 of {\em
  Lecture Notes in Computational Science and Engineering}, pages 129--178.
  Springer, Cham, 2016.

\bibitem{CoDPE:16}
B.~Cockburn, D.~A. Di~Pietro, and A.~Ern.
\newblock Bridging the {H}ybrid {H}igh-{O}rder and hybridizable discontinuous
  {G}alerkin methods.
\newblock {\em ESAIM Math. Model. Numer. Anal.}, 50(3):635--650, 2016.

\bibitem{CFHJSS:18}
B.~Cockburn, Z.~Fu, A.~Hungria, L.~Ji, M.~A. S\'{a}nchez, and F.-J. Sayas.
\newblock St{\"o}rmer-{N}umerov {HDG} methods for acoustic waves.
\newblock {\em J. Sci. Comput.}, 75(2):597--624, 2018.

\bibitem{CoGoL:09}
B.~Cockburn, J.~Gopalakrishnan, and R.~Lazarov.
\newblock Unified hybridization of discontinuous {G}alerkin, mixed, and
  continuous {G}alerkin methods for second order elliptic problems.
\newblock {\em SIAM J. Numer. Anal.}, 47(2):1319--1365, 2009.

\bibitem{CoGoS:10}
B.~Cockburn, J.~Gopalakrishnan, and F.-J. Sayas.
\newblock A projection-based error analysis of {HDG} methods.
\newblock {\em Math. Comp.}, 79(271):1351--1367, 2010.

\bibitem{CocQS:12}
B.~Cockburn, W.~Qiu, and K.~Shi.
\newblock Conditions for superconvergence of {HDG} methods for second-order
  elliptic problems.
\newblock {\em Math. Comp.}, 81(279):1327--1353, 2012.

\bibitem{CoQue:14}
B.~Cockburn and V.~Quenneville-B\'{e}lair.
\newblock Uniform-in-time superconvergence of the {HDG} methods for the
  acoustic wave equation.
\newblock {\em Math. Comp.}, 83(285):65--85, 2014.

\bibitem{CurnierA:88}
A.~Curnier and P.~Alart.
\newblock A generalized {N}ewton method for contact problems with friction.
\newblock {\em J. M\'ec. Th\'eor. Appl.}, 7(suppl. 1):67--82, 1988.

\bibitem{DabDe:20}
J.~Dabaghi and G.~Delay.
\newblock A unified framework for high-order numerical discretizations of
  variational inequalities.
\newblock hal-02969793, 2020.

\bibitem{Maso2006}
G.~Dal~Maso, A.~DeSimone, and M.~G. Mora.
\newblock {Q}uasistatic evolution problems for linearly elastic--perfectly
  plastic materials.
\newblock {\em Arch. Ration. Mech. Anal.}, 180(2):237--291, 2006.

\bibitem{Dauge_1988}
M.~Dauge.
\newblock {\em Elliptic boundary value problems on corner domains}, volume 1341
  of {\em Lecture Notes in Mathematics}.
\newblock Springer-Verlag, Berlin, 1988.

\bibitem{DiPDr:17}
D.~A. Di~Pietro and J.~Droniou.
\newblock A {H}ybrid {H}igh-{O}rder method for {L}eray-{L}ions elliptic
  equations on general meshes.
\newblock {\em Math. Comp.}, 86(307):2159--2191, 2017.

\bibitem{DiPDr:20}
D.~A. Di~Pietro and J.~Droniou.
\newblock {\em The {Hybrid High-Order} method for polytopal meshes}, volume~19
  of {\em Modeling, Simulation and Application}.
\newblock Springer, Cham, 2020.

\bibitem{DPDEr:15}
D.~A. Di~Pietro, J.~Droniou, and A.~Ern.
\newblock A discontinuous-skeletal method for advection-diffusion-reaction on
  general meshes.
\newblock {\em SIAM J. Numer. Anal.}, 53(5):2135--2157, 2015.

\bibitem{DiPDM:18}
D.~A. Di~Pietro, J.~Droniou, and G.~Manzini.
\newblock {D}iscontinuous {S}keletal {G}radient {D}iscretisation {M}ethods on
  polytopal meshes.
\newblock {\em J. Comput. Phys.}, 355:397--425, 2018.

\bibitem{DiPEr:12}
D.~A. Di~Pietro and A.~Ern.
\newblock {\em Mathematical {A}spects of {D}iscontinuous {G}alerkin {M}ethods},
  volume~69 of {\em Math\'ematiques \& Applications}.
\newblock Springer-Verlag, Berlin, 2012.

\bibitem{DiPEr:15}
D.~A. Di~Pietro and A.~Ern.
\newblock A {Hybrid High-Order} locking-free method for linear elasticity on
  general meshes.
\newblock {\em Comput. Meth. Appl. Mech. Engrg.}, 283:1--21, 2015.

\bibitem{DiPEr:17}
D.~A. Di~Pietro and A.~Ern.
\newblock Arbitrary-order mixed methods for heterogeneous anisotropic diffusion
  on general meshes.
\newblock {\em IMA J. Numer. Anal.}, 37(1):40--63, 2017.
\newblock Preprint originally available at hal-00918482v1 (2013).

\bibitem{DiPEL:14}
D.~A. Di~Pietro, A.~Ern, and S.~Lemaire.
\newblock An arbitrary-order and compact-stencil discretization of diffusion on
  general meshes based on local reconstruction operators.
\newblock {\em Comput. Meth. Appl. Math.}, 14(4):461--472, 2014.

\bibitem{DiPEL:16}
D.~A. Di~Pietro, A.~Ern, and S.~Lemaire.
\newblock A review of {H}ybrid {H}igh-{O}rder methods: formulations,
  computational aspects, comparison with other methods.
\newblock In {\em Building bridges: connections and challenges in modern
  approaches to numerical partial differential equations}, volume 114 of {\em
  Lect. Notes Comput. Sci. Eng.}, pages 205--236. Springer, Cham, 2016.

\bibitem{DPELS:16}
D.~A. Di~Pietro, A.~Ern, A.~Linke, and F.~Schieweck.
\newblock A discontinuous skeletal method for the viscosity-dependent {S}tokes
  problem.
\newblock {\em Comput. Methods Appl. Mech. Engrg.}, 306:175--195, 2016.

\bibitem{DiPKr:18}
D.~A. Di~Pietro and S.~Krell.
\newblock A {Hybrid High-Order} method for the steady incompressible
  {Navier--Stokes} problem.
\newblock {\em J. Sci. Comput.}, 74(3):1677--1705, 2018.

\bibitem{DiPSp:16}
D.~A. Di~Pietro and R.~Specogna.
\newblock An a posteriori-driven adaptive mixed high-order method with
  application to electrostatics.
\newblock {\em J. Comput. Phys.}, 326:35--55, 2016.

\bibitem{Djoko2007a}
J.~K. Djoko, F.~Ebobisse, A.~T. McBride, and B.~D. Reddy.
\newblock A discontinuous {G}alerkin formulation for classical and gradient
  plasticity. {I}. {F}ormulation and analysis.
\newblock {\em Comput. Methods Appl. Mech. Engrg.}, 196(37-40):3881--3897,
  2007.

\bibitem{DrEGH:10}
J.~Droniou, R.~Eymard, T.~Gallou\"{e}t, and R.~Herbin.
\newblock A unified approach to mimetic finite difference, hybrid finite volume
  and mixed finite volume methods.
\newblock {\em Math. Models Methods Appl. Sci.}, 20(2):265--295, 2010.

\bibitem{Droniou2015}
J.~Droniou and B.~P. Lamichhane.
\newblock Gradient schemes for linear and non-linear elasticity equations.
\newblock {\em Numer. Math.}, 129(2):251--277, 2015.

\bibitem{DuSay:19}
S.~Du and F.-J. Sayas.
\newblock {\em An invitation to the theory of the hybridizable discontinuous
  {G}alerkin method}.
\newblock SpringerBriefs in Mathematics. Springer, Cham, 2019.

\bibitem{Dunavant1985}
D.~A. Dunavant.
\newblock High degree efficient symmetrical {G}aussian quadrature rules for the
  triangle.
\newblock {\em Internat. J. Numer. Methods Engrg.}, 21(6):1129--1148, 1985.

\bibitem{Dupont:73}
T.~Dupont.
\newblock {$L\sp{2}$}-estimates for {G}alerkin methods for second order
  hyperbolic equations.
\newblock {\em SIAM J. Numer. Anal.}, 10:880--889, 1973.

\bibitem{codeaster}
{Electricit{\'e} de France}.
\newblock Finite element {\it code$\_$aster}, structures and thermomechanics
  analysis for studies and research.
\newblock Open source on www.code-aster.org, 1989--2019.

\bibitem{ErnGu:17}
A.~Ern and J.-L. Guermond.
\newblock Finite element quasi-interpolation and best approximation.
\newblock {\em M2AN Math. Model. Numer. Anal.}, 51(4):1367--1385, 2017.

\bibitem{ErnGu:21a}
A.~Ern and J.-L. Guermond.
\newblock {\em Finite Elements {I}: {A}pproximation and interpolation},
  volume~72 of {\em Texts in Applied Mathematics}.
\newblock Springer, Cham, 2021.

\bibitem{ErnGu:21b}
A.~Ern and J.-L. Guermond.
\newblock {\em Finite Elements {II}: {G}alerkin approximation, elliptic and
  mixed {PDE}s}, volume~73 of {\em Texts in Applied Mathematics}.
\newblock Springer, Cham, 2021.

\bibitem{ErnVo:20}
A.~Ern and M.~Vohral{\'{\i}}k.
\newblock Stable broken {$H^1$} and {${\boldsymbol H}(\mathrm{div})$}
  polynomial extensions for polynomial-degree-robust potential and flux
  reconstruction in three space dimensions.
\newblock {\em Math. Comp.}, 89(322):551--594, 2020.

\bibitem{ErnZa:20}
A.~Ern and P.~Zanotti.
\newblock A quasi-optimal variant of the hybrid high-order method for elliptic
  partial differential equations with {$H^{-1}$} loads.
\newblock {\em IMA J. Numer. Anal.}, 40:2163--2188, 2020.

\bibitem{Evans:98}
L.~C. Evans.
\newblock {\em Partial differential equations}, volume~19 of {\em Graduate
  Studies in Mathematics}.
\newblock American Mathematical Society, Providence, RI, 1998.

\bibitem{Eymard2010}
R.~{Eymard}, T.~{Gallou\"{e}t}, and R.~{Herbin}.
\newblock Discretization of heterogeneous and anisotropic diffusion problems on
  general nonconforming meshes {SUSHI}: a scheme using stabilization and hybrid
  interfaces.
\newblock {\em IMA J. Numer. Anal.}, 30(4):1009--1043, 2010.

\bibitem{FuCoS:15}
G.~Fu, B.~Cockburn, and H.~Stolarski.
\newblock Analysis of an {HDG} method for linear elasticity.
\newblock {\em Internat. J. Numer. Methods Engrg.}, 102(3-4):551--575, 2015.

\bibitem{Gr85}
P.~Grisvard.
\newblock {\em Elliptic problems in nonsmooth domains}, volume~24 of {\em
  Monographs and Studies in Mathematics}.
\newblock Pitman (Advanced Publishing Program), Boston, MA, 1985.

\bibitem{grundmol1978}
A.~Grundmann and H.~M. Moller.
\newblock Invariant integration formulas for the n-simplex by combinatorial
  methods.
\newblock {\em SIAM J. Numer. Analysis}, 15(2):282--290, 1978.

\bibitem{GuGuZ:18}
Q.~Guan, M.~Gunzburger, and W.~Zhao.
\newblock Weak-{G}alerkin finite element methods for a second-order elliptic
  variational inequality.
\newblock {\em Comput. Methods Appl. Mech. Engrg.}, 337:677--688, 2018.

\bibitem{Halphen1975}
B.~Halphen and Q.~Son~Nguyen.
\newblock Sur les mat{\'e}riaux standard g{\'e}n{\'e}ralis{\'e}s.
\newblock {\em J. Mecanique.}, 14:39--63, 1975.

\bibitem{Han2013}
W.~Han and B.~D. Reddy.
\newblock {\em {P}lasticity: {M}athematical {T}heory and {N}umerical
  {A}nalysis}.
\newblock Springer, New York, 2013.

\bibitem{Hansbo2010}
P.~Hansbo.
\newblock A discontinuous finite element method for elasto-plasticity.
\newblock {\em Int. J. Numer. Meth. Biomed. Eng.}, 26(6):780--789, 2010.

\bibitem{HarParVal13}
C.~Harder, D.~Paredes, and F.~Valentin.
\newblock A family of multiscale hybrid-mixed finite element methods for the
  {D}arcy equation with rough coefficients.
\newblock {\em J. Comput. Phys.}, 245:107--130, 2013.

\bibitem{HePiE:21}
F.~H{\'e}din, G.~Pichot, and A.~Ern.
\newblock A hybrid high-order method for flow simulations in discrete fracture
  networks.
\newblock In {\em Numerical mathematics and advanced applications---{ENUMATH}
  2019}, Lect. Notes Comput. Sci. Eng. Springer, Cham, 2021.

\bibitem{HesHa:08}
J.~S. Hesthaven and T.~Warburton.
\newblock {\em Nodal discontinuous {G}alerkin methods: {A}lgorithms, analysis,
  and applications}, volume~54 of {\em Texts in Applied Mathematics}.
\newblock Springer, New York, NY, 2008.

\bibitem{Horgan:95}
C.~O. Horgan.
\newblock {Korn}'s inequalities and their applications in continuum mechanics.
\newblock {\em SIAM Rev.}, 37:491--511, 1995.

\bibitem{JohLa:13}
A.~Johansson and M.~G. Larson.
\newblock A high order discontinuous {G}alerkin {N}itsche method for elliptic
  problems with fictitious boundary.
\newblock {\em Numer. Math.}, 123(4):607--628, 2013.

\bibitem{JoNeS:16}
L.~John, M.~Neilan, and I.~Smears.
\newblock Stable discontinuous {G}alerkin {FEM} without penalty parameters.
\newblock In {\em Numerical Mathematics and Advanced Applications ENUMATH
  2015}, Lecture Notes in Computational Science and Engineering, pages
  165--173. Springer, 2016.

\bibitem{KaLeC:2015}
H.~Kabaria, A.~J. Lew, and B.~Cockburn.
\newblock A hybridizable discontinuous {G}alerkin formulation for non-linear
  elasticity.
\newblock {\em Comput. Methods Appl. Mech. Engrg.}, 283:303--329, 2015.

\bibitem{Keast1986}
P.~Keast.
\newblock Moderate-degree tetrahedral quadrature formulas.
\newblock {\em Comput. Methods Appl. Mech. Engrg.}, 55(3):339 -- 348, 1986.

\bibitem{kikuchi-oden-88}
N.~Kikuchi and J.~T. Oden.
\newblock {\em Contact problems in elasticity: a study of variational
  inequalities and finite element methods}, volume~8 of {\em SIAM Studies in
  Applied Mathematics}.
\newblock Society for Industrial and Applied Mathematics (SIAM), Philadelphia,
  PA, 1988.

\bibitem{KrWWW:16}
J.~Kr{\"a}mer, C.~Wieners, B.~Wohlmuth, and L.~Wunderlich.
\newblock A hybrid weakly nonconforming discretization for linear elasticity.
\newblock {\em Proc. Appl. Math. Mech.}, 16(1):849--850, 2016.

\bibitem{Lehrenfeld:10}
C.~Lehrenfeld.
\newblock {\em Hybrid Discontinuous Galerkin methods for solving incompressible
  flow problems}.
\newblock PhD thesis, Rheinisch-Westf\"alische Technische Hochschule (RWTH)
  Aachen, Germany, 2010.

\bibitem{LehSc:16}
C.~Lehrenfeld and J.~Sch\"oberl.
\newblock High order exactly divergence-free hybrid discontinuous {G}alerkin
  methods for unsteady incompressible flows.
\newblock {\em Comput. Methods Appl. Mech. Engrg.}, 307:339--361, 2016.

\bibitem{Lemaire:21}
S.~Lemaire.
\newblock Bridging the hybrid high-order and virtual element methods.
\newblock {\em IMA J. Numer. Anal.}, 41(1):549--593, 2021.

\bibitem{Lemaitre1994}
J.~Lemaitre and J.-L. Chaboche.
\newblock {\em Mechanics of Solid Materials}.
\newblock University Press, Cambridge, 1994.

\bibitem{LipMa:14}
K.~Lipnikov and G.~Manzini.
\newblock A high-order mimetic method on unstructured polyhedral meshes for the
  diffusion equation.
\newblock {\em J. Comput. Phys.}, 272:360--385, 2014.

\bibitem{Liu2010}
R.~Liu, M.~F. Wheeler, C.~N. Dawson, and R.~H. Dean.
\newblock A fast convergent rate preserving discontinuous {G}alerkin framework
  for rate-independent plasticity problems.
\newblock {\em Comput. Methods Appl. Mech. Engrg.}, 199(49-52):3213--3226,
  2010.

\bibitem{Liu2013}
R.~Liu, M.~F. Wheeler, and I.~Yotov.
\newblock On the spatial formulation of discontinuous {G}alerkin methods for
  finite elastoplasticity.
\newblock {\em Comput. Methods Appl. Mech. Engrg.}, 253:219--236, 2013.

\bibitem{MacLean_2000}
W.~McLean.
\newblock {\em Strongly elliptic systems and boundary integral equations}.
\newblock Cambridge University Press, Cambridge, 2000.

\bibitem{Miehe2002}
C.~Miehe, N.~Apel, and M.~Lambrecht.
\newblock Anisotropic additive plasticity in the logarithmic strain space:
  modular kinematic formulation and implementation based on incremental
  minimization principles for standard materials.
\newblock {\em Comput. Methods Appl. Mech. Engrg.}, 191(47-48):5383--5425,
  2002.

\bibitem{MonSu:99}
P.~Monk and E.~S{\"u}li.
\newblock The adaptive computation of far-field patterns by a posteriori error
  estimation of linear functionals.
\newblock {\em SIAM J. Numer. Anal.}, 36(1):251--274, 1999.

\bibitem{MuWaY:15}
L.~Mu, J.~Wang, and X.~Ye.
\newblock A weak {G}alerkin finite element method with polynomial reduction.
\newblock {\em J. Comput. Appl. Math.}, 285:45--58, 2015.

\bibitem{Nguyen2012a}
N.~C. Nguyen and J.~Peraire.
\newblock Hybridizable discontinuous {G}alerkin methods for partial
  differential equations in continuum mechanics.
\newblock {\em J. Comput. Phys.}, 231(18):5955--5988, 2012.

\bibitem{NgPeC:11}
N.~C. Nguyen, J.~Peraire, and B.~Cockburn.
\newblock High-order implicit hybridizable discontinuous {G}alerkin methods for
  acoustics and elastodynamics.
\newblock {\em J. Comput. Phys.}, 230(10):3695--3718, 2011.

\bibitem{nit71}
J.~Nitsche.
\newblock \"{U}ber ein {V}ariationsprinzip zur {L}\"osung von
  {D}irichlet-{P}roblemen bei {V}erwendung von {T}eilr\"aumen, die keinen
  {R}andbedingungen unterworfen sind.
\newblock {\em Abh. Math. Sem. Univ. Hamburg}, 36:9--15, 1971.

\bibitem{Noels2006}
L.~Noels and R.~Radovitzky.
\newblock A general discontinuous {G}alerkin method for finite hyperelasticity.
  {F}ormulation and numerical applications.
\newblock {\em Internat. J. Numer. Methods Engrg.}, 68(1):64--97, 2006.

\bibitem{Ogden1997}
R.~W. Ogden.
\newblock {\em Non-linear elastic deformations}.
\newblock Dover Publications Inc., New York, NY, 1997.

\bibitem{Oikawa:15}
I.~Oikawa.
\newblock A hybridized discontinuous {G}alerkin method with reduced
  stabilization.
\newblock {\em J. Sci. Comput.}, 65(1):327--340, 2015.

\bibitem{PayWe:60}
L.~E. Payne and H.~F. Weinberger.
\newblock An optimal {P}oincar\'e inequality for convex domains.
\newblock {\em Arch. Rational Mech. Anal.}, 5:286--292, 1960.

\bibitem{SCNPC:17}
M.~A. S\'{a}nchez, C.~Ciuca, N.~C. Nguyen, J.~Peraire, and B.~Cockburn.
\newblock Symplectic {H}amiltonian {HDG} methods for wave propagation
  phenomena.
\newblock {\em J. Comput. Phys.}, 350:951--973, 2017.

\bibitem{Simo1992b}
J.~C. Simo.
\newblock Algorithms for static and dynamic multiplicative plasticity that
  preserve the classical return mapping schemes of the infinitesimal theory.
\newblock {\em Comput. Methods Appl. Mech. Engrg.}, 99:61--112, 1992.

\bibitem{Simo1998}
J.~C. Simo and T.~J.~R. Hughes.
\newblock {\em Computational Inelasticity}.
\newblock Springer, Berlin, 1998.

\bibitem{SomVi:07}
A.~Sommariva and M.~Vianello.
\newblock Product {G}auss cubature over polygons based on {G}reen's integration
  formula.
\newblock {\em BIT}, 47(2):441--453, 2007.

\bibitem{SomVi:09}
A.~Sommariva and M.~Vianello.
\newblock Gauss-{G}reen cubature and moment computation over arbitrary
  geometries.
\newblock {\em J. Comput. Appl. Math.}, 231(2):886--896, 2009.

\bibitem{Soon2008}
S.-C. Soon.
\newblock {\em Hybridizable Discontinuous Galerkin Method for Solid Mechanics}.
\newblock PhD thesis, University of Minnesota, MN, 2008.

\bibitem{SoCoS:09}
S.-C. Soon, B.~Cockburn, and H.~K. Stolarski.
\newblock A hybridizable discontinuous {G}alerkin method for linear elasticity.
\newblock {\em Internat. J. Numer. Methods Engrg.}, 80(8):1058--1092, 2009.

\bibitem{StNPC:16}
M.~Stanglmeier, N.~C. Nguyen, J.~Peraire, and B.~Cockburn.
\newblock An explicit hybridizable discontinuous {G}alerkin method for the
  acoustic wave equation.
\newblock {\em Comput. Methods Appl. Mech. Engrg.}, 300:748--769, 2016.

\bibitem{SudWa:13}
Y.~Sudhakar and W.~A. Wall.
\newblock Quadrature schemes for arbitrary convex/concave volumes and
  integration of weak form in enriched partition of unity methods.
\newblock {\em Comput. Methods Appl. Mech. Engrg.}, 258:39--54, 2013.

\bibitem{Eyck2008a}
A.~ten Eyck, F.~Celiker, and A.~Lew.
\newblock Adaptive stabilization of discontinuous {G}alerkin methods for
  nonlinear elasticity: analytical estimates.
\newblock {\em Comput. Methods Appl. Mech. Engrg.}, 197(33-40):2989--3000,
  2008.

\bibitem{Eyck2008}
A.~ten Eyck, F.~Celiker, and A.~Lew.
\newblock Adaptive stabilization of discontinuous {G}alerkin methods for
  nonlinear elasticity: motivation, formulation, and numerical examples.
\newblock {\em Comput. Methods Appl. Mech. Engrg.}, 197(45-48):3605--3622,
  2008.

\bibitem{Eyck2006}
A.~ten Eyck and A.~Lew.
\newblock Discontinuous {G}alerkin methods for non-linear elasticity.
\newblock {\em Internat. J. Numer. Methods Engrg.}, 67(9):1204--1243, 2006.

\bibitem{VeeVe:12}
A.~Veeser and R.~Verf{\"u}rth.
\newblock Poincar\'e constants for finite element stars.
\newblock {\em IMA J. Numer. Anal.}, 32(1):30--47, 2012.

\bibitem{WaWWZ:16}
C.~Wang, J.~Wang, R.~Wang, and R.~Zhang.
\newblock A locking-free weak {G}alerkin finite element method for elasticity
  problems in the primal formulation.
\newblock {\em J. Comput. Appl. Math.}, 307:346--366, 2016.

\bibitem{WangWei:18}
F.~Wang and H.~Wei.
\newblock Virtual element method for simplified friction problem.
\newblock {\em Appl. Math. Letters}, 85:125--131, 2018.

\bibitem{WangYe:13}
J.~Wang and X.~Ye.
\newblock A weak {G}alerkin finite element method for second-order elliptic
  problems.
\newblock {\em J. Comput. Appl. Math.}, 241:103--115, 2013.

\bibitem{WangYe:14}
J.~Wang and X.~Ye.
\newblock A weak {G}alerkin mixed finite element method for second order
  elliptic problems.
\newblock {\em Math. Comp.}, 83(289):2101--2126, 2014.

\bibitem{Wriggers2017a}
P.~Wriggers and B.~Hudobivnik.
\newblock A low order virtual element formulation for finite elasto-plastic
  deformations.
\newblock {\em Comput. Methods Appl. Mech. Engrg.}, 327:459--477, 2017.

\bibitem{Wriggers2017}
P.~Wriggers, B.~D. Reddy, W.~Rust, and B.~Hudobivnik.
\newblock Efficient virtual element formulations for compressible and
  incompressible finite deformations.
\newblock {\em Comput. Mech.}, 60(2):253--268, 2017.

\bibitem{WrRuR:16}
P.~Wriggers, W.~T. Rust, and B.~D. Reddy.
\newblock A virtual element method for contact.
\newblock {\em Comput. Mech.}, 58(6):1039--1050, 2016.

\bibitem{Wulfinghoff2017}
S.~Wulfinghoff, H.~R. Bayat, A.~Alipour, and S.~Reese.
\newblock A low-order locking-free hybrid discontinuous {G}alerkin element
  formulation for large deformations.
\newblock {\em Comput. Methods Appl. Mech. Engrg.}, 323:353--372, 2017.

\bibitem{YeZha:20}
X.~Ye and S.~Zhang.
\newblock A stabilizer-free weak {G}alerkin finite element method on polytopal
  meshes.
\newblock {\em J. Comput. Appl. Math.}, 371:112699, 9, 2020.

\bibitem{ZhWuX:19}
M.~Zhao, H.~Wu, and C.~Xiong.
\newblock Error analysis of {HDG} approximations for elliptic variational
  inequality: obstacle problem.
\newblock {\em Numer. Algorithms}, 81(2):445--463, 2019.

\end{thebibliography}



\end{document}